\documentclass[11pt, letterpaper]{amsart}

\usepackage{amssymb,latexsym, xcolor
}
\usepackage[dvips, hcentering,
includehead,width=14.2cm,
top=0.5cm,
height=21.6cm,
paperheight=23.8cm,paperwidth=17cm
]{geometry}
\usepackage{enumitem}
\usepackage{hyperref}
\usepackage{pict2e}
\usepackage{graphicx}
\usepackage{rotating}

\usepackage{ccaption}
\changecaptionwidth
\captionwidth{5.6in}

\renewcommand{\ref}{\hyperref}

\def\AA{{\mathbb A}}
\def\ZZ{{\mathbb Z}}
\def\RR{{\mathbb R}}
\def\QQ{{\mathbb Q}}

\def\CC{{\mathbb C}}

\def\PP{{\mathbb P}}

\def\eqdef{{\stackrel{\rm def}{=}}}
\def\muinf{{\mu(C,p,L^\infty)}}

\newcommand{\Disk}{{\mathbf{D}}}

\newcommand{\bn}{\boldsymbol{n}}
\newcommand{\VR}{V_{\RR}}
\newcommand{\CR}{C_{\RR}}
\newcommand{\Cm}{C_{(m)}}
\newcommand{\Cab}{C_{a,b}}
\newcommand{\Cabc}{C_{a,b,c}}

\newcommand{\nE}{\mathbf{n}(E)}

\DeclareMathOperator{\Span}{Span}
\DeclareMathOperator{\conv}{conv}
\DeclareMathOperator{\Fix}{Fix}
\DeclareMathOperator{\br}{Br}

\DeclareMathOperator{\Comp}{\mathbf{Comp}}

\DeclareMathOperator{\ord}{ord}
\DeclareMathOperator{\mt}{mult}
\DeclareMathOperator{\Sing}{\mathbf{Sing}}

\newtheorem{theorem}{Theorem}[section]
\newtheorem{proposition}[theorem]{Proposition}

\newtheorem{corollary}[theorem]{Corollary}
\newtheorem{lemma}[theorem]{Lemma}

\theoremstyle{definition}

\newtheorem{remark}[theorem]{Remark}

\newtheorem{example}[theorem]{Example}
\newtheorem{definition}[theorem]{Definition}
\newtheorem{problem}[theorem]{Problem}

\numberwithin{equation}{section}

\newcommand{\red}[1]{{\color{red} #1}}
\newcommand{\redsf}[1]{{#1}}
\newcommand{\blue}[1]{{\color{blue} #1}}
\newcommand{\green}[1]{{\color{green} #1}}
\newcommand{\cyan}[1]{{\color{cyan} #1}}
\newcommand{\gray}[1]{{\color{gray} #1}}

\newcommand{\eps}{{\varepsilon}}

\begin{document}

\title{Expressive curves}

\author{Sergey Fomin}
\address{Department of Mathematics, University of Michigan,
Ann Arbor, MI 48109, USA}
\email{fomin@umich.edu}

\author{Eugenii Shustin}
\address{School of Mathematical Sciences, Tel Aviv
University, 
Tel Aviv 69978,
Israel}
\email{shustin@tauex.tau.ac.il}


\thanks{The first author was supported by the NSF grants DMS-1664722 and DMS-2054231 and by
a Simons Fellowship.
The second author was supported by
the ISF grant  501/18 and the Bauer-Neuman Chair in Real and Complex Geometry.}

\date{August 20, 2023}

\begin{abstract}
We initiate the study of a class of real plane algebraic curves
which we call \emph{expressive}.
These are the curves whose defining polynomial has the smallest
number of critical  points allowed by the topology of the set of real points~of~a~curve.
This concept can be viewed as a global version of the notion of
a real morsification of an isolated plane curve singularity.

We prove that a plane curve~$C$ is expressive if (a)~each irreducible component
of~$C$~can be parametrized by real polynomials (either ordinary or trigonometric),
(b)~all singular points of $C$ in the affine plane are ordinary hyperbolic~nodes, and
(c)~the set of real points of $C$ in the affine plane is connected.
Conversely, an expressive curve with real irreducible components must satisfy conditions (a)--(c),
unless it exhibits some exotic behaviour at infinity.

We describe several constructions that produce expressive curves,
and discuss a large number of examples, including:
arrangements of lines, parabolas, and~circles;
Chebyshev and Lissajous curves;
hypotrochoids and epitrochoids; and much more.
\end{abstract}


\ \vspace{-0.2in}

\maketitle

\vspace{-.25in}

\tableofcontents

\vspace{-.35in}

\newpage

\section*{Introduction}

Let $g(x)\in\RR[x]$ be a polynomial of degree~$n$ whose $n$ roots are real and distinct.
Then $g$ has exactly $n-1$ critical points, all of them real, interlacing the roots of~$g$.

In this paper, we study the two-dimensional version of this phenomenon.
We~call~a bivariate real polynomial $G(x,y)\in\RR[x,y]$
(or the corresponding affine plane curve~$C$) 
\emph{expressive} if
\redsf{the locations of the critical points of~$G$ are
determined by the set of real points
$C_\RR=\{(x,y)\in\RR^2\mid G(x,y)=0\}$, as follows:
\begin{itemize}[leftmargin=.2in]
\item
there is precisely one extremum inside each bounded region of~$\RR^2\setminus\CR$;
\item
all other critical points of~$G$ are the saddles located at hyperbolic nodes of~$C$.
\end{itemize}
(Recall that a \emph{hyperbolic node} is an intersection of two smooth real local branches.)
In particular, all critical points of an expressive polynomial~$G$ are real.
}

\begin{figure}[ht]
\begin{center}
\includegraphics[scale=0.18, trim=12cm 15cm 15cm 10.5cm, clip]{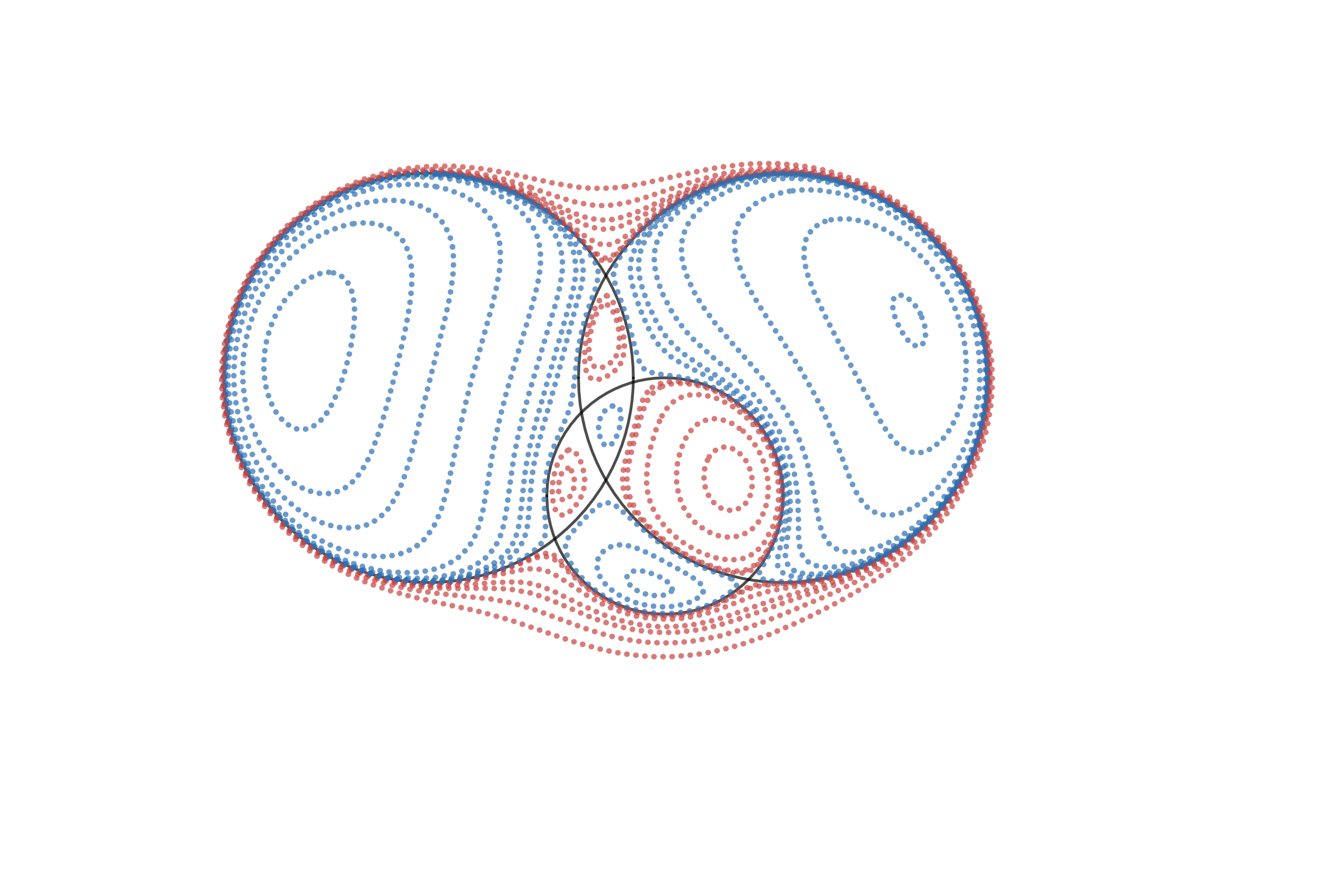}
\end{center}
\caption{
The expressive curve $C=\{G(x,y)=0\}$ in the picture is a union~of three circles, shown in solid black.
Dotted isolines represent level sets of~$G$.
The~polynomial~$G$ has 13 critical points:
6~saddles located at the double points (the hyperbolic nodes of~$C$)
plus 7 extrema, 
one in each bounded region of~$\RR^2\setminus C$.
}
\label{fig:3-circles}
\end{figure}

\begin{figure}[ht]
\begin{center}
\vspace{-.3in}
\includegraphics[scale=0.22, trim=15cm 24cm 15cm 13.5cm, clip]{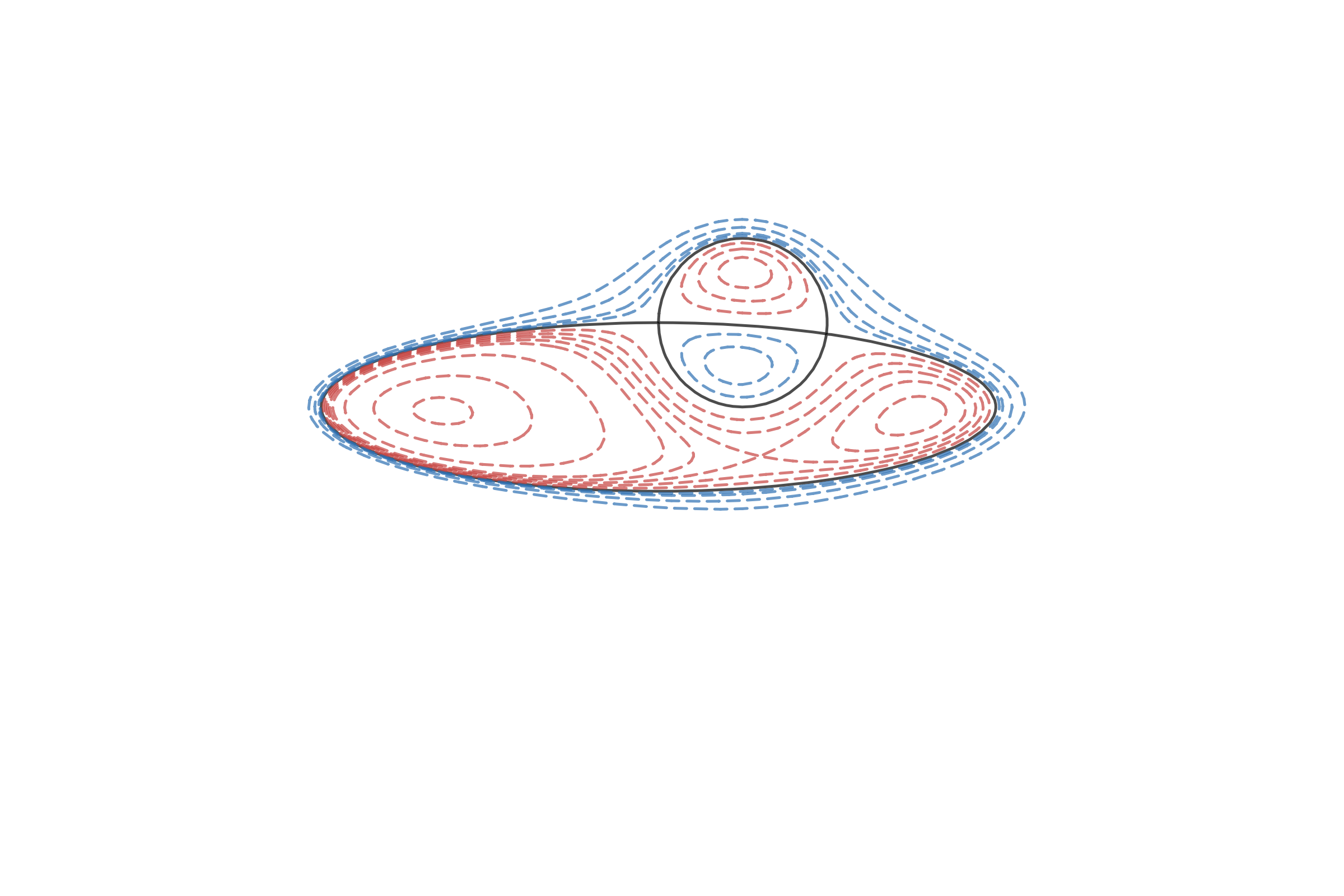}
\vspace{-.1in}
\end{center}
\caption{A non-expressive curve whose components are
a circle and an ellipse.
The bounded region at the bottom contains 3 critical points.}
\label{fig:isolines}
\end{figure}

An example is shown in Figure~\ref{fig:3-circles}.
For a non-example, see Figure~\ref{fig:isolines}.
\vspace{.1in}


\pagebreak

Our main result (Theorem~\ref{th:reg-expressive}) gives an explicit characterization of expressive curves,
subject to a mild requirement of ``$L^\infty$-regularity.'' 
(This requirement forbids some exotic behaviour of~$C$ at infinity.)
We prove that a plane algebraic curve~$C$
with real irreducible components is expressive and $L^\infty$-regular
if and only if
\begin{itemize}[leftmargin=.2in]
\item
each component of $C$ has a trigonometric or polynomial parametrization,
\item
all singular points of $C$ in the affine plane are real hyperbolic nodes, and
\item
the set of real points of $C$ in the affine plane is connected.
\end{itemize}

To illustrate, a union of circles is an expressive curve
provided any two of them intersect at two real points, as in Figure~\ref{fig:3-circles}.
On the other hand, the circle and the ellipse in Figure~\ref{fig:isolines}
intersect at four points, two of which are complex conjugate.
(In~the case of a pair of circles, those two points escape to infinity.)
\medskip

\begin{figure}[ht]
\begin{center}
\includegraphics[scale=0.232, trim=13.5cm 26cm 10.3cm 21cm, clip]{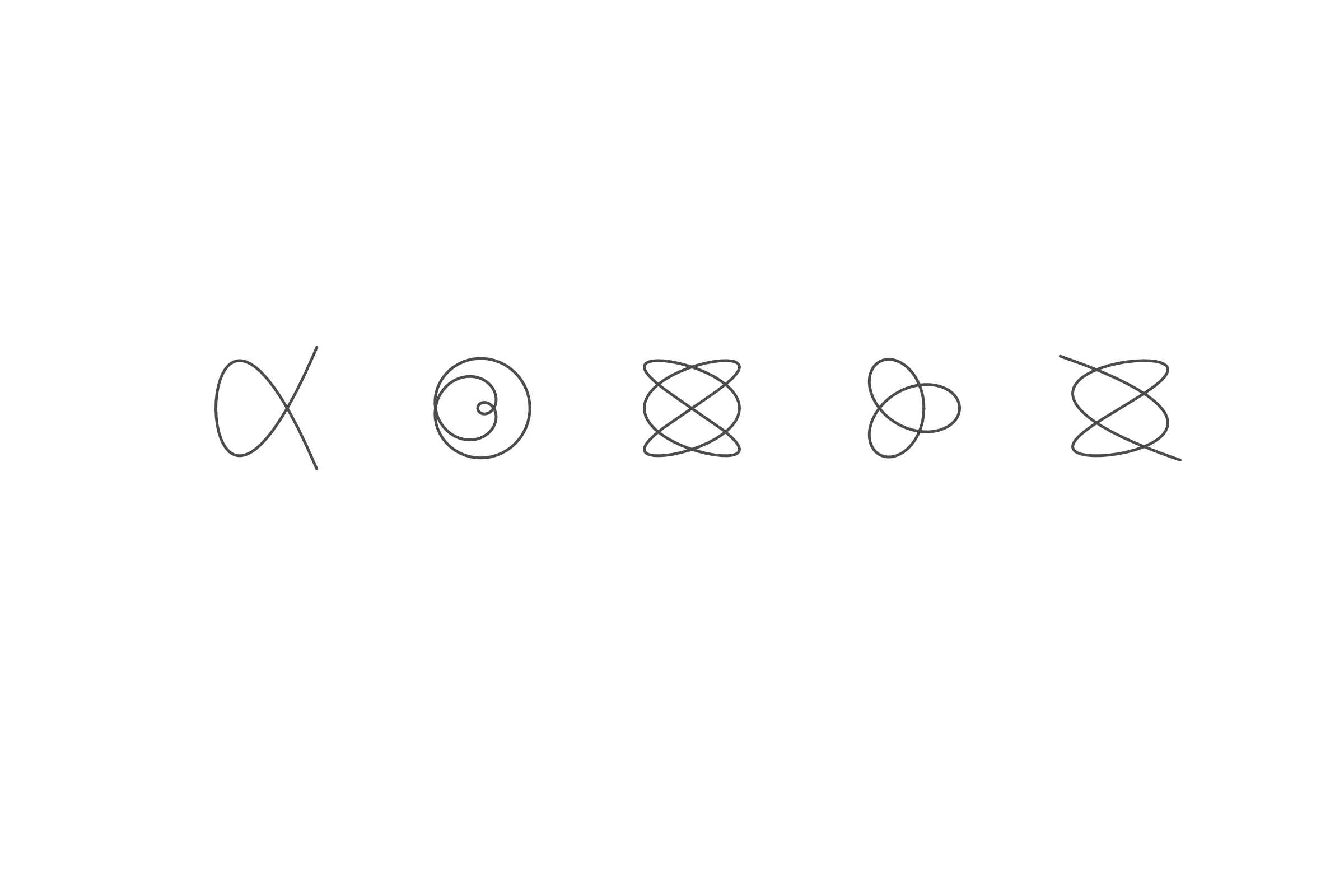} \\
\hspace{-.05in}
(a) \hspace{0.91in}
(b) \hspace{0.92in}
(c) \hspace{0.96in}
(d) \hspace{0.94in}
(e)  \\
\vspace{-.1in}
\end{center}
\caption{Irreducible expressive curves:
(a) a singular cubic;
(b) double lima\c con;
(c) $(2,3)$-Lissajous curve;
(d) $3$-petal hypotrochoid;
(e) $(3,5)$-Chebyshev curve.
}
\vspace{-.1in}
\label{fig:5-irr-expressive}
\end{figure}

\begin{figure}[ht]
\begin{center}
\vspace{-.15in}
\hspace{.2in}
\includegraphics[scale=0.237, trim=13.5cm 23cm 12cm 21cm, clip]{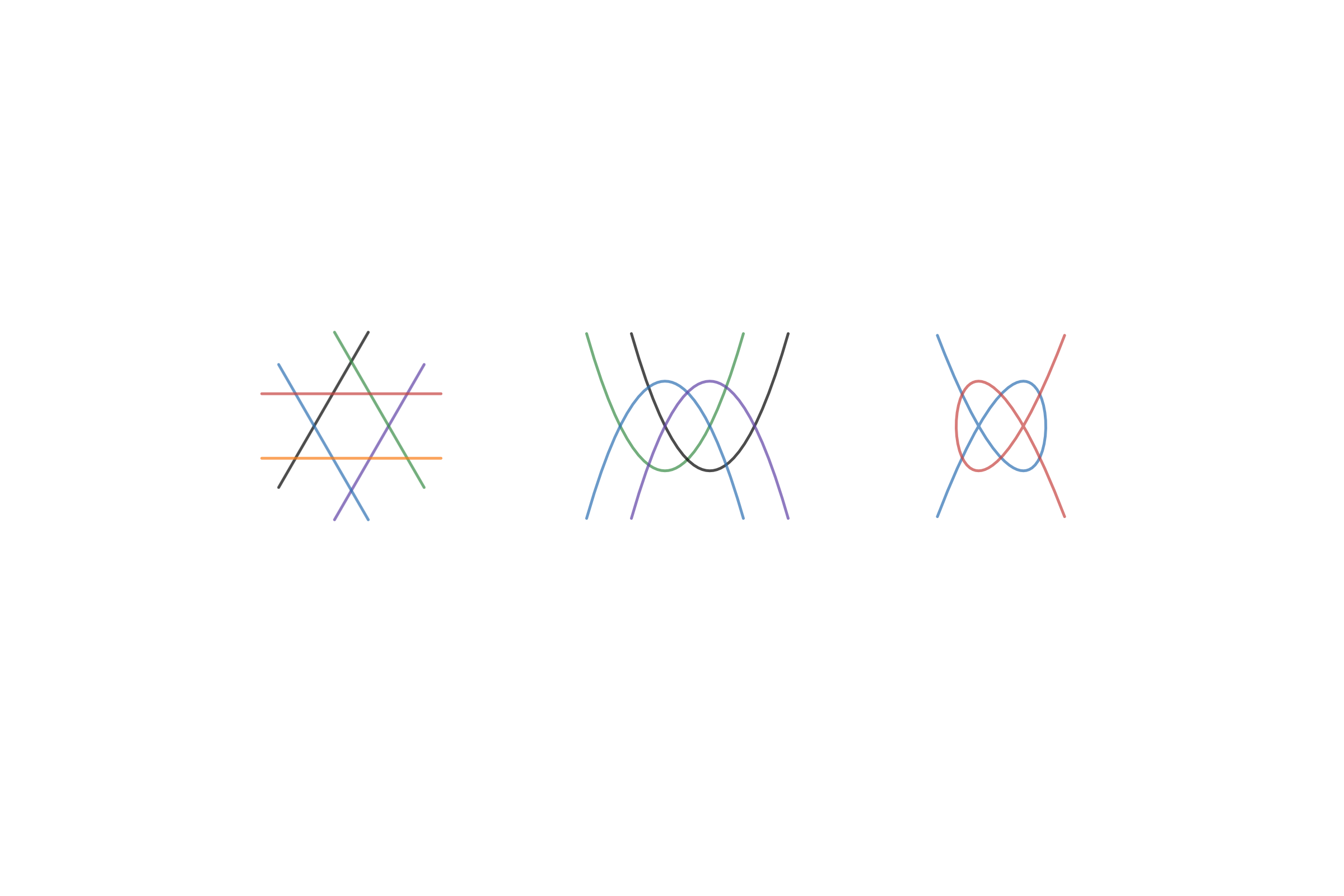} \\
\vspace{-.1in}
\end{center}
\caption{Reducible expressive curves: arrangements of
six lines, four parabolas, and two singular cubics.
}
\vspace{-.15in}
\label{fig:3-arrangements}
\end{figure}

The above characterization allows us to construct numerous examples of expressive plane curves,
including arrangements of lines, parabolas, circles, and singular cubics;
Chebyshev and Lissajous curves;
hypotrochoids and epitrochoids; and much more.
See Figures~\ref{fig:5-irr-expressive}--\ref{fig:3-arrangements} for
an assortment of examples;
many more are scattered throughout the paper.

On the face of it, expressivity is an analytic property
of a function~\redsf{$G:\RR^2\to\RR$}.
This is however an illusion: just like in the univariate~case,
in order to 
rule out \redsf{accidental} critical points,
we need~$G$ to be a polynomial of a certain kind.
Thus, expressivity is essentially an algebraic phenomenon.
Accordingly, its study requires tools of algebraic geometry and singularity theory.

For a real plane algebraic curve~$C$ to be expressive, one needs
\begin{align}
\notag
&\#\{\text{critical points of~$C$ in the complex affine plane}\} \\[-.05in]
\tag{$*$}
=\,&\#\{\text{double points in~$\CR$}\} +
\#\{\text{bounded components of $\RR^2\setminus\CR$}\} .
\end{align}
Since a generic plane curve of degree~$d$
has $(d-1)^2$ critical points, whereas expression~($*$) 
is typically smaller than~$(d-1)^2$,
we need all the remaining critical points to escape to infinity.
Our~analysis shows that this can only happen if each (real) irreducible component of~$C$
either has a unique point at infinity, or a pair of complex conjugate points;
moreover the components must intersect each other in the affine plane at real points,
specifically at hyperbolic nodes.
The requirement of having one or two points at infinity translates
into the condition of having a polynomial or trigonometric parametrization,
yielding the expressivity criterion formulated above.


As mentioned earlier, these results are established under the assumption of
\emph{$L^\infty$-regularity}, which concerns the behaviour of
the projective closure of~$C$ at the line at infinity~$L^\infty$.
This assumption ensures that the number of critical points accumulated at each point~$p\in C\cap L^\infty$
is determined in the expected way by the topology of~$C$ in the vicinity of~$p$
together with the intersection multiplicity $(C\cdot L^\infty)_p$.
All~polynomial and trigonometric curves are $L^\infty$-regular,
as are all expressive curves of degrees~$\le4$.

\vspace{.05in}

\noindent
\textbf{Section-by-section overview.}
Sections~\ref{sec:plane-curves}--\ref{sec:poly-trig} are devoted to algebraic geometry groundwork.
Section~\ref{sec:plane-curves} reviews basic background on plane algebraic curves, intersection numbers, and topological invariants of isolated singularities.
The number of critical points escaping to infinity is determined by the
intersection multiplicities of polar curves at infinity,
which are studied in Section~\ref{sec:polar-curves-at-infinity}.
Its main result is Proposition~\ref{prop:milnor-infinity},
which gives a lower bound for such a multiplicity in terms of the Milnor number
and the order of tangency between the curve and the line at infinity.
When this bound becomes an equality, a plane curve~$C$ is called $L^\infty$-regular.

In Section~\ref{sec:L-regular-curves}, we provide several criteria for $L^\infty$-regularity.
We also show (see Proposition~\ref{pr:hironaka-milnor})
that for an $L^\infty$-regular curve~$C\!=\!\{G(x,y)\!=\!0\}$
all of whose singular points in the affine plane are ordinary nodes,
the number of critical points of~$G$ is completely determined by
the number of those nodes,
the geometric genus of~$C$, and the number of local branches of~$C$ at infinity.
This statement relies on classical formulas due to H.~Hironaka~\cite{Hir} and J.~Milnor~\cite{Milnor}.

\redsf{The technical material of Sections~~\ref{sec:polar-curves-at-infinity}--\ref{sec:L-regular-curves}
can be safely skipped by the readers who are willing to treat the notion of $L^\infty$-regularity as a ``black-box'' genericity
condition that automatically holds for most, if not all, expressive curves that arise in applications.}

Section~\ref{sec:poly-trig} introduces polynomial and trigonometric curves,
the plane curves possessing a parametrization $t\mapsto (X(t),Y(t))$
in which both $X$ and $Y$ are polynomials, resp.\ trigonometric polynomials.
We review a number of examples of such curves, recall the classical result~\cite{Abhyankar-1988} characterizing polynomial curves as those
with a single place at infinity, and provide an analogous characterization for trigonometric curves.

Expressive curves are introduced in Section~\ref{sec:expressive}.
We formulate their basic properties and discuss a large number of examples,
which include an inventory of expressive curves of degrees~$\le 4$.

In Section~\ref{sec:divides}, we introduce 
\emph{divides}
and relate them to the notion of expressivity.

Section~\ref{sec:regular+expressive} contains our main results.
Using the aforementioned bounds and criteria, \linebreak[3]
we show (see Theorem~\ref{th:reg-expressive-irreducible})
that an irreducible real plane algebraic curve
is expressive and $L^\infty$-regular if and only if
it is either trigonometric or polynomial,
and moreover all its singular points in the complex affine plane are (real) hyperbolic nodes. \pagebreak[3]
This criterion is then extended (see Theorem~\ref{th:reg-expressive})
to general plane curves with real irreducible components.
Additional expressivity criteria are given in Section~\ref{sec:more-expressivity-criteria}.

\pagebreak[3]

As a byproduct, we obtain the following elementary statement
(see Corollary~\ref{cor:all-crit-pts-are-real}):
if $C=\{G(x,y)=0\}$ is a real polynomial or trigonometric affine plane curve
that intersects itself solely at hyperbolic nodes, then all critical points of~$G$ are real.

Multiple explicit constructions of expressive curves
are presented in \hbox{Sections~\ref{sec:constructions-irr}--\ref{sec:arrangements-trig}},
demonstrating the richness and wide applicability of the theory.
In Section~\ref{sec:constructions-irr}, we describe the procedures of
bending, doubling, and unfolding. Each of them can be
used to create new (more ``complicated'') expressive curves from existing ones.
Arrangements of lines, parabolas, and circles,
discussed in Section~\ref{sec:overlays}, provide another set of examples. \linebreak[3]
These examples are generalized in Section~\ref{sec:shifts+dilations}
to arrangements consisting of shifts, dilations and/or rotations of a given expressive curve.
Explicit versions of these constructions for polynomial (resp., trigonometric) curves
are presented in Section~\ref{sec:arrangements-poly}
(resp., Section~\ref{sec:arrangements-trig}).


In Section~\ref{sec:alternative-expressivity}, we briefly discuss alternative notions of expressivity:
a ``topological'' notion that treats real algebraic curves set-theoretically,
and an ``analytic'' notion that does not require the defining equation of a curve to be algebraic.

The class of divides which can arise from $L^\infty$-regular expressive curves
is studied in Section~\ref{sec:regular-expressive-divides}.
In particular, we show that a simple pseudoline arrangement belongs to this class
if and only if it is \emph{stretchable}.
In Section~\ref{sec:expressive-vs-algebraic}, we  compare this class to
the class of \emph{algebraic} divides studied in \cite{FPST}. 

\vspace{.05in}

\noindent
\textbf{Motivations and outlook.}
This work grew out of the desire to develop a global version
of the A'Campo--Guse\u\i n-Zade theory \cite{acampo-1975, acampo-1999, GZ1974, gusein-zade-2}
of morsifications of isolated singularities of plane curves.
The defining feature of such morsifications is a \emph{local} expressivity property,
which prescribes the locations (up to real isotopy) of the critical points of a morsified curve
in the vicinity of the original singularity.
In~this paper, expressivity is a \emph{global} property of a real plane algebraic curve,
prescribing  the locations of its critical points (again, up to real isotopy) on the entire affine plane.

In a forthcoming follow-up to this paper, we intend to develop a global analogue---in the setting of expressive curves---of A'Campo's theory of divides and their links.
As~shown in~\cite{FPST}, this theory has intimate connections to the combinatorics of
quivers, cluster mutations, and plabic graphs.

It would be interesting to explore the phenomenon of expressivity in higher dimensions,
and in particular find out which results of this paper generalize.

The concept of an expressive curve/hypersurface can be viewed as a generalization
of the notion of a line/hyperplane arrangement.
(Expressivity of such arrangements in arbitrary dimension
can be established by a log-concavity argument.)
This opens the possibility of extending the classical theory of hyperplane arrangements~\cite{Aguiar-Mahajan, Dimca, Stanley-arrangements}
to arrangements of expressive curves/surfaces.

\vspace{.05in}

\noindent
\textbf{Acknowledgments.}
We thank Michael Shapiro for stimulating \hbox{discussions}, and
Pavlo Pylyavskyy and Dylan Thurston
for the collaboration~\cite{FPST} which prompt\-ed our work on this project.
\redsf{We are grateful to the referee for multiple suggestions whose implementation improved the quality of the paper.}
We used \texttt{Sage} to compute resultants, and \texttt{Desmos} to draw curves.
While cataloguing expressive curves of degrees $\le 4$,
we made use of the classifications produced by A.~Korchagin and D.~Weinberg~\cite{Korchagin-Weinberg-2005}.

\pagebreak[3]

\section{Plane curves and their singularities}
\label{sec:plane-curves}

\begin{definition} 
Let $\PP^2$ denote the complex projective plane.
We fix homogeneous coordinates $x,y,z$ in~$\PP^2$.
Any homogeneous polynomial $F\in\CC[x,y,z]$ defines
a \emph{plane algebraic curve}~$C=Z(F)$ in $\PP^2$ given by
\[
C=Z(F)=\{F(x,y,z)=0\}.
\]
We understand the notion of a curve (and the notation~$Z(F)$) scheme-theoretically:
if the polynomial $F$ splits into factors, we count each component of the curve $C=Z(F)$ with the multiplicity of the corresponding factor.
\end{definition}

For two distinct points $p,q\in\PP^2$, we denote by $L_{pq}$
the line passing through $p$~and~$q$.
The \emph{line at infinity} $L^\infty\subset\PP^2$ is defined by $L^\infty=Z(z)$.

For $F$ a smooth function in $x,y,z$, we use the shorthand
\[
F_x=\tfrac{\partial F}{\partial x}\,\quad
F_y=\tfrac{\partial F}{\partial y}\,,\quad
F_z=\tfrac{\partial F}{\partial z}
\]
for the partial derivatives of~$F$.
The following elementary statement is well known, and easy to check.

\begin{lemma}[Euler's formula]
Let $F=F(x,y,z)$ be a homogeneous polynomial of degree~$d$.
Then
\begin{equation}
\label{eq:Euler}
d\cdot F=xF_x+yF_y+zF_z\,.
\end{equation}
\end{lemma}

\begin{definition} 
Let $F\in\CC[x,y,z]$ be a homogeneous polynomial in $x,y,z$.
For a point $q=(q_x,q_y,q_z)\in\PP^2$, we denote
\[
F_{(q)}=q_xF_x+q_yF_y+q_zF_z\,.
\]
The \emph{polar curve} $C_{(q)}$ associated with a curve $C=Z(F)$ and a point $q\in\PP^2$
is defined by 
$C_{(q)}=Z(F_{(q)})$.
In particular, for $q=(1,0,0)\in L^\infty$ (resp., $q=(0,1,0)\in L^\infty$), we get the polar curve $Z(F_x)$ (resp., $Z(F_y)$).
\end{definition}

For a point $p$ lying on two plane curves $C$ 
and $\tilde C$, 
we denote by 
$(C\cdot \tilde C)_p$
the \emph{intersection number} of these curves at~$p$.
We will also use this notation for analytic curves,
i.e., curves defined by analytic equations in a neighborhood of~$p$.

\begin{definition} 
\label{def:top-invariants}
Let $C=Z(F)$ be a plane algebraic curve,
and $p$ an isolated singular point of~$C$.
Let us recall the following topological invariants of the singularity $(C,p)$:
\begin{itemize}[leftmargin=.2in]
\item
the \emph{multiplicity} $\mt(C,p)=(C\cdot L)_p$, where $L$ is any line passing through~$p$ which is not tangent to the germ~$(C,p)$;

\item
the \emph{$\varkappa$-invariant}
$\varkappa(C,p)=(C\cdot C_{(q)})_p$, where $q\in \PP^2\setminus\{p\}$
    is such that the line~$L_{pq}$ 
    is not tangent to $(C,p)$;

\item
the \emph{number $\br(C,p)$ of local branches} (irreducible components) of the germ~$(C,p)$;

\item
the \emph{$\delta$-invariant} $\delta(C,p)$,
which can be determined from
\[
\varkappa(C,p)=2\delta(C,p)+\mt(C,p)-\br(C,p);
\]

\item
the \emph{Milnor number} $\mu(C,p)=(C_{(q')}\cdot C_{(q'')})_p$, where the points $q', q''\in\PP^2$ are chosen so that $p,q',q''$ are not collinear.
\end{itemize}
More generally, for any point $p\in C_{(q')}\cap C_{(q'')}$, not necessarily lying on the curve~$C$, we can define the Milnor number
\[
\mu(C,p)=(C_{(q')}\cdot C_{(q'')})_p
\]
(provided $p,q',q''$ are not collinear).
Note that for $p\notin L^\infty$, we can simply define
\begin{equation}
\label{eq:mu-x-y}
\mu(C,p)=(Z(F_x)\cdot Z(F_y))_p\,.
\end{equation}

See \cite[\S5 and~\S10]{Milnor} and \cite[Sections I.3.2 and~I.3.4]{GLS}
for additional details as well as basic properties of these invariants.
See also Remark~\ref{rem:invariants-informal} and Proposition~\ref{pr:invariants-identities}~below.
\end{definition}

\begin{remark}
\label{rem:invariants-informal}
All invariants listed in Definition~\ref{def:top-invariants} depend only on the topological type of the singularity at hand.
The Milnor number $\mu(C,p)$ measures the complexity of the singular point~$p$
viewed as a critical point of~$F$.
It is equal to the maximal number of critical points
that a small deformation of~$F$ may have in the vicinity~of~$p$.
The $\delta$-invariant is the maximal number of critical points lying on the deformed curve in a small deformation of the germ~$(C,p)$.
The $\varkappa$-invariant is the number of ramification points of a generic
projection onto a line of a generic deformation of the germ $(C,p)$.
\end{remark}

\begin{proposition}[{\cite[Propositions~3.35, 3.37, 3.38]{GLS}}]
\label{pr:invariants-identities}
Let $(C,p)$ be an isolated plane curve singularity as above. Then we have:
\begin{align}
\label{eq:milnor-formula}
\mu(C,p)&=2\delta(C,p)-\br(C,p)+1 \quad \text{(Milnor's formula)}; \\
\label{eq:C*Cq}
(C\cdot C_{(q)})_p&=\varkappa(C,p)+(C\cdot L_{pq})_p-\mt(C,p)
   \quad \text{for any $q\in\PP^2\setminus\{p\}$;}\\
\label{eq:kappa=mu+mult-1}
\varkappa(C,p)&=\mu(C,p)+\mt(C,p)-1.
\end{align}
\end{proposition}

\begin{example}
\label{ex:(x^2+z^2)(yz^2-x^3+x^2y)}
Consider the quintic curve $C=Z(F)$ defined by the polynomial
\[
F(x,y,z)=(x^2+z^2)(yx^2+yz^2-x^3)=(x+iz)(x-iz)(yx^2+yz^2-x^3).
\]
It has two points on the line at infinity~$L^\infty$,
namely $p_1=(0,1,0)$ and $p_2=(1,1,0)$.
At~$p_1$, the cubical component has an elliptic node,
and the two line components are the two tangents to the cubic at~$p_1$.
At~$p_2$, we have a smooth real local branch of the cubical component.
Direct computations show that
\begin{alignat*}{3}
\mt(C,p_1)&=4 &\qquad \mt(C,p_2)&=1\\
\varkappa(C,p_1)&=16&\qquad \varkappa(C,p_2)&=0 \\
\br(C,p_1)&=4 &\qquad \br(C,p_2)&=1 \\
\delta(C,p_1)&=8&\qquad \delta(C,p_2)&=0 \\
\mu(C,p_1)&=13&\qquad \mu(C,p_2)&=0 \\
(C\cdot L^\infty)_{p_1}&=4 &\qquad (C\cdot L^\infty)_{p_2}&=1\\
(Z(F_x)\cdot Z(F_y))_{p_1}&=16 &\qquad (Z(F_x)\cdot Z(F_y))_{p_2}&=0 \\
(C\cdot F_x)_{p_1}&=16 &\qquad (C\cdot F_x)_{p_2}&=0
\end{alignat*}
Note that \eqref{eq:mu-x-y} does not hold for $p=p_1$;
this is not a contradiction since $p_1\in L^\infty$.
\end{example}

\newpage

\section{Intersections of polar curves at infinity}
\label{sec:polar-curves-at-infinity}

In this section, we study the properties of intersection numbers of polar curves 
at their common points located at the line at infinity.

\begin{lemma}
\label{d1}
Let $F(x,y,z)\in\CC[x,y,z]$ be a non-constant homogeneous polynomial.
\begin{itemize}[leftmargin=.3in]
\item[{\rm (i)}]
The set $Z(F_{(q')})\cap Z(F_{(q'')})\cap L^\infty$
does not depend on the choice of
a pair of distinct points $q',q''\in L^\infty$.
Moreover this set is contained in~$C$.
\item[{\rm (ii)}]
For a point $p\in C\cap L^\infty$, the intersection multiplicity $(Z(F_{(q')})\cdot Z(F_{(q'')}))_p$ does not depend on the choice of a pair of distinct points $q',q''\in L^\infty$.
\end{itemize}
\end{lemma}

\begin{proof}
Any other pair $\hat q',\hat q''$ of distinct points in $L^\infty$ satisfies
\begin{equation}
\label{eq:hat-q}
\begin{array}{r}
\hat q'=a_{11}q'+a_{12}q'',\\
\hat q''=a_{21}q'+a_{22}q'',
\end{array}
\text{\quad with\quad}
\left| \begin{matrix} a_{11} & a_{12}\\
a_{21} & a_{22}
\end{matrix}\right|
\ne 0.
\end{equation}
Consequently
\begin{align*}
F_{(\hat q')}&=a_{11}F_{(q')}+a_{12}F_{(q'')},\\
F_{(\hat q'')}&=a_{21}F_{(q')}+a_{22}F_{(q'')}
\end{align*}
and the first claim in~(i) follows.
To establish the second claim, set $q'=(1,0,0)$ and $q''=(0,1,0)$
(i.e., take the polar curves $Z(F_x)$ and $Z(F_y)$),
and note that by Euler's formula~\eqref{eq:Euler},
$F$ vanishes as long as $F_z$, $F_y$ and $z$ vanish.

To prove (ii), factor the  nonsingular $2\times2$ matrix in \eqref{eq:hat-q}
into the product of an upper triangular and a lower triangular matrix;
then use that, for $bc\ne0$,
\begin{equation}
\label{e4e}
(Z(aG_1+bG_2)\cdot Z(cG_1))_p=(Z(G_2)\cdot Z(G_1))_p\,. \qedhere
\end{equation}
\end{proof}

\begin{definition}
\redsf{Let $C$ be a plane projective curve $C$ that does not contain the line at infinity~$L^\infty$ as a component.}
For a point $p\in L^\infty$, we denote
\[
\muinf \eqdef (Z(F_x)\cdot Z(F_y))_p\,.
\]
\redsf{Note that by Lemma~\ref {d1}(i)}, if $p$ lies on both $Z(F_x)$ and~$Z(F_y)$, then it necessarily
lies on~$C$, so $\muinf$ can only be nonzero at points $p\in C\cap L^\infty$.
\end{definition}

\begin{remark}
\label{r5}
For $p\in C\cap L^\infty$, the number $\muinf$ may differ from the Milnor number~$\mu(C,p)$ (cf.~\eqref{eq:mu-x-y}), since the points $p,q',q''$ lie on the same line~$L^\infty$.
Moreover, $\muinf$ is not determined by the topological type of the singularity~$(C,p)$,
as it also depends on its ``relative position'' with respect to the line~$L^\infty$.
The following example illustrates this phenomenon.
Consider the curves $Z(x^2y-z^3-xz^2)$ and $Z(x^2y^2-yz^3)$.
Each of them has an ordinary cusp (type~$A_2$) at the point $p=(0,1,0)$.
On the other hand, we have $\muinf=4$ in the former case versus
$\muinf=3$ in the latter.
Additional examples can be produced using Proposition~\ref{prop:milnor-infinity} below.
\end{remark}

\begin{remark}
\label{rem:Milnor-number-at-infinity}
The intersection number $\muinf$ is also different from ``the Milnor number at infinity''
(as defined, for instance, in~\cite{ABLMH, Siersma-Tibar})
since $\muinf$ depends on the choice of a point $p\in L^\infty$.
Moreover, $\muinf$ is not determined by the local topology of the configuration consisting of
the germ $(C,p)$ and the line~$L^\infty$.
\redsf{To see this, consider Example~\ref{e3.8} and Example~\ref{example:expressive-non-reg} with $p=p_1$.
In both cases, $(C,p)$ is an ordinary cusp transversal to~$L^\infty$.
In Example~\ref{e3.8}, we have $\mu(C,p,L^\infty)=4$, whereas
in Example~\ref{example:expressive-non-reg}, we get
$\mu(C,p,L^\infty)=\mu(C,p)+(C\cdot L^\infty)_p-1=3$
by Proposition~\ref{pr:nonmax-tangency})}.
\end{remark}

\pagebreak[3]

\begin{proposition}
\label{prop:milnor-infinity}
Let $C=Z(F)$ be an algebraic curve in~$\PP^2$. 
Let~$p
\in C\cap L^\infty$.
Then
\begin{equation}
\label{eq:e2}
\muinf
  \ge \mu(C,p)+(C\cdot L^\infty)_p-1.
\end{equation}
\end{proposition}

The proof of Proposition~\ref{prop:milnor-infinity} will rely on two lemmas,
one of them very~simple.

\begin{lemma}
\label{lem:F*Fq}
For any $q\in L^\infty\setminus\{p\}$, we have
\begin{equation*}
(C\cdot Z(F_{(q)}))_p=\mu(C,p)+(C\cdot L^\infty)_p-1.
\end{equation*}
\end{lemma}

\begin{proof}
Using \eqref{eq:C*Cq} and \eqref{eq:kappa=mu+mult-1}, we obtain:
\begin{align*}
(C\cdot Z(F_{(q)}))_p&=\varkappa(C,p)+(C\cdot L^\infty)_p-\mt(C,p) 
=\mu(C,p)+(C\cdot L^\infty)_p-1. 
\qedhere
\end{align*}
\end{proof}

\begin{lemma}
\label{lem:e1}
Let $Q$ be a local branch
(i.e., a reduced, irreducible component) of the germ of the curve $Z(F_y)$ at $p$.
Then
\begin{equation}
\label{eq:e1}
(Z(zF_z)\cdot Q)_p\ge(Z(F)\cdot Q)_p\ .
\end{equation}
\end{lemma}

 \begin{proof}
We will prove the inequality \eqref{eq:e1} inductively, by blowing up the point $p$. As a preparation step, we will apply a coordinate change intended to reduce the general case to a particular one, in which the blowing up procedure is easier to describe.

Without loss of generality, we assume that $p=(1,0,0)$. In a neighborhood of $p$, we can set $x=1$ and then work in the affine coordinates $y,z$. Abusing notation, for a homogeneous polynomial $G(x,y,z)$, we will write $G(y,z)$ instead of $G(1,y,z)$.

For any curve $Z(G)$,
the intersection multiplicity $(G\cdot Q)_p$ can be computed as follows.
Write $Q=\{f(y,z)=0\}$, where $f(y,z)$ is an irreducible element of the ring $\CC\{y,z\}$ of germs at $p$ of holomorphic functions in the variables $y$ and~$z$.
By \cite[Proposition~I.3.4]{GLS}, we have
\[
f(y,z)=u(y,z)\prod_{1\le i\le k}(y-\xi_i(z^{1/k})), \quad u(y,z)\in\CC\{y,z\},\quad u(0,0)\ne0,
\]
where each $\xi_i$ is a germ at zero of a holomorphic function vanishing at the origin.
Then \cite[Propostion I.3.10 (Halphen's formula)]{GLS} yields
\begin{equation}
(G\cdot Q)_p=\sum_{1\le i\le k}\ord_0G(\xi_i(z),z).
\label{e17}\end{equation}
It follows that the variable change
$(y,z)=\tau(y_1,z_1)\stackrel{\rm def}{=}(y_1,z_1^k)$ multiplies both sides of \eqref{eq:e1} by~$k$.
For $g\in\CC\{y,z\}$, let us denote $\tau^*g(y_1,z_1)=g\circ\tau(y_1,z_1))$.
Then $\tau^*Q$ splits into $k$ smooth branches
\[
Q_i=\{y_1-\xi_i(z_1)=0\}, \quad i=1,...,k,
\]
and it suffices to prove the inequality \eqref{eq:e1} with $Q$ replaced by each of the~$Q_i$'s.
Moreover, with respect to $Q_i$, the desired inequality is of the same type. 
Namely, $\tau^*(F_y)=(\tau^*F)_{y_1}$, and hence $Q_i$ is a local branch of the polar curve $Z((\tau^*F)_{y_1})$ of the curve $\tau^*C=Z(\tau^*F)$. Furthermore,
\begin{align*}
z_1(\tau^*F)_{z_1}(y_1,z_1)&=z_1\tfrac{\partial}{\partial z_1}[F(y_1,z_1^k)]\\
&=z_1F_y(y_1,z_1^k)+kz_1^kF_z(y_1,z_1^k)\\
&=z_1(\tau^*F_y)(y_1,z_1)+k\cdot(\tau^*(zF_z))(y_1,z_1),\end{align*}
which implies
\[
(Z(z_1(\tau^*F)_{z_1})\cdot Q_i)_p=(Z(\tau_*(zF_z))\cdot Q_i)_p,\quad i=1,...,k.
\]

We have thus reduced the proof of \eqref{eq:e1} to the case where $Q$ is a smooth curve germ transversal to the line~$L^\infty$.
To simplify notation, we henceforth write
$C,F,y,z$ instead of $\varphi^*C,\varphi^*F,y_1,z_1$, respectively.

We proceed by induction on $\mu(C,p)$. If $\mu(C,p)=0$, then $(C,p)$ is a smooth germ.
If $C$ intersects $L^\infty$ transversally, then $(C,p)$ is given by
\[
F(y,z)=ay+bz+\text{h.o.t.}, \quad a\ne0
\]
(hereinafter h.o.t.\ is a shorthand for ``higher order terms''),
implying \hbox{$F_y(0,0)=a\ne0$}.
Thus the polar curve $Z(F_y)$ does not pass through~$p$; consequently both sides of \eqref{eq:e1} vanish. If $C$ is tangent to $L^\infty$ at $p$, then it is transversal to $Q$ at $p$, so we have
\[
(C\cdot Q)_p=1=(Z(z)\cdot Q)_p\le(Z(zF_z)\cdot Q)_p\,.
\]

Suppose that $\mu(C,p)>0$. Then $m\stackrel{\rm def}{=}\mt(C,p)\ge2$.
If $C$ and $Q$ intersect transversally at $p$ (i.e., have no tangent in common), then
\[
(C\cdot Q)_p=\mt(C,p)\cdot\mt(Q,p)=m\cdot 1=m.
\]
On the other hand, $\mt(Z(F_z),p)\ge m-1$, and therefore
\[
(Z(zF_z)\cdot Q)\ge\mt(Z(zF_z),p)\cdot \mt(Q,p)\ge m\cdot 1=m=(C\cdot Q).
\]
If $C$ and $Q$ have a common tangent at $p$, we apply the blowing-up
$\pi:\widetilde\PP^2\to\PP^2$ of the plane at the point~$p$.
For a curve $D$ passing through~$p$, let $D^*$ denote its \emph{strict transform}, i.e.,  the closure of the preimage $\pi^{-1}(D\setminus\{p\})$ in~$\widetilde\PP^2$.
(For more details, see \cite[Section~I.3.3, p.~185]{GLS}.).
Since $Q$ is smooth, the strict transform $Q^*$ is smooth too, and intersects transversally the exceptional divisor $E$ at some point~$p^*$. More precisely, if $Q=Z(y-\eta z-\text{h.o.t.})$, then in the coordinates $(y_*,z_*)$ given by
$y=y_*z_*$, $z=z_*$, we have $E=Z(z_*)$ and $p^*=(\eta,0)$. Since $C$ and its polar curve $Z(F_y)$ are tangent to the line $Z(y-\eta z)$, the lowest homogeneous form of $F(y,z)$ is divisible by $(y-\eta z)^2$, while the lowest
homogeneous form of $F_y$ is divisible by $y-\eta z$ and, moreover, $\mt(Z(F_y),p)=m-1$.

We now recall some properties of the blowing-up.
For a curve $D\!=\!Z(G(y,z))$ pass\-ing
through~$p$, we have \cite[Prop.~I.3.21 and I.3.34, and computations on p.~186]{GLS}:
\begin{align*}
(D^*\cdot Q^*)_{p^*}&=(D\cdot Q)_p-\mt(D,p)\cdot\mt(Q,p)=(D\cdot Q)_p-\mt(D,p); \\
\textstyle\sum_{q\in D^*\cap E}\delta(D^*,q)&=\delta(D,p)-\tfrac12 \mt(D,p)(\mt(D,p)-1); \\
D^*&=Z(z_*^{-\mt(D,p)}G(y_*z_*,z_*)).
\end{align*}
We see that
\begin{align*}
F^*(y_*,z_*)&=z_*^{-m}F(y_*z_*,z_*),\\
(F^*)_{y_*}(y_*,z_*)&=
z_*^{1-m}F_y(y_*z_*,z_*)=(F_y)^*(y_*,z_*).
\end{align*}
Thus $Q^*$ is a local branch of the polar curve $(F^*)_{y_*}\!=\!0$ of the strict transform~$C^*$.
Hence after the blowing-up we come to the original setting. Furthermore,
\begin{align*}
\mu(C^*,p^*)&=2\delta(C^*,p^*)-\br(C^*,p^*)+1\\
&\le2\delta(C^*,p^*)\\
&\le2\textstyle\sum_{q\in C^*\cap E}\delta(C^*,q)\\
&=2\delta(C,p)-m(m-1)\\
&=\mu(C,p)+\br(C,p)-1-m(m-1)\\
&\le\mu(C,p)+m-1-m(m-1)\\
&<\mu(C,p).
\end{align*}
So we can apply the induction assumption.
Observe that $\mt(F_z,p)=m-1+r$ for some $r\ge0$.
It follows that
\begin{equation}
\begin{aligned}
(Z(F)\cdot Q)_p&=(Z(F^*)\cdot Q^*)_{p^*}+m,\\
(Z(zF_z)\cdot Q)_p&=1+(Z(F_z)\cdot Q)_p=(Z(F_z)^*)\cdot Q^*)_{p^*}+m+r.
\end{aligned}
\label{eq:e3}
\end{equation}
Now
\begin{align*}z_*(F^*)_{z_*}&=-mz_*^{-m}F(y_*z_*,z_*)+y_*z_*^{1-m}F_y(y_*z_*,z_*)+z_*^{1-m}F_z(y_*z_*,z_*)\\
&=-mF^*(y_*,z_*)+y_*(F^*)_{y_*}(y_*,z_*)+z_*^r(F_z)^*(y_*,z_*).
\end{align*}
Finally, the latter formula, the induction assumption, and (\ref{eq:e3}) imply
\begin{align*}
(Z(zF_z)\cdot Q)_p&=(Z((F_z)^*)\cdot Q^*)_{p^*}+m+r\\
&=(Z(z_*^r(F_z)^*)\cdot Q^*)_{p^*}+m\\
&\ge\min\{Z(z_*(F^*)_{z_*})\cdot Q^*)_{p^*},\, (Z(F^*)\cdot Q^*)_{p^*}\}+m\\
 &\!\!\!\stackrel{\eqref{eq:e1}}{=} (Z(F^*)\cdot Q^*)_{p^*}\}+m\\
&=(Z(F)\cdot Q)_p\,. \qedhere
\end{align*}
\end{proof}

\begin{proof}[Proof of Proposition~\ref{prop:milnor-infinity}]
We again assume $p=(1,0,0)$.
Set $d=\deg(F)$. In the local affine coordinates
$y, z$ (with $x=1$), Euler's formula~\eqref{eq:Euler} becomes
\[
d F(1,y,z)=F_x(1,y,z)+yF_y(1,y,z)+zF_z(1,y,z).
\]
Consequently $\muinf=(Z(F_x)\cdot Z(F_y))_p=(Z(d F-zF_z)\cdot Z(F_y))_p$.

Let ${\mathcal B}$ denote the set of local branches of the polar curve $F_y=0$ at~$p$.
Then
\begin{align*}
\muinf&=(Z(d F-zF_z)\cdot Z(F_y))_p\\
&=\textstyle\sum_{Q\in{\mathcal B}}
(Z(dF-zF_z)\cdot Q)_p\\
&\ge\textstyle\sum_{Q\in{\mathcal B}}\min\{(Z(F)\cdot Q)_p\,,(Z(zF_z)\cdot Q)_p\}\\
&=\textstyle\sum_{Q\in{\mathcal B}}(Z(F)\cdot Q)_p \qquad \text{(by \eqref{eq:e1})}\\
&=(Z(F)\cdot Z(F_y))_p \\
&=\mu(C,p)+(C\cdot L^\infty)_p-1 \quad \text{(by Lemma~\ref{lem:F*Fq}).}
\qedhere
\end{align*}
\end{proof}

\newpage

\section{\texorpdfstring{$L^\infty$}{L}-regular curves}
\label{sec:L-regular-curves}

\begin{definition}
\label{def:c-regular}
Let $C=Z(F(x,y,z))\subset\PP^2$ be a reduced plane algebraic curve
which does not contain the line at infinity~$L^\infty$ as a component.
The curve $C$ (or the polynomial~$F$)
is called \emph{$L^\infty$-regular} if
at each point $p\in C\cap L^\infty$, the formula~\eqref{eq:e2} becomes an equality:
\begin{equation}
\label{eq:c-regularity}
\muinf = \mu(C,p)+(C\cdot L^\infty)_p-1.
\end{equation}
\end{definition}

In the rest of this section, we provide $L^\infty$-regularity criteria
for large classes of plane curves.

\redsf{As mentioned in the Introduction, the technical material in this section can be skipped
if the reader is willing to take it on faith and to view the requirement of $L^\infty$-regularity as a genericity
condition that holds in all ``non-pathological'' examples arising in common applications.}

\begin{proposition}
\label{pr:regularity-via-milnor}
Let $C\!=\!Z(F(x,y,z))\!\subset\!\PP^2$ be a reduced algebraic curve of degree~$d$
which does not contain~$L^\infty$ as a component.
Assume that the polynomial $F(x,y,1)$ has $\xi<\infty$ critical points,
counted with multiplicities.
Then we have
\begin{equation}
\label{eq:xi-d-mu}
\xi  \le d^2-3d+1-\sum_{p\in C\cap L^\infty} (\mu(C,p)-1),
\end{equation}
with equality if and only if $C$ is $L^\infty$-regular.
\end{proposition}

\begin{proof}
In view of Proposition~\ref{prop:milnor-infinity}, we have
\begin{equation}
\label{eq:L-inequality}
\sum_{p\in C\cap L^\infty} \muinf
\ge d+\sum_{p\in C\cap L^\infty} (\mu(C,p)-1),
\end{equation}
with equality if and only if $C$ is $L^\infty$-regular.
Since $F$ has finitely many critical points,
B\'ezout's theorem for the polar curves $Z(F_x)$ and~$Z(F_y)$ applies,
yielding
\begin{equation}
\sum_{p\in C\cap L^\infty}\muinf = (d-1)^2 -\xi.
\end{equation}
The claim follows.
\end{proof}

\begin{example}
\label{ex:(x^2+z^2)(yz^2-x^3+yx^2)-regular}
As in Example~\ref{ex:(x^2+z^2)(yz^2-x^3+x^2y)},
consider the quintic curve $C=Z(F)$ defined by the polynomial
\[
F(x,y,z)=(x^2+z^2)(yx^2+yz^2-x^3)=(x+iz)(x-iz)(yx^2+yz^2-x^3).
\]
Set $G(x,y)=F(x,y,1)=(x^2+1)(yx^2+y-x^3)$. Then
\begin{align*}
G_x&= 2x(yx^2+y-x^3)+(x^2+1)(2xy-3x^2) =(x^2+1)(4xy-3x^2)-2x^4,  \\
G_y&=(x^2+1)^2,
\end{align*}
and we see that $G$ has no critical points in the complex $(x,y)$-plane; thus $\xi=0$.
Using the values of Milnor numbers computed in Example~\ref{ex:(x^2+z^2)(yz^2-x^3+x^2y)},
we obtain:
\[
d^2-3d+1-\sum_{p\in C\cap L^\infty} (\mu(C,p)-1)
=25-15+1-(13-1)-(0-1)
=0.
\]
It follows by Proposition~\ref{pr:regularity-via-milnor} that $C$ is $L^\infty$-regular.
Alternatively, one can check directly that the equality~\eqref{eq:c-regularity}
holds at $p_1$ and~$p_2$.
\end{example}

Recall that the \emph{geometric genus} of a plane curve~$C$ is defined by
\begin{equation}
\label{eq:g(C)}
g(C)=\sum_{C'\in\Comp(C)} (g(C')-1)+1,
\end{equation}
where $\Comp(C)$ is the set of irreducible components of~$C$,
and $g(C')$ denotes the genus of the normalization of a component~$C'$.

\begin{proposition}
\label{pr:hironaka-milnor}
Let $C\!=\!Z(F(x,y,z))\!\subset\!\PP^2$ be a reduced algebraic curve of degree~$d$.
Suppose that
\begin{itemize}[leftmargin=.2in]
\item
$C$ does not contain the line at infinity~$L^\infty$ as a component;
\item
all singular points of $C$ in the affine $(x,y)$-plane $\PP\setminus L^\infty$ are ordinary nodes;
\item
the polynomial $F(x,y,1)\in\CC[x,y]$ has finitely many critical points.
\end{itemize}
Let $\nu$ denote the number of nodes of~$C$ in the $(x,y)$-plane,
and let $\xi$ denote the number of critical points of the polynomial~$F(x,y,1)$,
counted with multiplicities. \linebreak[3]
Then we have
\begin{equation}
\label{eq:xi-le}
\xi  \le 2g(C)-1+2\nu+\sum_{p\in C\cap L^\infty}\br(C,p),
\end{equation}
with equality if and only if $C$ is $L^\infty$-regular.
\end{proposition}

\begin{proof}
By Hironaka's genus formula~\cite{Hir}
(cf.\ also \cite[Chapter II, (2.1.4.6)]{GLS1}), we~have
\begin{equation}
\label{eq:Hironaka}
g(C)=\tfrac{(d-1)(d-2)}{2}-\sum_{p\in\Sing(C)}\delta(C,p),
\end{equation}
where $\Sing(C)$ denotes the set of singular points of~$C$.
Combining this with Milnor's formula~\eqref{eq:milnor-formula}, we obtain:
\begin{align*}
\sum_{p\in C\cap L^\infty} (\mu(C,p)-1)
&=\sum_{p\in \Sing(C)} (\mu(C,p)-1) \\
&= \sum_{p\in \Sing(C)} (2\delta(C,p)-\br(C,p)) \\
&=(d-1)(d-2)-2g(C)-2\nu-\sum_{p\in C\cap L^\infty}\br(C,p).
\end{align*}
Therefore
\[
d^2-3d+1-\sum_{p\in C\cap L^\infty} (\mu(C,p)-1)
=2g(C)-1+2\nu+\sum_{p\in C\cap L^\infty}\br(C,p),
\]
and the claim follows from Proposition~\ref{pr:regularity-via-milnor}.
\end{proof}

Our next result (Proposition~\ref{pr:nonmax-tangency} below)
shows that equation \eqref{eq:c-regularity} holds under certain rather mild local conditions,
To state these conditions, we will need to recall some terminology and notation.

\begin{definition}[{\cite[Definitions I.2.14--I.2.15]{GLS}}]
\label{def:Newton-diagram}
We denote by $\Gamma(G)$ the \emph{Newton diagram} of a bivariate polynomial~$G$,
i.e., the union of the edges of the Newton polygon of~$G$
which are visible from the origin.
The \emph{truncation} of~$G$ along an edge~$e$ of~$\Gamma(G)$ is the sum of all monomials in~$G$ corresponding to the integer points in~$e$.

An isolated singularity of an affine plane curve $\{G(y,z)=0\}$ at the origin is called \emph{Newton nondegenerate} (with respect to the local affine coordinates~$(y,z)$)
if the Newton diagram $\Gamma(G)$ intersects each of the coordinate axes,
and the truncation of~$G$ along any edge of the Newton diagram is a quasihomogeneous polynomial without critical points in~$(\CC^*)^2$.
\end{definition}

\begin{proposition}
\label{pr:nonmax-tangency}
Let $C=Z(F(x,y,z))\subset\PP^2$ be a reduced curve not containing the line at infinity $L^\infty=Z(z)$ as a component.
Let $p\in C\cap L^\infty$.
\redsf{Without loss of generality, assume that $p=(1,0,0)$.}
Suppose that $p$ is either a smooth point of~$C$,
or a singular point of $C$ such that
\begin{align}
\label{eq:newton-nondeg}
&\!\!\text{the singularity $(C,p)$ is Newton nondegenerate, in the local coordinates~$y,z$;}\\
\label{eqeq}
&\!\!
(Z(y)\cdot C)_p<\deg C=\deg F.
\end{align}
Then
\begin{equation}
\label{eq:local-regularity}
\muinf = \mu(C,p)+(C\cdot L^\infty)_p-1.
\end{equation}
Thus, if conditions~\eqref{eq:newton-nondeg}--\eqref{eqeq} hold at every point $p\in C\cap L^\infty$,
then $C$ is $L^\infty$-regular.
\end{proposition}

\begin{remark}
Although we did not find Proposition~\ref{pr:nonmax-tangency} in the literature,
similar results---proved using similar tools---appeared before, see for example~\cite{ABLMH}.
\end{remark}

\begin{proof}
If $p$ is a smooth point of $C$ with the tangent $L\ne L^\infty$, then
\[
F(1,y,z)=ay+bz+\text{h.o.t.} \quad (a\ne0),
\]
implying that $p\not\in Z(F_y)$.
The $L^\infty$-regularity follows:
\[
(Z(F_x)\cdot Z(F_y))_p=0=\mu(C,p)+(C\cdot L^\infty)_p-1.
\]
If $C$ is smooth at $p$ with the tangent line $L^\infty$, then
\[
F(1,y,z)=ay^n+bz+\text{h.o.t.} \quad (ab\ne0,\ n>1),
\]
which implies the Newton nondegeneracy as well as condition~\eqref{eqeq}:
\[(Z(y)\cdot C)_p=1<n\le\deg C\ .\]
Thus, this situation can be viewed as a particular case
of the general setting where we have a singular point~$p$ satisfying conditions~\eqref{eq:newton-nondeg}--\eqref{eqeq}.
We next turn to the treatment of this setting.

We proceed in two steps.
We first consider semi-quasihomogeneous singular points, and then move to the general case.
We set $x=1$ in a neighborhood of $p$ and write $F(y,z)$ as a shorthand for~$F(1,y,z)$.

\smallskip

(1) Assume that $\Gamma(F)$ is a segment with endpoints $(m,0)$ and~$(0,n)$.
By the assumptions of the lemma, $m\le d=\deg F$ and $n<d$. The Newton nondegeneracy condition means that the truncation $F^{\Gamma(F)}$ of $F$ on $\Gamma(F)$ is a square-free quasihomogeneous polynomial.

Assuming that $s=\gcd\{m,n\}$, $m=m_1s$, $n=n_1s$, we can write
\[
F^{\Gamma(F)}(y,z)=\sum_{k=0}^sa_ky^{m_1k}z^{n_1(s-k)},\quad\text{where}\ a_0a_s\ne0 .
\]
Consider the family of polynomials
\[
F_t(y,z)=t^{-mn}F(yt^n,zt^m)=F^{\Gamma(F)}(y,z)+\sum_{in+jm>mn}c_{ij}t^{in+jm-mn}\redsf{y^iz^j}, \quad t\in[0,1].
\]
Note that $F_0=F^{\Gamma(F)}$ and the polynomials $F$ and $F_t$, $0<t<1$, differ by a linear change of the variables. This together with the lower semicontinuity of the intersection multiplicity implies
\begin{equation}\muinf=(Z(F_y)\cdot Z(dF-zF_z))_p\le(Z(F^{\Gamma(F)}_y)\cdot Z(dF^{\Gamma(F)}-zF^{\Gamma(F)}_z))_p\ .\label{e10}\end{equation}
Here
\begin{align}
F^{\Gamma(F)}_y&=\sum_{k=1}^sm_1ka_ky^{m_1k-1}z^{n_1(s-k)},\nonumber\\
dF^{\Gamma(F)}-zF^{\Gamma(F)}_z&=
\sum_{k=0}^s(d-n_1(s-k))a_ky^{m_1k}z^{n_1(s-k)}.
\label{e12}
\end{align}
Since $a_s\ne0$ and $n_1s=n<d$, these are nonzero polynomials, and moreover $F^{\Gamma(F)}_y$ splits into $l_1\ge0$ factors of type $z^{n_1}+\alpha y^{m_1}$, $\alpha\ne0$ and the factor $y^{m-1-l_1m_1}$, while
$zF^{\Gamma(F)}_z-dF^{\Gamma(F)}$ splits into $l_2\ge0$ factors of type $z^{n_1}+\beta y^{m_1}$, $\beta\ne0$, and the factor~$z^{n-l_2n_1}$.
Observe that the polynomials $F^{\Gamma(F)}_y$ and $zF^{\Gamma(F)}_z-F^{\Gamma(F)}$ are coprime.
(Otherwise, they would have a common factor $z^{n_1}+\gamma y^{m_1}$ with $\gamma\ne0$,
which would also be a divisor of the polynomial
\[
n_1yF^{\Gamma(F)}+m_1(zF^{\Gamma(F)}_z-dF^{\Gamma(F)})=m_1(n-d)F^{\Gamma(F)}.
\]
\redsf{Then $F^{\Gamma(F)}$ and its derivative $F_y^{\Gamma(F)}$ would have a common factor,}
contradicting the square-freeness of~$F^{\Gamma(F)}$.)
Since
\begin{align*}
(Z(y)\cdot Z(z))_p&=1, \\
(Z(y)\cdot Z(z^{n_1}+\beta y^{m_1}))_p&=n_1, \\
(Z(z)\cdot Z(z^{n_1}+\alpha y^{m_1}))_p&=m_1\quad \text{as}\quad \alpha\ne0, \\
(Z(z^{n_1}+\alpha y^{m_1})\cdot Z(z^{n_1}+\beta y^{m_1}))_p&=m_1n_1\quad\text{as}\quad\alpha\ne\beta,
\end{align*}
the right-hand side of (\ref{e10}) equals
\begin{align*}
&l_1l_2\cdot(Z(z^{n_1}+\alpha y^{m_1})\cdot Z(z^{n_1}+\beta y^{m_1}))_p+
l_1(n-l_2n_1)\cdot(Z(z)\cdot Z(z^{n_1}+\alpha y^{m_1}))_p \\
&\quad +(m-1-l_1m_1)l_2\cdot(Z(y)\cdot Z(z^{n_1}+\beta y^{m_1}))_p
+(m-1-l_1m_1)\cdot(Z(y)\cdot Z(z))_p \\
=&(m-1)n\\
=&(m-1)(n-1)+m-1\\
=&\mu(C,p)+(C\cdot L^\infty)_p-1\ ,
\end{align*}
which together with (\ref{eq:e2}) yields the desired equality.

\smallskip
(2) Suppose that $\Gamma(F)$ consists of $r\ge2$ edges $\sigma^{(1)},...,\sigma^{(r)}$ successively ordered so that $\sigma^{(1)}$ touches the axis of exponents of $y$ at the point $(m,0)$, where $m=(C\cdot L^\infty)_p$, and $\sigma^{(r)}$ touches the axis of exponents of $z$ at the point $(0,n)$, where $n=(C\cdot Z(y))_p<d=\deg F$. By the hypotheses of the lemma, for any edge $\sigma=\sigma^{(i)}$, the truncation $F^\sigma(y,z)$ is the product of $y^az^b$, $a,b\ge0$, and of a quasihomogeneous, square-free polynomial $F_0^\sigma$, whose Newton polygon $\Delta(F_0^\sigma)$ is the segment $\sigma_0$ with endpoints on the coordinate axes, obtained from $\sigma$ by translation along the vector $(-a,-b)$.

Note that the minimal exponent of $z$ in the polynomial $dF(0,z)-zF_z(0,z)$ is $n$, and hence
\begin{equation}(Z(F_y)\cdot Z(dF-zF_z))_p=(Z(yF_y)\cdot Z(dF-zF_z))_p-n\ .\label{e11}\end{equation}
Next, we note that the Newton diagram $\Gamma(yF_y)$ contains entire edges $\sigma^{(1)},...,\sigma^{(r-1)}$ and some part of the edge $\sigma^{(r)}$, while
$\Gamma(dF-zF_z)=\Gamma(F)$, since the monomials of $dF-zF_z$ and of $F$ on the Newton diagram $\Gamma(F)$ are in bijective correspondence, and the corresponding monomials differ by a nonzero constant factor, cf. (\ref{e12}).

By \cite[Proposition~I.3.4]{GLS}, we can split the polynomial $F$ inside the ring $\CC\{y,z\}$ into the product
$$F=\varphi_1...\varphi_r,\quad\Gamma(\varphi_i)=\sigma^{(i)}_0,\ \varphi_i^{\sigma^{(i)}_0}=
F_0^{\sigma^{(i)}},\ i=1,...,r,$$
and similarly
$$dF-zF_z=\psi_1...\psi_r,\quad\Gamma(\psi_i)=\sigma^{(i)}_0,\ \psi_i^{\sigma^{(i)}_0}=
(dF-zF_z)_0^{\sigma^{(i)}},\ i=1,...,r,$$
$$yF^y=\theta_1...\psi_r,\quad \Gamma(\theta_i)=\sigma^{(i)}_0,\ \theta_i^{\sigma^{(i)}_0}=
(yF_y)_0^{\sigma^{(i)}},\ i=1,...,r-1,$$ while
$$\theta_r^{\sigma_0^{(r)}}=(yF_y)^{\sigma^{(r)}}\cdot z^{-c}$$ for $c$ the minimal exponent of $z$
in $(yF_y)^{\sigma^{(r)}}$. Thus,
\begin{equation}(Z(yF_y)\cdot Z(dF-zF_z))_p=\sum_{1\le i,j\le r}(Z(\psi_i)\cdot Z(\theta_j))_p\ .\label{e14}\end{equation}
We claim that
\begin{equation}(Z(\psi_i)\cdot Z(\theta_j))_p=(Z(\varphi_i)\cdot Z(\theta_j))_p\quad \text{for all}\ 1\le i,j\le r\ .\label{e15}\end{equation}
Having this claim proven, we derive from (\ref{e14}) that
\begin{align*}(Z(yF_y)\cdot Z(dF-zF_z))_p&=\sum_{1\le i,j\le r}(Z(\varphi_i)\cdot Z(\theta_j))_p\\
&=(Z(yF_y)\cdot Z(F))_p=(Z(y)\cdot Z(F))_p+(Z(F_y)\cdot Z(F))_p\\
&=n+\varkappa(C,p)+(C\cdot L^\infty)_p-\mt(C,p)\\
&=n+(\mu(C,p)+\mt(C,p)-1)+(C\cdot L^\infty)_p-\mt(C,p)\\
&=n+\mu(C,p)+(C\cdot L^\infty)_p-1\ ,\end{align*}
which completes the proof in view of (\ref{e11}).

The equality (\ref{e15}) follows from the fact that both sides of the relation depend only on the geometry of the segments $\sigma_0^{(i)}$ and $\sigma_0^{(j)}$.
Suppose that $1\le i<j\le r$.
We have  $\sigma^{(i)}=[(m',0),(0,n')]$, $\sigma^{(j)}=[(m'',0),(0,n'')]$ with $\frac{m'}{n'}>\frac{m''}{n''}$. By  the Newton-Puiseux algorithm \cite[pp.~165--170]{GLS}, the function $\varphi_i(y,z)$ (or $\psi_i(y,z)$) splits into the product of $u(y,z)\in\CC\{y,z\}$, $u(0,0)\ne0$, and $n'$ factors of the form
$$z-\alpha y^{m'/n'}+\text{h.o.t.},\quad \alpha\ne0,$$ and then by (\ref{e17}) we get
$$(Z(\varphi_i)\cdot Z(\theta_j))_p=(Z(\psi_i)\cdot Z(\theta_i))_p=m''n'.$$
This holds even for $j=r$ since the only monomial of $\theta_r(y,z)$ that comes into play is $\beta y^{m''}$, $\beta\ne0$.
The case $1\le j<i\le r$ is settled analogously.
Now suppose \hbox{$1\le i=j\le r$}. As we observed in the first part of the proof, the pairs of quasihomogeneous polynomials $\varphi_i^{\sigma^{(i)}_0}$ and
$\theta_i^{\sigma^{(i)}_0}$, $\psi_i^{\sigma^{(i)}_o}$ and $\theta_i^{\sigma^{(i)}_0}$ are coprime (even for $i=r$!). Then the above computation yields
$$(Z(\varphi_i)\cdot Z(\theta_i))_p=(Z(\psi_i)\cdot Z(\theta_i))_p=m'n'\ ,$$ where
$\sigma^{(i)}_0=[(m',0),(0,n')]$.
\end{proof}

The following example shows that condition \eqref{eqeq} in Proposition~\ref{pr:nonmax-tangency} cannot be removed.

\begin{example}\label{e3.8}
Consider the cubic $C=Z(xy^2-z^3-yz^2)$ (cf.\ Remark~\ref{r5}).
It has a Newton nondegenerate singular point $p=(1,0,0)\in L^\infty$,
an ordinary cusp with the tangent line $L=Z(y)\ne L^\infty$.
Thus $(C\cdot L)=3=\deg(C)$, so~\eqref{eqeq} fails.
Also,
\begin{align*}
\muinf=(Z(y^2) \cdot Z(2xy-z^2))_p&=4, \\
\mu(C,p)+(C\cdot L^\infty)_p-1=2+2-1&=3.
\end{align*}
so \eqref{eq:local-regularity} fails as well.
\end{example}

For another (more complicated) instance of this phenomenon, see Example~\ref{example:expressive-non-reg}.

\begin{corollary}
Let $C=Z(F(x,y,z))\subset\PP^2$ be a reduced curve not containing
the line at infinity $L^\infty=Z(z)$ as a component.
Let $G(x,y)=F(x,y,1)$.
Assume that the Newton polygon $\Delta(G)$ intersects each of the coordinate axes in points different from the origin,
and the truncation of~$G$ along any edge of the boundary~$\partial\Delta(G)$ not visible from the origin 
is a square-free polynomial, except possibly for factors of the form $x^i$ or~$y^j$.
Then $C$ is $L^\infty$-regular.
\end{corollary}

\begin{proof}
The intersection $C\cap L^\infty$ is determined by the top degree form of $G(x,y)$,
which has the form $x^iy^jf(x,y)$, with $f(x,y)$ a square-free homogeneous polynomial. Hence $C\cap L^\infty$ consists of the points $(0,1,0)$ (if $i>0$), $(1,0,0)$ (if $j>0$),
and $\deg(f)$ other points, at which $C$ is smooth and $L^\infty$-regular
(see Proposition~\ref{pr:nonmax-tangency}).
The Newton diagram of each of the points $(0,1,0)$ and $(1,0,0)$ consists of some edges of $\partial\Delta(G)$ mentioned in the lemma,
and therefore these points (if they lie on~$C$) satisfy the requirements of Proposition \ref{pr:nonmax-tangency}.
\end{proof}

\begin{remark}
\label{rem:Linf-reducible}
A reducible plane curve~$C=Z(F)$ with $L^\infty$-regular components~does not have to be $L^\infty$-regular.
For example, take $F=(xy-1)(xy-2)$: each factor is $L^\infty$-regular, but the
product is not, since $F_x$ and $F_y$ have a common divisor $2xy-3$.

Conversely, a curve may be $L^\infty$-regular even when one of its components is not.
\redsf{For example, let  $F=x^2y+x^3+x^2z+xz^2+z^3$.
Then $Z((x-y)F)$ is $L^\infty$-regular by Proposition~\ref{pr:nonmax-tangency}.
On the other hand, $L^\infty$-regularity of $Z(F)$ fails at $p=(0,1,0)$:
direct computation yields $(Z(F_x)\cdot Z(F_y))_p=4$,
whereas $\mu(Z(F),p)+(Z(F)\cdot L^\infty)_p-1=3$.}
\end{remark}

\newpage

\section{Polynomial and trigonometric curves}
\label{sec:poly-trig}

In Sections~\ref{sec:poly-trig} and~\ref{sec:expressive}, we introduce several classes of affine
(rather than projective) plane curves.
Before we begin, let us clarify what we mean by a real affine plane~curve.
There are two different notions here: an algebraic and a topological one:

\begin{definition}
\label{def:curves}
As usual, a reduced real \emph{algebraic curve}~$C$ in the complex affine plane
$\AA^2\cong\CC^2$ is the vanishing set
\[
C=V(G)=\{(x,y)\in\CC^2\mid G(x,y)=0\}
\]
of a squarefree bivariate polynomial $G(x,y)\in\RR[x,y]\subset\CC[x,y]$.
We view $C$ as a subset of~$\AA^2$, and implicitly identify it with the polynomial~$G(x,y)$
(viewed up to a constant nonzero factor),
or with the principal ideal generated by it, the ideal of polynomials vanishing on~$C$.

Alternatively, one can consider a ``topological curve''~$\CR$ in the real affine $(x,y)$-plane~$\RR^2$,
defined as the set of real points of an algebraic curve~$C$ as above:
\[
\CR=\VR(G)=\{(x,y)\in\RR^2\mid G(x,y)=0\}.
\]
In contrast to the real algebraic curve $V(G)\subset\CC^2$,
the real algebraic set~$\VR(G)$---even when it is one-dimensional---does not determine the polynomial~$G(x,y)$ up to a scalar factor.
In other words, an algebraic curve~$C$ is not determined by the set of its real points~$\CR$,
even when $\CR$ ``looks like'' an algebraic curve.
(Roughly speaking, this is because $C$ can have ``invisible'' components
which either have no real points at all, or all such points are isolated in~$\RR^2$.)
There is however a canonical choice, provided $\CR$ is nonempty and without isolated points:
we can let $C$ be the Zariski closure of~$\CR$,
or equivalently let $G(x,y)$ be the \emph{minimal polynomial} of~$\CR$,
a real polynomial of the smallest possible degree
satisfying \hbox{$\VR(G)=\CR$}. (The~minimal polynomial is defined up to a nonzero real factor.)
\end{definition}

Throughout this paper, we switch back and forth between a projective curve $Z(F(x,y,z))$
and its affine counterpart $V(G(x,y))$, where $G(x,y)=F(x,y,1)$. \linebreak[3]
(Remember that the line at infinity $L^\infty=Z(z)=\PP^2\setminus\AA^2$ is fixed throughout.)
For example, we say that $V(G)$ is $L^\infty$-regular if and only if $Z(F)$ is $L^\infty$-regular,
cf.\ Definition~\ref{def:c-regular}.

\begin{definition}
\label{def:poly-curve}
Let $C$ be a complex curve in the affine $(x,y)$-plane.
We say that
$C$ is a \emph{polynomial curve}
if it has a polynomial para\-metrization,
i.e., if there exist polynomials $X(t), Y(t)\in\CC[t]$ such that the map
$t\mapsto (X(t),Y(t))$
is a (birational, i.e., generically one-to-one) parametrization of~$C$.

A projective algebraic curve~$C=\{F(x,y,z)=0\}\subset\PP^2$ is called \emph{polynomial}
if $C$ does not contain the line at infinity $L^\infty=\{z=0\}$,
and the portion of $C$ contained in the affine $(x,y)$-plane (i.e., the curve $\{F(x,y,1)=0\}$)
is an affine polynomial curve.
\end{definition}

\begin{remark}
Not every polynomial map defines a polynomial parametrization.
For example, $t\mapsto (t^2-t,t^4-2t^3+t)$ is not a polynomial (or~birational)
parametrization, since it is not generically one-to-one:
$(X(t),Y(t))=(X(1-t),Y(1-t))$.
\end{remark}

\begin{example}
\label{ex:poly-with-elliptic-node}
The cubic $y^2=x^2(x-1)$
is a real polynomial curve, with a polynomial parametrization $t\mapsto (t^2+1,t(t^2+1))$.
Note that this curve has an elliptic node $(0,0)$, attained for imaginary parameter values $t=\pm\sqrt{-1}$.
\end{example}

\begin{example}
The ``witch of Agnesi'' cubic $x^2y+y-1$ is rational but not polynomial.
Indeed, if $X(t)$ is a positive-degree polynomial in~$t$, then $Y(t)=\frac{1}{(X(t))^2+1}$ is not.
\end{example}

\begin{figure}[ht]
\begin{center}
\vspace{-.1in}
\includegraphics[scale=0.25, trim=21cm 26.5cm 23cm 23cm, clip]{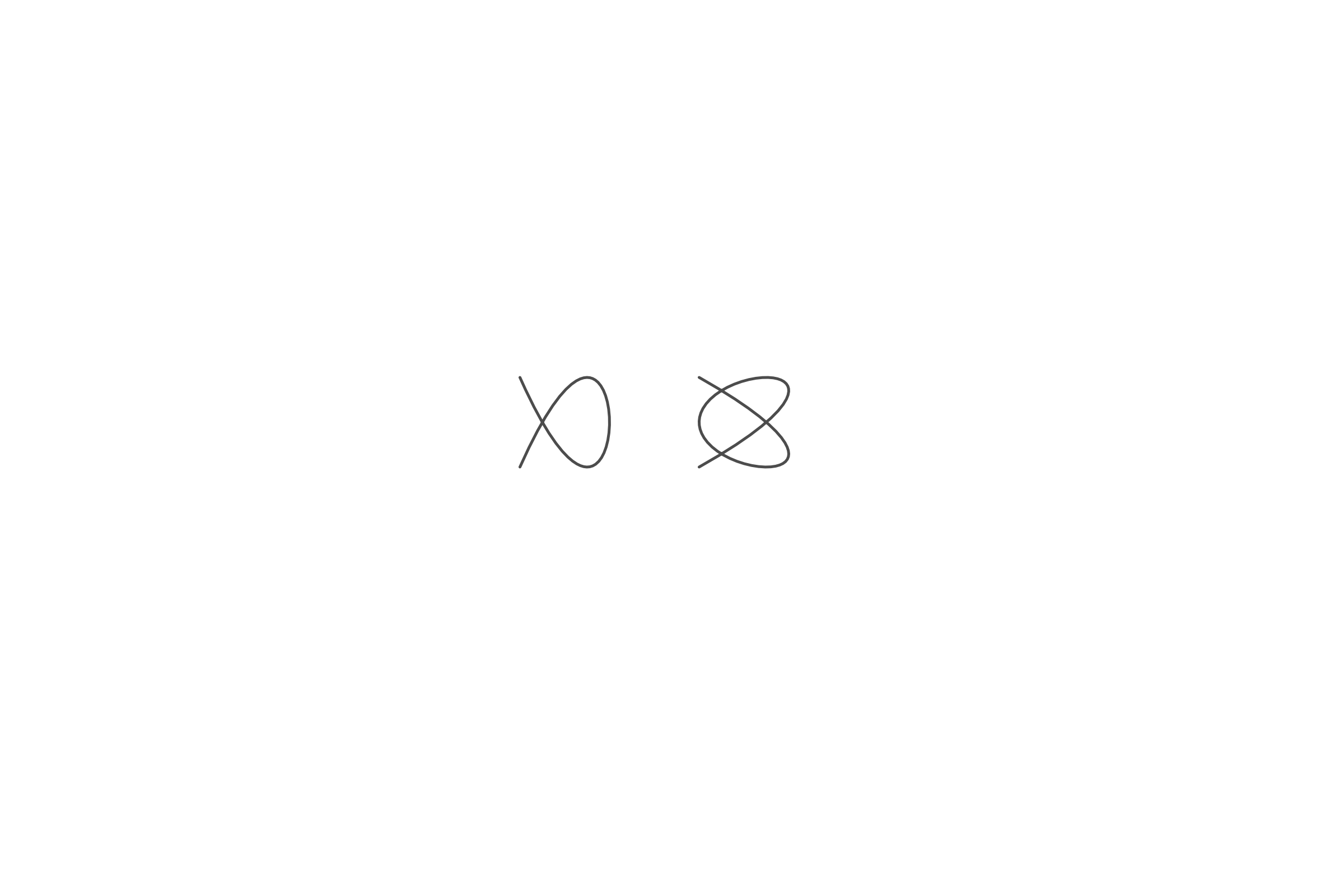}
\vspace{-.1in}
\end{center}
\caption{The Chebyshev curves with parameters $(3,2)$ and $(3,4)$.}
\vspace{-.2in}
\label{fig:chebyshev-34-32}
\end{figure}

\begin{example}[Irreducible Chebyshev curves]
Recall that the \emph{Chebyshev polynomials} of the first kind
are the univariate polynomials~$T_a(x)$
(here $a\in\ZZ_{>0}$) defined~by
\begin{equation}
\label{eq:chebyshev-poly}
T_a(\cos\varphi)=\cos(a\varphi).
\end{equation}
The polynomial $T_a(x)$ has integer coefficients, and is an even (resp., odd) function of~$x$ when $a$ is even (resp., odd).
We note that $T_a(T_b(t))=T_b(T_a(t))=T_{ab}(t)$.

Let $a$ and $b$ be coprime positive integers.
The \emph{Chebyshev curve} with parameters $(a,b)$
is given by the equation
\begin{equation}
\label{eq:Chebyshev-curve}
T_a(x)+T_b(y)=0.
\end{equation}
It is not hard to see that this curve is polynomial; let us briefly sketch why.
(For a detailed exposition, see \cite[Section~3.9]{Fischer}.)
Without loss of generality, let us assume that $a$ is odd.
Then the $(a,b)$-Chebyshev curve has a polynomial parametrization
$t\mapsto (-T_b(t), T_a(t))$.
Indeed, $T_a(-T_b(t))+T_b(T_a(t))=-T_{ab}(t)+T_{ab}(t)=0$.

To illustrate, consider the Chebyshev curves with parameters $(3,2)$ and $(3,4)$
shown in Figure~\ref{fig:chebyshev-34-32}.
The $(3,2)$-Chebyshev curve is a nodal Weierstrass cubic
\begin{equation}
\label{eq:32-Chebyshev}
4x^{3}-3x+2y^{2}-1=0,
\end{equation}
or parametrically $t\mapsto (-2t^2+1,4t^3-3t)$.
The $(3,4)$-Chebyshev curve (see Figure~\ref{fig:chebyshev-34-32})
is a quartic given by the equation
\begin{equation}
\label{eq:34-Chebyshev}
4x^{3}-3x+8y^{4}-8y^{2}+1=0,
\end{equation}
or by the polynomial parametrization
$t\mapsto (-8t^4+8t^2-1, 4t^3-3t)$.
\end{example}

\begin{lemma}
\label{lem:poly-curve}
For a real plane algebraic curve~$C$, the following are equivalent:
\begin{itemize}[leftmargin=.2in]
\item[{\rm (1)}]
$C$ is polynomial;
\item[{\rm (2)}]
$C$ has a parametrization $t\mapsto (X(t),Y(t))$
with $X(t),Y(t)\in\RR[t]$;
\item[{\rm (3)}]
$C$ is rational, with a unique local branch at infinity.
\end{itemize}
\end{lemma}

\begin{proof}
The equivalence (1)$\Leftrightarrow$(3) (for complex curves) is well known;
see, e.g.,~\cite{Abhyankar-1988}.
The implication (2)$\Rightarrow$(1) is obvious. It remains to show that (3)$\Rightarrow$(2).
The local branch of $C$ at infinity must by real, since otherwise complex conjugation would yield another such branch.
Consequently, the set of real points of~$C$ has a one-dimensional
connected component which contains the unique point $p\in C\cap L^\infty$.
It follows that the normalization map $\bn:\PP^1\to C\hookrightarrow\PP^2$ (which is nothing but a rational parametrization of $C$) pulls back the complex conjugation on $C$ to the antiholomorphic involution $c:\PP^1\to\PP^1$ which is a reflection with a fixed point set $\Fix(c)\simeq S^1$,
and $p$ lifts to a point in~$\Fix(c)$.
(There is another possible antiholomorphic involution on $\PP^1$, the antipodal one,
corresponding to real plane curves with a finite real point~set.)
Thus, we can choose coordinates $(t_0,t_1)$ on~$\PP^1$ so that $c(t_0,t_1)=(\overline t_0,\overline t_1)$ and the preimage of $p$ is $(0,1)$. Hence the map $\bn$ can be expressed~as
\[
x=X(t_0,t_1),\quad y=Y(t_0,t_1),\quad z=t_0^d\,,
\]
where $(t_0,t_1)\in\PP^1$, $X$~and $Y$ are bivariate homogeneous polynomials of degree $d=\deg C$, and by construction
\[
\overline x=X(\overline t_0,\overline t_1),\quad \overline y=Y(\overline t_0,\overline t_1),
\]
which means that $X$ and $Y$ have real coefficients.
\end{proof}

Recall that a \emph{trigonometric polynomial} is a finite linear combination
of functions of the form $t\mapsto \sin(kt)$ and/or $t\mapsto \cos(kt)$, with $k\in\ZZ_{\ge 0}$.

\begin{definition}
\label{def:trig-curve}
We say that a real \redsf{algebraic} curve $C$ in the affine $(x,y)$-plane is a \emph{trigonometric curve}
if there exist real trigonometric polynomials $X(t)$ and~$Y(t)$
such that $t\mapsto (X(t),Y(t))$ is a parametrization of~$\CR$, the set of real points of~$C$,
generically one-to-one for $t\in [0,2\pi)$.

A projective real algebraic curve~$C=\{F(x,y,z)=0\}\subset\PP^2$ is called \emph{trigonometric}
if $C$ does not contain the line at infinity $L^\infty=\{z=0\}$,
and the portion of $C$ contained in the affine $(x,y)$-plane 
is an affine trigonometric curve.
\end{definition}

\begin{remark}
Not every trigonometric map gives a trigonometric parametrization.
For example, $t\mapsto (\cos(t),\cos(2t))$ is not a trigonometric parametrization
of its image (a segment of the parabola $y\!=\!2x^2\!-\!1$), since it is not generically one-to-one
on~$[0,2\pi)$.
\end{remark}

\begin{example}
The most basic example of a trigonometric curve is the circle $t\mapsto (\cos(t),\sin(t))$,
or more generally an ellipse
\[
t\mapsto (A\cos(t),B\sin(t)) \quad (A, B\in\RR_{>0}).
\]
\end{example}

\begin{figure}[ht]
\begin{center}
\includegraphics[scale=0.2, trim=11cm 25cm 16cm 23cm, clip]{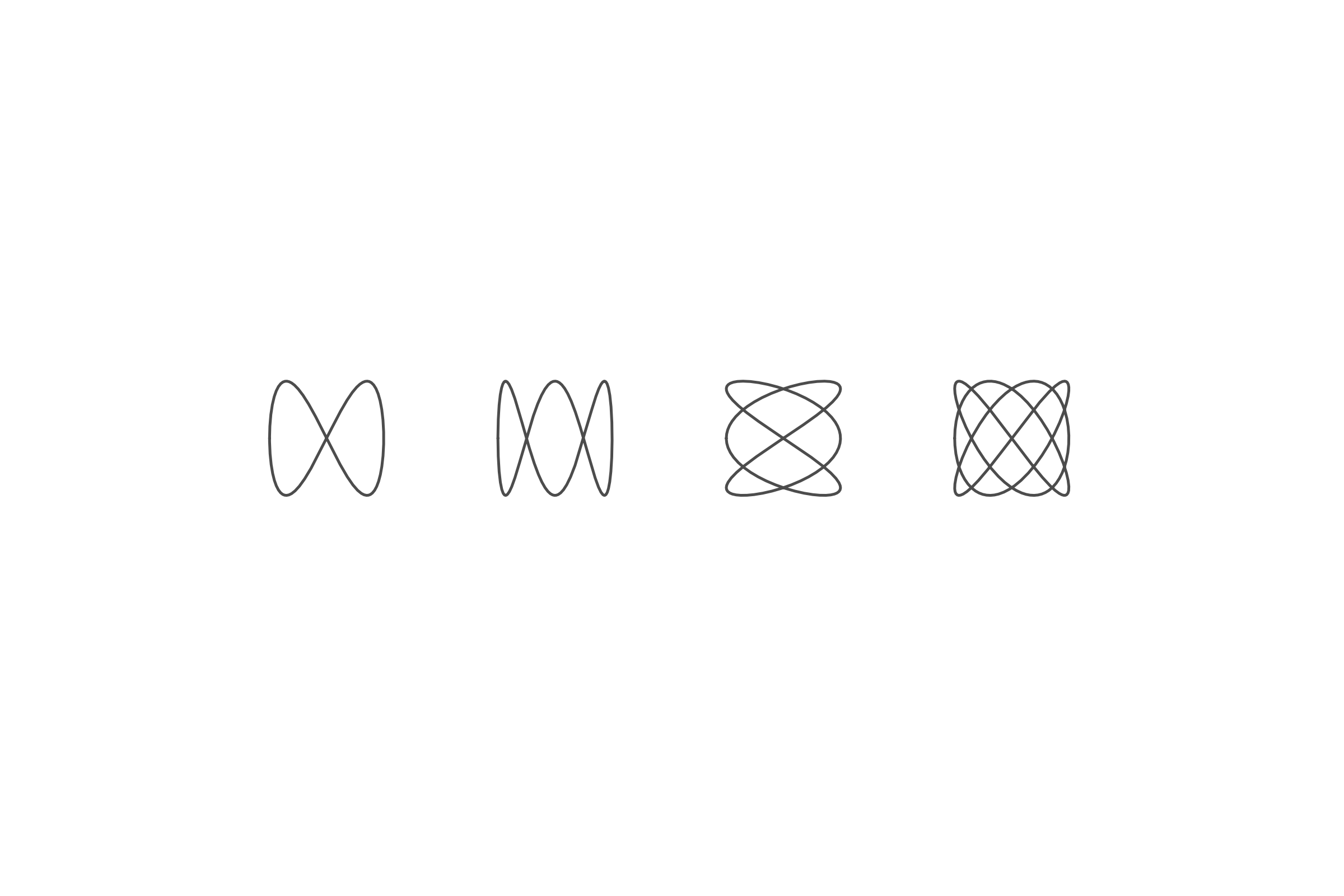}
\end{center}
\caption{The Lissajous curves with parameters $(2,1)$, $(3,1)$, $(2,3)$, and $(4,3)$.}
\vspace{-.2in}
\label{fig:lissajous-12-13-32-34}
\end{figure}

\begin{example}[Lissajous curves] 
\label{ex:Lissajous}
Let $k$ and $\ell$ be coprime positive integers, with $\ell$~odd.
The \emph{Lissajous curve} with parameters $(k,\ell)$
is a trigonometric curve defined by the parametrization
\[
t\mapsto(\cos(\ell t),\sin(kt)).
\]
The algebraic equation for this curve is
\begin{equation}
\label{eq:Lissajous-algebraic}
T_{2k}(x)+T_{2\ell }(y)=0,
\end{equation}
cf.~\eqref{eq:chebyshev-poly}.
(Indeed, $T_{2k}(\cos(\ell t))+T_{2\ell}(\cos(\frac{\pi}{2}-kt))=\cos(2k\ell t)+\cos(\ell\pi-2k\ell t)=0$.)
Note that \eqref{eq:Lissajous-algebraic} looks exactly like~\eqref{eq:Chebyshev-curve},
except that now the indices $2k$ and $2\ell$ are not coprime
(although $k$ and~$\ell$ are).

To illustrate, the $(2,3)$-Lissajous curve is given by the equation
\begin{equation}
\label{eq:32-Lissajous}
8x^4-8x^2+1+32y^6-48y^4+18y^2-1=0,
\end{equation}
or by the trigonometric parametrization
\begin{equation}
\label{eq:32-Lissajous-parametric}
t\mapsto (\cos(3t),\sin(2t)).
\end{equation}
Several Lissajous curves, including this one, are shown in Figure~\ref{fig:lissajous-12-13-32-34}.
\end{example}

\begin{figure}[ht]
\begin{center}
\hspace{-0.7in}\includegraphics[scale=0.2, trim=1cm 22cm 5cm 22cm, clip]{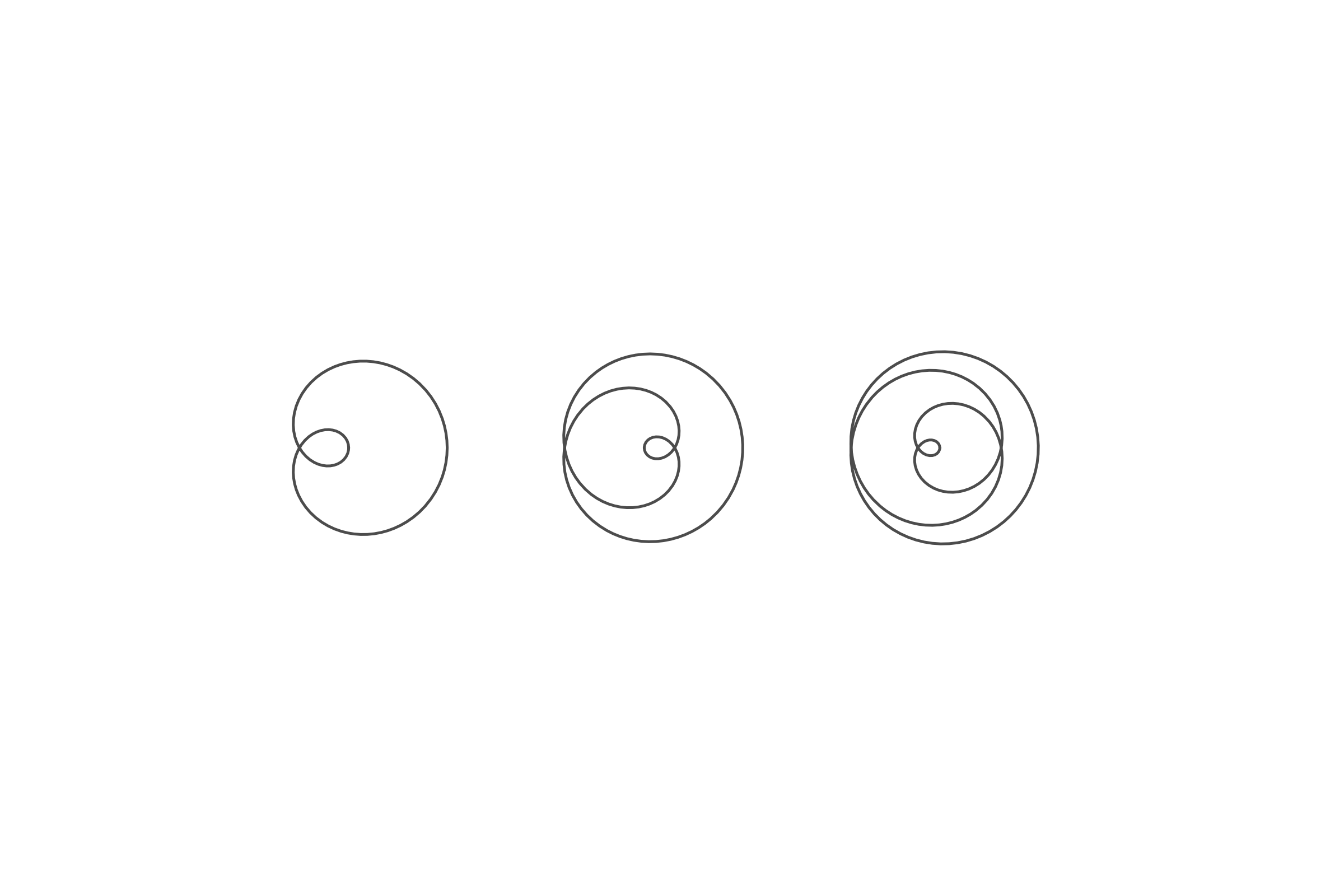}
\end{center}
\vspace{-.1in}
\caption{The rose curves $r=\cos\bigl(\tfrac{\theta}{3}\bigr)$ (cf.~\eqref{eq:limacon}),
$r=\cos\bigl(\tfrac{\theta}{5}\bigr)$, and $r=\cos\bigl(\tfrac{\theta}{7}\bigr)$.}
\label{fig:multi-limacons}
\end{figure}

\begin{example}[Rose curves]
\label{ex:rose-curves}
A \emph{rose curve} with parameter~$q=\frac{a}{b}\in\QQ_{>0}$ is defined in polar coordinates by
$r=\cos(q\theta)$, $\,\theta\in [0,2\pi b)$.
While a general rose curve has a complicated singularity at the origin,
it becomes nodal when $q\!=\!\frac{1}{2k+1}$, with $k\!\in\!\ZZ_{>0}$.
In that case, we get a ``multi-lima\c con,'' defined in polar coordinates by
$r=\cos\bigl(\tfrac{\theta}{2k+1}\bigr)$,
or equivalently by
\begin{equation}
\label{eq:multi-limacon}
r T_{2k+1}(r)-x=0.
\end{equation}
Note that the left-hand side of~\eqref{eq:multi-limacon} is a polynomial in $r^2=x^2+y^2$, so
it is an algebraic equation in $x$ and~$y$.
This is a trigonometric curve, with a trigonometric parametrization given by
\[
t\mapsto \Bigl(\frac{\cos(kt)+\cos((k+1)t)}{2}, \frac{\sin(kt)+\sin((k+1)t)}{2}\Bigr).
\]
The cases $k=1,2,3$ are shown in Figure~\ref{fig:multi-limacons}.

In the special case $k=1$, we get $rT_3(r)=4r^4-3r^2$,
and the equation~\eqref{eq:multi-limacon}~becomes
\begin{equation}
\label{eq:limacon}
4(x^2+y^2)^2-3(x^2+y^2)-x=0.
\end{equation}
This quartic curve is one of the incarnations of the \emph{lima\c con} of \'Etienne Pascal.
\end{example}

\vspace{-.2in}


\begin{lemma}
\label{lem:trig-curve}
For a real plane algebraic curve~$C$, the following are equivalent:
\begin{itemize}[leftmargin=.2in]
\item[{\rm (1)}]
$C$ is trigonometric;
\item[{\rm (2)}]
there exist polynomials $P(\varphi),Q(\varphi)\in\CC[\varphi]$ such that the map $\CC^*\to\AA^2$
given by
\begin{equation}
\label{eq:param-PQ}
\varphi\longmapsto (P(\varphi)+\overline P(\varphi^{-1}),Q(\varphi)+\overline Q(\varphi^{-1}))
\end{equation}
is a birational parametrization of~$C$.
\item[{\rm (3)}]
$C$ is rational, with two complex conjugate local branches at infinity and with an infinite real point set.\end{itemize}
\end{lemma}

\pagebreak[3]

\begin{proof}
$\boxed{(1)\Rightarrow(2)}$
The correspondence
$t\leftrightarrow \varphi=\exp(t\sqrt{-1})$
establishes a bianalytic isomorphism
between $(0,2\pi)$ and $S^1\setminus\{1\}$.
(Here $S^1=\{|\varphi|=1\}\subset\CC$.)
Under this correspondence, we have
\begin{equation}
\label{e4.5}
a\cos(kt)+b\sin(kt)=\tfrac{a-b\sqrt{-1}}{2}\varphi^k+\tfrac{a+b\sqrt{-1}}{2}\varphi^{-k}
\quad (a,b\in\RR),
\end{equation}
so any trigonometric polynomial in~$t$ transforms into a Laurent polynomial in~$\varphi$
of the form $P(\varphi)+\overline P(\varphi^{-1})$.
Thus a trigonometric parametrization of a curve~$C$ yields its
parametrization of the form~\eqref{eq:param-PQ};
this parametrization is generically one-to-one along $S^1$ and therefore extends to a birational map $\PP^1\to C$.

$\boxed{(2)\Rightarrow(1)}$
A parametrization \eqref{eq:param-PQ} of a curve~$C$
sends the circle~$S^1$ generically one-to-one to~$\CR$,
the real point set of~$C$, see the formulas~\eqref{e4.5}.
The same formulas \eqref{e4.5} convert the parametrization \eqref{eq:param-PQ} restricted to the circle~$S^1$ into a trigonometric parametrization $t\in[0,2\pi)\mapsto(X(t),Y(t))$ of~$\CR$.

$\boxed{(2)\Rightarrow(3)}$
A parametrization (\ref{eq:param-PQ}) intertwines the standard real structure in~$\PP^2$
and the real structure defined by the involution $c(\varphi)=\overline\varphi^{-1}$ on~$\CC\setminus\{0\}$.
Thus, it takes the set $S^1=\Fix(c)$ to the set~$\CR$ of real points of~$C$,
while the conjugate points $\varphi=0$ and $\varphi=\infty$ of~$\PP^1$ go to the points of~$C$ at infinity  determining two complex conjugate local branches at infinity.

$\boxed{(3)\Rightarrow(2)}$
Assuming~(3), the normalization map $\bn:\PP^1\to C$ pulls back the~standard complex conjugation in $\PP^2$ to the standard complex conjugation on~$\PP^1$,
while the circle $\RR\PP^1$ maps to the one-dimensional connected component of~$\CR$, and some complex conjugate points $\alpha,\overline\alpha\in\PP^1$ go to infinity. The automorphism of $\PP^1$ defined~by
\[
\varphi=\frac{s-\alpha}{s-\overline\alpha}
\]
takes the points $s=\alpha$ and $s=\overline\alpha$ to $0$ and $\infty$, respectively, the circle $\RR\PP^1$ to the circle~$S^1$, and the standard complex conjugation to the involution~$c$, see above. Hence the parametrization $\bn:\PP^1\to C$ goes to a parametrization (\ref{eq:param-PQ}).
\end{proof}

\begin{lemma}
\label{lem:real-points}
Let $C$ be a real polynomial (resp., trigonometric) nodal plane curve,
with a parametrization $\{(X(t),Y(t))\}$ as in Definition~\ref{def:poly-curve}
(resp., Definition~\ref{def:trig-curve}).
Assume that $C$ has no elliptic nodes.
Then $\CR=\{(X(t),Y(t))\mid t\in\RR\}$.
\end{lemma}

\begin{proof}
This lemma follows from the well known fact (see, e.g., \cite[Proposition~1.9]{IMR})
that the real point set of a real nodal rational curve~$C$ in $\RR\PP^2$
is the disjoint union of a circle $\RR\PP^1$ generically immersed in $\RR\PP^2$  and a finite set of elliptic nodes.
\end{proof}

If we allow elliptic nodes, the conclusion of Lemma~\ref{lem:real-points} can fail,
cf.\ Example~\ref{ex:poly-with-elliptic-node}.

Recall that an irreducible nodal plane curve of degree~$d$ has at most $\frac{(d-1)(d-2)}{2}$ nodes,
with the upper bound only attained for rational curves.
\redsf{In the case of trigonometric or polynomial curves, maximizing the number of nodes
has direct geometric consequences:}

\begin{proposition}[{\rm G.~Ishikawa \cite[Proposition~1.4]{Ishikawa-trigonometric}}]
\label{pr:trig-no-flex}
Let $C$ be a trigonometric curve of degree~$d$
with $\frac{(d-1)(d-2)}{2}$ real hyperbolic nodes.
Then $C$ has no inflection points.
\end{proposition}

\begin{proposition}
\label{pr:poly-no-flex}
A \redsf{(complex)} plane polynomial curve of degree $d$ with $\frac{(d-1)(d-2)}{2}$ nodes has no inflection points.
\end{proposition}

\begin{proof}
By Hironaka's genus formula (\ref{eq:Hironaka}), the projective closure $\hat C$ of $C$ has a single smooth point~$p$ on~$L^\infty$, with $(\hat C\cdot L^\infty)_p=d$.
We then determine the number of inflection points of~$C$ (in the affine plane)
using Pl\"ucker's formula (see, e.g., \cite[Chapter IV, Sections 6.2--6.3]{Walker}):
\[
2d(d-2)-(d-2)-6\cdot\tfrac{(d-1)(d-2)}{2}=0. \qedhere
\]
\end{proof}

\begin{figure}[ht]
\begin{center}
\vspace{-.2in}
\includegraphics[scale=0.14, trim=1cm 15cm 2cm 13cm, clip]{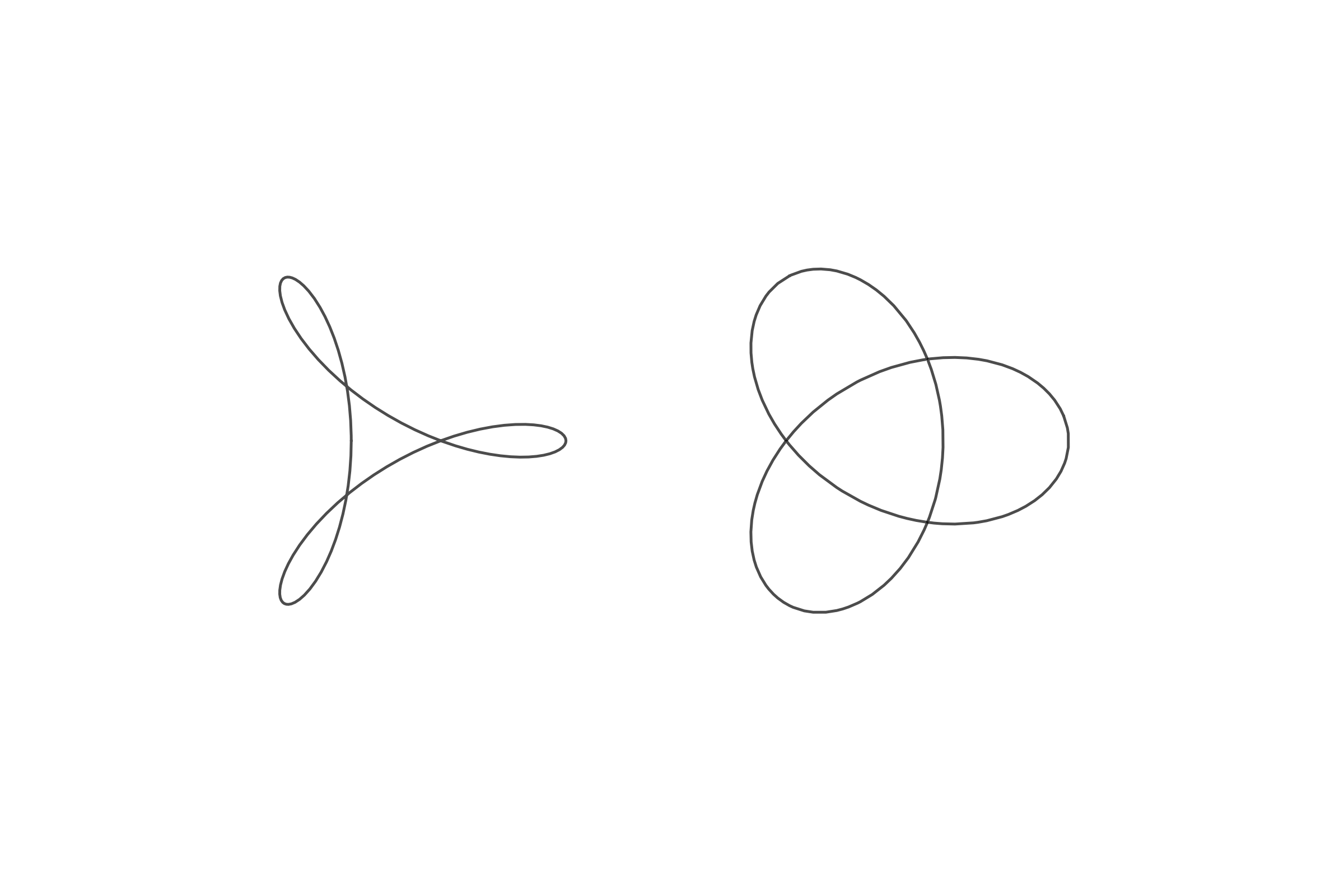}
\end{center}
\vspace{-.2in}
\caption{Three-petal hypotrochoids \eqref{eq:3-petal-param},
with $a=\frac34$ (left) and $a=2$ (right).
}
\label{fig:3-petal}
\end{figure}
\begin{figure}[ht]
\vspace{-.2in}
\begin{center}
\includegraphics[scale=0.2, trim=7cm 20cm 7cm 20cm, clip]{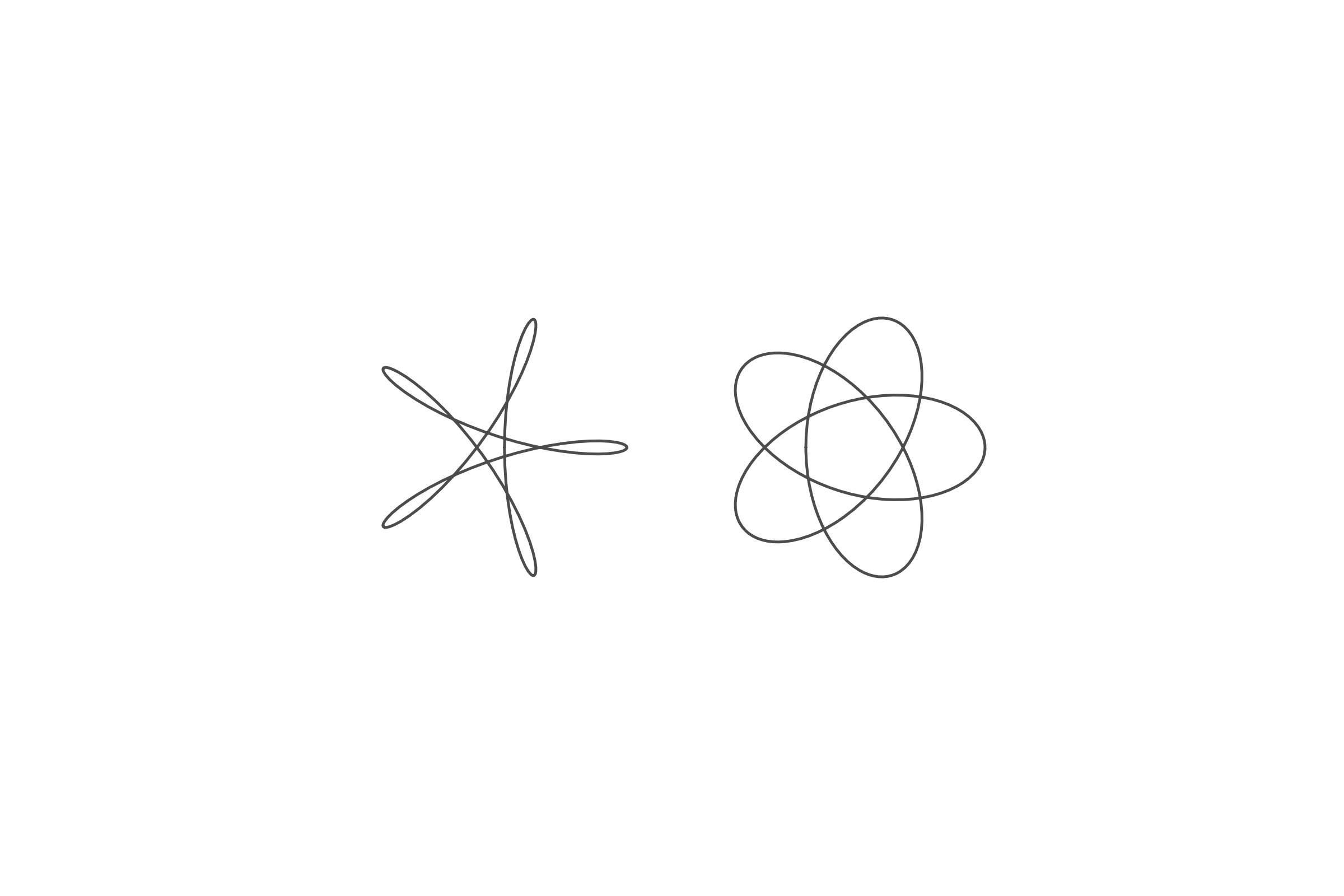}
\end{center}
\caption{Five-petal hypotrochoids $\{(\cos(2t)+a\cos(3t),\sin(2t)-a\sin(3t))\}$,
with $a=\frac56$ (left) and $a=2$ (right).
Each is a trigonometric curve of degree~\hbox{$d\!=\!6$}, with $\frac{(d-1)(d-2)}{2}=
10$~nodes.}
\label{fig:5-petal}
\end{figure}

\begin{example}
\label{ex:hypotrochoids}
The curve
\begin{equation}
\label{eq:hypotrochoid}
t\mapsto (\cos((k-1)t)+a\cos(kt), \sin((k-1)t)-a\sin(kt))
\end{equation}
is a trigonometric curve of degree~$d=2k$.
(It is a special kind of \emph{hypotrochoid}, cf.\ Definition~\ref{def:hypotrochoids} below.)
For \redsf{suitably} chosen real values of~$a$,
this curve has $\frac{(d-1)(d-2)}{2}=(k-1)(2k-1)$ real hyperbolic nodes,
as in Proposition~\ref{pr:trig-no-flex}.
See Figures~\ref{fig:3-petal} and~\ref{fig:5-petal}.

In the special case $k=1$ illustrated in Figure~\ref{fig:3-petal},
we get a three-petal hypotrochoid, a quartic trigonometric curve with 3~nodes
given by the parametrization
\begin{equation}
\label{eq:3-petal-param}
t\mapsto (\cos(t)+a\cos(2t),\sin(t)-a\sin(2t)),
\end{equation}
or by the algebraic equation
\begin{equation}
\label{eq:3-petal-alg}
a^2(x^2+y^2)^2+(-2a^4+a^2+1)(x^2+y^2)+(a^2-1)^3-2ax^3+6axy^2=0.
\end{equation}
\end{example}

\clearpage
\newpage

\section{Expressive curves and polynomials}
\label{sec:expressive}

\begin{definition}
\label{def:expressive-poly}
Let $G(x,y)\!\in\!\RR[x,y]\!\subset\!\CC[x,y]$ be a polynomial with real coefficients. \linebreak[3]
Let $C=V(G)$ be the corresponding affine algebraic curve,
and let $\CR=\VR(G)$ be the set of its real points, see Definition~\ref{def:curves}.
We say that $G(x,y)$ is an \emph{expressive polynomial}
(resp., $C$~is an \emph{expressive curve})~if
\begin{itemize}[leftmargin=.2in]
\item
all critical points of $G$ (viewed as a polynomial in $\CC[x,y]$) are real;
\item
all critical points of $G$ are Morse (i.e., have nondegenerate Hessians);
\item
each bounded component of $\RR^2\setminus \CR$
contains exactly one critical point of~$G$;
\item
each unbounded component of $\RR^2\setminus \CR$
contains no critical points;
\item
$\CR$ is connected, and contains at least two (hence infinitely many) points.
\end{itemize}
\end{definition}

\begin{remark}
Let $G(x,y)$ be a real polynomial with real Morse critical points.
Then each double point of  $\VR(G)$ must be a critical point of~$G$ (a~saddle).
Also, each bounded connected component of $\RR^2\setminus \VR(G)$
must contain at least one critical point (an extremum).
Thus, for~$G$ to be expressive,
it must have the smallest possible number of complex critical points
that is allowed by the topology of~$\VR(G)$:
a saddle at each double point,
one extremum within each bounded component of $\RR^2\setminus \VR(G)$, and nothing else.
\end{remark}

\begin{example}
\label{ex:quadratic}
The following quadratic polynomials are expressive:
\begin{itemize}[leftmargin=.2in]
\item
$G(x,y)=x^2-y$ has no critical points;
\item
$G(x,y)=x^2+y^2-1$ has one critical point $(0,0)$ (a minimum) lying inside
the unique bounded component of $\RR^2\setminus \VR(G)$;
\item
$G(x,y)=x^2-y^2$ has one critical point $(0,0)$ (a saddle), a hyperbolic node.
\end{itemize}
The following quadratic polynomials are not expressive:
\begin{itemize}[leftmargin=.2in]
\item
$G(x,y)=x^2-y^2-1$ has a critical point $(0,0)$
in an unbounded component of $\RR^2\setminus \VR(G)$;
besides, $\VR(G)$ is not connected;
\item
$G(x,y)=x^2+y^2$ has $\VR(G)$ consisting of a single point;
\item
$G(x,y)=x^2+y^2+1$ has $\VR(G)=\varnothing$;
\item
$G(x,y)=x^2-1$ and $G(x,y)=x^2$ have non-Morse critical points.
\end{itemize}
\end{example}

\begin{lemma}
\label{lem:expressive-poly}
Let $G(x,y)$ be an expressive polynomial.
Then:
\begin{itemize}[leftmargin=.2in]
\item
$G$ is squarefree (i.e, not divisible by a square of a non-scalar polynomial);
\item
$G$ has finitely many critical points, \redsf{all of them real};
\item
each critical point of $G$ is either a saddle or an extremal point;
\item
all saddle points of $G$ lie on $\VR(G)$;
they are precisely the singular points of~$V(G)$;
\item
each bounded connected component of $\RR^2\setminus \VR(G)$
 is simply connected.
 \end{itemize}
\end{lemma}

\begin{proof}
As the critical points of $G$ are real and Morse, each of them is either a saddle
or a local (strict) extremum of~$G$, viewed as a function $\RR\to\RR$.
The extrema must be located outside~$\CR$,
one per bounded connected component of $\RR^{2}\setminus\CR$.
The saddles must lie on~$\CR$, so they are precisely the double points of~it.
We conclude that $G$ has finitely many critical points.
Consequently $G$ is squarefree.
Finally, since $\CR$ is connected, each bounded component of $\RR^{2}\setminus\CR$
must be simply connected.
\end{proof}

\pagebreak[3]

Definition~\ref{def:expressive-poly}
naturally extends to homogeneous polynomials in three variables,
and to algebraic curves in the projective plane:

\begin{definition}
\label{def:expressive-xyz}
Let $F(x,y,z)\!\in\!\RR[x,y,z]\!\subset\!\CC[x,y,z]$ be a homogeneous polynomial \linebreak
with real coefficients, and $C\!=\!Z(F)$ the corresponding projective algebraic~curve.
Assume that $F(x,y,z)$ is not divisible by~$z$.
(In other words, $C$ does not contain the line at infinity~$L^\infty$.)
We~call $F$ and~$C$ \emph{expressive} if the bivariate polynomial
$F(x,y,1)$ is expressive in the sense of Definition~\ref{def:expressive-poly},
or equivalently the affine curve $C\setminus L^\infty\subset\AA^2=\PP^2\setminus L^\infty$
is expressive.
\end{definition}

In the rest of this section, we examine examples of expressive and non-expressive curves and polynomials.

\begin{example}[Conics]
Among real conics (cf.\ Example~\ref{ex:quadratic}),
a parabola, an ellipse, and a pair of crossing real lines are expressive,
whereas a hyperbola, a pair of parallel (or identical) lines,
and a pair of complex conjugate lines (an elliptic node) are not.
\end{example}

\begin{example}
\label{ex:two-graphs}
Let $f_1(x),f_2(x)\in\RR[x]$ be two distinct real univariate polynomials of degrees~$\le d$
such that $f_1-f_2$ has $d$ distinct real roots.
(Thus, at least one of $f_1,f_2$ has degree~$d$.)
We claim that the polynomial
\[
G(x,y)=(f_1(x)-y)(f_2(x)-y)
\]
is expressive.
To prove this, we first introduce some notation.
Let  $x_1,\dots,x_d$ be the roots of $f_1-f_2$, and let $z_1,\dots,z_{d-1}$ be the roots
of its derivative $f_1'-f_2'$; they are also real and distinct, by Rolle's theorem.
The critical points of~$G$ satisfy
\begin{align*}
G_x(x,y)&=f'_1(x)(f_2(x)-y) + (f_1(x)-y)f'_2(x)=0, \\
G_y(x,y)&=-f_1(x)-f_2(x)+2y=0.
\end{align*}
It is straightforward to see that these equations have $2d-1$ solutions:
$d$~hyperbolic nodes $(x_k,f_1(x_k))=(x_k,f_2(x_k))$, for $k=1,\dots,d$,
as well as $d-1$ extrema at the points $(z_k,\frac12(f_1(z_k)+f_2(z_k))$, for $k=1,\dots,d-1$.
The conditions of Definition~\ref{def:expressive-poly} are now easily verified.
(Alternatively, use the coordinate change $(x,\tilde y)=(x,y-f_1(x))$
to reduce the problem to the easy case when one of the two polynomials is~$0$.)
\end{example}

\begin{example}[Lemniscates]
\label{ex:lemniscates}
The \emph{lemniscate of Huygens} (or Gerono)
is given by the equation
$y^2+4x^4-4x^2=0$, or by the parametrization $t\mapsto (\cos(t),\sin(2t))$.
This curve, shown in Figure~\ref{fig:lissajous-12-13-32-34} on the far left,
is a Lissajous curve with parameters~$(2,1)$, cf.\ Example~\ref{ex:Lissajous}.
The polynomial $G(x,y)=y^2+4x^4-4x^2$
has three real critical points:
a saddle at the hyperbolic node~$(0,0)$,
plus two extrema $(\pm \frac{1}{\sqrt{2}},0)$ inside the two bounded connected components
of $\RR^2\setminus \VR(G)$.
Thus, this lemniscate is expressive.

By contrast, another quartic curve with a similar name (and a similar-looking set of real points),
the \emph{lemniscate of Bernoulli}
\begin{equation}
\label{eq:lemniscate-Bernoulli}
(x^2+y^2)^2-2x^2+2y^2=0,
\end{equation}
is not expressive, as the polynomial $G(x,y)=(x^2+y^2)^2-2x^2+2y^2$
has critical points $(0,\pm i)$ outside~$\RR^2$.
\end{example}

\newpage

Many more examples of expressive and non-expressive polynomials (or curves)
are given in Tables~\ref{table:conics+cubics} and~\ref{table:conics+cubics+quartics},
and later in the paper.
\bigskip

\begin{table}[ht]
\begin{tabular}{|p{1.2in}|p{1.9in}|p{2.0in}|}
\hline
\multicolumn{3}{|c|}{\textbf{expressive conics}} \\
  \hline
 $G(x,y)$ & real curve $\VR(G)$ & critical points \\
\hline
$x^2-y$ & parabola & none \\
$x^2-y^2$ & two lines & saddle\\
$x^2+y^2-1$ & ellipse & extremum \\
\hline
\multicolumn{3}{|c|}{\textbf{non-expressive conics}} \\
\hline
$G(x,y)$ & real point set $\VR(G)$ & why not expressive? \\
\hline
$x^2-y^2-1$ & hyperbola & saddle in an unbounded region 
\\
$x^2+y^2+1$ & imaginary ellipse & $\VR(G)$ is empty
\\
$x^2+y^2$ & elliptic node & $\VR(G)$ is a single point  \\
$x^2+a$ \ ($a\in\RR$) & two parallel lines & critical points are not Morse
\\
\hline
%
%
\multicolumn{3}{|c|}{\textbf{expressive cubics}} \\
\hline
$G(x,y)$ & real point set $\VR(G)$ & critical points \\
\hline
$x^3-y$ & cubic parabola & none \\
$(x^2-y)x$ & parabola and its axis & single saddle \\
$x^3-3x+2-y^2$ & nodal Weierstrass cubic & one saddle, one extremum \\
$(x-1)xy$ & two parallel lines + line
& two saddles 
\\
$(x^2-y)(y-1)$ & parabola + line & two saddles, one extremum \\
$(x+y-1)xy$ & three lines & three saddles, one extremum \\
$(x^2+y^2-1)x$ & ellipse + line & two saddles, two extrema \\
\hline
\multicolumn{3}{|c|}{\textbf{non-expressive cubics}} \\
\hline
$G(x,y)$ & real point set $\VR(G)$ & why not expressive? \\
\hline
$x^3-y^2$ & semicubic parabola 
   & critical point is not Morse \\
$x^3-3x-y^2$ & two-component elliptic curve & saddle 
in an unbounded region
\\
$x^3-3x+3-y^2$ & one-component elliptic curve & saddle 
in an unbounded region
\\
$x^3+3x-y^2$ & one-component elliptic curve & 
two non-real critical points \\
$x^3+xy^2+4xy+y^2$ & oblique strophoid & two non-real critical points \\
$x^2y-x^2+2y^2$ & Newton's species \#54 & 
two non-real critical points \\
$yx^2-y^2-xy+x^2$ & Newton's species \#51 & two non-real critical points \\
$x^3+y^3+1$ & Fermat cubic & critical point is not Morse \\
$x^3+y^3-3xy$ & folium of Descartes & two non-real critical points \\
$x^2y+y-x$ & serpentine curve & 
two non-real critical points \\
$x^2y+y-1$ & witch of Agnesi & 
two non-real critical points \\
$x^2y$ & double line + line & critical points are not Morse
\\
$x(x^2-y^2)$ & three concurrent lines & critical point 
is not Morse \\
$(x^2-y^2-1)y$ & hyperbola + line & two non-real critical points \\
$(xy-1)x$ & hyperbola + asymptote & $\VR(G)$ is not connected \\
\hline
\end{tabular}
\caption{Expressive and non-expressive conics and cubics.
Unless specified otherwise, lines are placed so as to maximize the number of crossings. }
\label{table:conics+cubics}
\end{table}

\clearpage
\newpage

\begin{table}[ht]
\begin{tabular}{|c|c|p{1.9in}|p{2.95in}|c|}
\hline
$d$ & $\xi$ & $G(x,y)$ & real point set $\VR(G)$  \\
\hline
\multicolumn{4}{|c|}{} \\[-.18in]
\hline
$1$ & 0 & $x$ & line  \\
\hline
$2$ & 0 & $x^2-y$ & parabola  \\
$2$ & 1 & $x^2-y^2$ & two crossing lines \\
$2$ & 1 & $x^2+y^2-1$ & ellipse  \\
\hline
$3$ & 0 & $x^3-y$ & cubic parabola  \\
$3$ & 1 & $(x^2-y)x$ & parabola and its axis  \\
$3$ & 2 & $x^3-3x+2-y^2$ & nodal Weierstrass cubic  \\
$3$ & 2 & $(x-1)xy$ & three lines, two of them parallel  \\
$3$ & 3 & $(x^2-y)(y-1)$ & parabola crossed by a line  \\
$3$ & 4 & $(x+y-1)xy$ & three generic lines  \\
$3$ & 4 & $(x^2+y^2-1)x$ & ellipse crossed by a line   \\
\hline
$4$ & 0 & $y-x^4$ & quartic parabola \\
$4$ & 0 & $(y-x^2)^2-x$ & \\
$4$ & 1 & $(y-x^2)^2-x^2$ & two aligned parabolas \\
$4$ & 1 & $(y-x^2)^2+x^2-1$ & \\
$4$ & 1 & $(y-x^3)x$ & cubic parabola and its axis \\
$4$ & 2 & $(y-x^2)^2-xy$& \\
$4$ & 2 & $(y-x^2)x(x-1)$& parabola + two lines parallel to its axis \\
$4$ & 3 & $y^2-(x^2-1)^2$ & co-oriented parabolas crossing at two points\\
$4$ & 3 & $y^2+4x^4-4x^2$ & lemniscate of Huygens \\
$4$ & 3 & $x(x-1)(x+1)y$ & three parallel lines crossed by a fourth\\
$4$ & 3 & 
                 $4(x^2+y^2)^2-3(x^2+y^2)-x$ & lima{\c c}on \\
$4$ & 5 & $x(x-1)y(y-1)$ & two pairs of parallel lines\\
$4$ & 5 & $(x^2+y^2-1)(x^2-2x+y^2)$ & two circles crossing at two points\\
$4$ & 5 & $(y-x^3+x)y$ & cubic parabola + line \\
$4$ & 5 & $(x^3-3x+2-y^2)(x-a)$ & nodal Weierstrass cubic $\!+\!$ line crossing it at~$\infty$\\
$4$ & 6 & 
$4x^3-3x+8y^4-8y^2+1$ & $(3,4)$-Chebyshev curve\\
$4$ & 6 & $(y^2+4x-6)^2 - x^3 - 3x^2 $ & \\ 
$4$ & 6 & $(y-x^2)(x-a)(y-1) $ & parabola + line + line parallel to the axis \\
$4$ & 6 & $(y-x^2)(y-1)(y-2) $ & parabola + two parallel lines \\
$4$ & 7 & $(y-x^2+1)(x-y^2+1) $ & two parabolas crossing at four points\\
$4$ & 7 & see \eqref{eq:3-petal-alg} & three-petal hypotrochoid \\
$4$ & 7 & $(x^{3}-3x+2-y^{2})(x+y-a)$& nodal Weierstrass cubic + line \\
$4$ & 7 & $(x+y)(x-y)(x-1)(x-a)$ & line + line + two parallel lines \\
$4$ & 7 & $(x^2+y^2-1)x(x-a)$  & ellipse + two parallel lines\\
$4$ & 8 & $(4y-x^2)(x+y)(ax+y+1) $& parabola + two lines \\
$4$ & 8 & $(x^2+y^2-1)(y-4x^2+2)$ & ellipse and parabola crossing at four points \\
$4$ & 9 & $x(y+1)(x-y)(x+y-1)$ & four lines \\
$4$ & 9 & $(x^2+y^2-1)x(x+y-a)$& ellipse + two lines\\
$4$ & 9 & $(x^2 + 4y^2)(4x^2 + y^2)$ & two ellipses crossing at four points \\
\hline
\end{tabular}
\caption{Expressive curves $C\!=\!\{G\!=\!0\}$ of degrees $d\le 4$.
Unless noted otherwise, lines are drawn to maximize the number of crossings.
The value~$a\!\in\!\RR$ 
is to be chosen appropriately.
We denote by~$\xi$ the number of critical points~of~$G$.
}
\label{table:conics+cubics+quartics}
\end{table}

\clearpage
\newpage

\begin{remark}
Definition~\ref{def:expressive-poly} can be generalized to allow arbitrary ``hyperbolic''
singular points, i.e., isolated real singular points all of whose local branches are real.
\end{remark}

\begin{remark}
In can be verified by an exhaustive case-by-case analysis that every expressive curve
of degree $d\le 4$ is $L^\infty$-regular.
Starting with $d\ge 5$, this is no longer the case,
cf.\ Example~\ref{example:expressive-non-reg} below:
an expressive curve~$C$ need not be $L^\infty$-regular, even when $C$ is rational.
Still, examples like this one are rare.
\end{remark}

\begin{example}
\label{example:expressive-non-reg}
%
Consider the real rational quintic curve~$C$
parametrized by
\[
x=t^2, \ 
y=t^{-1}+t^{-2}-t^{-3}.
\]
The set of its real points in the $(x,y)$-plane
consists of two interval components corresponding to the negative and
positive values of~$t$, respectively. These components intersect at the point $(1,1)$,
attained for $t=1$ and $t=-1$. The algebraic equation of~$C$ is obtained as follows:
\begin{align*}
&y-x^{-1}=t^{-1}-t^{-3} ,\\
&(y-x^{-1})^2=t^{-2}-2t^{-4}+t^{-6}=x^{-1}-2x^{-2}+x^{-3}, \\
&x^3y^2-2x^2y-x^2+3x-1=0.
\end{align*}
The Newton triangle of $C$ is $\conv\{(0,0),(3,2),(2,0)\}$, which means that
\begin{itemize}[leftmargin=.2in]
\item
at the point $p_1=(1,0,0)$, the curve $C$ has a type~$A_2$ (ordinary cusp) singularity,
tangent to the axis $y=0$;
\item
at the point $p_2=(0,1,0)$, the curve $C$ has a type~$E_8$ singularity, tangent to the axis $x=0$.
\end{itemize}
In projective coordinates, we have $C=Z(F)$ where
\begin{align*}
F&=x^3y^2-2x^2yz^2-x^2z^3+3xz^4-z^5,\\
F_x&=3x^2y^2-4xyz^2-2xz^3+3z^4,\\
F_y&=2x^3y-2x^2z^2.
\end{align*}
It is now easy to check that the only critical point of $F(x,y,1)$
is the hyperbolic node $(1,1)$ discussed above.
It follows that the curve $C$ is expressive.

We next examine the behavior of the polar curves at infinity.
In a neighborhood of the point $p_2=(0,1,0)$, we set $y=1$ and obtain
\begin{align*}
F_x&=3x^2-4xz^2-2xz^3+3z^4,\\
F_y&=2x^3-2x^2z^2.
\end{align*}
The curve $Z(F_x)$ has two smooth local branches at~$p_2$, quadratically
tangent to the axis $Z(x)$;
the curve $Z(F_y)$ has the double component $Z(x)$
and a smooth local branch quadratically tangent to~$Z(x)$.
Thus the
intersection multiplicity~is
\[
(Z(F_x)\cdot Z(F_y))_{p_2}=12>\mu(C,p_2)+(C\cdot L^\infty)_{p_2}-1
=8+3-1=10,
\]
so $L^\infty$-regularity fails at~$p_2$.
(The intersection at $p_1$ is regular by
Proposition~\ref{pr:nonmax-tangency}.)
\end{example}

\begin{example}[Lissajous-Chebyshev curves]
\label{ex:lissajous-chebyshev}
Let $a$ and $b$ be positive integers. 
The \emph{Lissajous-Chebyshev curve}
with parameters~$(a,b)$ is given by the equation
\begin{equation}
\label{eq:Ta+Tb}
T_a(x)+T_b(y)=0.
\end{equation}
When  $a$ and~$b$ are coprime, with $a$ odd, we recover the polynomial Chebyshev curve
with parameters $(a,b)$, see~\eqref{eq:Chebyshev-curve}.
When both $a$ and $b$ are even, with $\frac{a}{2}$ and $\frac{b}{2}$ coprime,
we recover the trigonometric Lissajous curve
with parameters $(\frac{a}{2},\frac{b}{2})$, see~\eqref{eq:Lissajous-algebraic}.
This explains our use of the term ``Lissajous-Chebyshev curve.''
The general construction appears in the work of S.~Guse{\u\i}n-Zade~\cite{GZ1974},
who observed (without using this terminology) that the Lissajous-Chebyshev curve with parameters~$(a,b)$
provides a morsification of an isolated quasihomogeneous singularity of type~$(a,b)$.

There is also a variant of the Lissajous-Chebyshev curve defined by
\begin{equation}
\label{eq:Ta-Tb}
T_a(x)-T_b(y)=0.
\end{equation}
When $a$ (resp.,~$b$) is odd, this curve is a mirror image of the Lissajous-Chebyshev curve~\eqref{eq:Ta+Tb},
under the substitution $x:=-x$ (resp., $y:=-y$).
However, when both $a$ and $b$ are even, the two curves differ.
For example, for $(a,b)=(4,2)$, the curve defined by~\eqref{eq:Ta+Tb}
is the lemniscate of Huygens (see Example~\ref{ex:lemniscates}),
whereas the curve defined by~\eqref{eq:Ta-Tb} is a union of two parabolas.

It is not hard to verify that every Lissajous-Chebyshev curve
(and every curve $V(T_a(x)-T_b(y))$) is expressive.
The critical points of 
$T_a(x)\pm T_b(y)$ are found from the equations
\[
T_a'(x)=T_b'(y)=0,
\]
so they are of the form $(x_i,y_j)$ where $x_1,\dots,x_{a-1}$ (resp., $y_1,\dots,y_{b-1}$)
are the (distinct) roots of~$T_a$ (resp.,~$T_b$).
Since the total number of nodes and bounded components of $\RR^2\setminus\CR$
is easily seen to be exactly $(a-1)(b-1)$, the claim follows.
\end{example}

\begin{figure}[ht]
\begin{center}
\includegraphics[scale=0.25, trim=18cm 19.5cm 11cm 20cm, clip]{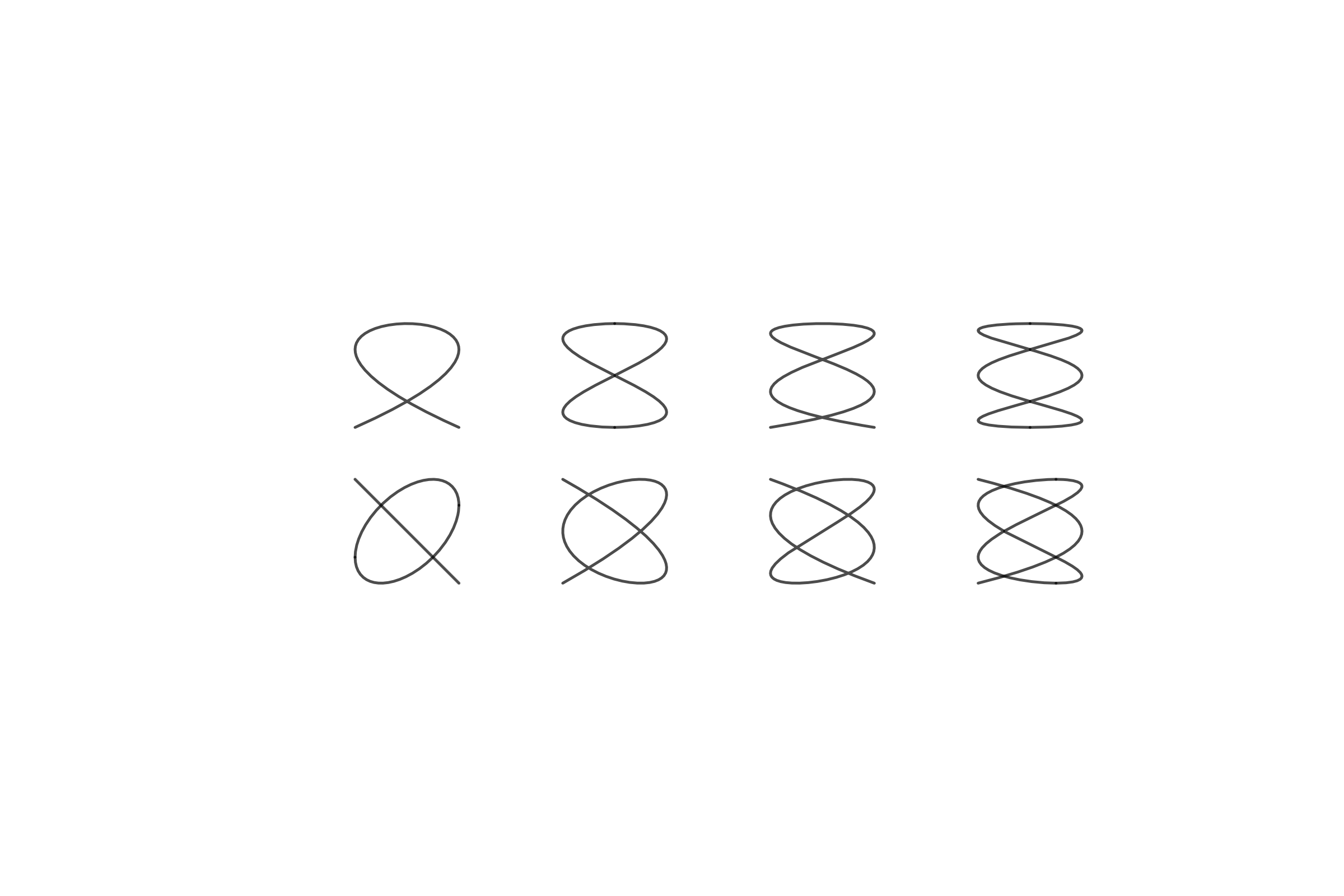}
\end{center}
\caption{Lissajous-Chebyshev curves with parameters $(2,3)$, $(2,4)$, $(2,5)$, $(2,6)$
(top row) and $(3,3)$, $(3,4)$, $(3,5)$, $(3,6)$ (bottom row).}
\vspace{-.2in}
\label{fig:lissajous-chebyshev}
\end{figure}

\begin{example}[Multi-lima\c cons]
\label{ex:multi-limacons}
Recall that the multi-lima\c con with parameter~$k$ is a trigonometric curve~$C$ of degree~$2k+2$
given by the equation~\eqref{eq:multi-limacon}.
It is straightforward to verify that the corresponding polynomial has $2k+1$ critical points,
all of them located on the $x$ axis.
Comparing this to the shape of the curve~$C_\RR$, we conclude that $C$ is expressive.
\end{example}

\newpage

\section{Divides}
\label{sec:divides}

The notion of a divide was first introduced and studied by
N.~A'Campo \cite{acampo-ihes, ishikawa}.
The version of this notion that we use in this paper differs slightly from
A'Campo's, and from the version used in~\cite{FPST}.

\begin{figure}[ht]
\begin{center}
\setlength{\unitlength}{1.6pt}
\begin{picture}(80,100)(10,-50)
\thicklines
\put(40,-51.5){\cyan{\line(0,1){102.5}}}
\put(0,0){\red{\qbezier(0,20)(40,-10)(60,-10)}}
\put(0,0){\red{\qbezier(60,-10)(80,-10)(80,0)}}
\put(0,0){\red{\qbezier(60,10)(80,10)(80,0)}}
\put(0,0){\red{\qbezier(0,-20)(40,10)(60,10)}}
\put(-40,0){\qbezier(135,20)(95,-10)(75,-10)}
\put(-40,0){\qbezier(75,-10)(55,-10)(55,0)}
\put(-40,0){\qbezier(75,10)(55,10)(55,0)}
\put(-40,0){\qbezier(135,-20)(95,10)(75,10)}
\put(40,0){\blue{\circle{35}}}
\thicklines
\put(63,23){{\makebox(0,0){\footnotesize region}}}
\put(63,19){\gray{\vector(0,-1){15}}}
\put(76.3,-25.5){{\makebox(0,0){\footnotesize node}}}
\put(76.3,-22.5){\gray{\vector(0,1){15}}}
\put(28.3,-33.5){{\makebox(0,0){\footnotesize branch}}}
\put(29.9,-30.5){\gray{\vector(0,1){15}}}
\put(47.5,0){\green{\circle{104}}}
\end{picture}
\end{center}
\caption{A divide, its branches, regions, and nodes.}
\vspace{-.2in}
\label{fig:divide}
\end{figure}
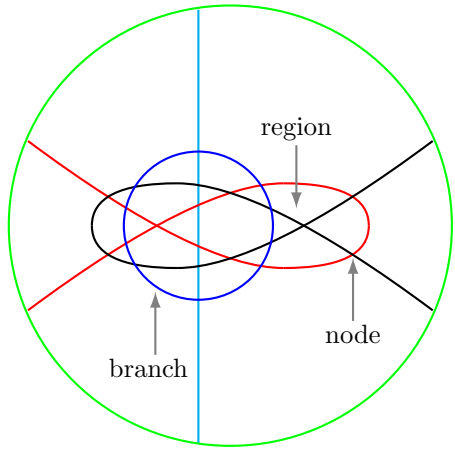

\begin{definition}
\label{def:divide}
Let 
$\Disk$ be a
disk in the real~plane~$\RR^2$.
A 
\emph{divide}~$D$ in~$\Disk$ is the image of a generic relative immersion of a finite set of intervals and circles into~$\Disk$ satisfying the conditions listed below.
The images of these immersed intervals and circles are called the \emph{branches} of~$D$.
They must satisfy the following conditions:
\begin{itemize}[leftmargin=.2in]
\item
the immersed circles do not intersect the boundary~$\partial\Disk$;
\item
the endpoints of the immersed intervals lie on~$\partial\Disk$, and are pairwise distinct;
\item
these immersed intervals intersect $\partial D$ transversally;
\item
all intersections and self-intersections of the branches
are transversal.
\end{itemize}
We view divides as topological objects, i.e., we do not distinguish between divides related by a diffeomorphism between their respective ambient disks.

A divide is called \emph{connected} if the union of its branches is connected.

The connected components of the complement $\Disk\setminus D$
which are disjoint from $\partial\Disk$ are the \emph{regions} of~$D$.
If $D$ is connected, then each region of~$D$ is simply connected.
We refer to the singular points of~$D$ as its \emph{nodes}.
See Figure~\ref{fig:divide}.
\end{definition}

\begin{remark}
Although all divides of interest to us are connected, we forego any connectivity requirements in the above definition of a divide,
which therefore is slightly more~general than \cite[Definition~2.1]{FPST}.
\end{remark}

The main focus of~\cite{FPST} was on the class of \emph{algebraic} divides
coming from real~\emph{morsi\-fi\-cations} of isolated plane curve singularities,
see \cite[Definition~2.3]{FPST}.
Here~we~study a different (albeit related) class of divides which arise from real algebraic curves:

\begin{definition}
\label{def:divide-polynomial}
Let $G(x,y)\in\RR[x,y]$ be a real polynomial such that
each real singular point of the curve $V(G)\subset\CC^2$
is a \emph{hyperbolic node} (an intersection of two smooth real local branches).
Then the portion of $\VR(G)$ contained in a sufficiently large disk
\begin{equation}
\label{eq:Disk_R}
\Disk_R=\{(x,y)\in\RR^2\mid x^2+y^2\le R^2\}
\end{equation}
gives a divide in~$\Disk_R$.
Moreover this divide does not depend (up to homeomorphism) on the choice of~$R\gg0$. We denote this divide~by~$D_G$.
\end{definition}

%

%
%

\pagebreak[3]

%
%
%
%
%

\begin{proposition}
\label{pr:expressivity-criterion}
Let $G(x,y)\in\RR[x,y]$ be a real polynomial with $\xi<\infty$ critical points.
Assume that the real algebraic set $\VR(G)=\{G=0\}\subset\RR^2$ is nonempty,
and each singular point of $\VR(G)$ is a hyperbolic node.
Let~$\nu$ be the number of such nodes, and let $\iota$ be the number of
interval branches of the divide~$D_G$.
Then
\begin{equation}
\label{eq:expressivity-xi}
\xi \ge 2\nu-\iota+1,
\end{equation}
with equality if and only if the polynomial $G$ is expressive.
\end{proposition}

\begin{proof}
Let $K$ denote the union of the divide~$D_G$ and all its regions,
viewed as a closed subset of~$\RR^2$.
Let $\sigma$ be the number of connected components of~$K$.
Since all these components are simply connected,
the Euler characteristic of~$K$ is equal~to~$\sigma$.
On the other hand, $K$ can be split into 0-dimensional cells
(the nodes, plus the ends of interval branches),
1-dimensional cells (curve segments of $D$ connecting nodes),
ovals (smooth closed components of~$D_G$),
and regions.
The number~$\beta$ of 1-dimensional cells satisfies $2\beta\!=\!2\iota+4\nu$
(by counting endpoints), implying $\beta\!=\!\iota+2\nu$.
We thus~have
\[
\sigma = \chi(K)=(\nu+2\iota)-(\iota+2\nu)+\rho-h = \iota-\nu+\rho-h,
\]
where $\rho$ is the number of regions in~$D_G$,
and $h$ denotes the total number of holes (the sum of first Betti numbers) over all regions.

Denote $u=\xi-\nu-\rho$.
The set of critical points of~$G$ contains all of the nodes,
plus at least one extremum per region.
Thus $u=0$ if $G$ has no other critical points, and $u>0$ otherwise.

Putting everything together, we obtain:
\begin{align*}
\xi=\nu+\rho+u
=\nu+\sigma-\iota+\nu+h+u
=2\nu-\iota+\sigma+h+u.
\end{align*}
Since $\sigma\ge 1$ and $h,u\ge 0$, we get~\eqref{eq:expressivity-xi}.
Moreover $\xi = 2\nu-\iota+1$ if and only if $K$ is connected,
all regions are simply connected, and $G$ has exactly $\nu+\rho$ critical~points.
All these conditions are satisfied if $G$ is expressive;
conversely, they ensure expressivity.
($\VR(G)$~is connected if $K$ is connected and each region is simply connected.)
\end{proof}

\begin{remark}
Table~\ref{table:conics+cubics+quartics} does not list any expressive quartic
with $\xi=4$.
We~can now explain why.
\redsf{First, by Proposition \ref{pr:expressivity-criterion}, $\xi=4$ and $d=4$ would imply that $2\nu=\iota+3$.

\textbf{Case~1}:
$\nu=3$, $\iota=3$.
Then our quartic $C$ has three real irreducible components,
\emph{viz.}, a conic (necessarily a parabola) and two lines crossing it, with three hyperbolic nodes total,
and no other nodes.
Such a configuration is impossible:
if both lines are parallel to the axis of the parabola, then we get two nodes;
otherwise, at least four.

\textbf{Case~2}:
$\nu=2$, $\iota=1$.
Then $C$ is irreducible.
(Otherwise, $C$~would split into~two conics---a parabola and an ellipse---intersecting
at four real points, contrary to $\nu=2$.)
An irreducible quartic~$C$ with $\iota=1$ has one real local branch $Q$ centered at some real point $p\in L^\infty$,
plus perhaps a pair of complex conjugate local branches centered on~$L^\infty$. 
Let us list the possibilities. 

\textbf{Case~2A}:
$Q$ is smooth, $(Q\cdot L^\infty)_p=4$, with no other local branches centered~on~$L^\infty$.

\textbf{Case~2B}:
$Q$ is smooth, $(Q\cdot L^\infty)_p=2$, with two smooth complex conjugate branches transversal to~$L^\infty$. 
These two local branches either cross $L^\infty$ at different points, or at the same point~$p'$. 
By the genus formula~\eqref{eq:Hironaka} and due to $\nu=2$, we have $p'\ne p$, with a nodal singularity at~$p'$. 

\textbf{Case~2C}:
$Q$ is singular. In this case, again using $\nu=2$ and the genus formula~\eqref{eq:Hironaka}, 
we conclude that $Q$ is of type~$A_2$. 
The intersection multiplicity of a cusp and a line is either $2$ or~$3$.  
In our setting, $(Q\cdot L^\infty)_p=2$, 
with two additional smooth complex conjugate local branches of $C$ centered at distinct complex conjugate points on~$L^\infty$. 

All cases 2A--2C are subject to Proposition~\ref{pr:nonmax-tangency}, so $C$ is $L^\infty$-regular. 
A direct computation then yields that $\sum_{q\in C\cap L^\infty}\mu(C,q,L^\infty)\le3$ in every case. 
This, however, is in contradiction with
$\sum_{q\in C\cap L^\infty}\mu(C,q,L^\infty)=(4-1)^2-\xi=5$.}
\end{remark}


\section{\texorpdfstring{$L^\infty$}{L}-regular expressive curves}
\label{sec:regular+expressive}

\begin{proposition}
\label{pr:expressive+regular}
Let $C=Z(F)\subset\PP^2$ be a reduced algebraic curve defined
by a real homogeneous polynomial~$F(x,y,z)\in\RR[x,y,z]$.
Assume that
\begin{itemize}[leftmargin=.35in]
\item[{\rm (a)}]
all irreducible components of $C$ are real;
\item[{\rm (b)}]
$C$ does not contain the line at infinity~$L^\infty$ as a component;
\item[{\rm (c)}]
all singular points of $C$ in the affine $(x,y)$-plane are real hyperbolic nodes;
\item[{\rm (d)}]
the polynomial $F(x,y,1)\in\CC[x,y]$ has finitely many critical points;
\item[{\rm (e)}]
the set of real points $\{F(x,y,1)=0\}\subset\RR^2$ is nonempty.
\end{itemize}
Then the following are equivalent:
\begin{itemize}[leftmargin=.35in]
\item[{\rm (i)}]
the curve $C$ is expressive and $L^\infty$-regular;
\item[{\rm (ii)}]
each irreducible component of $C$ is rational,
with a set of local branches at infinity consisting of
either  a unique (necessarily real) local branch,
or a pair of complex conjugate local branches, possibly based at the same real point.
\end{itemize}
\end{proposition}

\begin{proof}
The proof is based on Propositions~\ref{pr:hironaka-milnor} and~\ref{pr:expressivity-criterion}.
Let us recall the relevant notation, and introduce additional one:
\begin{align*}
d&=\deg(C),\\
\iota&=\text{number of interval branches of~$D_G$}\\
s&=|\Comp(C)|=\text{number of irreducible components of~$C$},\\
s_1&= \text{number of components of $C$ with a real local branch at infinity}, \\
s_2&= \text{number of components of $C$ with a pair of complex conjugate}\\
&\hspace{.225in} \text{local branches at infinity}.
\end{align*}
Combining Propositions~\ref{pr:hironaka-milnor} and~\ref{pr:expressivity-criterion},
we conclude that
\begin{equation*}
2g(C)-2+\iota+\sum_{p\in C\cap L^\infty}\br(C,p)
\ge 0,
\end{equation*}
or equivalently (see \eqref{eq:g(C)})
\begin{equation}
\label{eq:restated-reg-expressive}
2\sum_{C'\in\Comp(C)} g(C')+\iota+\sum_{p\in C\cap L^\infty}\br(C,p)
\ge 2s,
\end{equation}
with equality if and only if $C$ is both expressive and $L^\infty$-regular.

On the other hand, we have the inequalities
\begin{align}
\label{eq:sum-g(C')}
\sum_{C'\in\Comp(C)} g(C')&\ge 0, \\[-.05in]
\iota &\ge s_1\,, \\
\sum_{p\in C\cap L^\infty}\br(C,p) &\ge s_1+2s_2\,, \\
\label{eq:s+s}
2s_1 + 2s_2 &\ge 2s,
\end{align}
whose sum yields~\eqref{eq:restated-reg-expressive}.
Therefore we have equality in~\eqref{eq:restated-reg-expressive} if and only if
each of \eqref{eq:sum-g(C')}--\eqref{eq:s+s} is an equality.
This is precisely statement~(ii).
\end{proof}

\begin{example}
The curve $C$ from Example~\ref{example:expressive-non-reg}
satisfies requirements (a)--(e) of Proposition~\ref{pr:expressive+regular}.
Since $C$ is not $L^\infty$-regular, condition~(ii) 
must~fail.
Indeed, while $C$ is rational, it has two real local branches at infinity,
centered at the points $p_1$ and~$p_2$.
\end{example}

\begin{proposition}
\label{prop:nonempty}
Let $C$ be an expressive $L^\infty$-regular plane curve whose irreducible components are all real.
Then each component of $C$ is either trigonometric or polynomial.
\end{proposition}

\begin{proof}
Since $C$ is $L^\infty$-regular and expressive, with real components,
the requirements (a)--(e) of Proposition~\ref{pr:expressive+regular}
are automatically satisfied, cf.\ Lemma~\ref{lem:expressive-poly}.
Consequently the statement (ii) of Proposition~\ref{pr:expressive+regular} holds.
By Lemmas~\ref{lem:poly-curve} and~\ref{lem:trig-curve}, this would imply
that each component of $C$ is either trigonometric or polynomial,
provided the real point set of each component is infinite.
It thus suffices to show that each component $B$ of $C$ has an infinite real point set in~$\AA^2$.
In fact, it is enough to show that this set is nonempty,
for if it were finite and nonempty, then~$B$---hence~$C$---would have an elliptic node,
contradicting the expressivity of~$C$ (cf.\ Lemma~\ref{lem:expressive-poly}).

It remains to prove that each component of $C$ has a nonempty real point set in~$\AA^2$.

We argue by contradiction.
Let $C\!=\!Z(F)$.
Suppose that $B\!=\!Z(G)$ is a component of~$C$ without real points in~$\AA^2$.
We claim that the rest of~$C$ is given by a polynomial of the form~$H(G,z)$,
where $H\in\RR[u,v]$ is a bivariate polynomial.
Once we establish this claim, it will follow that the polynomials $F_x$ and $F_y$
have a non-trivial common factor $\frac{\partial}{\partial u}(uH(u,v))\big|_{u=G,v=z}$,
contradicting the finiteness of the intersection $Z(F_x)\cap Z(F_y)$.

Let us make a few preliminary observations. First, the degree $d=\deg G=\deg B$ must be even.
Second, by Proposition~\ref{pr:expressive+regular}, $B$~has two complex conjugate branches centered on~$L^\infty$.
Third, in view of the expressivity of $C$, the affine curve $B\cap\AA^2$ is disjoint from any other component $B'$ of $C$, implying that
\begin{equation}
B\cap B'\subset L^\infty.
\label{el1}
\end{equation}

Consider two possibilities.

\textbf{Case~1}: $B\cap L^\infty$ consists of two complex conjugate points $p$ and $\overline p$.
Let $Q$ and $\overline Q$ be the local branches of $B$ centered at $p$ and $\overline p$, respectively.

Let $B'=Z(G')$ be some other component of $C$, of degree $d'=\deg B'=\deg G'$.
Since $B'$ is real and satisfies statement~(ii) of Proposition~\ref{pr:expressive+regular}
as well as \eqref{el1}, we~get
\[
B'\cap B=B'\cap L^\infty=B\cap L^\infty=\{p,\overline p\};
\]
moreover $B'$ has a unique local branch $R$ (resp.,~$\overline R$)
at the point $p$ (resp.~$\overline p$).
We have
\begin{align*}
(R\cdot Q)_p=(\overline R\cdot\overline Q)_{\overline p}=\tfrac{d'd}{2},\qquad
((L^\infty)^{d'}\cdot Q)_p=((L^\infty)^{d'}\cdot\overline Q)_{\overline p}&=\tfrac{d'd}{2}.
\end{align*}
It follows that any curve $\hat B$ in the pencil
$\Span\{B',(L^\infty)^{d'}\}$ satisfies
\begin{equation}
(\hat B\cdot Q)_p\ge\tfrac{d'd}{2},\quad (B'\cdot \overline Q)_{\overline p}\ge\tfrac{d'd}{2}.
\label{el2}
\end{equation}
Pick a point $q\in B\cap\AA^2$.
Since $q\not\in B'\cup L^\infty$, there exists a curve $\tilde B\in\Span\{B',(L^\infty)^{d'}\}$ containing~$q$.
It then follows by B\'ezout and by \eqref{el2} that $\tilde B$ must contain $B$ as a component.
In particular, $d'\ge d$. By symmetry, the same argument yields $d\ge d'$.
Hence $d'=d$, $B\in\Span\{B',(L^\infty)^d\}$, and the polynomial $G'$ defining the curve $B'$ satisfies $G'=\alpha G+\beta z^d$ for some $\alpha,\beta\in\CC\setminus\{0\}$.
The desired claim follows.

\pagebreak[3]

\textbf{Case~2}: $B\cap L^\infty$ consists of one (real) point~$p$.
Then $B$ has two complex conjugate branches $Q$ and $\overline Q$ centered at~$p$.

Let $B'=Z(G')$ be a component of $C$ different from $B$, and let $B'$ have a unique (real) local branch~$R$, necessarily centered at~$p$.
Denote $d'=\deg B'=\deg G'$. Then
\[
(R\cdot Q)_p=(R\cdot \overline Q)_p=\tfrac{d'd}{2}, \qquad
((L^\infty)^{d'}\cdot Q)_p=((L^\infty)^{d'}\cdot\overline Q)_p=\tfrac{d'd}{2},
\]
which implies (cf.\ Case~1) that any curve $\hat B'\in\Span\{B',(L^\infty)^{d'}\}$ satisfies
\[
(\hat B'\cdot q)_p\ge\tfrac{d'd}{2},\quad (\hat B'\cdot \overline Q)_p\ge\tfrac{d'd}{2},
\]
and then we conclude---as above---that there exists a curve $\tilde B'\in\Span\{B',(L^\infty)^{d'}\}$ containing $B$ as a component. On the other hand,
\[
(B\cdot R)_p=d'd,\quad ((L^\infty)^d\cdot R)_p=d'd,
\]
which in a similar manner implies that there exists a curve $\tilde B\in\Span\{B,(L^\infty)^d\}$
containing $B'$ as a component.
We conclude that $d'=d$, $B\in\Span\{B',(L^\infty)^d\}$,
and finally, $G'=\alpha G+\beta z^d$ for $\alpha,\beta\in\CC\setminus\{0\}$, as desired.

Now let $B'=Z(G')$ be a component of $C$ different from $B$, and let it have a couple of complex conjugate local branches $R$ and~$\overline R$ centered at $p$. Since $B'$ is real, we have
\[
(B'\cdot Q)_p=(B'\cdot \overline Q)_p=\tfrac{d'd}{2},\quad
((L^\infty)^{d'}\cdot Q)_p=((L^\infty)^{d'}\cdot\overline Q)_p=\tfrac{d}{2},
\]
and since $B$ is real, we have
\[
(B\cdot R)_p=(B\cdot\overline R)_p=\tfrac{d'd}{2},\quad ((L^\infty)^d\cdot R)_p=((L^\infty)^d\cdot\overline R)_p=\tfrac{d'd}{2}.
\]
Thus the above reasoning applies again, yielding $d'=d$ and $B'\in\Span\{B,(L^\infty)^d\}$.
Hence $G'=\alpha G+\beta z^d$ for $\alpha,\beta\in\CC\setminus\{0\}$, and we are done.
\end{proof}

The following example shows that in Proposition~\ref{prop:nonempty},
the requirement that all components are real cannot be dropped.

\begin{example}
\label{ex:reducible-quintic}
The quintic curve $C=Z(F)$ defined by the polynomial
\[
F(x,y,z)=(x^2+z^2)(yx^2+yz^2-x^3)
\]
has two non-real components $Z(x\pm z\sqrt{-1})$.
In Example~\ref{ex:(x^2+z^2)(yz^2-x^3+yx^2)-regular},
we verified that this curve is $L^\infty$-regular.
It is also expressive, because the polynomial $G(x,y)=F(x,y,1)$
has no critical points in the complex affine plane (see Example~\ref{ex:(x^2+z^2)(yz^2-x^3+yx^2)-regular}) 
and the set 
\begin{equation}
\label{eq:y=x^3/...}
\VR(G)=\{(x,y)\in\RR^2\mid y=\tfrac{x^3}{x^2+1}\}
\end{equation}
is connected.
On the other hand, the real irreducible component
$\widetilde C=Z(yx^2+yz^2-x^3)$
is neither trigonometric nor polynomial because it has two real points at infinity,
see Example~\ref{ex:(x^2+z^2)(yz^2-x^3+x^2y)}.
Furthermore, $\widetilde C$ is not expressive,
since the polynomial
$yx^2+y-x^3$
has two critical points $(\pm \sqrt{-1}, \pm \frac32 \sqrt{-1})$ outside~$\RR^2$.
\end{example}

Proposition~\ref{prop:nonempty} immediately implies the following statement.

\begin{corollary}
\label{cor:reg-expressive-irreducible}
An irreducible $L^\infty$-regular expressive curve is either trigonometric or polynomial.
\end{corollary}

We note that such a curve also needs to be immersed, meeting itself transversally at real hyperbolic nodes (thus, no cusps, tacnodes, or triple points).

\medskip

We next provide a partial converse to Corollary~\ref{cor:reg-expressive-irreducible}.

\begin{proposition}
\label{pr:irr-expressive}
Let $C$ be a real polynomial or trigonometric curve
whose singular set in the affine plane~$\AA^2=\PP^2\setminus L^\infty$
consists solely of hyperbolic nodes.
Then $C$ is expressive and $L^\infty$-regular.
\end{proposition}

\begin{proof}
By Lemmas~\ref{lem:poly-curve} and~\ref{lem:trig-curve},
a real polynomial or trigonometric curve~$C=Z(F)$ is a
real rational curve with one real or two complex conjugate local branches at infinity,
and with a nonempty set of real points in~$\AA^2$.
Thus, conditions (ii), (a), (b), (c), and~(e) of Proposition \ref{pr:expressive+regular}
are satisfied, so in order to obtain~(i), we only need to establish~(d).
That is, we need to show that the polynomial $F$ has finitely many critical points in the affine plane $\AA^2$.
We will prove this by contradiction.

Denote $d=\deg C=\deg F$. Suppose that $Z(F_x)\cap Z(F_y)$ contains a (real, possibly reducible) curve $B$ of a positive degree $d'<d$.

We first observe that $B\cap C\cap\AA^2=\emptyset$.
Assume not. If $q\in B\cap C\cap \AA^2$, then $q\in\Sing(C)\cap\AA^2$,
so $q$ must be a hyperbolic node of~$C$,
implying $(Z(F_x)\cdot Z(F_y))_q=\mu(C,q)=1$;
but since $q$ lies on a common component of $Z(F_x)$ and $Z(F_y)$,
we must have $(Z(F_x)\cdot Z(F_y))_q=\infty$.
Since $B$ is real, and $C$ has either one real or two complex conjugate local branches at infinity,
it follows that $B\cap C=L^\infty\cap C$.

\textbf{Case~1}: $C$ has a unique (real) branch $Q$ at a point $p\in L^\infty$.
Then $(B\cdot Q)_p=d'd$. On the other hand, $((L^\infty)^{d'}\cdot Q)_p=d'd$.
Hence any curve $\hat B\in\Span\{B,(L^\infty)^{d'}\}$ satisfies $(\hat B\cdot Q)_p\ge d'd$.
For any point $q\in C\setminus\{p\}$ there is a curve $\tilde B\in\Span\{B,(L^\infty)^{d'}\}$ passing through~$q$.
Hence $C$ is a component of~$\tilde B$, in contradiction with $d'<d$.

\textbf{Case~2}: $C$ has two complex conjugate branches $Q$ and~$\overline Q$
centered at (possibly coinciding) points $p$ and~$\overline p$ on $L^\infty$, respectively.
Since $B$ is real, we have
\[
(B\cdot Q)_p=(B\cdot\overline Q)_{\overline p}=((L^\infty)^{d'}\cdot Q)_p=((L^\infty)^{d'}\cdot\overline Q)_{\overline p}=d'd/2.
\]
This implies that any curve $\hat B\in\Span\{B,(L^\infty)^{d'}\}$ satisfies both
$(\hat B\cdot Q)_p\ge d'd/2$ and $(\hat Q\cdot\overline Q)_{\overline p}\ge d'd/2$,
resulting in a contradiction as in Case~1.
\end{proof}

\begin{remark}
In Proposition~\ref{pr:irr-expressive}, the requirement concerning the singular set 
cannot be dropped: a real polynomial curve~$C$ may have elliptic nodes
(see Example~\ref{ex:poly-with-elliptic-node}) or cusps
(consider the cubic $(t^2,t^3)$), preventing~$C$ from being expressive.
\end{remark}

\begin{remark}
\label{rem:poly-trig-expr=>reg}
Proposition~\ref{pr:irr-expressive} implies that if a real curve $C$ is polynomial or trigonometric,
and also expressive, then it is necessarily $L^\infty$-regular.
Indeed, expressivity in particular means that all singular points of~$C$ are hyperbolic nodes.
\end{remark}

\pagebreak[3]

\begin{example}
\label{ex:hypotrochoid-expressive}
Recall from Example~\ref{ex:hypotrochoids} that for
appropriately chosen values of the real parameter~$a$,
the hypotrochoid~$C$ given by the equation~\eqref{eq:hypotrochoid}
is a nodal trigonometric curve of degree~$d=2k$
with the maximal possible number of real hyperbolic nodes,
namely $\frac{(d-1)(d-2)}{2}=(k-1)(2k-1)$.
Thus $C$ has no other singular points in the affine plane,
and consequently is expressive by Proposition~\ref{pr:irr-expressive}.
\end{example}

Combining Corollary~\ref{cor:reg-expressive-irreducible} and Proposition~\ref{pr:irr-expressive},
we obtain:

\begin{theorem}
\label{th:reg-expressive-irreducible}
For a real plane algebraic curve~$C$, the following are equivalent:
\begin{itemize}[leftmargin=.2in]
\item
$C$ is irreducible, expressive and $L^\infty$-regular;
\item
$C$ is either trigonometric or polynomial, and all its singular points in the affine plane~$\AA^2$ are hyperbolic nodes.
\end{itemize}
\end{theorem}

\begin{example}
Theorem~\ref{th:reg-expressive-irreducible} is illustrated in Table~\ref{table:irreducible-expressive}.
Each real curve in this table is either polynomial or trigonometric.
We briefly explain why all these curves are expressive (hence $L^\infty$-regular, see Remark~\ref{rem:poly-trig-expr=>reg}).

Expressivity of Chebyshev and Lissajous curves was established in Example~\ref{ex:lissajous-chebyshev}.
%
The lima\c con of Pascal was discussed in Example~\ref{ex:multi-limacons}.
In Example~\ref{ex:hypotrochoid-expressive}, we saw~that
a hypotrochoid~\eqref{eq:hypotrochoid} is expressive for suitably chosen values of~$a$.
For a more general statement, see Proposition~\ref{pr:epi-hypo} below.
As to the remaining curves in Table~\ref{table:irreducible-expressive},
all we need to check is that each of their singular points in the affine plane is a hyperbolic node.
The curves $V(x^d-y)$ are smooth, so there is nothing to prove.
Ditto for the ellipse, as well as the curves $V((y\!-\!x^2)^2\!-x)$ and $V((y\!-\!x^2)^2\!+x^2\!-\!1)$.
Finally, $V((y\!-\!x^2)^2\!-\!xy)$ has a single singular point in~$\AA^2$, a hyperbolic node at the origin.
\end{example}


\begin{table}[ht]
\begin{center}
\begin{tabular}{|c|c|p{1.7in}|p{1.8in}|p{1.25in}|}
\hline
$d$ & $\xi$ & $G(x,y)$
& $(X(t), Y(t))$ & $V(G)$
\\
\hline
\multicolumn{4}{|c|}{} \\[-.18in]
\hline
$1$ & 0 & $x-y$
& $(t, t)$ & line
\\
\hline
$2$ & 0 & $x^2-y$
& $(t, t^2)$ & parabola
\\
$2$ & 1 & $x^2+y^2-1$
& $(\cos(t),\sin(t))$ & ellipse
\\
\hline
$3$ & 0 & $x^3-y$
& $(t,t^3)$ & cubic parabola
\\
$3$ & 2 &
$4x^{3}-3x+2y^{2}-1$
&
$(-2t^2+1,4t^3-3t)$
& $(3,2)$-Chebyshev
\\
\hline
$4$ & 0 & $x^4-y$ & $(t,t^4)$ & quartic parabola \\
$4$ & 0 & $(y-x^2)^2-x$ & $(t^2,t^4+t)$ & \\
$4$ & 1 & $(y-x^2)^2+x^2-1$ & $(\cos(t), \sin(t)\!+\!\cos^2(t))$ & \\
$4$ & 2 & $(y-x^2)^2-xy$ & $(t^2-t,t^4-t^3)$ & \\
$4$ & 3 & $y^2+4x^4-4x^2$ & $(\cos(t), \sin(2t))$ & (1,2)-Lissajous
\\
$4$ & 3 & 
$4(x^2+y^2)^2\!-\!3(x^2+y^2)\!-\!x$
& 
$(\cos(t)\!\cos(3t),\cos(t)\!\sin(3t)$
& lima{\c c}on
       \\
$4$ & 6 & $4x^3-3x+8y^4-8y^2+1$  & $(-8t^4+8t^2-1, 4t^3-3t)$
& $(3,4)$-Chebyshev
\\
$4$ & 7 & see \eqref{eq:3-petal-alg} & see \eqref{eq:3-petal-param}
& 
(2,1)-hypotrochoid
\\
\hline
\end{tabular}
\vspace{-.1in}
\end{center}
\caption{Irreducible expressive curves $C=V(G)$  of degrees $d\le 4$.
All curves are $L^\infty$-regular.
For~each curve, a trigonometric or polynomial parametrization $(X(t),Y(t))$ is shown.
We denote by~$\xi$ the number of critical points~of~$G$.
}
\label{table:irreducible-expressive}
\end{table}

\pagebreak[3]

The following construction provides a rich source of examples of expressive curves.

\begin{definition}[Epitrochoids and hypotrochoids]
\label{def:hypotrochoids}
Let $b$ and $c$ be coprime nonzero integers, with $b>|c|$.
Let $u$ and $v$ be nonzero reals.
The trigonometric curve
\begin{align}
\label{def:epi-hypo1}
x&=u\cos bt+v\cos ct, \\
\label{def:epi-hypo2}
y&=u\sin bt-v\sin ct.
\end{align}
is called a \emph{hypotrochoid} if $c>0$, and an \emph{epitrochoid} if $c<0$.
It is a rational curve of degree~$2b$.
\end{definition}

\begin{figure}[ht]
\begin{center}
\includegraphics[scale=0.28, trim=16.5cm 15cm 17.5cm 18.5cm, clip]{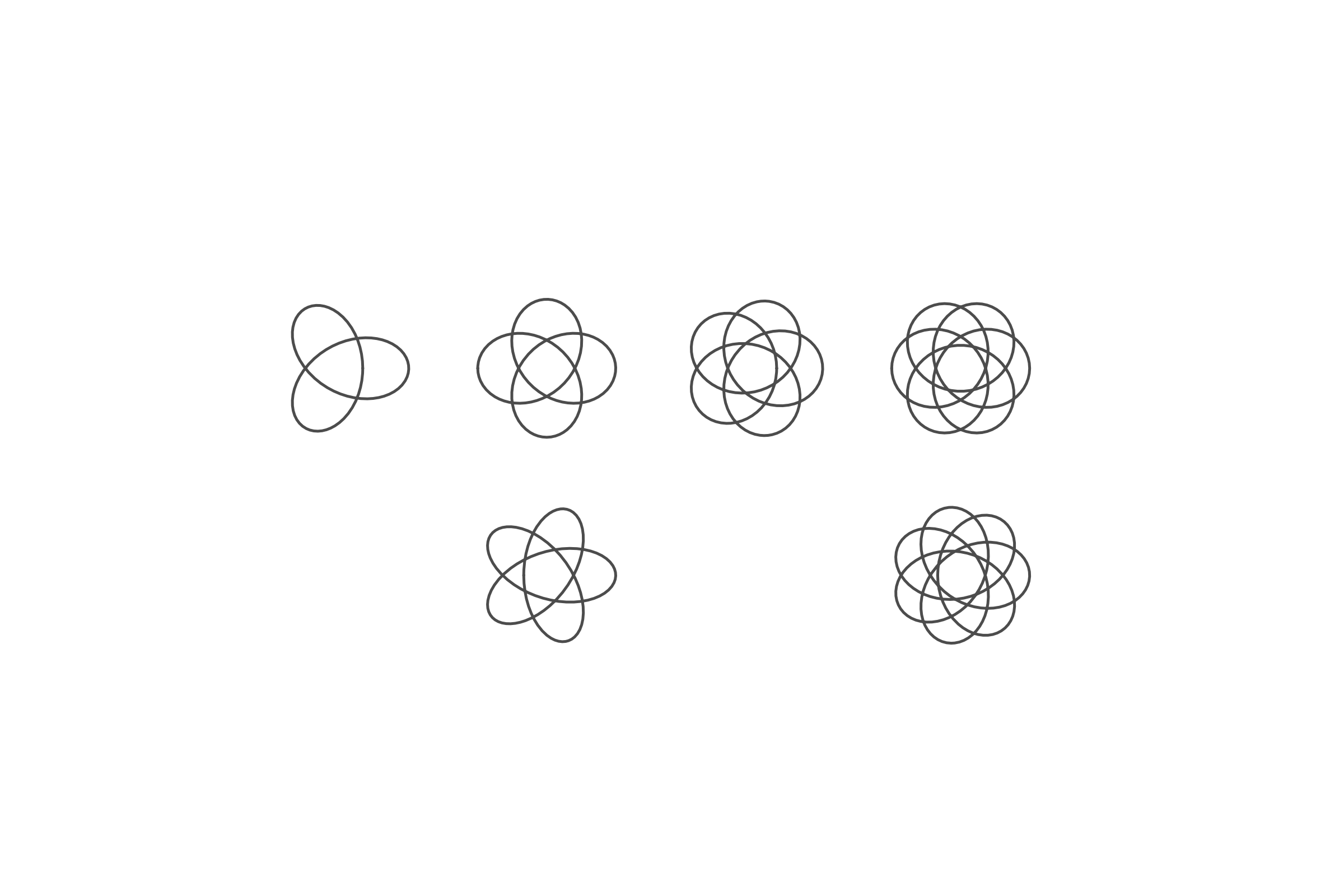}
\end{center}
\caption{\redsf{Expressive} hypotrochoids defined by the parametric equations
\eqref{def:epi-hypo1}--\eqref{def:epi-hypo2}.
Top row: $(b,c)=(2,1), (3,1), (4,1), (5,1)$.
Bottom row:
$(b,c)=(3,2), (5,2)$.
}
\vspace{-.2in}
\label{fig:6-hypotrochoids}
\end{figure}

\begin{figure}[ht]
\begin{center}
\includegraphics[scale=0.28, trim=17cm 28cm 14cm 33.5cm, clip]{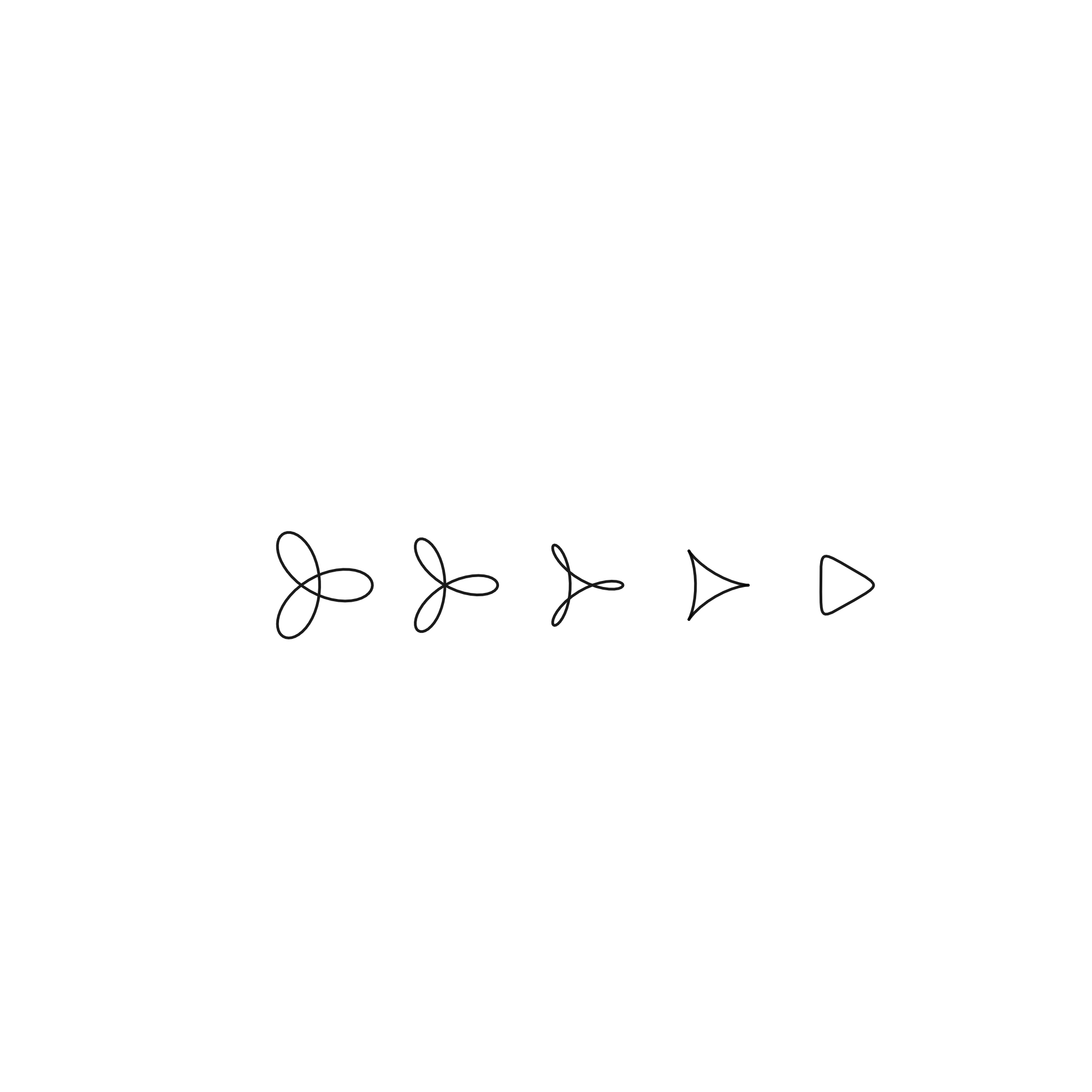}
\end{center}
\caption{\redsf{Hypotrochoids defined by the parametric equations
\eqref{def:epi-hypo1}--\eqref{def:epi-hypo2}
with $b=2$, $c=1$, $v=1$, and $u\in\{1.25, 1, 0.75, 0.5, 0.25\}$, shown left-to-right in this order.
The~first and the third hypotrochoids (with $u=1.25$ and $u=0.75$, respectively)
are expressive; the remaining three are not.}
}
\vspace{-.2in}
\label{fig:5-general-hypotrochoids}
\end{figure}

\begin{example}
\label{ex:hypotrochoids-again}
A hypotrochoid with (coprime) parameters $(b,c)$ and \redsf{suitably chosen ratio~$\frac{v}{u}$}
has $b+c$ ``petals,''
see Figure~\ref{fig:6-hypotrochoids}.  
\redsf{The number of petals can change as $\frac{v}{u}$ changes, cf.\ Figure~\ref{fig:5-general-hypotrochoids}.}

When $b=c+1$, we recover Example~\ref{ex:hypotrochoids},
cf.\ Figure~\ref{fig:3-petal} ($(b,c)=(2,1)$) and
Figure~\ref{fig:5-petal} ($(b,c)\!=\!(3,2)$).
\end{example}

\begin{figure}[ht]
\vspace{-.1in}
\begin{center}
\includegraphics[scale=0.25, trim=14cm 15.5cm 14cm 16cm, clip]{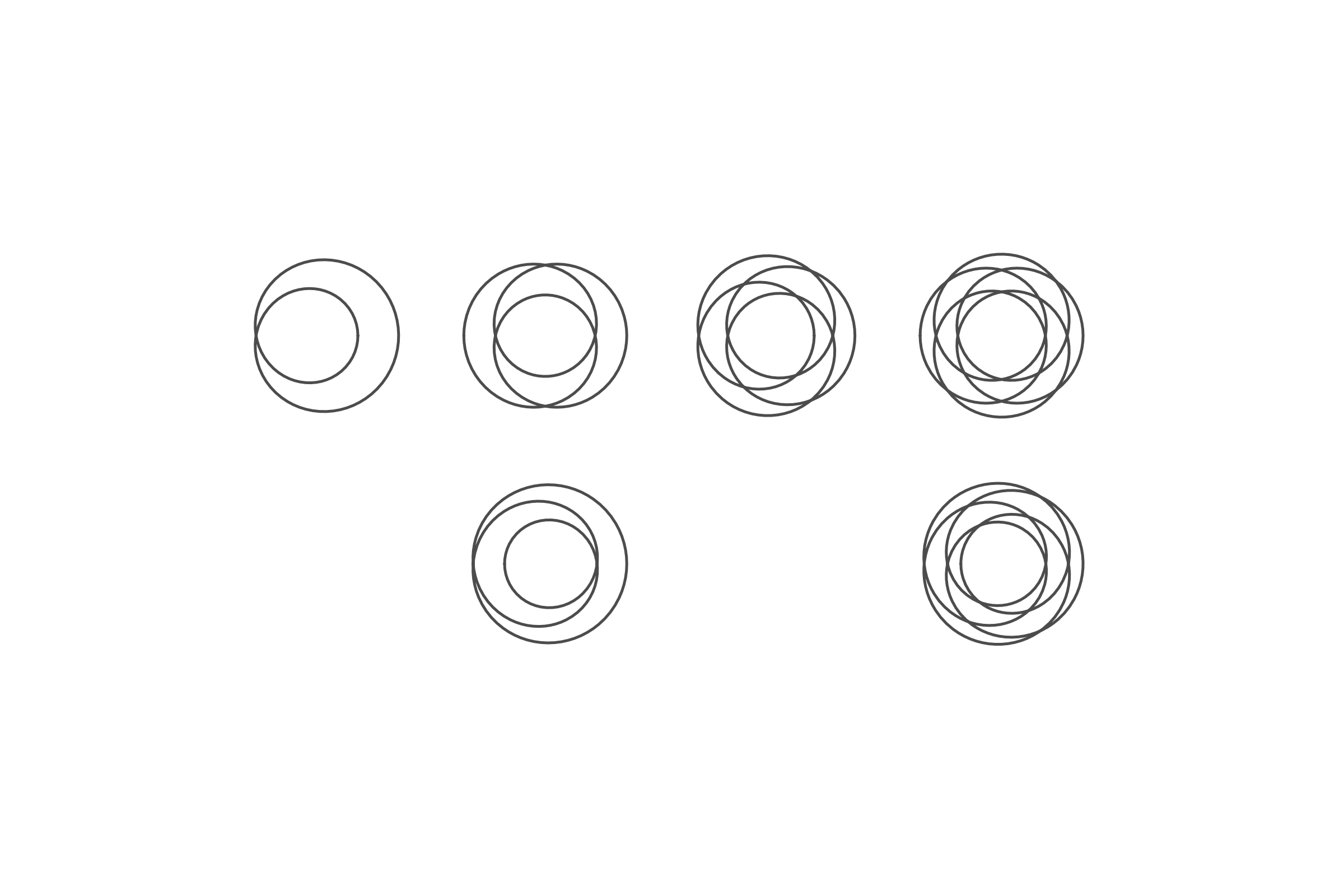}
\end{center}
\caption{Epitrochoids defined by the parametric equations
\eqref{def:epi-hypo1}--\eqref{def:epi-hypo2}.
Top row: $(b,c)=(2,-1), (3,-1), (4,-1), (5,-1)$.
Bottom row: $(b,c)=(3,2), (5,2)$.}
\vspace{-.2in}
\label{fig:6-epitrochoids}
\end{figure}

\begin{example}
An epitrochoid with (coprime) parameters $b$ and $c$ (here $b>-c>0$)
has $b+c$ inward-pointing ``petals,'' \linebreak[3]
see Figure~\ref{fig:6-epitrochoids}.

When $b=-c+1$, we recover the ``multi-lima\c cons'' of Example~\ref{ex:rose-curves}
and Figure~\ref{fig:multi-limacons} (with $b\in\{2,3,4\}$, $c=1-b$).  
\end{example}

\begin{proposition}
\label{pr:epi-hypo}
Let $C$ be an epitrochoid or hypotrochoid
given by \eqref{def:epi-hypo1}--\eqref{def:epi-hypo2}
(As~in Definition~\ref{def:hypotrochoids},
$b$ and $c$ are coprime integers, with $b>|c|$;
and \hbox{$u,v\in\RR^*$}.) \linebreak[3]
Then $C$ is expressive and $L^\infty$-regular 
if it has $(b+c)(b-1)$ hyperbolic nodes~in~$\AA^2$.
In~particular,  this holds if $\frac{|v|}{|u|}$ is sufficiently small.
\end{proposition}

\begin{proof}
Setting $\tau=e^{it}$, we convert the trigonometric parametrization \eqref{def:epi-hypo1}--\eqref{def:epi-hypo2} of~$C$ into a rational one:
\begin{align*}
x&=\ \ \,\tfrac{1}{2}(u\tau^{2b}+v\tau^{b+c}+v\tau^{b-c}+u), \\
y&=-\tfrac{i}{2}(u\tau^{2b}-v\tau^{b+c}+v\tau^{b-c}-u), \\
z&=\tau^b.
\end{align*}
This curve has two points at infinity, namely $(1,i,0)$ and $(1,-i,0)$,
corresponding to $\tau=0$ and $\tau=\infty$, respectively.
Noting that the line $x+iy=0$ passes through $(1,i,0)$,
we change the coordinates by replacing $y$ by
\begin{align*}
x+iy&=u\tau^{2b}+v\tau^{b-c}.
\end{align*}
In a neighborhood of $(1,i,0)$, this yields the following parametrization:
$$x=1,\quad x+iy=\tfrac{2v}{u}\tau^{b-c}+\text{h.o.t.},\quad z=\tfrac{2}{u}\tau^b+\text{h.o.t.}$$
Thus, at the point $(1,i,0)$ we have a semi-quasihomogeneous singularity of weight $(b-c,b)$,
and similarly at~$(1,-i,0)$.
The $\delta$-invariant at each of the two points is equal to $\frac{1}{2}(b-c-1)(b-1)$.
Hence in the affine plane~$\AA^2$, there remain
\[
\tfrac{1}{2}(2b-1)(2b-2)-(b-c-1)(b-1)=(b+c)(b-1)
\]
nodes. So if all these nodes are real hyperbolic (and distinct), then $C$ is expressive
and $L^\infty$-regular by Proposition~\ref{pr:irr-expressive}
(or Theorem~\ref{th:reg-expressive-irreducible}).

It remains to show that $C$ has $(b+c)(b-1)$ hyperbolic nodes
as long as $u,v\in\RR^*$ are chosen appropriately
(in particular, if $\frac{|v|}{|u|}$ is sufficiently small).
Consider the rational parametrization of $C$ given by
\begin{align*}
x+iy&=u\tau^b+v\tau^{-c},\\
x-iy&=u\tau^{-b}+v\tau^c, \\
z&=1.
\end{align*}
We need to show that all solutions $(\tau,\sigma)\in(\CC^*)^2$ of
\begin{align}
\label{eq:u(tau-theta1)}
&u\tau^b+v\tau^{-c}=u\sigma^b+v\sigma^{-c},\\
\label{eq:u(tau-theta2)}
&u\tau^{-b}+v\tau^c=u\sigma^{-b}+v\sigma^c,\\
\label{eq:tau-neq-sigma}
&\tau\neq\sigma
\end{align}
satisfy $|\tau|=|\sigma|=1$.
Let us rewrite \eqref{eq:u(tau-theta1)}--\eqref{eq:u(tau-theta2)} as
\begin{align*}
u(\tau^b-\sigma^b)-v(\tau\sigma)^{-c}(\tau^c-\sigma^c)&=0,\\
-u(\tau\sigma)^{-b}(\tau^b-\sigma^b)+v(\tau^c-\sigma^c)&=0.
\end{align*}
The condition $\gcd(b,c)=1$ implies that at least one of $\tau^b-\sigma^b$
and $\tau^c-\sigma^c$ is nonzero (or else $\tau=\sigma$).
Consequently
\[
\det\left(\begin{matrix} u & -v(\tau\sigma)^{-c} \\
-u(\tau\sigma)^{-b} & v\end{matrix}\right)
=uv(1-(\tau\sigma)^{-b-c})
=0,
\]
meaning that $\omega=\tau\sigma$ must be a root of unity: $\omega^{b+c}=1$.
With respect to~$\tau$ and~$\omega$,
the conditions \eqref{eq:u(tau-theta1)}--\eqref{eq:tau-neq-sigma}
become:
\begin{align}
\label{eq:omega}
&\omega^{b+c}=1, \\
\label{eq:tau-eps}
&u(\tau^b-\omega^b\tau^{-b})+v(\tau^{-c}-\omega^{-c}\tau^c)=0, \\
&
\label{eq:tau^2}
\tau^2\neq\omega.
\end{align}
Let $\lambda$ be a square root of~$\omega$, i.e., $\lambda^2=\omega$.
Then \eqref{eq:tau^2} states that $\tau\notin\{\lambda,-\lambda\}$.

For each of the $b+c$ possible roots of unity~$\omega$,
\eqref{eq:tau-eps}~is an algebraic equation of degree~$2b$ in~$\tau$.
We claim that if $u,v\in\RR^*$ are suitably chosen,
then all $2b$ solutions of this equation lie on $S^1=\{|\tau|=1\}$.
It is easy to see that this set of solutions contains the two values
$\tau=\pm\lambda$ which we need to exclude, leaving us with $2(b-1)$ solutions
for each~$\omega$ satisfying~\eqref{eq:omega}.
We also claim that all these $2(b+c)(b-1)$ solutions are distinct,
thereby yielding $(b+c)(b-1)$ hyperbolic nodes of the curve, as desired.

Let us establish these claims.
Denote $\eps=\lambda^{b+c}\in\{-1,1\}$.
Replacing $\tau$ by $\rho=\tau\lambda^{-1}$
(thus $\tau=\lambda\rho$), we transform~\eqref{eq:tau-eps} into
\[
u\eps(\rho^b-\rho^{-b})=v(\rho^c-\rho^{-c}).
\]
Making the substitution $\rho=e^{i\alpha}$, we translate this into
\begin{equation}
\label{eq:two-sinusoids}
u\eps \sin(b\alpha)=v\sin(c\alpha).
\end{equation}
If $|\frac{v}{u}|$ is sufficiently small, then the equation~\eqref{eq:two-sinusoids}
clearly has $2b$ distinct real solutions in the interval $[0,2\pi)$, as claimed.
Finally, it is not hard to see that all resulting $2(b+c)(b-1)$ values of~$\tau$ are distinct.
\end{proof}

We return to the general treatment of (potentially reducible) expressive $L^\infty$-regular curves.
First, a generalization of Proposition~\ref{pr:irr-expressive}:

\begin{proposition}
\label{pr:sufficient-expressive}
Let $C$ be a reduced real plane curve such that
\begin{itemize}[leftmargin=.2in]
\item
each component of $C$ 
is real, and either polynomial or trigonometric;
\item
the singular set of $C$ in the affine plane $\AA^2$ consists solely of hyperbolic nodes;
\item
the set of real points of $C$ in the affine plane is connected.
\end{itemize}
Then $C$ is expressive and $L^\infty$-regular.
\end{proposition}

\begin{proof}
The proof utilizes the approach used in the proof of Proposition~\ref{pr:irr-expressive}.
Arguing exactly as at the beginning of the latter proof,
we conclude that all we need to show is that the polynomial $F$ defining $C$ has finitely many critical points in the affine plane~$\AA^2$.
Once again, we argue by contradiction. Assuming that $Z(F_x)\cap Z(F_y)$ contains a real curve $B$ of a positive degree $d'<d=\deg(C)$, and reasoning as in the earlier proof, we conclude that for any component $C'$ of~$C$, there exists a curve $\hat B_{C'}\in\Span\{B,(L^\infty)^{d'}\}$ containing $C'$ as a component.
In view of the connectedness of~$\CR$ and the fact that different members of the pencil $\Span\{B,(L^\infty)^{d'}\}$ are disjoint from each other in the affine plane $\AA^2$,
we establish that there is just one curve $\hat B\in\Span\{B,(L^\infty)^{d'}\}$ that contains all the components of $C$---but this contradicts the inequality $d'<d$.
\end{proof}

The following theorem is the main result of this paper.

\begin{theorem}
\label{th:reg-expressive}
Let $C$ be a reduced real plane algebraic curve, with all irreducible components real.
The following are equivalent:
\begin{itemize}[leftmargin=.2in]
\item
$C$ is expressive and $L^\infty$-regular;
\item
each irreducible component of~$C$ is either trigonometric or polynomial, \\
all singular points of $C$ in the affine plane~$\AA^2$ are hyperbolic nodes, and \\
the set of real points of $C$ in the affine plane is connected.
\end{itemize}
\end{theorem}

\begin{proof}
Follows from Propositions~\ref{prop:nonempty} and~\ref{pr:sufficient-expressive}.
\end{proof}

\begin{corollary}
Let $C$ be an $L^\infty$-regular expressive plane curve
whose irreducible components are all real.
Let $C'$ be a subcurve of~$C$, i.e., a union of a subset of irreducible components.
If the set of real points of~$C'$ in the affine plane is connected,
then $C'$ is also $L^\infty$-regular and expressive.
\end{corollary}

\begin{proof}
Follows from Theorem~\ref{th:reg-expressive}.
\end{proof}

We conclude this section by a corollary whose statement is entirely elementary,
and in particular does not involve the notion of expressivity.

\begin{corollary}
\label{cor:all-crit-pts-are-real}
Let $C=V(G(x,y))$ be a real polynomial or trigonometric affine plane curve
which intersects itself solely at hyperbolic nodes.
Then all critical points of the polynomial $G(x,y)$ are real.
\end{corollary}

\begin{proof}
Immediate from Proposition~\ref{pr:irr-expressive} (or Theorem~\ref{th:reg-expressive-irreducible}).
\end{proof}


\section{More expressivity criteria}
\label{sec:more-expressivity-criteria}

We first discuss the irreducible case.
By Theorem~\ref{th:reg-expressive-irreducible}, an irreducible plane curve is
expressive and $L^\infty$-regular if and only if it is
is either trigonometric or polynomial, and all its singular points outside~$L^\infty$
are real hyperbolic nodes.
The last condition is the usually the trickiest to verify.

One simple case is when the number of hyperbolic nodes attains its maximum:

\begin{corollary}
\label{cor:max-number-of-nodes}
Let $C$ be a real polynomial or trigonometric curve of degree~$d$
with $\frac{(d-1)(d-2)}{2}$ hyperbolic nodes.
Then $C$ is expressive and $L^\infty$-regular.
\end{corollary}

\begin{proof}
In view of Hironaka's formula \eqref{eq:Hironaka}, the curve~$C$ has no other singular points
besides the given hyperbolic nodes.
The claim follows by Proposition~\ref{pr:irr-expressive}.
\end{proof}

Examples illustrating Corollary~\ref{cor:max-number-of-nodes}
include Lissajous-Chebyshev curves~\eqref{eq:Ta+Tb} with parameters $(d,d-1)$
as well as hypotrochoids with parameters \hbox{$(k,k-1)$} (cf.~\eqref{eq:hypotrochoid},
\redsf{with suitably chosen value of~$a$}).

\begin{remark}
In Corollary~\ref{cor:max-number-of-nodes},
the requirement that the curve $C$ is polynomial or trigonometric cannot be dropped.
For example, there exists an irreducible real quadric with three hyperbolic nodes which is not expressive.
\end{remark}

Corollary~\ref{cor:max-number-of-nodes} can be generalized as follows.

\begin{corollary}
\label{cor:nodes-Newton}
Let $C=V(G(x,y))$ be a real polynomial or trigonometric curve with $\nu$ hyperbolic nodes.
Suppose that the Newton polygon of $G(x,y)$ has $\nu$ interior integer points.
Then $C$ is expressive and $L^\infty$-regular.
\end{corollary}

\begin{proof}
It is well known (see \cite{Baker1893} or \cite[Section~4.4]{Fulton})
that the maximal possible number of nodes of an irreducible plane curve
with a given Newton polygon
(equivalently, the arithmetic genus of a curve in the linear system
spanned by the monomials in the Newton polygon, on the associated toric surface)
equals the number of interior integer points in the Newton polygon.
The claim then follows by Proposition~\ref{pr:irr-expressive}.
\end{proof}

Applications of Corollary~\ref{cor:nodes-Newton} include arbitrary irreducible
Lissajous-Chebyshev curves~\eqref{eq:Ta+Tb}.

\medskip

It is natural to seek an algorithm for verifying whether a given immersed
real polynomial or trigonometric curve~$C$, say one given by an explicit parametrization,
is expressive.
By Theorem~\ref{th:reg-expressive-irreducible}, this amounts to checking that
each point of self-intersection of~$C$ in the affine plane~$\AA^2$
corresponds to two real values of the parameter.
In the case of a polynomial curve
\[
t\mapsto (P(t),Q(t)),
\]
this translates into requiring that
\begin{itemize}[leftmargin=.2in]
\item
the resultant
(with respect to either variable $s$ or $t$) of the polynomials
\[
\widehat P(t,s)=\frac{P(t)-P(s)}{t-s} \quad \text{and}\quad
\widehat Q(t,s)=\frac{Q(t)-Q(s)}{t-s}
\]
has \redsf{simple} real roots and
\item
\redsf{the corresponding points $(P(t),Q(t))$ are simple (hyperbolic) nodes of~$C$}. 
\end{itemize}

\begin{figure}[ht]
\begin{center}
\includegraphics[scale=0.18, trim=18cm 24cm 20cm 18cm, clip]{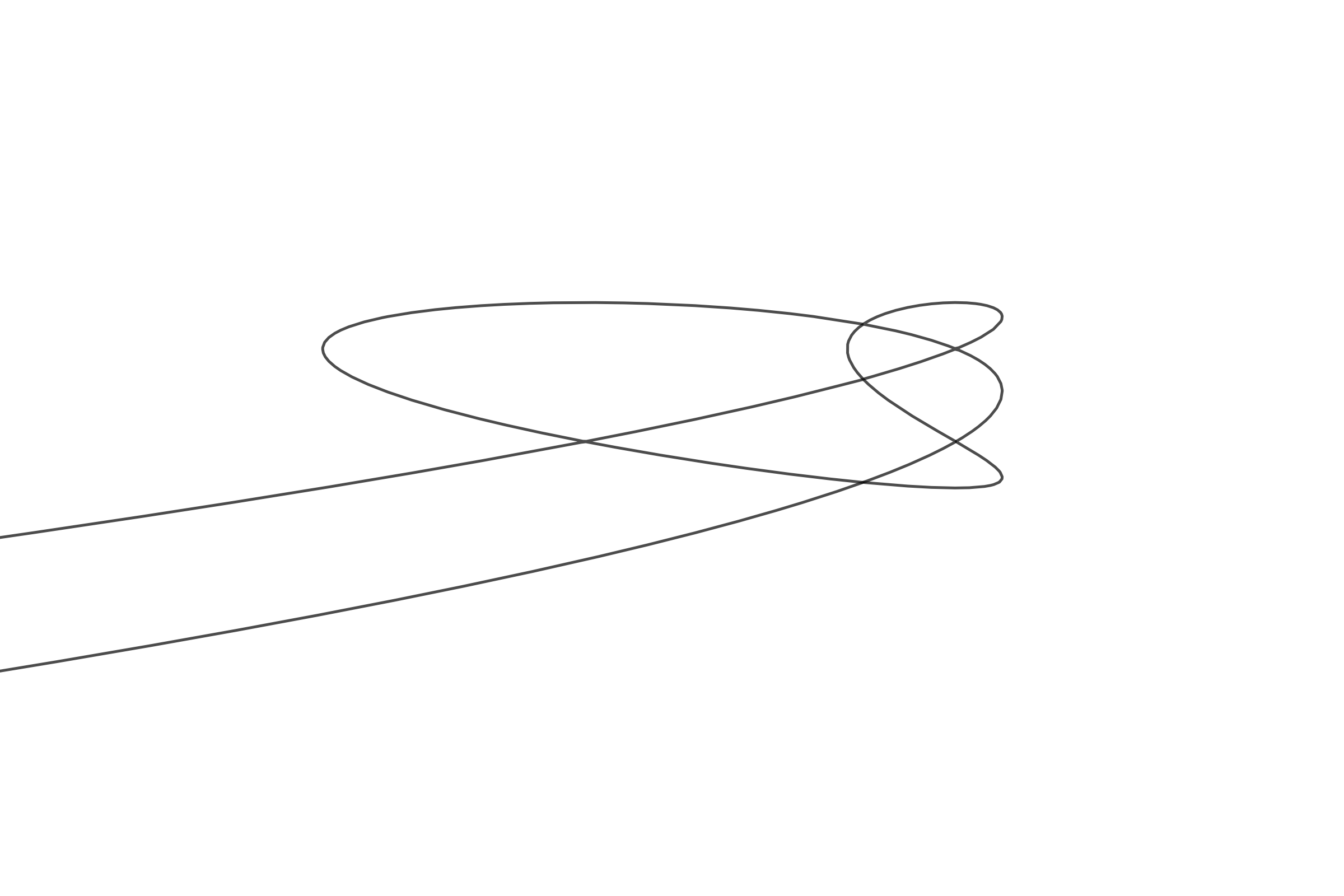}
\vspace{-.1in}
\end{center}
\caption{The curve $C=V(2x^3+3x^2-1+(4y^3-3y+\tfrac{x^2}{2})^2)$.}
\vspace{-.1in}
\label{fig:newton-degenerate}
\end{figure}

\begin{example}
Consider the sextic curve
\[
t\mapsto (-8t^{6}+24t^{4}+4t^{3}-18t^{2}-6t+1,-2t^{4}+4t^{2}-1).
\]
In this case,
\begin{align*}
\widehat P(t,s)&=-2 (2s^2 + 2st + 2t^2 - 3) (2s^3 + 2t^3 - 3s - 3t - 1),\\
\widehat Q(t,s)&=-2(s+t)(s^2+t^2-2).
\end{align*}
The resultant of $\widehat P$ and $\widehat Q$ (which we computed using \texttt{Sage})
is equal to
\[
-256 (2s^2 - 3)(2s^2 - 2s - 1)(2s^2 + 2s - 1)(8s^6 - 24s^4 - 4s^3 + 18s^2 + 6s - 1).
\]
All its 12 roots are real, so the curve is expressive, with 6 nodes.
See Figure~\ref{fig:newton-degenerate}.
\end{example}

The case of a trigonometric curve can be treated in a similar way.
Let $C$ be a trigonometric curve with a Laurent parametrization
\[
t\mapsto (x(t),y(t))=(P(t)+\overline P(t^{-1}),Q(t)+\overline Q(t^{-1}))
\]
(cf.\ \eqref{eq:param-PQ}), with $P(t)=\sum_k \alpha_k t^k$, $Q(t)=\sum_k \beta_k t^k$.
We write down the differences 
\begin{align*}
\frac{x(t)-x(s)}{t-s}
=\sum_k(t^{k-1}+t^{k-2}s+...+s^{k-1})\Bigl(\alpha_k-\frac{\overline\alpha_k}{t^k s^k}\Bigr),\\[-.05in]
\frac{y(t)-y(s)}{t-s}=\sum_k(t^{k-1}+t^{k-2}s+...+s^{k-1})\Bigl(\beta_k-\frac{\overline\beta_k}{t^k s^k}\Bigr),
\end{align*}
clear the denominators by multiplying by an appropriate power of~$ts$,
and require that all values of $t$ and $s$ for which the resulting polynomials vanish
(equivalently, all roots of their resultants) lie on the unit circle.

\begin{proposition}
Let $C\!=\!V(G)$ be an $L^\infty$-regular curve, with $G(x,y)\in\RR[x,y]$.
Assume that $\CR$ is connected, and each singular point of it
is a hyperbolic node, as in Definition~\ref{def:divide-polynomial}.
Let $\nu$ denote the number of such nodes,
and let $\rho$ be the number of regions of the divide~$D_G$.
Then $C$ is expressive if and~only~if
\begin{equation}
\label{eq:nu+rho}
\nu+\rho=(d-1)^2-\sum_{p\in C\cap L^\infty} \mu(C,p,L^\infty).
\end{equation}
\end{proposition}

\begin{proof}
By B\'ezout's theorem, the right-hand side of~\eqref{eq:nu+rho}
is the number of critical points of~$G$ in the affine plane, counted with multiplicities.
Since each region of $D_G$ contains at least one (real) critical point,
and each node of~$\CR$ is a critical point,
the only way for~\eqref{eq:nu+rho} to hold is for $C$ to be expressive.
\end{proof}

\newpage

\section{Bending, doubling, and unfolding}
\label{sec:constructions-irr}

In this section, we describe several transformations
which can be used to construct new examples of expressive curves from existing ones.
The simplest transformation of this kind is the ``bending'' procedure
based on the following observation:

\begin{proposition}
\label{pr:bending}
Let $f(x,y),g(x,y)\in\RR[x,y]$ be such that the map
\begin{equation}
\label{eq:real-biregular}
(x,y)\mapsto (f(x,y),g(x,y))
\end{equation}
is a biregular automorphism~of~$\AA^2$.
If the curve $C=V(G(x,y))$ is expressive, then so is the curve 
\[
\widetilde C=V(G(f(x,y),g(x,y))).
\]
If, in addition, $C$ is $L^\infty$-regular, with real components, then so is~$\widetilde C$.
\end{proposition}

\begin{figure}[ht]
\begin{center}
\includegraphics[scale=0.18, trim=18cm 19.5cm 22cm 21cm, clip]{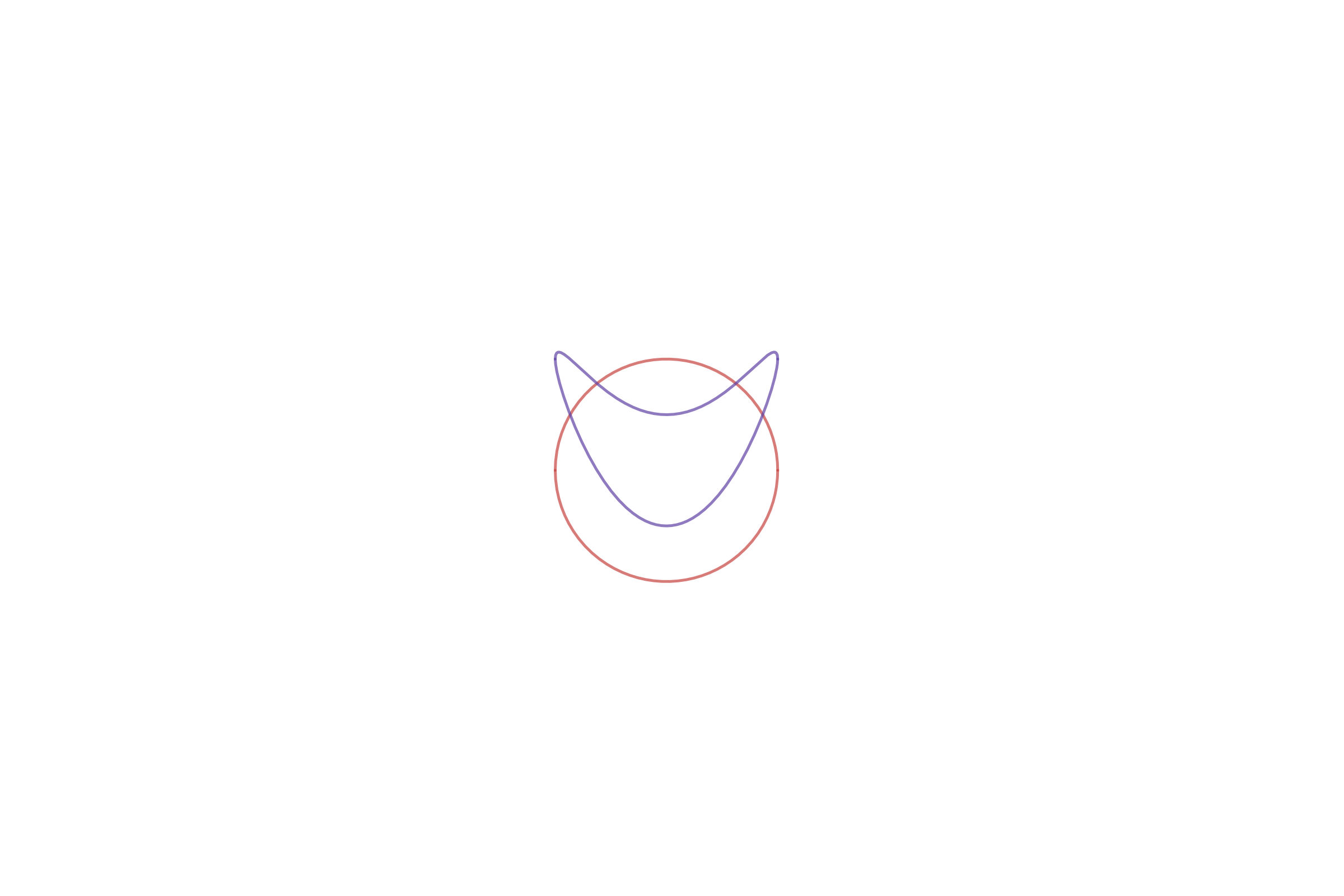}
\end{center}
\vspace{-.2in}
\caption{The curves $C=V(x^2+y^2-1)$ and $\widetilde C=V(x^2+(2y-2x^2)^2-1)$.}
\vspace{-.2in}
\label{fig:bending}
\end{figure}

Proposition~\ref{pr:bending} is illustrated in Figure~\ref{fig:bending}.

\begin{proof}
The automorphism~\eqref{eq:real-biregular} is an invertible change of variables that restricts
to a diffeomorphism of the real plane~$\RR^2$.
As such, it sends (real) critical points to (real) critical points,
does not change the divide of the curve,
and preserves expressivity.

Let us now assume that $C$ is expressive and $L^\infty$-regular, with real components.
In view of Theorem~\ref{th:reg-expressive}, all we need to show is that all components of~$\widetilde C$
are polynomial or trigonometric.
Geometrically, a polynomial (resp., trigonometric) component
is a Riemann sphere punctured at one real point (resp., two complex conjugate points)
and equivariantly immersed into the plane.
This property is preserved under real biregular automorphisms of~$\AA^2$.
\end{proof}

\begin{remark}
It is well known~\cite{vanderKulk} that the group of automorphisms of the affine plane
is generated by affine transformations together with the transformations of the form
$(x,y)\mapsto (x,y+P(x))$, for $P$ a polynomial.
This holds over any field of characteristic zero, in particular over the reals.
\end{remark}

\begin{example}
Several examples of bending can be extracted from Table~\ref{table:irreducible-expressive}.
Applying the automorphism $(x,y)\mapsto (x,y+x-x^m)$ to the line $V(x-y)$, we get
$V(x^m-y)$.
The automorphism $(x,y)\mapsto (x,y-x^2)$ transforms
the parabola $V(y^2-x)$ into $V((y-x^2)^2-x)$,
the ellipse $V(x^2+y^2-1)$ into $V(x^2+(y-x^2)^2-1)$,
and the nodal cubic $V(y^2-xy-x^3)$ into $V((y-x^2)^2-xy)$.
(In turn, $V(y^2-xy-x^3)$ and $V(4x^{3}-3x+2y^{2}-1)$ are related to each other
by an affine transformation.)
\end{example}

\begin{figure}[ht]
\begin{center}
\includegraphics[scale=0.2, trim=24cm 14cm 25cm 8cm, clip]{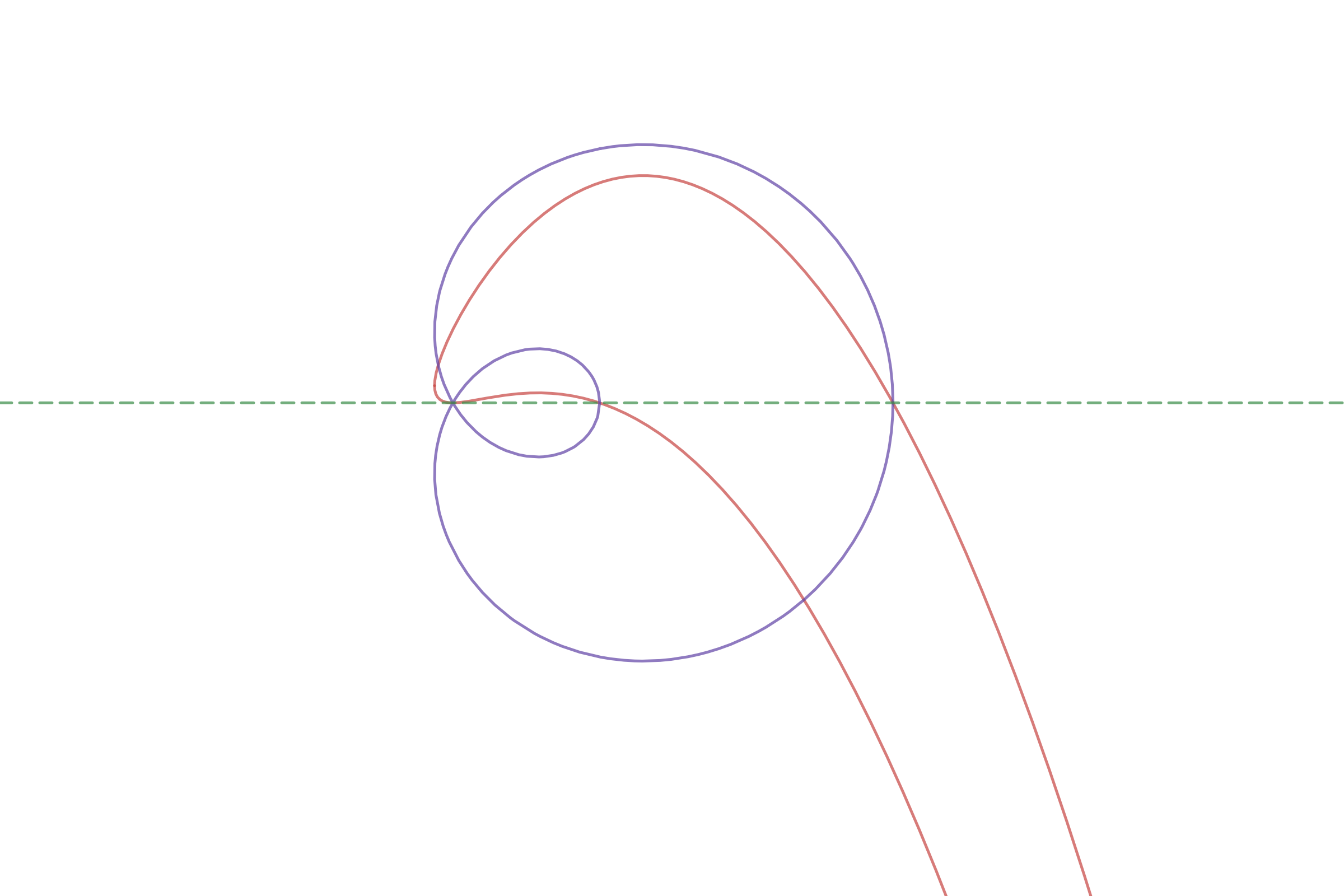}
\qquad\qquad
\includegraphics[scale=0.2, trim=30cm 15cm 29cm 2cm, clip]{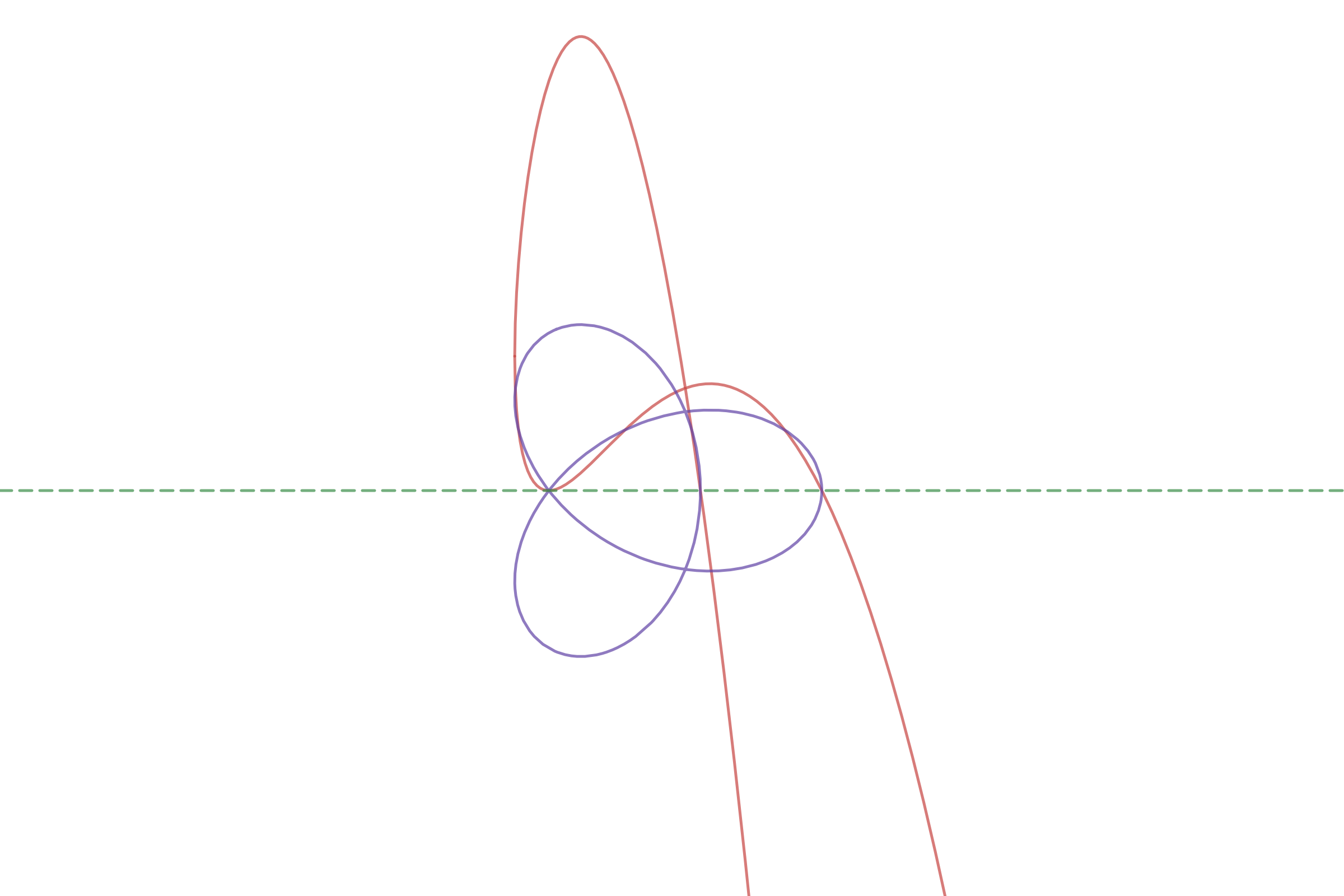}
\end{center}
\caption{A curve $C=V(G(x,y))$ and the ``doubled'' curve $\widetilde C=V(G(x,y^2))$.
On the left: $G(x,y^2)$ is the left-hand side of the equation~\eqref{eq:limacon}.
Cf.\ Figure~\ref{fig:multi-limacons}.
On the right: $G(x,y^2)$ is the left-hand side of~\eqref{eq:3-petal-alg}, with $a=2$.
Cf.\ Figure~\ref{fig:3-petal}.}
\label{fig:doubling-to-limacon-and-hypotrochoid}
\vspace{-.2in}
\end{figure}

\begin{example}
The polynomial expressive curves shown in Figure~\ref{fig:doubling-to-limacon-and-hypotrochoid}
(in red) are obtained by bending a parabola and
a nodal cubic.
\end{example}

We next discuss the ``doubling'' construction which transforms a plane curve \linebreak
\hbox{$C=V(G(x,y))$} into
a new curve $\widetilde C=V(G(x,y^2))$.
Proposition~\ref{pr:doubling} below shows~that, under certain conditions,
this procedure preserves expressivity.
See Figure~\ref{fig:doubling-to-limacon-and-hypotrochoid}.

\begin{proposition}
\label{pr:doubling}
Let $C$ be an expressive $L^\infty$-regular curve whose components are all real
(hence polynomial or trigonometric, see Proposition~\ref{prop:nonempty}).
Suppose that
\begin{itemize}[leftmargin=.2in]
\item
each component $B=V(G(x,y))$ of $C$, say with $\deg_x(G)=d$,
intersects $Z(y)$ in $d$ real points (counting multiplicities),
all of which are smooth points of~$C$; \redsf{moreover,
\begin{itemize}[leftmargin=.15in]
\item[{\rm $\circ$}]
if $B$ is trigonometric, then it intersects $Z(y)$ at $\frac{d}{2}$ points of quadratic tangency;
\item[{\rm $\circ$}]
if $B$ is polynomial, then it intersects $Z(y)$ at $\frac{d}{2}$, $\frac{d-1}{2}$ or $\frac{d-2}{2}$
points of quadratic tangency, with 0, 1, or 2 points of transverse intersection, respectively;
\end{itemize}
} 
\item
all nodes of $C$ lie in the real half-plane $\{y>0\}$;
\item
at each point of quadratic tangency between $C$ and $Z(y)$,
the local real branch of $C$ lies in the upper half-plane $\{y>0\}$;
\item
each connected component of the set $\{(x,y)\in C\cap\RR^2\mid y<0\}$
is unbounded.
\end{itemize}
Then the curve $\widetilde C=V(G(x,y^2))$ is expressive and $L^\infty$-regular.
\end{proposition}

\begin{proof}
The curve $\widetilde C$ is nodal by construction.
It is not hard to see that the set $\widetilde C_\RR$ is connected,
and all the nodes of $\widetilde C$ are hyperbolic.
(See the proof of Proposition~\ref{pr:accordion} for a more involved version of the argument.)
In view of Theorem \ref{th:reg-expressive}, it remains to show that
all components of $\widetilde C$ are real polynomial or trigonometric.

We only treat the trigonometric case, since the argument in the polynomial case is similar.
(Cf.\ also the proof of Lemma~\ref{lem:accordion2-poly}, utilizing such an argument in a more
complicated context.)
The natural map $\widetilde C\!\to \!C$ lifts to a two-sheeted ramified covering
$\rho:\widetilde C^{\vee}\!\to\! C^\vee$ between respective normalizations.
The restriction of $\rho$ to  $\widetilde C^{\vee}\!\setminus\!\rho^{-1}(Z(y))$
is an unramified two-sheeted covering,
and each component of $\widetilde C$ contains a one-dimensional
fragment of the real point set, hence is real.
In fact, $\rho$ is not ramified at all, since
each point in $C\cap Z(y)$ lifts to a node,
and hence to two preimages in $\widetilde C^{\vee}$.
Since $C^\vee=\CC^*$,
it follows that $\widetilde C^{\vee}$ is a union of at most two disjoint copies of~$\CC^*$.
We conclude that $\widetilde C$ consists of one or two trigonometric components.
\end{proof}

\begin{figure}[ht]
\begin{center}
\includegraphics[scale=0.13, trim=18cm 7cm 18cm 3cm, clip]{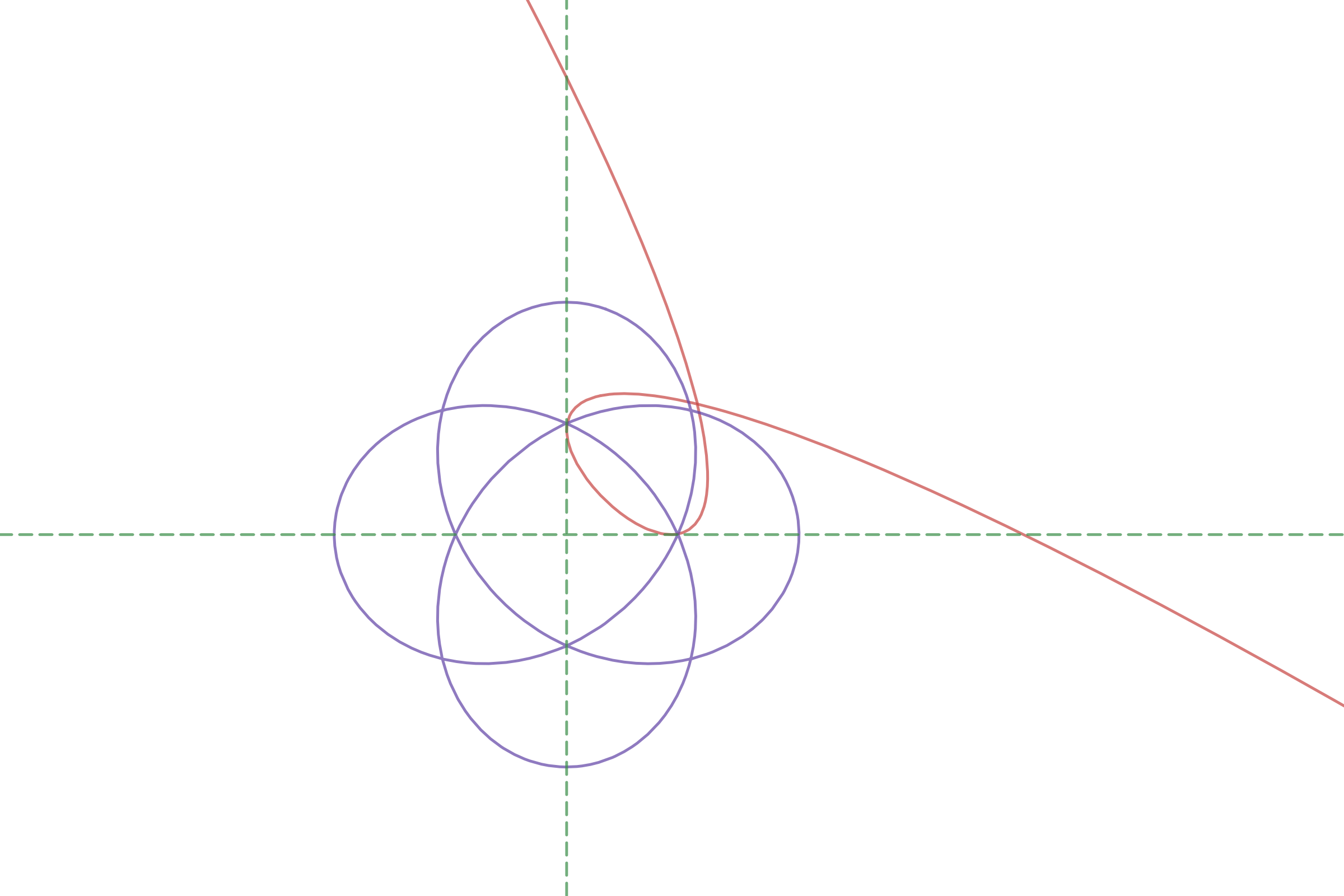}
\qquad
\includegraphics[scale=0.13, trim=18cm 16cm 18cm 3cm, clip]{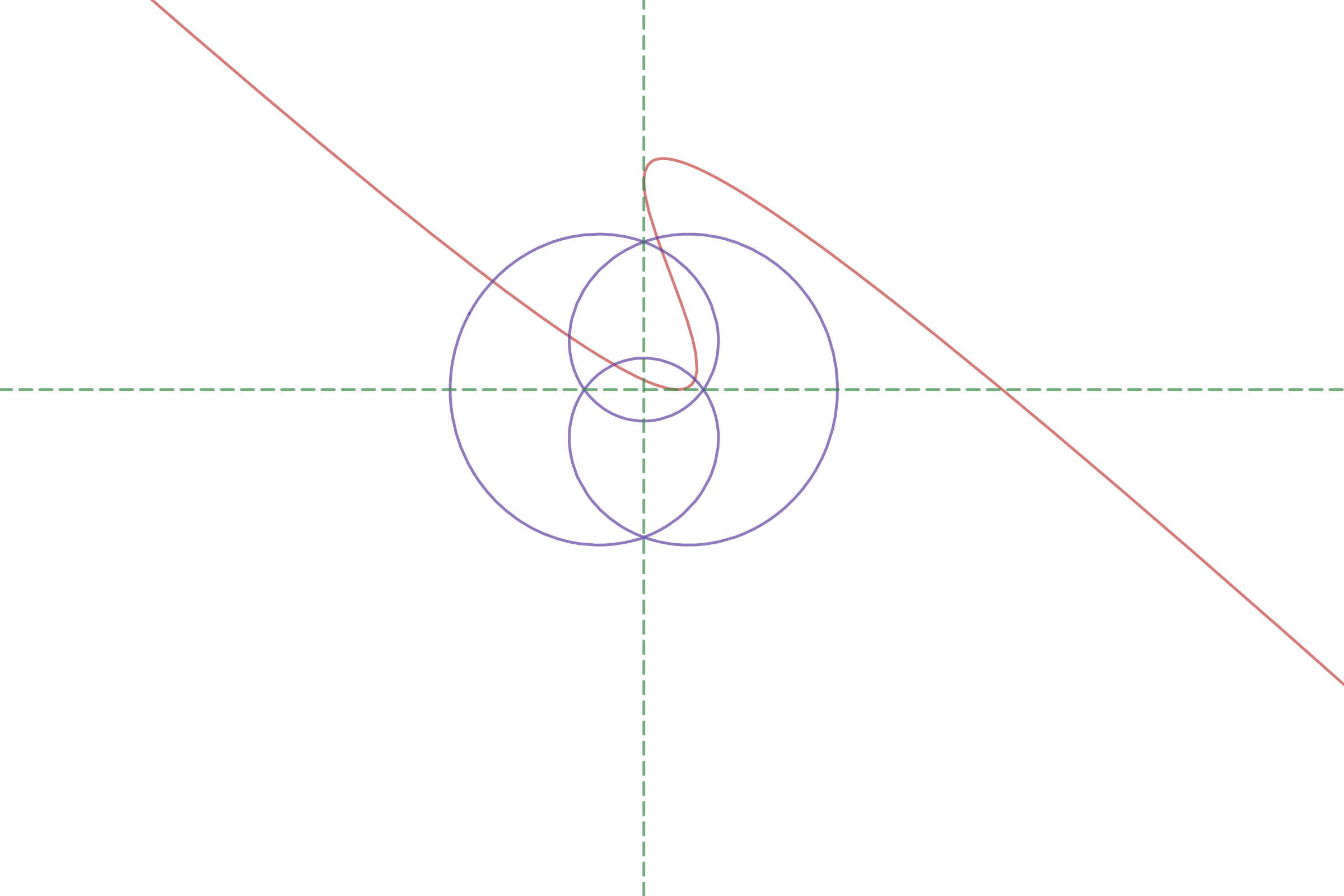}
\vspace{-.1in}
\end{center}
\caption{A curve $C=V(G(x,y))$ and its ``double-double'' $C'=V(G(x^2,y^2))$.
Left: $C$ is a nodal cubic tangent to both axes, $\widetilde C$ is a 4-petal hypotrochoid,
cf.\ Definition~\ref{def:hypotrochoids}, with $b=3$, $c=1$.
Right: $C$ is a cubic parabola tangent to both axes, $\widetilde C$ is an epitrochoid
with $b=3$, $c=-1$.
}
\vspace{-.1in}
\label{fig:double-reflection}
\end{figure}

If a curve $C=V(G(x,y))$ satisfies the conditions in Proposition~\ref{pr:doubling}
with respect to each of the coordinate axes $Z(x)$ and~$Z(y)$,
with all points in $C\cap Z(x)$ (resp.,~$C\cap Z(y)$)
located on the positive ray $\{x=0, y>0\}$ (resp., $\{y=0, x>0\}$),
the one can apply the doubling transformation twice,
obtaining an expressive curve $\widetilde C=V(G(x^2,y^2))$.
A~couple of examples are shown in Figure~\ref{fig:double-reflection}.

\begin{figure}[ht]
\begin{center}
\includegraphics[scale=0.25, trim=16cm 25.5cm 18cm 18cm, clip]{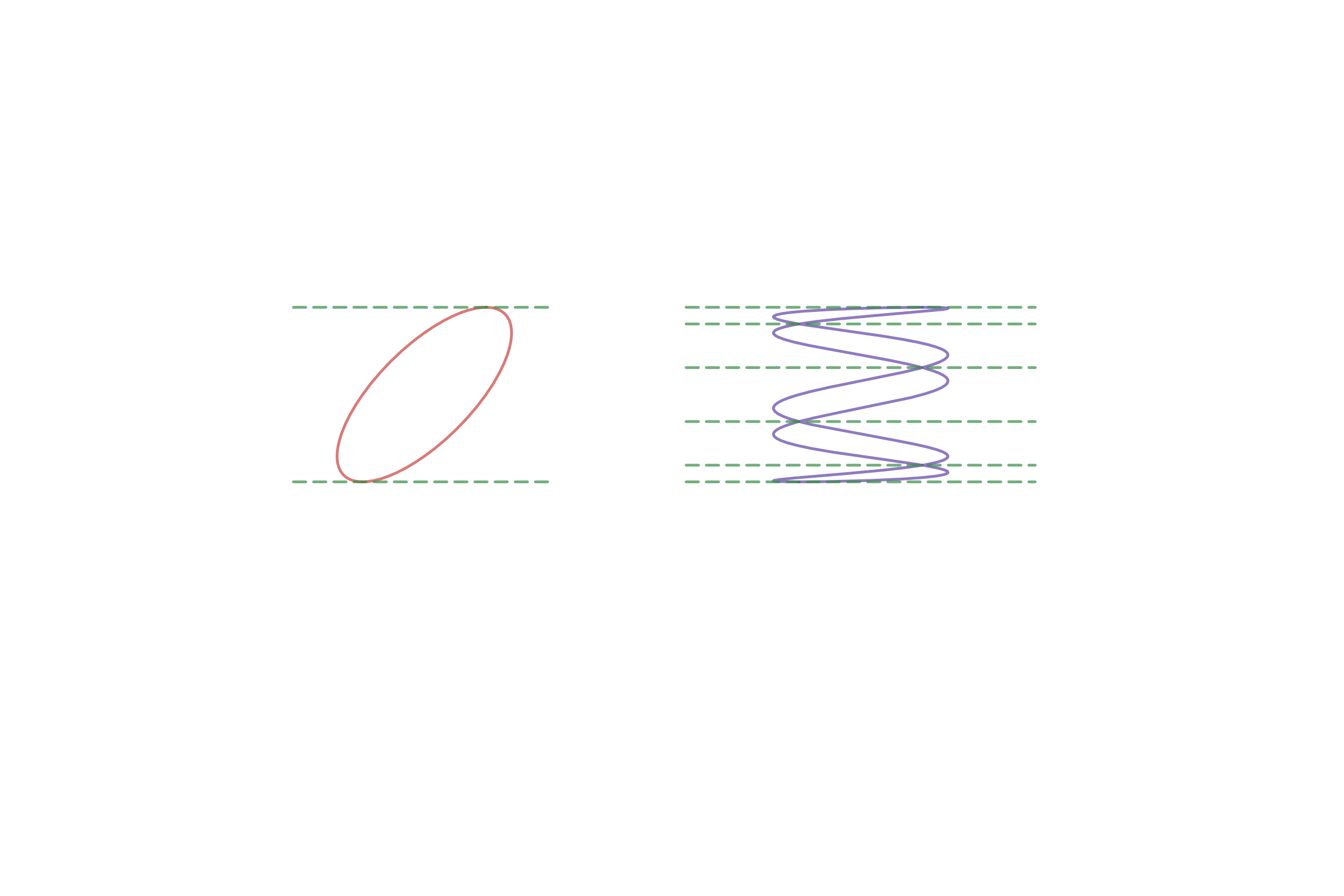}
\end{center}
\caption{An ellipse $C=V(G(x,y))$ and its unfolding $C=V(G(x,T_5(y)))$.
Here $G(x,y)=y^2-\sqrt{2}xy+x^{2}-\frac12$.
The green dashed lines are given by the equations $y\!=\!\pm 1$ (on the left)
and $T_5(y)\!=\!16y^5-20y^3+5y\!=\!\pm1$ (on the right).}
\vspace{-.2in}
\label{fig:ellipse-infolding}
\end{figure}

\begin{figure}[ht]
\begin{center}
\includegraphics[scale=0.3, trim=19cm 30cm 19cm 18cm, clip]{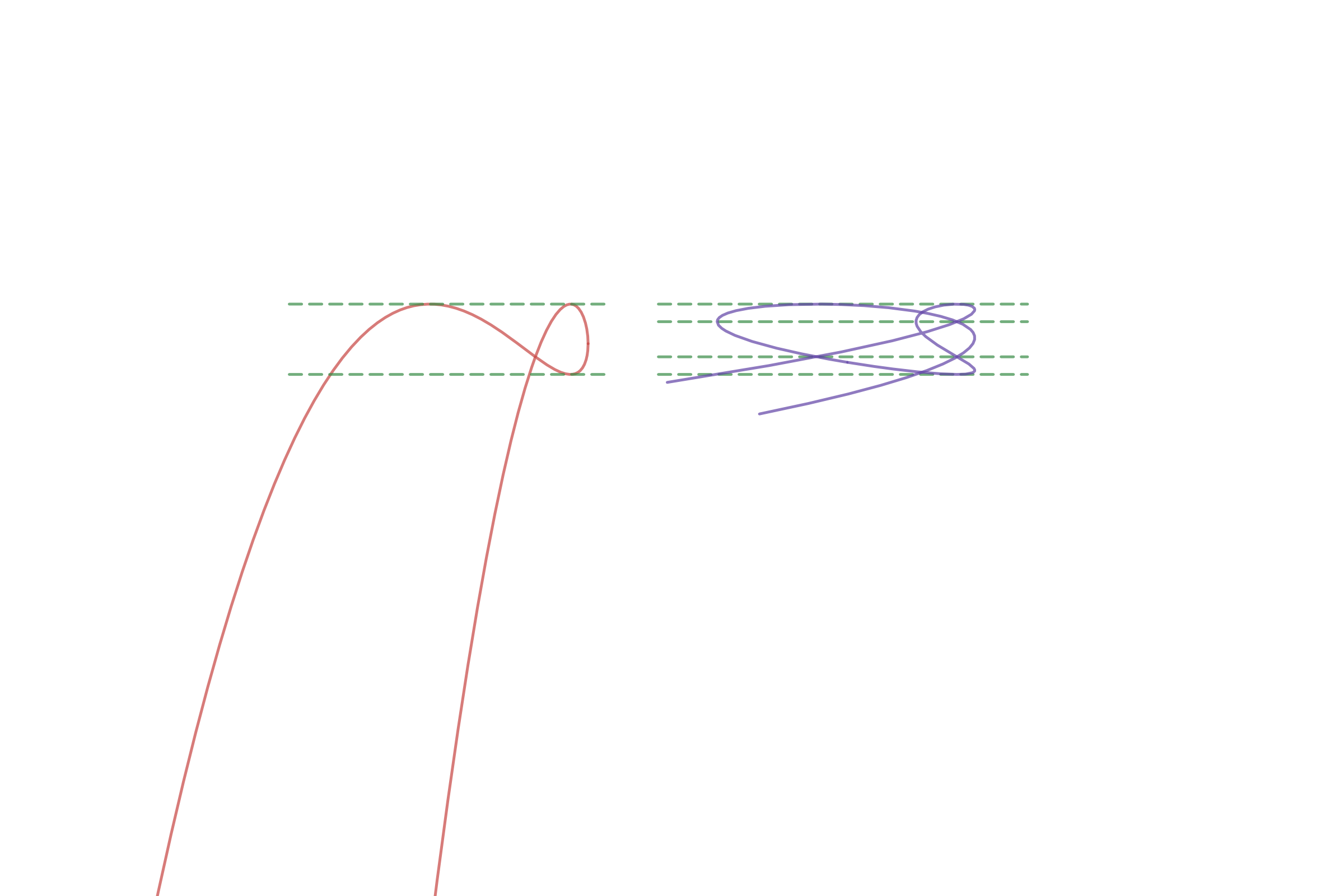}
\end{center}
\caption{A curve $C=V(G(x,y))$ and its triple unfolding $C=V(G(x,T_3(y)))$.
Here $G(x,y)=8x^{3}-12x^{2}+(2y+(x-1)^{2})^{2}$.
The green dashed lines are given by the equations $y=\pm 1$ (on the left)
and $T_3(y)=4y^3-3y=\pm1$. Cf.\ Figure~\ref{fig:newton-degenerate}.
}
\vspace{-.2in}
\label{fig:accordion-exotic}
\end{figure}

The remainder of this section is devoted to the discussion of ``unfolding.''
This is a transformation of algebraic curves that utilizes the coordinate change
\begin{equation}
(x,y)=(x,T_m(u)).
\label{eq:accordion}
\end{equation}
(As before, $T_m$ denotes the $m$th Chebyshev polynomial of the first kind,
see~\eqref{eq:chebyshev-poly}.)
A~precise description of unfolding is given in Proposition~\ref{pr:accordion} below.
To get a general idea of how this construction works,
take a look at the  examples in Figures~\ref{fig:ellipse-infolding}--\ref{fig:accordion-exotic}.
As these examples illustrate, the bulk of the unfolded curve (viewed up to an isotopy of the real plane)
is obtained by stitching together  $m$ copies of the input curve~$C$, or more precisely
the portion of~$C$ contained in the strip $\{-1\le y\le 1\}$.

\begin{lemma}
\label{lem:accordion1}
Let $C=V(F(x,y))$ be a trigonometric curve.
Suppose that
\begin{itemize}[leftmargin=.2in]
\item
the strip $\{-1<y<1\}$ contains a one-dimensional fragment of $\CR$;
\item
$C$ intersects each of the lines $y=\pm 1$ in $d/2$ points;
\item
all of these points are smooth points of
quadratic tangency between $C$ and $Z(y^2-1)$.
\end{itemize}
Then for any~$m\in\ZZ_{>0}\,$, the curve $\Cm$ defined by
\begin{equation}
\label{eq:Cm}
\Cm=V(F(x,T_m(u)))
\end{equation}
is a union of trigonometric curves.
\end{lemma}

\begin{proof}
The natural map $\Cm\to C$ given by~\eqref{eq:accordion} lifts to the $m$-sheeted
ramified covering map $\rho:\Cm^{\vee}\to C^\vee$ between the normalizations.
The restriction
\[
\rho:\Cm^{\vee}\setminus\rho^{-1}(Z(y^2-1))\to C^\vee\setminus Z(y^2-1)
\]
is an unramified $m$-sheeted covering, and each of the components of $\Cm$ contains a one-dimensional
fragment of the real point set, hence is real.
Let us show that $\rho$ is not ramified at all.
If $p\in C\cap Z(y-1)$ and $T_m^{-1}(1)$ consists of $a$ simple roots and $b$ double roots, then
$p$ lifts to $a$ smooth points (where $\Cm$ is quadratically tangent to the lines $u=\lambda$ with $\lambda$ running over
all these simple roots) and $b$ nodes, totaling $a+2b=m$ preimages in $\Cm^{\vee}$.

As $C$ is trigonometric, $C^\vee\!=\!\CC^*$.
Since a cylinder can only be covered by a cylinder,
$\Cm^{\vee}$~is a union of disjoint copies of $\CC^*$ (not necessarily $m$ of them),
so $\Cm$ is a union of trigonometric components.
\end{proof}

\pagebreak[3]

\begin{lemma}
\label{lem:accordion2-poly}
Let $C=V(F(x,y))$ be a real polynomial curve.
Let $d=\deg_x F(x,y)$.
Suppose that
\begin{itemize}[leftmargin=.2in]
\item
the strip $\{-1<y<1\}$ contains a one-dimensional fragment of $\CR$;
\item
$C$ intersects $Z(y^2-1)$ in $2d$ points
(counting multiplicities), all of which are smooth points of~$C$;
\item
these points include $d-1$ quadratic tangencies and two transverse intersections.
\end{itemize}
Then for any $m\in\RR_{>0}\,$, the curve $\Cm$ defined by~\eqref{eq:Cm}
is a union of polynomial or trigonometric components.
\end{lemma}

\begin{proof}
It is not hard to see that $C$ has
a polynomial parametrization $t\mapsto (P(t),Q(t))$ with \hbox{$Q(t)=\pm T_d(t)$}.
This follows from ``Chebyshev's equioscillation theorem'' of approximation theory
(due to E.~Borel and A.~Markov, see, e.g., \cite[Section~1.1]{Sodin-Yuditskii}
or~\cite[Theorem~3.4]{Mason-Handscomb}).
In~the rest of the proof, we assume that $Q(t)=T_d(t)$, as the case $Q(t)=-T_d(t)$
is completely similar.

We note that one can slightly vary the coefficients of $P$ while
keeping the intersection properties of $C$ with $Z(y^2-1)$
(and maintaining expressivity if $C$ has this property).
Thus, we can assume that $P(t)$ is a generic polynomial
(in particular, with respect to $Q(t)=T_d(t)$).

\textbf{Case~1}: $\gcd(d,m)=1$. Observe that
\begin{equation}
\label{eq:Cm-param-coprime}
\tau\mapsto (x,u)=(P(T_m(\tau)),T_d(\tau))
\end{equation}
is a parametrization of a polynomial curve lying inside~$\Cm$.
Indeed,
\[
F(P(T_m(\tau)),T_m(T_d(\tau)))=F(P(T_m(\tau)),T_d(T_m(\tau)))=0
\]
because $F(P(t),Q(t))=0$.

Since $\frac{dx}{d\tau}$ and $\frac{du}{d\tau}$ never vanish simultaneously
(thanks to the genericity of~$P$ and the coprimeness of $d$ and~$m$),
the map~\eqref{eq:Cm-param-coprime} is an immersion of $\CC$ into the affine plane.
Since $u=0$ at $d$ points, while
$\deg_x F(x,T_m(u))=d$, the image of this immersion is the entire curve~$\Cm$,
which is therefore polynomial.

\textbf{Case~2}: $c=\gcd(m,d)>1$.
Let $m=cr$, $d=cs$.
The curve $\Cm$ is given in the affine $(x,u)$-plane by the equations
\begin{equation*}
x=P(t),\quad T_m(u)=T_d(t)
\label{eq:Cm-parametric}
\end{equation*}
involving an implicit parameter~$t$.
Setting $t=T_r(\tau)$ we rewrite this as
$$x=P(T_r(\tau)),\quad T_m(u)
=T_m(T_s(\tau))
$$
(since $T_d(T_r(\tau))=T_m(T_s(\tau))$).
The equation $T_m(u)=T_m(u')$ has solutions $u=u'$, $u=-u'$ (for $m$ even) as well as
$$\arccos u=\pm \arccos u'+2\pi\tfrac{k}{m},\quad k=0,\dots,m-1, \quad u,u'\in[-1,1].
$$
From this, we obtain the following components of~$\Cm$,
all of which turn out to be either polynomial or trigonometric.
The polynomial components are:
\begin{align*}
&x=P(T_r(\tau)),\quad u=T_s(\tau), \\
&x=P(T_r(\tau)),\quad u=-T_s(\tau) \quad (m\in 2\ZZ).
\end{align*}
The trigonometric components are (here we set $\tau=\cos\theta$):
$$x=P(\cos(r\theta)),\quad
u=\cos(s\theta\pm2\pi\tfrac{k}{m}),
\quad 0<k<\tfrac{m}{2}.
$$
One can sort out which of these components are distinct by taking into account
that $x$ is invariant with respect to the substitutions $\theta\mapsto-\theta$
and $\theta\mapsto\theta+2\pi\frac{j}{r}$.

\redsf{Finally, we note that the curve $C_{(m)}$ has no multiple components. 
To see that, take a line $x=x_0$ such that the polynomial $F(x_0,y)$ has only simple roots $y_1,...,y_r$ which differ from $\pm1$. 
(Such a line $x=x_0$ does exist, since otherwise $C$ would have multiple components or contain one of the lines $y=\pm1$.) Then the polynomial $F(x_0,T_m(u))$ has only simple roots which all come from the equations $T_m(u)=y_i$, $i=1,...,r$.}
\end{proof}


\begin{proposition}
\label{pr:accordion}
Let $C=V(G(x,y))$ be an expressive $L^\infty$-regular curve all of whose components are real
(hence polynomial or trigonometric, cf.\ Proposition~\ref{prop:nonempty}).
Suppose~that
\begin{itemize}[leftmargin=.2in]
\item
each component $V(F(x,y))$ of $C$, say with $\deg_x(F)=d$,
intersects $Z(y^2-1)$ in $2d$ real points (counting multiplicities),
all of which are smooth points of~$C$;
\item
moreover, these intersections occur in one of the following two ways:
\begin{itemize}[leftmargin=.2in]
\item[{\rm $\circ$}]
$d$ points of quadratic tangency if the component is trigonometric, or
\item[{\rm $\circ$}]
$d-1$ points of quadratic tangency and two points of transverse intersection\\
 if the component is polynomial;
\end{itemize}
\item
all nodes of $C$ lie in the strip $\{-1<y<1\}$;
\item
at each point of quadratic tangency between $C$ and $Z(y^2-1)$,
the local real branch of $C$ lies in the strip $\{-1\le y\le 1\}$.
\end{itemize}
Then the curve $\Cm=V(G(x,T_m(y)))$ is expressive and $L^\infty$-regular.
\end{proposition}

\begin{proof}
By construction, the curve $\Cm$ is nodal.
By Lemmas \ref{lem:accordion1} and~\ref{lem:accordion2-poly}, the components of $\Cm$ are real polynomial or trigonometric.
In view of Theorem \ref{th:reg-expressive},
it remains to show that the set $C_{(m),\RR}$ is connected,
and that all the nodes of $\Cm$ are hyperbolic.

The set of real points $C_\RR$ is connected by Definition~\ref{def:expressive-poly}.
Since all the nodes of~$C$ lie inside the strip $\{-1<y<1\}$, it follows that
the set $C_\RR\cap\{-1\le y\le1\}$ is connected.

Let $\{a_1<\cdots<a_{m+1}\}$ be the set of roots of the equation $T_m^2(a)=1$.
Then each~set
\begin{equation}
\label{eq:substips}
C_{(m),\RR}\cap\{a_j\le y\le a_{j+1}\}, \quad
j=1,...,m,
\end{equation}
is an image of $C_\RR\cap\{-1\le y\le1\}$ under a homeomorphism of the strip
$\{-1\le y\le1\}$ onto the strip $\{a_j\le y\le a_{j+1}\}$.
Furthermore, for each $j=2,...,m$, the sets
\[
C_{(m),\RR}\cap\{a_{j-1}\le y\le a_j\}
\quad \text{and}\quad
C_{(m),\RR}\cap\{a_j\le y\le a_{j+1}\}
\]
are attached to each other at their common points along the line $Z(y-a_j)$.
(This set of common points is nonempty since it includes the images of the intersection of $C_\RR$
with one of the two lines $Z(y\pm 1)$.)
To obtain the entire set $C_{(m),\RR}\,$, we attach to the (connected) union of the $m$ sets~\eqref{eq:substips}
the diffeomorphic images of the intervals forming the set $C_\RR\setminus\{-1\le y\le 1\}$ (if any).
We conclude that $C_{(m),\RR}$ is connected.

Regarding the nodes of $\Cm$, we observe that they come in two flavours.
First, as one of the $m$ preimages of a node of $C$ contained in the strip $\{-1<y<1\}$;
all these preimages are real, hence hyperbolic.
Second, as a preimage of a tangency point between $C$ and $Z(y^2-1)$;
this again yields a hyperbolic node.
\end{proof}

\newpage

\section{Arrangements of lines, parabolas, and circles}
\label{sec:overlays}

We next discuss ways of putting together
several expressive curves to create a new (reducible) expressive curve.
Our key tool is the following corollary. 

\begin{corollary}
\label{cor:overlays}
Let $C_1,\dots,C_k$  be expressive and $L^\infty$-regular plane curves such that
\begin{itemize}[leftmargin=.2in]
\item
each pair $C_i$ and $C_j$ intersect each other in~$\AA^2$ at (distinct) hyperbolic nodes,  and
\item
the set $C_\RR$ of real points of the curve $C=C_1\cup \cdots C_k$ is connected.
\end{itemize}
Then $C$ is expressive and $L^\infty$-regular.
\end{corollary}

\begin{proof}
Follows from Theorem~\ref{th:reg-expressive}.
\end{proof}

One easy consequence of Corollary~\ref{cor:overlays} is the following construction.

\begin{corollary}[cf.\ Example~\ref{ex:two-graphs}]
\label{cor:overlays-graphs}
Let $f_1(x),\dots,f_k(x)\in\RR[x]$.
Assume that each poly\-nomial $f_i(x)-f_j(x)$ has real roots,
and all such roots (over all pairs $\{i,j\}$) are pairwise distinct.
Then the curve $V(\prod_i (y-f_i(x)))$ is expressive and $L^\infty$-regular.
\end{corollary}

\begin{figure}[ht]
\begin{center}
\includegraphics[scale=0.25, trim=15cm 23cm 14cm 12cm, clip]{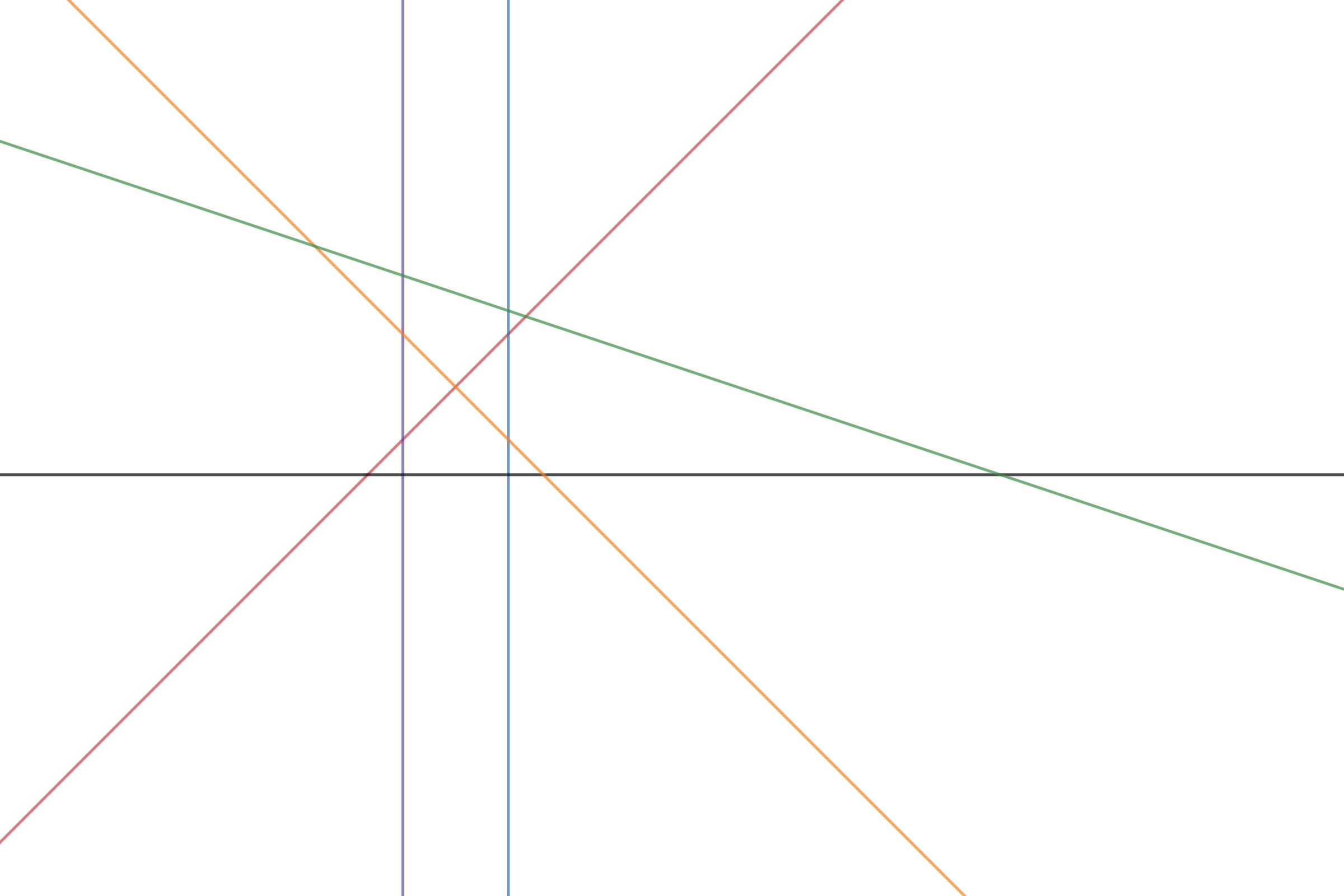}
\vspace{-.1in}
\end{center}
\caption{A line arrangement. 
}
\vspace{-.2in}
\label{fig:six-lines}
\end{figure}

\begin{example}[Line arrangements]
\label{ex:line-arrangements}
A \redsf{connected} arrangement of distinct real lines in the plane forms an expressive and $L^\infty$-regular curve,
as long as no three lines intersect at a point.
Parallel lines are allowed. See Figure~\ref{fig:six-lines}.
\end{example}

\begin{figure}[ht]
\begin{center}
\medskip
\includegraphics[scale=0.32, trim=21cm 18cm 21cm 24cm, clip]{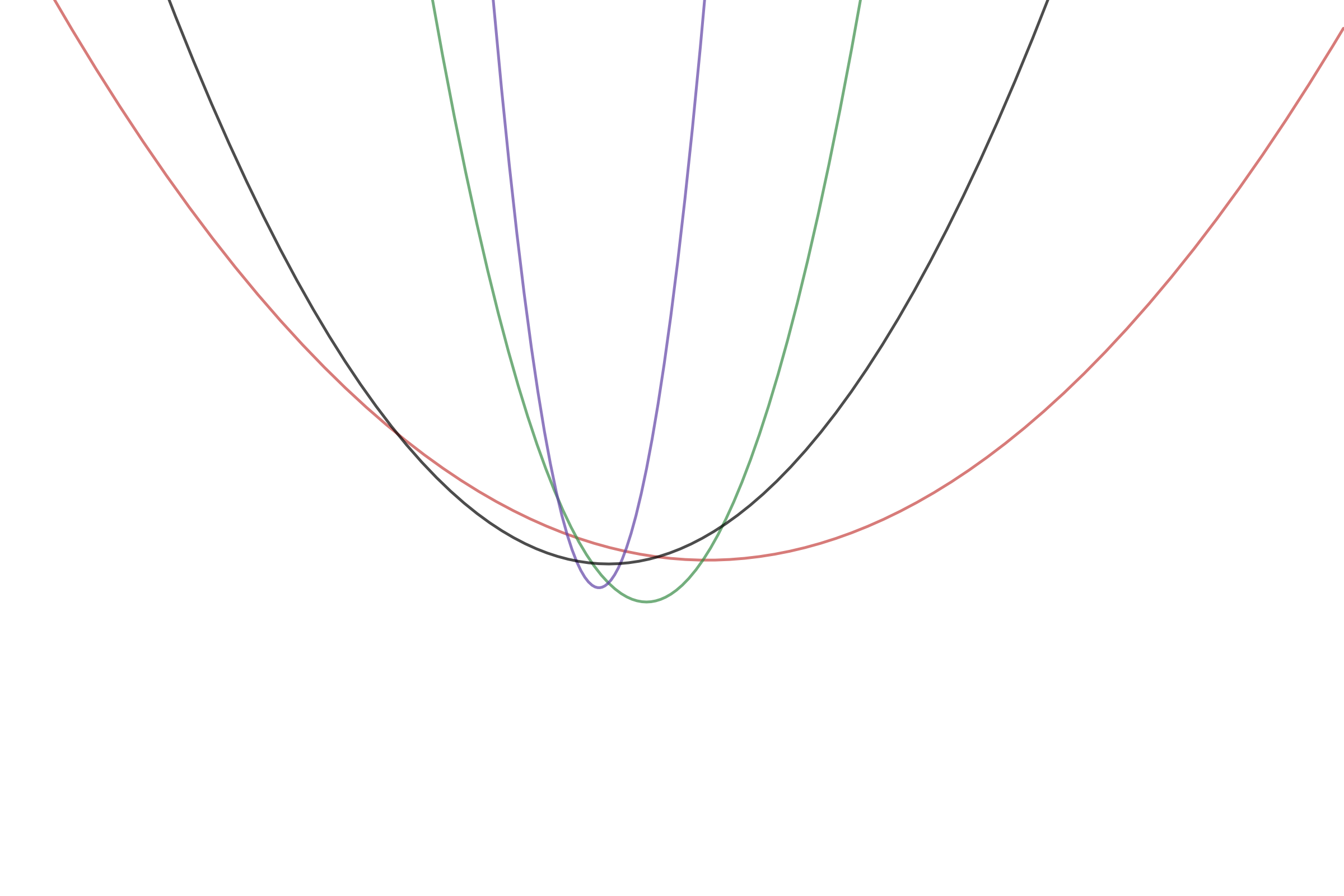}
\vspace{-.1in}
\end{center}
\caption{An arrangement of four co-oriented parabolas.
Each pair of parabolas intersect trans\-versally at two points.
}
\label{fig:four-parabolas}
\end{figure}

\begin{figure}[ht]
\begin{center}
\medskip
\includegraphics[scale=0.3, trim=25cm 29cm 17cm 11cm, clip]{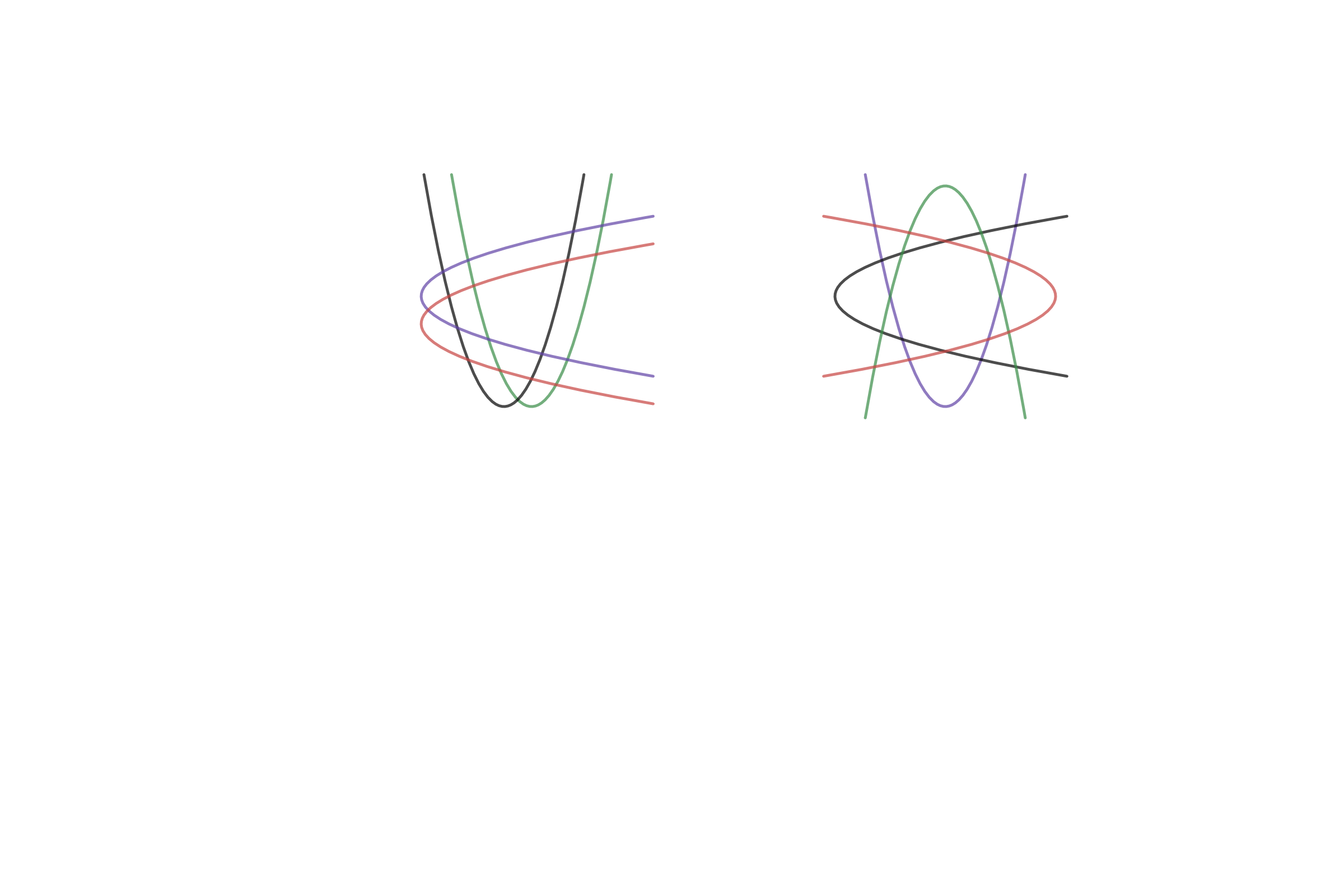}
\vspace{-.1in}
\end{center}
\caption{Two arrangements of four parabolas.
Each pair of parabolas intersect trans\-versally at one, two, or four points, all of them real.
}
\label{fig:two-4-parabolas}
\end{figure}

\begin{example}[Arrangements of parabolas]
\label{ex:arrangements-of-parabolas}
Let $C$ be a union of distinct real parabolas in the affine plane.
Then $C$ is an expressive and $L^\infty$-regular curve provided
\begin{itemize}[leftmargin=.2in]
\item
the set of real points of $C$ is connected;
\item
no three parabolas intersect at a point;
\item
all intersections between parabolas are transverse;
\item
for each pair of parabolas $P_1$ and $P_2$, one of the following options holds:
\begin{itemize}[leftmargin=.2in]
\item
$P_1$ and $P_2$ differ by a shift of the plane, or
\item
$P_1$ and $P_2$ have parallel (or identical) axes, and intersect at 2~points, or
\item
$P_1$ and $P_2$ intersect at 4 points.
\end{itemize}
\end{itemize}
See Figures~\ref{fig:four-parabolas} and~\ref{fig:two-4-parabolas}.
\end{example}

Examples~\ref{ex:line-arrangements} and~\ref{ex:arrangements-of-parabolas}
have a common generalization:

\begin{example}[Arrangements of lines and parabolas]
\label{ex:lines+parabolas}
Let $C$ be a union of distinct real lines and parabolas in the affine plane.
Then $C$ is an expressive and $L^\infty$-regular curve provided
\begin{itemize}[leftmargin=.2in]
\item
the set of real points of $C$ is connected;
\item
no three of these curves intersect at a point;
\item
all intersections between these lines and parabolas are transverse;
\item
each line intersects every parabola at one or two points;
\item
each pair of parabolas intersect in one of the ways listed in Example~\ref{ex:arrangements-of-parabolas}.
\end{itemize}
\end{example}

\begin{figure}[ht]
\begin{center}
\includegraphics[scale=0.26, trim=21cm 17.5cm 21cm 11cm, clip]{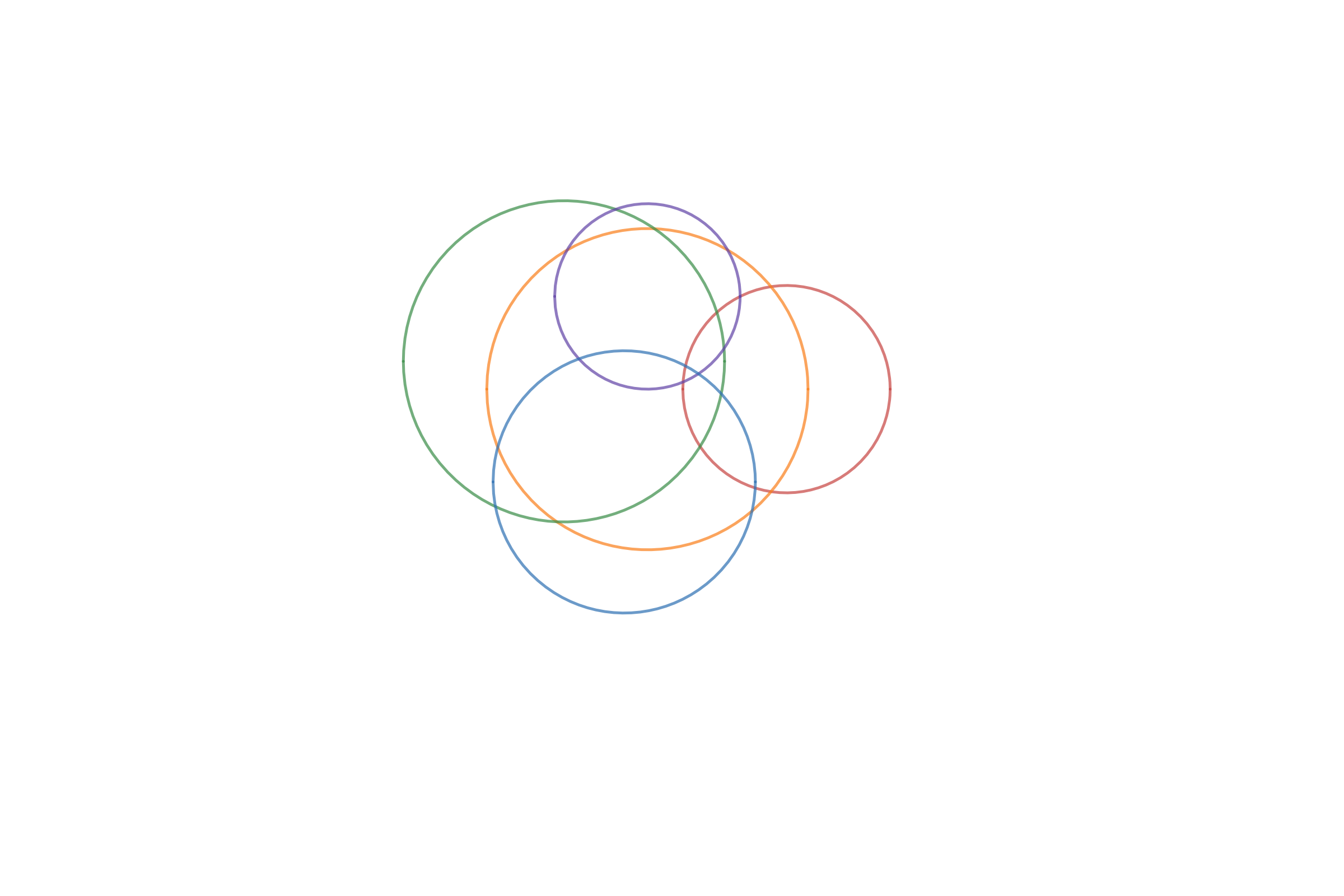}
\vspace{-.1in}
\end{center}
\caption{A circle arrangement.}
\vspace{-.1in}
\label{fig:five-circles}
\end{figure}

Another elegant application of Corollary~\ref{cor:overlays} involves
arrangements of circles:

\begin{example}[Circle arrangements]
\label{ex:circle-arrangements}
Let $\{C_i\}$ be a collection of circles on the real affine plane
such that each pair of circles intersect transversally at two real points,
with no triple intersections.
Then the curve $\bigcup C_i$ is expressive and $L^\infty$-regular.
See Figure~\ref{fig:five-circles}.
\end{example}

\begin{figure}[ht]
\begin{center}
\includegraphics[scale=0.26, trim=21cm 17.5cm 21cm 10cm, clip]{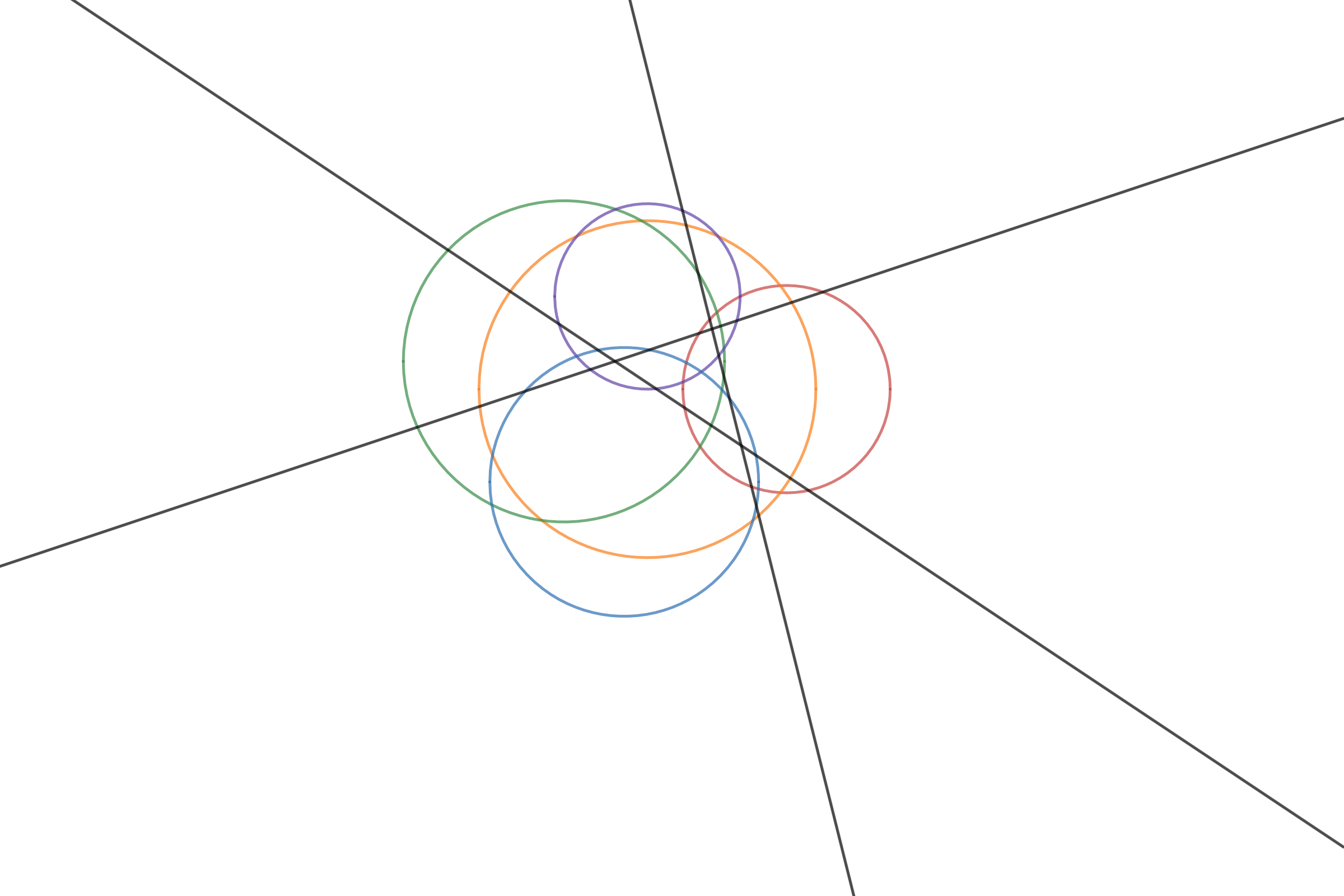}
\end{center}
\caption{An arrangement of circles and lines.}
\vspace{-.2in}
\label{fig:five-circles+three-lines}
\end{figure}

Here is a common generalization of Examples~\ref{ex:line-arrangements}
and~\ref{ex:circle-arrangements}:

\begin{example}[Arrangements of lines and circles]
Let $\{C_i\}$ be a collection of lines and circles on the real affine plane
such that each circle intersects every line
(resp., every other circle) transversally at two points,
with no triple intersections.
Then the curve $\bigcup C_i$ is expressive and $L^\infty$-regular.
See Figure~\ref{fig:five-circles+three-lines}.
\end{example}

\newpage

\section{Shifts and dilations}
\label{sec:shifts+dilations}

In this section, we obtain lower bounds for an intersection multiplicity
(at a point $p\in L^\infty$)
between a plane curve~$C$
and another curve obtained from~$C$ by a shift or dilation.
In Sections~\ref{sec:arrangements-poly}--\ref{sec:arrangements-trig},
we will use these estimates to derive expressivity criteria for unions of polynomial or trigonometric curves.

Without loss of generality, we assume that $p=(1,0,0)$ throughout this section.


To state our bounds, we will need some notation involving Newton diagrams:

\begin{definition}
\label{def:Newton-Gamma}
Let $C=Z(F(x,y,z))$ be a plane projective curve that contains neither of the lines $Z(z)=L^\infty$
and $Z(y)$ as a component.
Furthermore assume that $p=(1,0,0)\in C\cap L^\infty$.
We denote by $\Gamma(C,p)$ the Newton diagram of the polynomial
\begin{equation}
\label{eq:F1}
F_1(y,z)=F(1,y,z)\in\CC[y,z]
\end{equation}
at the point $(0,0)$, see Definition~\ref{def:Newton-diagram}.
Since $F_1$ is not divisible by $y$ or~$z$,
the Newton diagram $\Gamma(C,p)$ touches both coordinate axes.
We denote by $S^-(\Gamma(C,p))$ the area of the domain bounded by $\Gamma(C,p)$ and these axes.
\end{definition}

\begin{proposition}
\label{pr:gen-dilation}
Let $C=Z(F(x,y,z))$ be a projective curve that contains neither $Z(z)=L^\infty$ nor~$Z(y)$ as a component.
Let $c\in\CC^*$, and let $C_c$ denote the dilated~curve
\begin{equation}
\label{eq:dilX}
C_c=Z(F(cx,cy,z)).
\end{equation}
Assume that $p=(1,0,0)\in C\cap L^\infty$.
Then
\begin{equation}
(C\cdot C_c)_p\ge2S^-(\Gamma(C,p)).
\label{est-dilation}
\end{equation}
\end{proposition}

\begin{proof}
In the coordinates $(y,z)$, the dilation $C\leadsto C_c$ can be regarded as a transition from the polynomial $F_1$ (see~\eqref{eq:F1}) to the polynomial $F_c(y,z)=F(c,cy,z)$.
The polynomials $F_1$ and $F_c$ have the same Newton diagram at~$p$
(resp., Newton polygon).
Furthermore, let $G(y,z)$ be a polynomial with the same Newton polygon
and with generic coefficients.
By the lower semicontinuity of the intersection number, we have
$$(C\cdot C_c)_p\ge(C\cdot Z(G(y,z)))_p.$$
By Kouchnirenko's theorem \cite[1.18, Th\'eor\`eme~III${}'$]{Kouchnirenko},
the total intersection multiplicity of the curves $C$ and $Z(G)$ in the torus $(\CC^*)^2$
equals $2S(P)$, twice the area of the Newton polygon $P$ of~$F_1$.
Let us now deform $F_1$ and $G$ by adding all monomials underneath the Newton diagram
$\Gamma(C,p)$, with sufficiently small generic coefficients.
Again by Kouchnirenko's theorem, the total intersection multiplicity
of the deformed curves in $(\CC^*)^2$ equals $2S(P)+2S^-(\Gamma(C,p))$.
To establish the bound (\ref{est-dilation}), it remains to notice
that the extra term $2S^-(\Gamma(C,p))$ occurring in the latter intersection multiplicity
geometrically comes from simple intersection points in a neighborhood of~$p$,
obtained by breaking up the (complicated) intersection of $C$ and $Z(G)$ at~$p$.
\end{proof}

To state the analogue of Proposition~\ref{pr:gen-dilation} for shifted curves,
we need to recall some basic facts about the \emph{Newton-Puiseux algorithm}
\cite[Algorithm~I.3.6]{GLS}.
This algorithm assigns each local branch~$Q$ of a curve $C$
at the point $p=(1,0,0)\in C\cap L^\infty$
to an edge $E\!=\!E(Q)$ of the Newton diagram~$\Gamma(C,p)$,
cf.\ Definition~\ref{def:Newton-Gamma}.
We denote~by
\[
\nE=(n_y,n_z)\in\ZZ_{>0}^2
\]
the primitive integral normal vector to~$E$, with positive coordinates.

\pagebreak[3]


\begin{lemma}[{cf.\ \cite[Section~I.3.1]{GLS}}]
\label{l:wzwy}
Let $Q$ be a local branch of $C$ at \hbox{$p\!=\!(1,0,0)\!\in\! C$}.
Assume that $Q$ is tangent to~$L^\infty$.
Let $E=E(Q)$.
Then
\begin{equation}
\label{eq:nE=}
\nE=(n_y,n_z)=\tfrac{1}{r} (m,d)
\end{equation}
where
\begin{align}
\label{eq:d(Q)}
d=d(Q)&=(Q\cdot L^\infty)_p\,,\\
\label{eq:m(Q)}
m=m(Q)&=\mt(Q)<d,\\
\label{eq:r(Q)}
r=r(Q)&=\gcd(m,d).
\end{align}
\end{lemma}

We denote by $F^E(y,z)=F^{E(Q)}(y,z)$ the truncation of $F_1$ (see~\ref{eq:F1})
along the edge $E=E(Q)$.
The polynomial $F^E$ is quasihomogeneous
with respect to the weighting of $y$ and $z$ by $n_y$ and $n_z$, respectively.
We denote
\begin{align}
\label{eq:rho(Q)}
\rho=\rho(Q)&=\lim_{\substack{q=(1,y,z)\in Q\\ q\to p}} \frac{z^{n_y}}{y^{n_z}}\in\CC^*,\\
\label{eq:eta(Q)}
\eta(Q)&=\text{multiplicity of $(z^{n_y}-\rho y^{n_z})$ as a factor of $F^E(y,z)$.}
\end{align}
It is not hard to see that $\rho(Q)$ is well defined, and that $\eta(Q)\ge 1$.

\begin{proposition}
\label{pr:gen-shift}
Let $C=Z(F(x,y,z))$ be a projective curve not containing the line $Z(z)=L^\infty$
as a component.
Assume that $p=(1,0,0)\in C\cap L^\infty$, and that $C$ is not tangent to the line~$Z(y)$ at~$p$.
Let $a,b\in\CC$, and let $\Cab$ denote the shifted curve
\begin{equation}
\label{eq:shiftX}
\Cab=Z(F(x+az,y+bz,z)).
\end{equation}
Then
\begin{align}
\label{est-shift}
(C\cdot C_{a,b})_p\ge 2S^-(\Gamma(C,p))&-\mt(C,p)+(C\cdot L^\infty)_p +\sum_Q \min(r(Q), \eta(Q)-1),
\end{align}
where the sum is over all local branches $Q$ of~$C$ at~$p$ which are tangent to~$L^\infty$.
\end{proposition}

(For the definitions of $S^-(\Gamma(C,p))$, $r(Q)$ and $\eta(Q)$,
see Definition~\ref{def:Newton-Gamma}, \eqref{eq:d(Q)}--\eqref{eq:r(Q)} and \eqref{eq:rho(Q)}--\eqref{eq:eta(Q)}, respectively.)

\begin{proof}
Since the intersection multiplicity is lower semicontinuous, we may assume that $a$ and $b$
are generic complex numbers.

Let us denote
\begin{align*}
\mathcal{Q}(C,p) &= \text{the set of local branches of $C$ at $p$;}\\
\mathcal{Q}_0(C,p)&= \text{the set of local branches tangent to~$L^\infty$;}\\
\mathcal{Q}_1(C,p)&= \text{the set of local branches transversal to~$L^\infty$.}
\end{align*}
The local branches in $\mathcal{Q}_0(C,p)$ correspond to
the edges of $\Gamma(C,p)$ such that \hbox{$n_y<n_z$}, cf.\ Lemma~\ref{l:wzwy}.
Let $\Gamma_0(C,p)$ denote the union of these edges.
The local branches in $\mathcal{Q}_1(C,p)$, if~any, correspond to
the unique edge $E_{(1,1)}\!\subset\!\Gamma(C,p)$ with \hbox{$\mathbf{n}(E_{(1,1)})\!=\!(1,1)$}.
Figure~\ref{fig:Gamma(C,p)} illustrates the case where the edge $E_{(1,1)}$ is present;
equivalently, some local branches are transversal to~$L^\infty$.

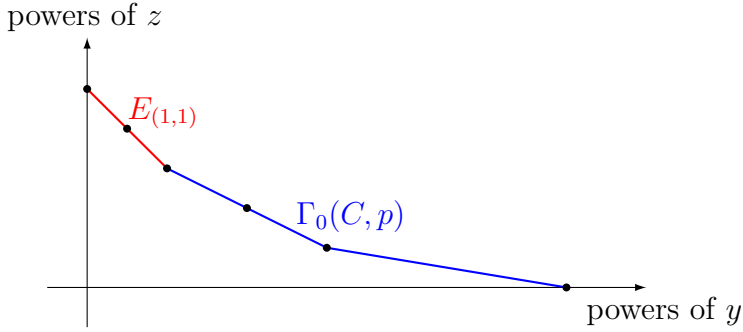
\begin{figure}[ht]
\begin{center}
\setlength{\unitlength}{1.5pt}
\begin{picture}(150,77)(0,0)
\thinlines
\put(0,10){\vector(1,0){150}}
\put(10,0){\vector(0,1){73}}
\thicklines
\put(10,60){\red{\line(1,-1){20}}}
\put(30,40){\blue{\line(2,-1){40}}}
\put(70,20){\blue{\line(6,-1){60}}}
\put(-10,76){powers of~$z$}
\put(135,2){powers of~$y$}
\put(29,54){\makebox(0,0){\red{$E_{(1,1)}$}}}
\put(76,28){\makebox(0,0){\blue{$\Gamma_0(C,p)$}}}

\put(10,60){\circle*{2}}
\put(20,50){\circle*{2}}
\put(30,40){\circle*{2}}
\put(50,30){\circle*{2}}
\put(70,20){\circle*{2}}
\put(130,10){\circle*{2}}
\end{picture}
\vspace{-.15in}
\end{center}
\caption{The Newton diagram $\Gamma(C,p)$.
The edge $E_{(1,1)}$ corresponds to the local branches transversal to~$L^\infty$.
The remaining edges of $\Gamma(C,p)$ form the subdiagram~$\Gamma_0(C,p)$;
they correspond to the local branches tangent to~$L^\infty$.}
\vspace{-.15in}
\label{fig:Gamma(C,p)}
\end{figure}

Let us consider the family of curves $\{C_{\lambda a,\lambda b}\}_{0\le\lambda\le1}$
interpolating between \hbox{$C=C_{0,0}$} and~$C_{a,b}\,$.
Being equisingular at the point~$p$,
this deformation preserves the number of local branches at~$p$
as well as their topological characteristics.
It therefore descends to families of individual local branches at~$p$, yielding
an equisingular bijection
\begin{align*}
\mathcal{Q}(C,p)&\longrightarrow\mathcal{Q}(C_{a,b},p)\\
Q&\longmapsto Q_{a,b}\,.
\end{align*}
Since the fixed point set of the shift $(x,y,z)\mapsto(x+az,y+bz,z)$ is the line $L^\infty$, this bijection restricts to bijections $\mathcal{Q}_0(C,p)\to\mathcal{Q}_0(C_{a,b},p)$
and $\mathcal{Q}_1(C,p)\to\mathcal{Q}_1(C_{a,b},p)$.

To obtain the desired lower bound on $(C\cdot C_{a,b})_p\,$, we will exploit the decomposition
\begin{equation}
(C\cdot C_{a,b})_p=\sum_{Q,Q'\in\mathcal{Q}(C,p)}(Q\cdot Q'_{a,b})_p
= \Sigma_{00}+\Sigma_{01}+\Sigma_{10}+\Sigma_{11}\,,
\label{eq-branch}
\end{equation}
where we use the notation
\[
\Sigma_{\varepsilon\delta}=\sum_{\substack{
Q\in\mathcal{Q}_\varepsilon(C,p)\\
Q'\in\mathcal{Q}_\delta(C,p)
}}
(Q\cdot Q'_{a,b})_p\,, \quad \text{for $\varepsilon,\delta\in\{0,1\}$.}
\]

\begin{lemma}
\label{lem:sum-QQ}
Suppose $\mathcal{Q}_1(C,p)\neq\varnothing$. Then
\begin{equation}
\label{eq:three-sigmas}
\Sigma_{01}+\Sigma_{10}+\Sigma_{11}=2S^-(E_{(1,1)}),
\end{equation}
where $S^-(E_{(1,1)})$ denotes the area of the trapezoid
bounded by the edge~$E_{(1,1)}$,
the coordinate axes, and the vertical line through the rightmost endpoint of~$E_{(1,1)}$.
\end{lemma}

\begin{proof}
For $Q\in\mathcal{Q}_1(C,p)$,
the tangent lines to $Q$ and $Q_{a,b}$ differ from each other.
The foregoing discussion implies that,
for any $Q\in\mathcal{Q}_1(C,p)$ and $Q'\in\mathcal{Q}(C,p)$, we have
\[
(Q\cdot Q'_{a,b})_p=(Q_{a,b}\cdot Q')_p=\mt Q\cdot \mt Q'. 
\]
We then observe that
\begin{align*}
\sum_{Q\in{\mathcal Q}_0(C,p)} \mt Q=\ell_0 &\,\eqdef\, \text{length of the projection of $\Gamma_0(C,p)$ to the vertical axis}, \\
\sum_{Q\in{\mathcal Q}_1(C,p)} \mt Q=\ell_1 &\,\eqdef\, \text{length of the projection of $E_{(1,1)}$ to either of the axes},
\end{align*}
implying
\[
\Sigma_{01}+\Sigma_{10}+\Sigma_{11}=\ell_0\ell_1+\ell_1\ell_0+\ell_1^2=2S^-(E_{(1,1)}).
\qedhere
\]
\end{proof}


In light of \eqref{eq-branch} and \eqref{eq:three-sigmas}, it remains to obtain the desired
lower bound for the summand~$\Sigma_{00}\,$.
To simplify notation, we will pretend, for the time being, that all local branches are tangent to~$L^\infty$,
so that $\mathcal{Q}(C,p)=\mathcal{Q}_0(C,p)$ and $\Gamma(C,p)=\Gamma_0(C,p)$.

Let $\mathcal{Q}(E)$ denote the set of local branches of $C$ at $p$
associated with an edge~$E$ of~$\Gamma(C,p)$.
Equivalently, $\mathcal{Q}(E)=\{Q\in\mathcal{Q}(C,p)\mid E(Q)=E\}$.

\begin{lemma}
\label{lem:E1E2}
Let $E_1$ and $E_2$ be two distinct edges of the Newton diagram~$\Gamma(C,p)$.
Assume that $E_2$ is located above and to the left of~$E_1$, so that
\[
E_1\!=\![(i_1,j_1),(i'_1,j'_1)], \quad
E_2\!=\![(i_2,j_2),(i'_2,j'_2)], \quad
i_2'\!<\!i_2\!\le\! i_1'\!<\!i_1\,,\quad
j_1\!<\!j_1'\!\le\! j_2\!<\!j_2'\,.
\]
Then
\begin{equation}
\label{el11.6}
\sum_{Q\in\mathcal{Q}(E_1)}\sum_{Q'\in\mathcal{Q}(E_2)}(Q\cdot Q'_{a,b})_p=
\sum_{Q\in\mathcal{Q}(E_1)}\sum_{Q'\in\mathcal{Q}(E_2)}(Q_{a,b}\cdot Q')_p
=(j'_1-j_1)(i_2-i'_2).
\end{equation}
\end{lemma}

Note that $(j'_1-j_1)(i_2-i'_2)$ is the area of the rectangle
formed by the intersections of the horizontal lines passing through~$E_1$
with the vertical lines passing through~$E_2$.
See Figure~\ref{fig:E1E2}.

\begin{figure}[ht]
\begin{center}
\setlength{\unitlength}{1.5pt}
\begin{picture}(130,57)(20,0)
\thinlines
\put(30,10){\vector(1,0){120}}
\put(30,10){\vector(0,1){43}}
\thicklines
\put(30,40){\blue{\line(2,-1){40}}}
\put(70,20){\blue{\line(6,-1){60}}}
\put(10,56){powers of~$z$}
\put(150,2){powers of~$y$}
\put(55,35){\makebox(0,0){\blue{$E_2$}}}
\put(101,20){\makebox(0,0){\blue{$E_1$}}}

\put(24,40){\makebox(0,0){\blue{$j_2'$}}}
\put(17,20){\makebox(0,0){\blue{$j_2\!=\!j_1'$}}}
\put(24,10){\makebox(0,0){\blue{$j_1$}}}

\put(31,4){\makebox(0,0){\blue{$i_2'$}}}
\put(70,4){\makebox(0,0){\blue{$i_2\!=\!i_1'$}}}
\put(130,4){\makebox(0,0){\blue{$i_1$}}}

\put(30,40){\circle*{2}}
\put(30,20){\circle*{2}}
\put(30,10){\circle*{2}}
\put(70,10){\circle*{2}}
\put(50,30){\circle*{2}}
\put(70,20){\circle*{2}}
\put(130,10){\circle*{2}}
\end{picture}
\vspace{-.15in}
\end{center}
\caption{Two edges in the Newton diagram $\Gamma(C,p)$.
Here $\Gamma(C,p)\!=\!\Gamma_0(C,p)$.}
\vspace{-.15in}
\label{fig:E1E2}
\end{figure}
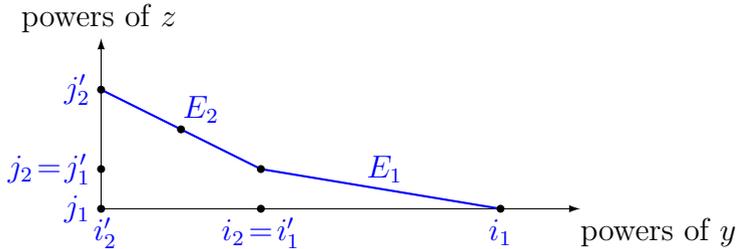

\redsf{We provide two proofs of Lemma~\ref{lem:E1E2}. 
The first proof, based on the Bernstein-Kouchnirenko mixed volume formula, 
relies on a genericity assumption whose justification we do not provide. 
(Note however that this assumption is not needed for the proof of the ``$\ge$'' inequality in~\eqref{el11.6}.
This inequality will be sufficient for the upcoming proof of Proposition \ref{pr:gen-shift}.) 
The second proof is rigorous but more technical. 

\begin{proof}[Proof~1 of Lemma~\ref{lem:E1E2} (sketch)]
Assume that the truncations of $F(1,y,z)$ along the edges $E_1$ and $E_2$ are square-free. 
Construct the right triangles $\tau_1$ and~$\tau_2$ as shown in Figure~\ref{fig-BK}(a). 
By Bernstein's theorem~\cite{DBernstein}, the left and middle terms in \eqref{el11.6} 
can be bounded from below by the difference between the mixed area of $E_1, E_2$ 
and the mixed area of $\tau_1, \tau_2$ (see Figure~\ref{fig-BK}(b,c)). The claim follows.
\end{proof}
}

\begin{figure}[ht]
\begin{center}
\setlength{\unitlength}{1.5pt}
\begin{picture}(115,48)(0,55)
\thinlines
\put(0,65){\vector(1,0){113}}
\put(0,65){\vector(0,1){40}}
\thicklines
\put(0,95){\blue{\line(2,-1){40}}}\put(0,75){\blue{\line(1,0){40}}}
\put(40,75){\blue{\line(6,-1){60}}}\put(40,75){\blue{\line(0,-1){10}}}
\put(25,90){\makebox(0,0){\blue{$E_2$}}}\put(15,80){\makebox(0,0){\blue{$\tau_2$}}}
\put(71,75){\makebox(0,0){\blue{$E_1$}}}\put(52,68){\makebox(0,0){\blue{$\tau_1$}}}
%
%

\put(0,95){\circle*{2}}
\put(0,75){\circle*{2}}
\put(40,65){\circle*{2}}
\put(40,75){\circle*{2}}
\put(100,65){\circle*{2}}

\put(45,55){\rm(a)}
\end{picture}
\\[.1in]

\setlength{\unitlength}{1.5pt}
\begin{picture}(115,50)(125,60)
\thinlines
\put(125,65){\vector(1,0){113}}
\put(125,65){\vector(0,1){40}}
\thicklines
\put(125,95){\blue{\line(2,-1){40}}}\put(125,95){\blue{\line(6,-1){60}}}
\put(165,75){\blue{\line(6,-1){60}}}\put(185,85){\blue{\line(2,-1){40}}}

\put(125,95){\circle*{2}}\put(185,85){\circle*{2}}
\put(165,75){\circle*{2}}
\put(225,65){\circle*{2}}

\put(170,55){\rm(b)}
\linethickness{0.1pt}
\multiput(165,75)(3,-0.5){20}{\line(-2,1){40}}
\end{picture}
\qquad
\setlength{\unitlength}{1.5pt}
\begin{picture}(115,50)(60,5)
\thinlines
\put(60,10){\vector(1,0){113}}
\put(60,10){\vector(0,1){40}}
\thicklines
\put(60,40){\blue{\line(2,-1){40}}}\put(60,40){\blue{\line(6,-1){60}}}\put(60,20){\blue{\line(1,0){40}}}
\put(100,20){\blue{\line(6,-1){60}}}\put(120,30){\blue{\line(2,-1){40}}}\put(100,20){\blue{\line(0,-1){10}}}

\put(60,40){\circle*{2}}
\put(60,20){\circle*{2}}
\put(60,10){\circle*{2}}
\put(100,10){\circle*{2}}
\put(120,30){\circle*{2}}
\put(100,20){\circle*{2}}
\put(160,10){\circle*{2}}

\put(105,0){\rm(c)}
\linethickness{0.1pt}
\multiput(100,20)(3,-0.5){20}{\line(-2,1){40}}
\multiput(100,20)(0,-1){10}{\line(-1,0){40}}
\end{picture}

\end{center}
\caption{
\redsf{
Proof 1 of Lemma \ref{lem:E1E2}. (a) Right triangles $\tau_1, \tau_2$.
(b) Mixed area of $E_1, E_2$. (c) Mixed area of $\tau_1,\tau_2$.}
}
\vspace{-.15in}
\label{fig-BK}
\end{figure}

\begin{proof}[Proof~2 of Lemma~\ref{lem:E1E2}]
In the coordinates $(1,y,z)$, the shift transformation
\[
(1,y,z)\mapsto(1+az,y+bz,z)
\]
converts the polynomial $F_1(y,z)=F(1,y,z)$ defining the affine curve $C\setminus Z(x)$
into the polynomial $F_{a,b}(y,z)=F(1+az,y+bz,z)$ defining~$C_{a,b}\setminus Z(x)$.
Under this transformation, each monomial $y^iz^j$ of $F_1$ becomes
\begin{equation}y^iz^j+\sum_{\renewcommand{\arraystretch}{0.6}
\begin{array}{c}
\scriptstyle{i',j'\ge0}\\
\scriptstyle{i'+j'>0}\end{array}}c_{i'j'}y^{i-i'}z^{j+i'+j'}.
\label{eq3}
\end{equation}
In view of the assumption $\mathcal{Q}(C,p)=\mathcal{Q}_0(C,p)$,
the monomials appearing in the sum above correspond to integer points lying strictly above the Newton diagram~$\Gamma(C,p)$ of~$F_1$.
In particular, $F_{a,b}$ has the same Newton diagram~$\Gamma(C,p)$,
and the same truncations to its edges.
Furthermore, each local branch $Q\in\mathcal{Q}(C,p)$ and its counterpart $Q_{a,b}\in\mathcal{Q}(C_{a,b},p)$ are associated with the same edge of~$\Gamma(C,p)$.
The union of the local branches of $C_{a,b}$ associated with the edge~$E_1$
can be defined by an analytic equation $f(y,z)=0$
whose Newton diagram at the origin
is the line segment $[(i_1-i'_1,0),(0,j'_1-j_1)]$
(cf.\ the Newton-Puiseux algorithm \cite[Algorithm~I.3.6]{GLS}).
For the same reason, a local branch $Q$ of $C$ at $p$ associated with the edge $E_2$ has a parametrization
\[
\begin{array}{l}
y=\varphi(t)=t^m,\\[.03in]
z=\psi(t)=\alpha t^d+O(t^{d+1}),
\end{array}
\quad
\alpha\!\ne\!0,\quad
|t| \! \ll \!1,\quad
d\!=\!(Q\cdot L^\infty)_p\,,\quad
m\!=\!\mt(Q)
\]
(cf.\ \eqref{eq:d(Q)}--\eqref{eq:m(Q)})
where $\frac{d}{m}=\frac{j'_2-j_2}{i_2-i'_2}$.
Since $\frac{j'_2-j_2}{i_2-i'_2}>\frac{j'_1-j_1}{i_1-i'_1}$,
we obtain
\[
f(\varphi(t),\psi(t))=t^{m(i_1-i_1')}(\beta+O(t)),\quad\beta\ne0.
\]
The statement of the lemma now follows from the fact that the total multiplicity
of the local branches of $C$ at $p$ associated with the edge $E_2$ equals $j_2'-j_2\,$.
\end{proof}

We are now left with the task of computing
\begin{equation}
\sum_E \sum_{Q,Q'\in\mathcal{Q}(E)}(Q\cdot Q'_{a,b})_p\,,
\label{eq:last-sum}
\end{equation}
where the first sum runs over the edges $E$ of the Newton diagram.
In this part of the proof, we continue to assume, for the sake of simplifying the exposition,
that all local branches are tangent to~$L^\infty$.
Since the shift $(a,y,z)\mapsto(1+az,y+bz,z)$ acts independently on the analytic factors of~$F_1(y,z)$,
we furthermore assume, in our computation of
$\sum_{Q,Q'\in\mathcal{Q}(E)}(Q\cdot Q'_{a,b})_p$ (see~\eqref{eq:last-sum}),
that the Newton diagram $\Gamma(C,p)$ consists of
a single edge $E=[(0,M),(D,0)]$, with $M<D$ and $\nE=(n_y,n_z)$.

Pick a local branch $Q\in\mathcal{Q}(E)$.
It admits an analytic parametrization of the form
\begin{equation}
\begin{array}{l}
x=1, \quad \\
y=\varphi(t)=t^m,\quad \\
z=\psi(t)=\alpha t^d+O(t^{d+1}),
\end{array}
\label{eq2}
\end{equation}
where $t$ ranges over a small disk in $\CC$ centered at zero,
$m=r(Q)\cdot n_y$, $d=r(Q)\cdot n_z$, and $r(Q)=\gcd(d,m)$,
cf.\ \eqref{eq:nE=} and~\eqref{eq:r(Q)}.

\begin{lemma}
\label{lem:Step3}
We have
\begin{equation}
\label{eq:sum-Q'}
\sum_{Q'\in\mathcal{Q}(E)}(Q\cdot Q'_{a,b})_p=dM-m+d+\min(r(Q),\eta(Q)-1).
\end{equation}
\end{lemma}

\redsf{Before providing a proof of Lemma~\ref{lem:Step3}, we note that the weaker inequality
\begin{equation}
\label{eq:sum-Q'-ineq}
\sum_{Q'\in\mathcal{Q}(E)}(Q\cdot Q'_{a,b})_p \ge dM-m+d 
\end{equation}
can be deduced directly from the Bernstein-Kouchnirenko formula, 
similarly to Proof~1 of Lemma~\ref{lem:E1E2} provided above. 
As in the case of Lemma~\ref{lem:E1E2}, this lower bound 
would be sufficient to complete the upcoming proof of Proposition~\ref{pr:gen-shift},
restricted to the case of a Newton non-degenerate singularity~$(C,p)$.
}

\begin{proof}
The left-hand side of~\eqref{eq:sum-Q'} is the minimal exponent of $t$
appearing in the expansion of $F_{a,b}(\varphi(t),\psi(t))$ into a power series in~$t$.
Since $F_1(\varphi(t),\psi(t))=0$, we may instead substitute \eqref{eq2}
into the difference $F_{a,b}(y,z)-F_1(y,z)$,
or equivalently into the monomials of the second summand in~\eqref{eq3}
(corresponding to individual monomials $y^iz^j$ of~$F_1(y,z)$).
Evaluating $y^{i-i'}z^{j+i'+j'}$ at $y=\varphi(t)$, $z=\psi(t)$, we obtain
\[
(t^m)^{i-i'}(t^d(\alpha+O(t)))^{j+i'+j'}=t^{mi+dj+(d-m)i'+dj'}(\alpha^{j+i'+j'}+O(t)).
\]
To get the minimal value of the exponent $mi+dj+(d-m)i'+dj'$,
we need to minimize $mi+dj$ (which is achieved for $(i,j)\in E$) and take $i'=1$ and $j'=0$.
Developing $F_{a,b}(y,z)=F(1+az,y+bz,z)$ into a power series in $a$ and~$b$,
we see that the corresponding monomials $y^{i-i'}z^{j+i'+j'}=y^{i-1}z^{j+1}$
in $F_{a,b}(y,z)-F(1,y,z)$ add up to $bzF_y(1,y,z)$.
We conclude that the desired minimal exponent of~$t$ occurs when we substitute
$(y,z)=(\varphi(t),\psi(t))$ either into $bzF^E_y(y,z)$
or into a monomial $y^{i-1}z^{j+1}$ such that $(i,j)$ is one of the  integral points
closest to the edge~$E$ and lying above~$E$.
The latter condition reads $n_yi+n_zj=n_zM+1$.

The truncation $F^E(y,z)$ of $F_1(y,z)$ 
has the form
\begin{equation}
F^E(y,z)=
\prod_{k=1}^n(z^{n_y}-\beta_ky^{n_z})^{r_k},
\label{eq:trun}
\end{equation}
where $\beta_1,\dots,\beta_n\in\CC$ are distinct, and $r_1,\dots,r_n\in\ZZ_{>0}$.
Developing $F^E(y,z)$ into~a power series in~$t$, we see that the monomials of $F^E$
yield the minimal exponent of~$t$.
Since $F_1(\varphi(t),\psi(t))=0$, these minimal powers of $t$ must cancel out,
implying that, for some $k_0\in\{1,\dots,n\}$, we have (cf.\ \eqref{eq:rho(Q)}, \eqref{eq:eta(Q)}):
\begin{align*}
\rho=\rho(Q)&=\lim_{t\to 0}\frac{(\alpha t^d)^{n_y}}{(t^m)^{n_z}}=\alpha^{n_y}=\beta_{k_0}\,,\\
\eta=\eta(Q)&=r_{k_0}\,.
\end{align*}

The factorization formula \eqref{eq:trun} implies that
\[
bzF^E_y(y,z)=bn_zy^{n_z-1}z\sum_{k=1}^n\Bigl((-\beta_k)^{r_k}(z^{n_y}-\beta_ky^{n_z})^{r_k-1}\prod_{l\ne k}(z^{n_y}-\beta_ly^{n_z})^{r_l}\Bigr).
\]
We then compute the minimal exponent for $bzF^E_y(y,z)$:
\begin{align*}bn_zy^{n_z-1}z\big|_{y=\varphi(t),z=\psi(t)}&=O(t^{mn_z-m+d}),\\
(z^{n_y}-\beta_ky^{n_z})^{r_k-1}\big|_{y=\varphi(t),z=\psi(t)}&=O(t^{dn_y(r_k-1)}),\quad k\ne k_0,\\
(z^{n_y}-\rho y^{n_z})^{\eta-1}\big|_{y=\varphi(t),z=\psi(t)}&=O(t^{(dn_y+1)(\eta-1)}),\\
(z^{n_y}-\beta_ly^{n_z})^{r_l}\big|_{y=\varphi(t),z=\psi(t)}&=O(t^{dn_yr_l}),\quad l\ne k_0,\\
(z^{n_y}-\rho y^{n_z})^{\eta}\big|_{y=\varphi(t),z=\psi(t)}&=O(t^{(dn_y+1)\eta}), \\
bzF^E_y(y,z)\big|_{y=\varphi(t),z=\psi(t)}&=O(t^{dM-m+d+\eta-1}).
\end{align*}
Also, for $n_yi+n_zj=n_zM+1$, we have
\[
y^{i-1}z^{j+1}\big|_{y=\varphi(t),z=\psi(t)}=O(t^{dM-m+d+r(Q)}).
\]
Consequently
\begin{align*}
\sum_{Q'\in\mathcal{Q}(E)}(Q\cdot Q'_{a,b})_p
&=\min(dM-m+d+\eta-1,dM-m+d+r(Q)), 
\end{align*}
as desired.
\end{proof}

We are now ready to complete the proof of Proposition~\ref{pr:gen-shift}.
We first note that in Lemma~\ref{lem:Step3},
\[
m=\mt Q,\quad d=(Q\cdot L^\infty)_p\,,
\]
and moreover
\begin{equation}
\label{eq:M,D}
\sum_{Q\in\mathcal{Q}(E)}\mt Q=M,\quad \sum_{Q\in\mathcal{Q}(E)}(Q\cdot L^\infty)_p=D.
\end{equation}
Therefore
\begin{align}
\notag
\sum_{Q,Q'\in\mathcal{Q}(E)}(Q\cdot Q'_{a,b})_p
&=\sum_Q (dM-m+d+\min(r(Q),\eta(Q)-1)) \\
\label{eq:DM-M+D...}
&=DM-M+D+\sum_Q \min(r(Q),\eta(Q)-1).
\end{align}
Note that $DM$ is twice the area of the right triangle with hypotenuse~$E$.

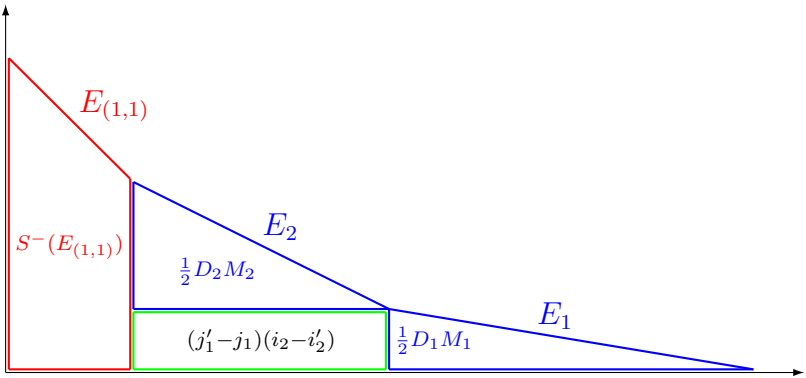
\begin{figure}[ht]
\begin{center}
\setlength{\unitlength}{3.2pt}
\begin{picture}(125,60)(10,10)
\thinlines
\put(10,10){\vector(1,0){125}}
\put(10,10){\vector(0,1){58}}
\thicklines
\put(10.5,59.5){\red{\line(1,-1){19}}}
\put(30,40){\blue{\line(2,-1){40}}}
\put(70,20){\blue{\line(6,-1){57}}}
\put(27,52){\makebox(0,0){\red{$E_{(1,1)}$}}}
\put(53,33){\makebox(0,0){\blue{$E_2$}}}
\put(96,19){\makebox(0,0){\blue{$E_1$}}}

\put(10.5,10.5){\red{\line(0,1){49}}}
\put(10.5,10.5){\red{\line(1,0){19}}}
\put(29.5,10.5){\red{\line(0,1){30}}}
\put(20,30){\makebox(0,0){\red{$ S^-(E_{(1,1)})$}}}
\put(43,26){\makebox(0,0){\blue{$ \frac12 D_2 M_2$}}}
\put(77,15){\makebox(0,0){\blue{$ \frac12 D_1 M_1$}}}
\put(50,15){\makebox(0,0){$ (j'_1-j_1)(i_2-i'_2)$}}
\put(30,10.5){\green{\line(1,0){39.5}}}
\put(30,19.5){\green{\line(1,0){39.5}}}
\put(30,10.5){\green{\line(0,1){9}}}
\put(69.5,10.5){\green{\line(0,1){9}}}
\put(30,40){\blue{\line(0,-1){20}}}
\put(30,20){\blue{\line(1,0){40}}}
\put(70,20){\blue{\line(0,-1){9.5}}}
\put(70,10.5){\blue{\line(1,0){57}}}
\end{picture}
\end{center}
\caption{The area $S^-(\Gamma(C,p))$ under the Newton diagram $\Gamma(C,p)$ is obtained by
adding three kinds of contributions:
\red{(a)} the area $S^-(E_{(1,1)})$ of the trapezoid under\-neath the edge~$E_{(1,1)}$,
cf.\ Lemma~\ref{lem:sum-QQ};
\green{(b)} the areas $(j'_1\!-\!j_1)(i_2\!-\!i'_2)$ of rectangles obtained from pairs of edges
$E_1,E_2$ of~$\Gamma_0(C,p)$, cf.\ Lemma~\ref{lem:E1E2}; and
\blue{(c)} the areas $\frac12 DM$ of right triangles adjacent to the edges in~$\Gamma_0(C,p)$,
cf.~\eqref{eq:DM-M+D...}.}
\label{fig:2*area}
\end{figure}

As illustrated in Figure~\ref{fig:2*area}, adding up the contributions $DM$ from all edges~$E$
of~$\Gamma_0(C,p)$,
together with the contributions coming from Lemmas~\ref{lem:sum-QQ} and~\ref{lem:E1E2},
we obtain $2S^-(\Gamma(C,p))$, cf.\ Definition~\ref{def:Newton-Gamma}.

Finally, in view of~\eqref{eq:M,D}, we have
\begin{align*}
\sum_{E\subset \Gamma_0(C,p)} (-M+D)
&=\sum_{Q\in \mathcal{Q}_0(C,p)} (-\mt(Q)+(Q\cdot L^\infty)_p)\\
&=\sum_{Q\in \mathcal{Q}(C,p)} (-\mt(Q)+(Q\cdot L^\infty)_p)\\
&=-\mt(C,p)+(C\cdot L^\infty)_p\,.
\end{align*}
Putting everything together, we obtain~\eqref{est-shift}.
\end{proof}

\newpage

\section{Arrangements of polynomial curves}
\label{sec:arrangements-poly}

In this section, we generalize Example~\ref{ex:lines+parabolas} to arrangements of
polynomial curves obtained from a given curve by shifts, dilations, and/or rotations.
We start by obtaining upper bounds on the number of intersection points of two polynomial curves
related by one of these transformations.
These bounds lead to expressivity criteria for arrangements consisting of such curves.

Recall 
that for $a,b\in\CC$ and $c\in\CC^*$, we denote by
\begin{align}
\label{eq:Cab}
\Cab&=Z(F(x+az,y+bz,z)), \\
C_c&=Z(F(cx,cy,z))
\end{align}
the curves obtained from a plane curve $C=Z(F(x,y,z))$ by a shift and a dilation, respectively.
We will also use the notation
\begin{equation}
\label{eq:shift+dilation}
\Cabc=Z(F(cx+az,cy+bz,z))
\end{equation}
for a curve obtained from~$C$ by a combination of a shift and a dilation.

As before, we identify a projective curve $C=Z(F(x,y,z))$ with its restriction to the affine plane
$\AA^2=\PP^2\setminus L^\infty$ given by $C=V(G(x,y))$, where $G(x,y)=F(x,y,1)$.
Under this identification, we have
\begin{align}
\Cab&=V(G(x+a,y+b)), \\
C_c&=V(G(cx,cy)),\\
\label{eq:Cabs=V(G)}
\Cabc&=V(G(cx+a,cy+b)).
\end{align}

\begin{remark}
\label{rem:shift+dilation}
Unless $c=1$ (the case of a pure shift),
the transformation $C\leadsto \Cabc$ can be viewed as a pure dilation centered
at some point $o\in\CC$ (where $o$ may be different from~$0$).
\end{remark}

\begin{corollary}
\label{cor:shift-poly}
Let $C$ be a real polynomial curve of degree~$d$.
Let $C\cap L^\infty=\{p\}$  and $m=\mt(C,p)$.
Then we have, for $a,b\in\CC$ and $c\in\CC^*$:
\begin{align}
\label{eq:dilation-poly}
(C\cdot \Cabc)_p &\ge dm ;\\
\label{eq:shift-poly}
(C\cdot \Cab)_p &\ge (d-1)(m+1)+\gcd(d,m).
\end{align}
\end{corollary}

\begin{proof}
In view of Remark~\ref{rem:shift+dilation} and the inequality
\[
(d-1)(m+1)+\gcd(d,m)\ge dm,
\]
it suffices to establish \eqref{eq:dilation-poly} in the case $a=b=0$,
i.e., with $\Cabc$ replaced by~$C_c$.

The bounds \eqref{eq:dilation-poly}--\eqref{eq:shift-poly} are obtained
by applying Propositions~\ref{pr:gen-dilation} and~\ref{pr:gen-shift} to
the case of a polynomial curve~$C$, while noting that in this case,
\begin{align*}
S^-(\Gamma(C,p))&=\tfrac12 md,\\
\mt(C,p)&=m,\\
(C\cdot L^\infty)_p&=d,\\
r(Q)=\eta(Q)&=\gcd(d,m). \qedhere
\end{align*}
\end{proof}

\begin{proposition}
\label{pr:dilations-poly}
Let $C$ be a real polynomial curve with a parametrization
\[
t\mapsto (P(t),Q(t)), \quad \deg(P)=d, \quad \deg(Q)=d'\neq d.
\]
Let $a,b,c\in\CC$, $c\neq 0$, and let $\Cabc$ be the shifted and dilated curve given by~\eqref{eq:shift+dilation}.
Assume that $C$ and $\Cabc$ have $N$ intersection points in $\AA^2=\PP^2\setminus L^\infty$,
\redsf{counting with multiplicity}. 
Then $N\le dd'$.
\end{proposition}

\begin{proof}
Without loss of generality, assume that $d>d'$.
By B\'ezout's theorem, the inter\-section multiplicity $(C\cdot\Cabc)_p$
at the point $p=(1,0,0)\in L^\infty$ is at most~\hbox{$d^2-N$}.
Applying \eqref{eq:dilation-poly}, with $m=d-d'$, we obtain $d^2-N\ge d(d-d')$.
The claim follows.
\end{proof}

Proposition~\ref{pr:dilations-poly} and Theorem~\ref{th:reg-expressive} imply
that if $C$ and $\Cabc$ intersect at $dd'$ hyperbolic nodes in the real affine plane,
then the union $C\cup\Cabc$ is expressive.

\begin{figure}[ht]
\begin{center}
\includegraphics[scale=0.16, trim=20cm 24.5cm 25cm 13cm, clip]{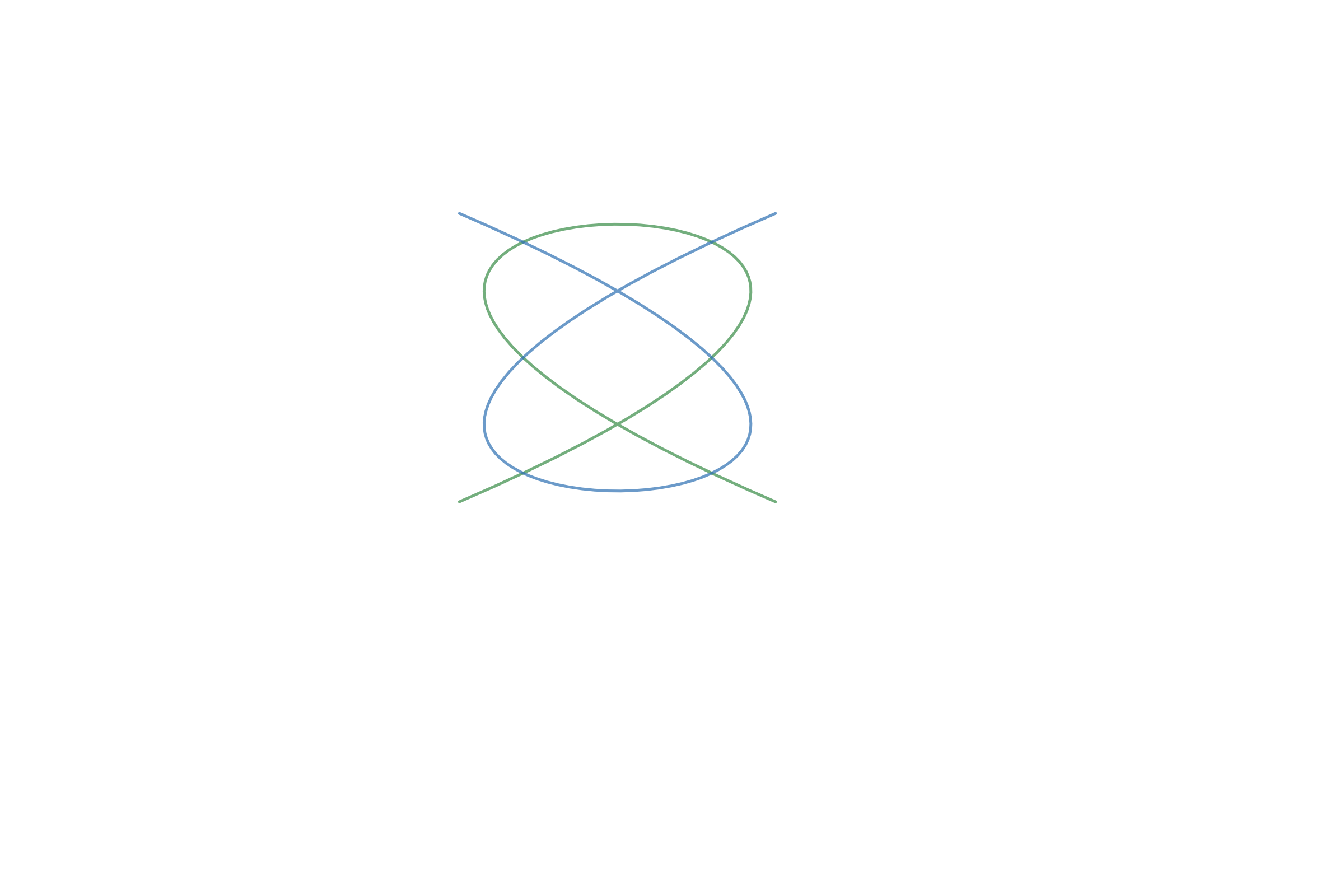}
\vspace{-.1in}
\end{center}
\caption{An expressive cubic and its reflection, intersecting at 6 real points.
The~resulting two-component curve is expressive.}
\vspace{-.1in}
\label{fig:chebyshev23-dilation}
\end{figure}

\begin{figure}[ht]
\begin{center}
\includegraphics[scale=0.28, trim=20cm 22.5cm 22cm 18cm, clip]{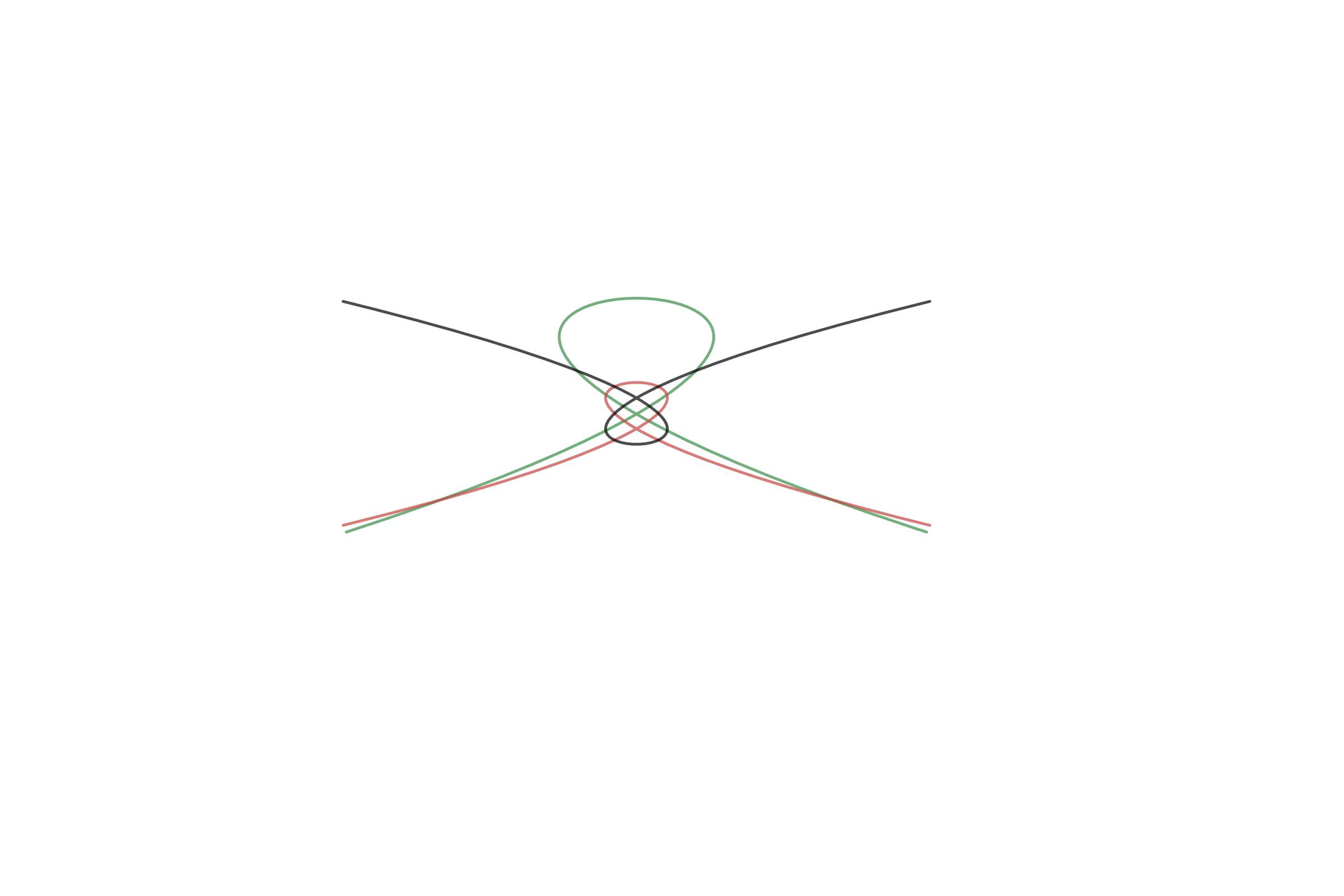}
\end{center}
\caption{An expressive curve whose components are singular cubics related to each other
by dilations. Each pair of them intersect at 6 hyperbolic nodes.}
\vspace{-.1in}
\label{fig:three-dilations}
\end{figure}

\begin{example}
\label{ex:23-dilation}
Let $C$ be the $(2,3)$-Chebyshev curve, the singular cubic given by
\begin{equation}
\label{eq:23-Chebyshev}
2x^2-1+4y^3-3y=0
\end{equation}
or parametrically by
\[
t\mapsto (4t^3-3t, -2t^2+1),
\]
cf.\ \eqref{eq:32-Chebyshev}.
Applying Proposition~\ref{pr:dilations-poly} (with $d=3$ and $d'=2$),
we see that the curve $C$ and its dilation~$\Cabc$ ($c\neq 1$)
can intersect in the real affine plane in at most 6~points.
When this bound is attained, the union $C\cup\Cabc$ is expressive.
Figure~\ref{fig:chebyshev23-dilation} shows one such example,
with $a=b=0$ and $c=-1$ (so $\Cabc$ is a reflection of~$C$).
Cf.\ also Figure~\ref{fig:three-dilations}.
\end{example}

\begin{proposition}
\label{pr:shifts-poly}
Let $C$ be a real polynomial curve with a parametrization
\[
t\mapsto (P(t),Q(t)), \quad \deg(P)=d, \quad \deg(Q)=d'<d.
\]
Let $a,b\in\CC$, and let $\Cab$ be the shifted curve given by~\eqref{eq:Cab}.
Assume that $C$ and $\Cab$ have $N$ intersection points in $\AA^2=\PP^2\setminus L^\infty$,
\redsf{counting with multiplicity}. 
Then
\begin{equation}
\label{eq:N-shift}
N\le dd'-d'-\gcd(d,d')+1.
\end{equation}
\end{proposition}

\begin{proof}
We use the same arguments as in the proof of Proposition~\ref{pr:dilations-poly} above,
with the lower bound~\eqref{eq:dilation-poly} replaced by~\eqref{eq:shift-poly}:
\begin{align*}
N\le d^2-(C\cdot \Cab)_p &\le d^2-(d-1)(m+1)-\gcd(d,m) \\
                                      &=d^2-(d-1)(d-d'+1)-\gcd(d,d') \\
                                      &=dd'-d'+1-\gcd(d,d'). \qedhere
\end{align*}
\end{proof}

For special choices of shifts, the bound \eqref{eq:N-shift} can be strengthened.
Here is one example, involving the Lissajous-Chebyshev curves,
see Example~\ref{ex:lissajous-chebyshev}.
(Note that such a curve does not have to be polynomial:
it could be trigonometric or reducible.)

\begin{proposition}
\label{pr:shifts-lissajous}
Let $C$ be a Lissajous-Chebyshev curve given by
\begin{equation}
\label{eq:lissajous-kl}
T_k(x)+T_\ell(y)=0.
\end{equation}
Let $a\!\in\!\CC$, and let $C_{a,0}$ (resp.,~$C_{0,a}$)
be the shift of~$C$ in the $x$ (resp.,~$y$) direction.
Assume that $C$ and $C_{a,0}$ (resp.,~$C_{0,a}$) intersect at $N_x$ (resp.,~$N_y$) points in~$\AA^2$.
Then
\begin{align}
\label{eq:x-shift-lissajous}
N_x &\le (k-1)\ell, \\
\label{eq:y-shift-lissajous}
N_y &\le k(\ell-1).
\end{align}
\end{proposition}

We note that in the case when $\gcd(k,\ell)=1$ and $k<\ell$, the bound \eqref{eq:y-shift-lissajous} matches~\eqref{eq:N-shift},
with $d=k$ and $d'=\ell$, whereas \eqref{eq:x-shift-lissajous} gives a stronger bound.

\begin{proof}
Due to symmetry, it suffices to prove~\eqref{eq:y-shift-lissajous}.
The intersection of the curves $C$ and~$C_{0,a}$ is given by
\begin{align*}
\begin{cases}
T_k(x)+T_\ell(y)=0 \\
T_k(x)+T_\ell(y+a)=0
\end{cases}
\Longleftrightarrow{\ }
\begin{cases}
T_k(x)+T_\ell(y)=0 \\
T_\ell(y+a)-T_\ell(y)=0.
\end{cases}
\end{align*}
The equation $T_\ell(y+a)-T_\ell(y)=0$ has at most $\ell-1$ roots;
each of these values of~$y$ then gives at most $k$ possible values of~$x$.
\end{proof}

\begin{figure}[ht]
\begin{center}
\includegraphics[scale=0.25, trim=20cm 24.5cm 20cm 16.5cm, clip]{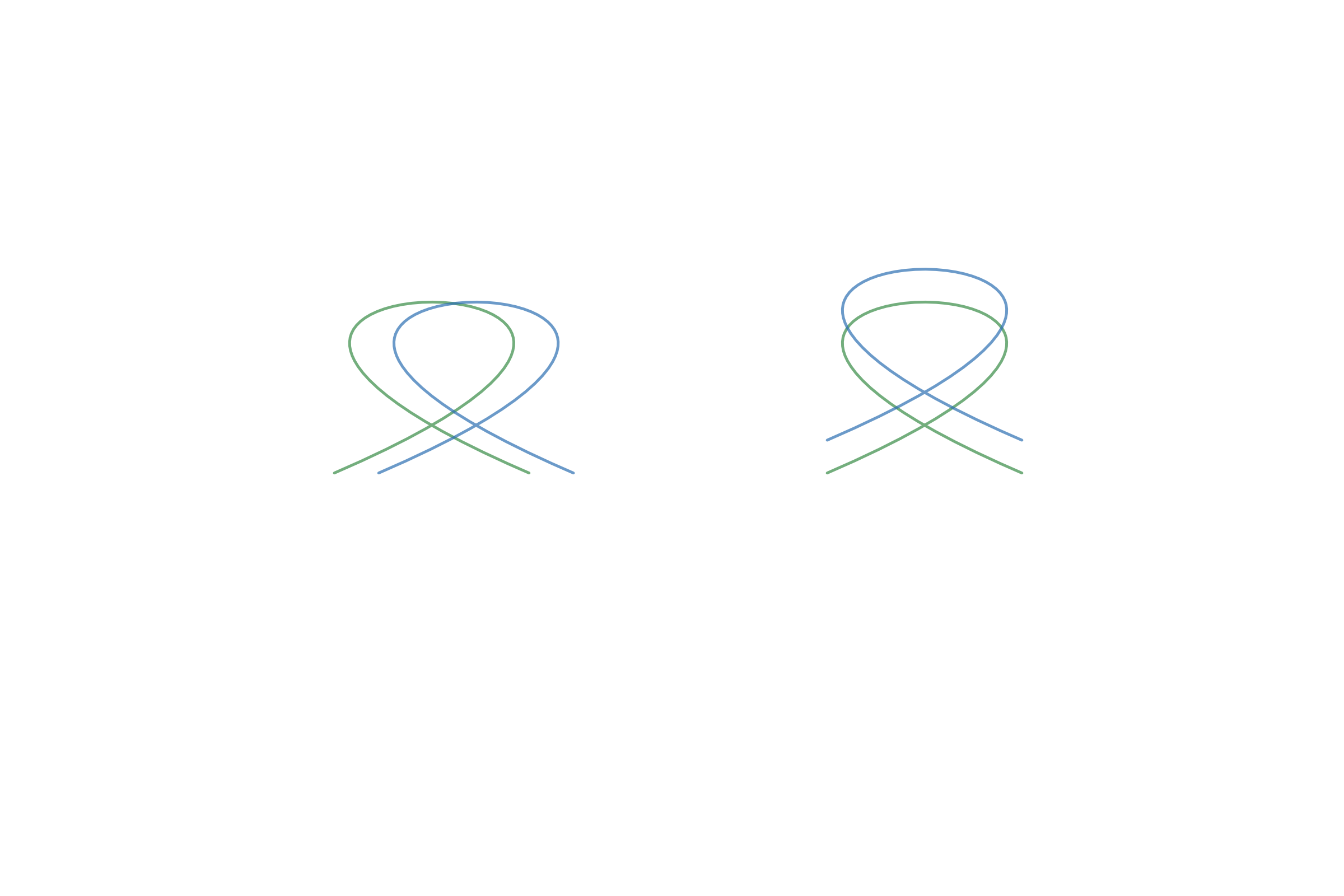}
\vspace{-.1in}
\end{center}
\caption{An expressive cubic and its shift, intersecting at 3 or 4 real points,
depending on the direction of the shift.
The~resulting two-component curve is expressive.}
\vspace{-.1in}
\label{fig:chebyshev23-shifts}
\end{figure}

\begin{example}
\label{ex:23-shifts}
As in Example~\ref{ex:23-dilation},
let $C$ be the $(2,3)$-Chebyshev curve~\eqref{eq:23-Chebyshev}.
Applying Proposition~\ref{pr:shifts-lissajous} (with $k=2$ and $\ell=3$),
we see that the curve $C$ and its vertical shift~$C_{0,b}$ ($b\neq 0$)
can intersect in the affine plane in at most 4~points.
More generally, Proposition~\ref{pr:shifts-poly} (with $d=3$ and $d'=2$)
gives the upper bound of~4 for the number of intersection points
between $C$ and its nontrivial shift~$\Cab$.
On the other hand, in the case of a horizontal shift, we get at most $3$~points of intersection.
When these bounds are attained, with all intersection points real, the union $C\cup\Cab$ is expressive.
See Figure~\ref{fig:chebyshev23-shifts}.
\end{example}

In addition to shifts and dilations, we can consider other linear changes of variables
that can be used to construct new expressive curves.
Let us illustrate one such construction using the example of Lissajous-Chebyshev curves:

\begin{proposition}
\label{pr:rescale-lissajous}
Let $C$ be the $(k,\ell)$-Lissajous-Chebyshev curve given by~\eqref{eq:lissajous-kl},
with $\ell>k\ge 2$.
For $q\in\CC^*$, $q\neq 1$, let $C_{[q]}$ denote the curve defined by
\begin{equation}
\label{eq:lissajous-kl-q}
T_k(\tfrac{x}{q^\ell})+T_\ell(\tfrac{y}{q^k})=0.
\end{equation}
Assume that $C$ and $C_{[q]}$ intersect at $N$ points in~$\AA^2$.
Then
\begin{equation}
\label{eq:x-rescale-lissajous}
N \le k(\ell-2).
\end{equation}
\end{proposition}

\begin{proof}
Since
$T_k(x)=2^{k-1}x^k+O(x^{k-1})$ and $T_\ell(y)=2^{\ell-1}y^\ell+O(y^{\ell-1})$,
the equations defining $C$ and $C_{[q]}$ can be written as
\begin{align}
\label{eq:leading-lissajous}
2^{k-1}x^k+O(x^{k-2})+2^{\ell-1}y^\ell+O(y^{\ell-2})&=0, \\
\label{eq:leading-lissajous-q}
2^{k-1}q^{-k\ell}x^k+O(x^{k-2})+2^{\ell-1}q^{-k\ell}y^\ell+O(y^{\ell-2})&=0.
\end{align}
Multiplying \eqref{eq:leading-lissajous} by $q^{k\ell}$ and subtracting~\eqref{eq:leading-lissajous},
we get an equation of the form
\begin{equation}
\label{eq:leading-lissajous-q-1}
O(x^{k-2})+O(y^{\ell-2})=0.
\end{equation}
We thus obtain a system of two algebraic equations of the form \eqref{eq:leading-lissajous} and \eqref{eq:leading-lissajous-q-1}.
Their Newton polygons are contained in the triangles with vertices $(0,0), (k,0), (0,\ell)$
and  $(0,0), (k-2,0), (0,\ell-2)$, respectively.
The mixed area of these two triangles is equal to $k(\ell-2)$
(here we use that $\ell>k\ge 2$ and therefore $\frac{\ell-2}{k-2}>\frac{\ell}{k}$).
By Bernstein's theorem~\cite{DBernstein}, this system of equations has at most $k(\ell-2)$ solutions.
\end{proof}

\begin{figure}[ht]
\begin{center}
\vspace{-.1in}
\includegraphics[scale=0.28, trim=20cm 25cm 25cm 20cm, clip]{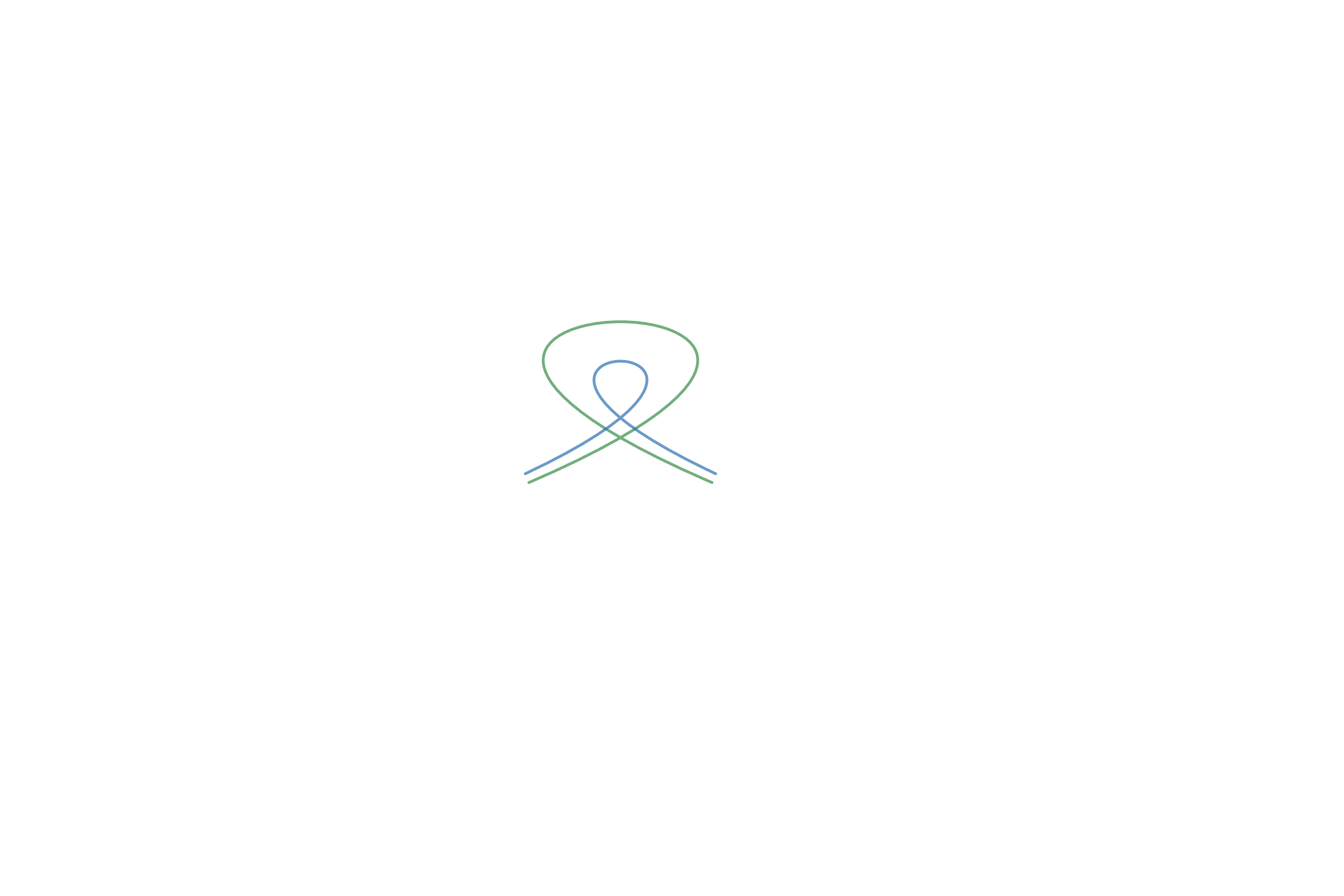}
\vspace{-.15in}
\end{center}
\caption{An expressive cubic and its rescaling~\eqref{eq:lissajous-kl-q}, intersecting at 2 real points.
The~resulting two-component curve is expressive.}
\vspace{-.15in}
\label{fig:chebyshev23-q}
\end{figure}

\begin{example}
\label{ex:23-q}
Let $C$ be the $(2,3)$-Chebyshev curve~\eqref{eq:23-Chebyshev},
as in Examples~\ref{ex:23-dilation} and~\ref{ex:23-shifts}.
Applying Proposition~\ref{pr:rescale-lissajous} (with $k=2$ and $\ell=3$),
we see that the curve $C$ and the rescaled curve~$C_{[q]}$ defined by~\eqref{eq:lissajous-kl-q}
can intersect in the affine plane in at most 2~points.
When they do intersect at 2 real points, the union $C\cup\Cab$ is expressive.
See Figure~\ref{fig:chebyshev23-q}.
\end{example}

\begin{figure}[ht]
\begin{center}
\vspace{-.1in}
\includegraphics[scale=0.28, trim=20cm 25cm 25cm 19cm, clip]{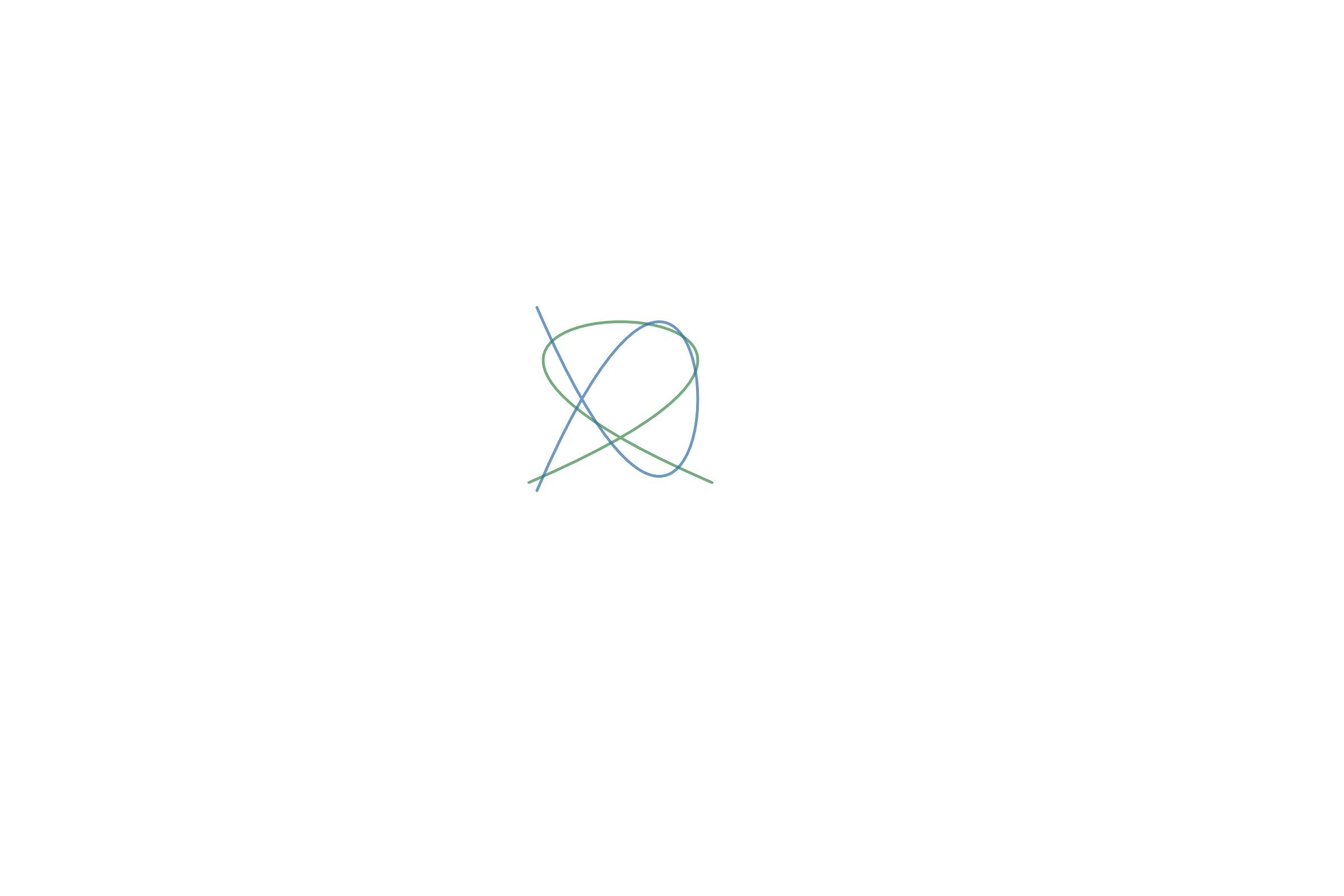}
\vspace{-.15in}
\end{center}
\caption{Two expressive cubics related by a $90^\circ$ rotation,
and intersecting at 9 real points.
The~resulting two-component curve is expressive.}
\vspace{-.15in}
\label{fig:chebyshev23-x}
\end{figure}

\begin{remark}
\label{rem:generic-crossings}
Let $C'$ be a curve obtained from a plane curve~$C=V(G(x,y))$ of degree~$d$
by an arbitrary affine change of variables:
\begin{equation*}
C'=V(G(c_{11}x + c_{12}y+a,c_{21}x+c_{22}y+b)).
\end{equation*}
Then $C$ and $C'$ intersect in~$\AA^2$ in at most~$d^2$ points.
Thus, if they intersect at $d^2$ hyperbolic nodes, then $C\cup C'$ is expressive.
See Figure~\ref{fig:chebyshev23-x}.
\end{remark}

\begin{figure}[ht]
\begin{center}
\includegraphics[scale=0.3, trim=19cm 26cm 19.5cm 24cm, clip]{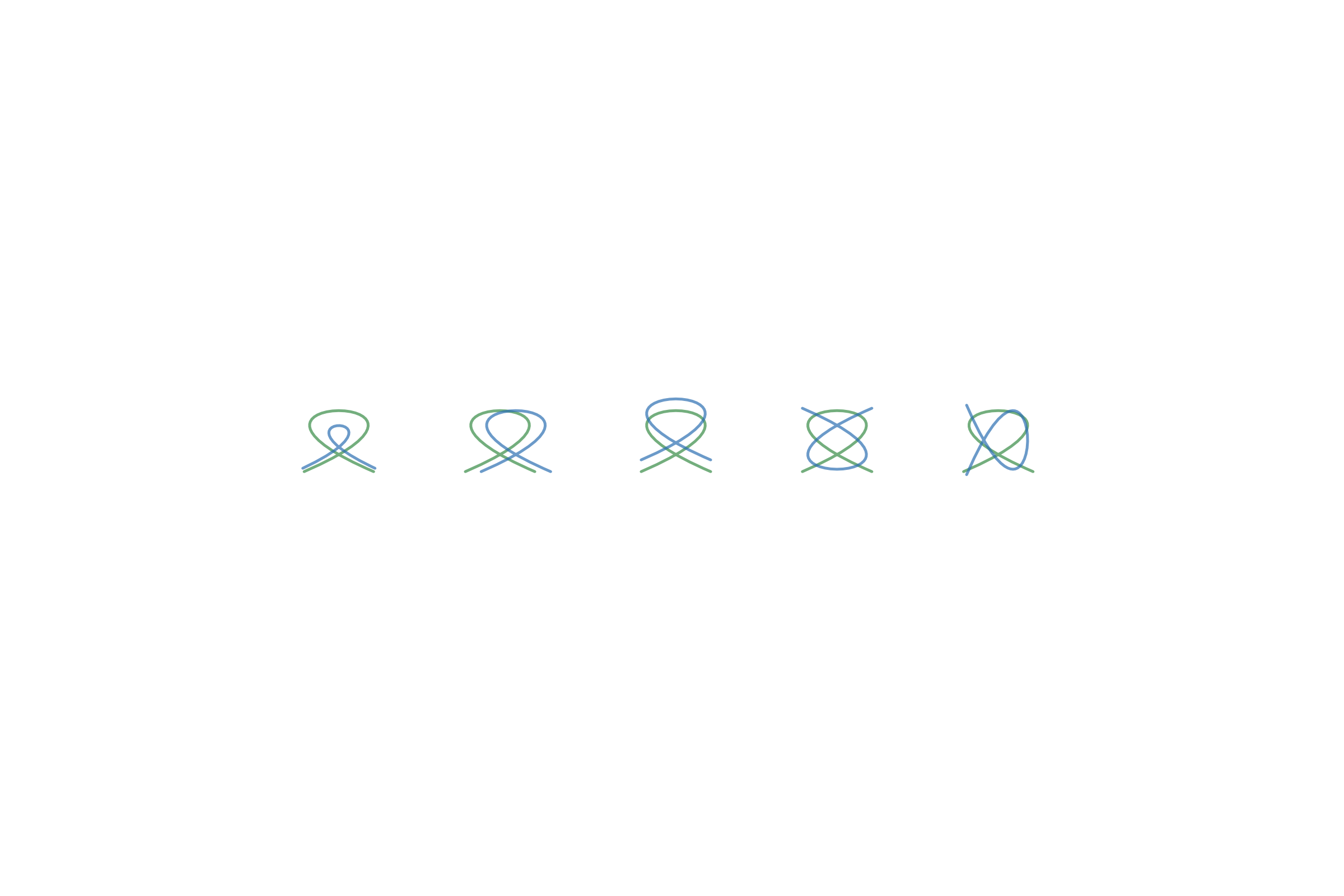} \\
\hspace{-.05in}
(a) \hspace{0.97in}
(b) \hspace{0.98in}
(c) \hspace{0.92in}
(d) \hspace{0.92in}
(e)  \\
\vspace{-.1in}
\end{center}
\caption{Two expressive cubics forming a two-component expressive curve.
The two components intersect at 2, 3, 4, 6, and 9 real points, respectively.
}
\vspace{-.1in}
\label{fig:chebyshev23-5}
\end{figure}

\begin{figure}[ht]
\begin{center}
\includegraphics[scale=0.3, trim=16.5cm 33cm 21cm 17cm, clip]{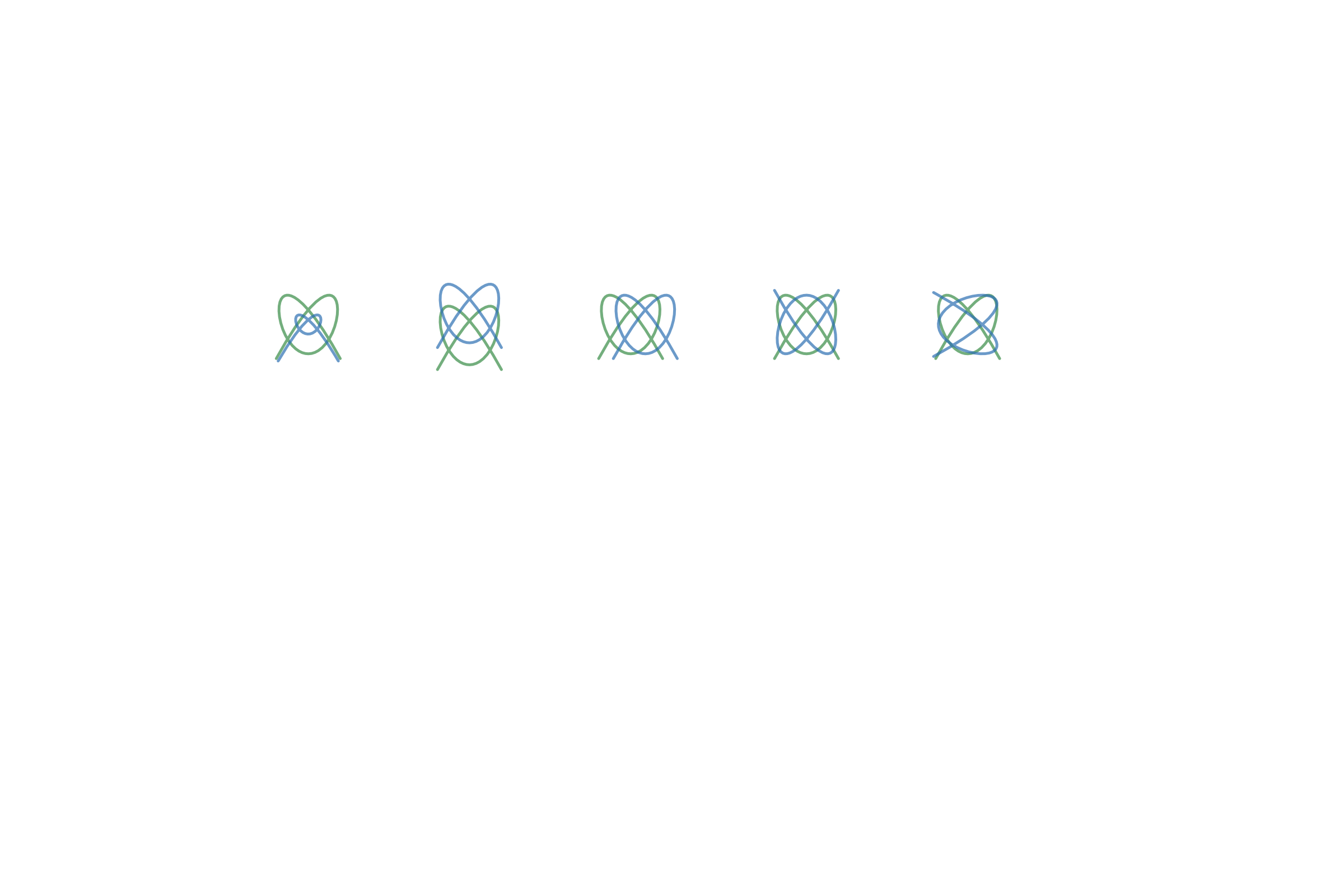} \\
\hspace{-.05in}
(a) \hspace{0.93in}
(b) \hspace{0.98in}
(c) \hspace{0.99in}
(d) \hspace{0.92in}
(e)  \\
\vspace{-.1in}
\end{center}
\caption{Two expressive quartics forming a two-component expressive curve.
The two components intersect at 6, 8, 10, 12, and 16 real points, respectively.
}
\vspace{-.1in}
\label{fig:chebyshev34-5}
\end{figure}

\begin{example}
\label{ex:5-chebyshevs}
Figures~\ref{fig:chebyshev23-5}--\ref{fig:chebyshev34-5}
show five different ways to arrange two Chebyshev curves
(the $(2,3)$-Chebyshev cubic and the $(3,4)$-Chebyshev quartic, respectively)
related to each other by an affine transformation of the plane~$\AA^2$
so that the resulting two-component curve is expressive.
These pictures illustrate:
\begin{itemize} 
\item[(a)]
Proposition~\ref{pr:rescale-lissajous}, cf.\ Example~\ref{ex:23-q};
\item[(b,c)]
Proposition~\ref{pr:shifts-lissajous}, cf.\ Example~\ref{ex:23-shifts};
\item[(d)]
Proposition~\ref{pr:dilations-poly}, cf.\ Example~\ref{ex:23-dilation}; and
\item[(e)]
Remark~\ref{rem:generic-crossings}.
\end{itemize}
\end{example}

\begin{figure}[ht]
\begin{center}
\includegraphics[scale=0.28, trim=20cm 23.5cm 22cm 15.5cm, clip]{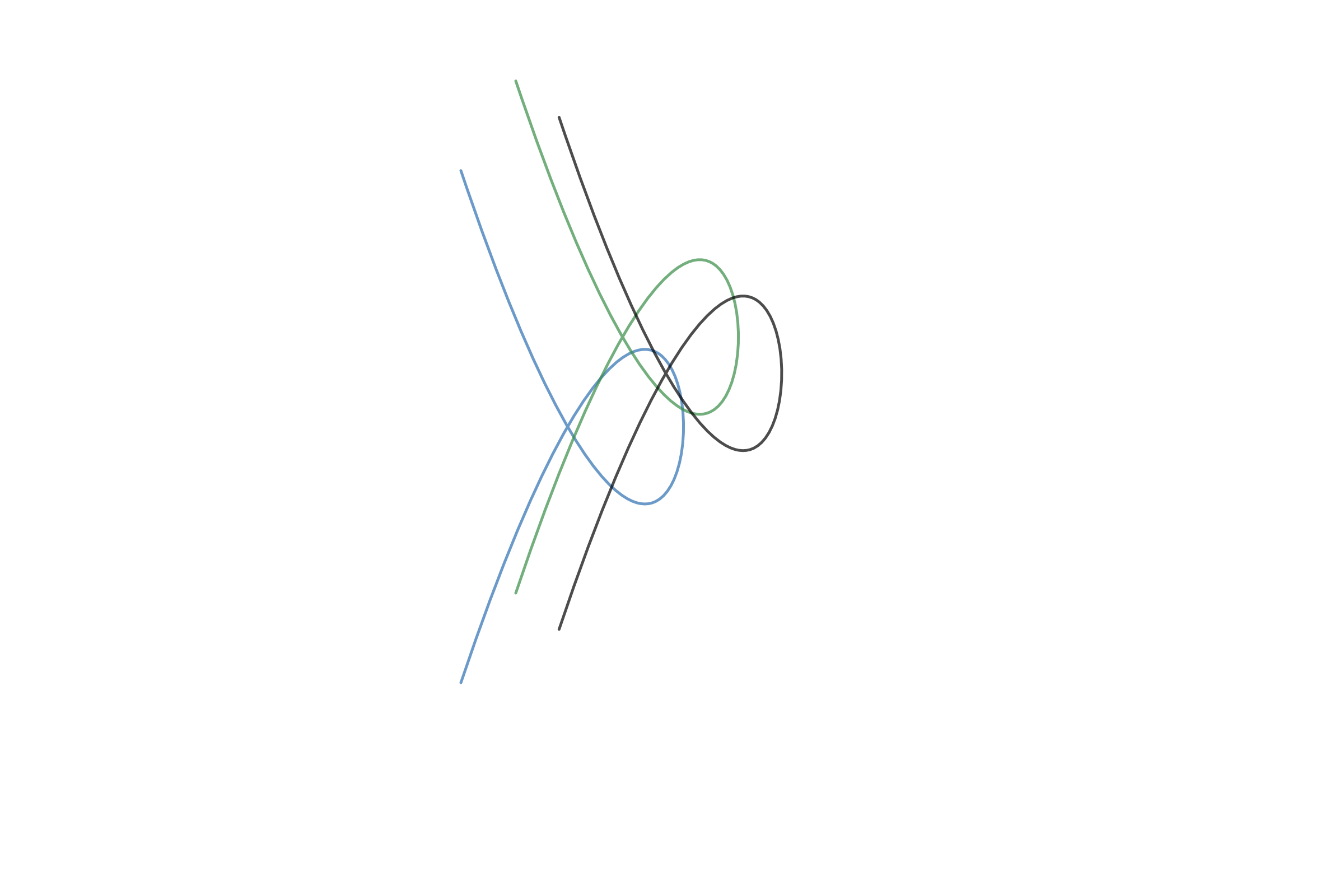}
\end{center}
\caption{An expressive curve whose three components are translations of the same nodal cubic.
Each pair of components intersect at 4 hyperbolic nodes.}
\vspace{-.1in}
\label{fig:three-nodal-cubics}
\end{figure}

\begin{figure}[ht]
\begin{center}
\includegraphics[scale=0.28, trim=20cm 22cm 22cm 16.5cm, clip]{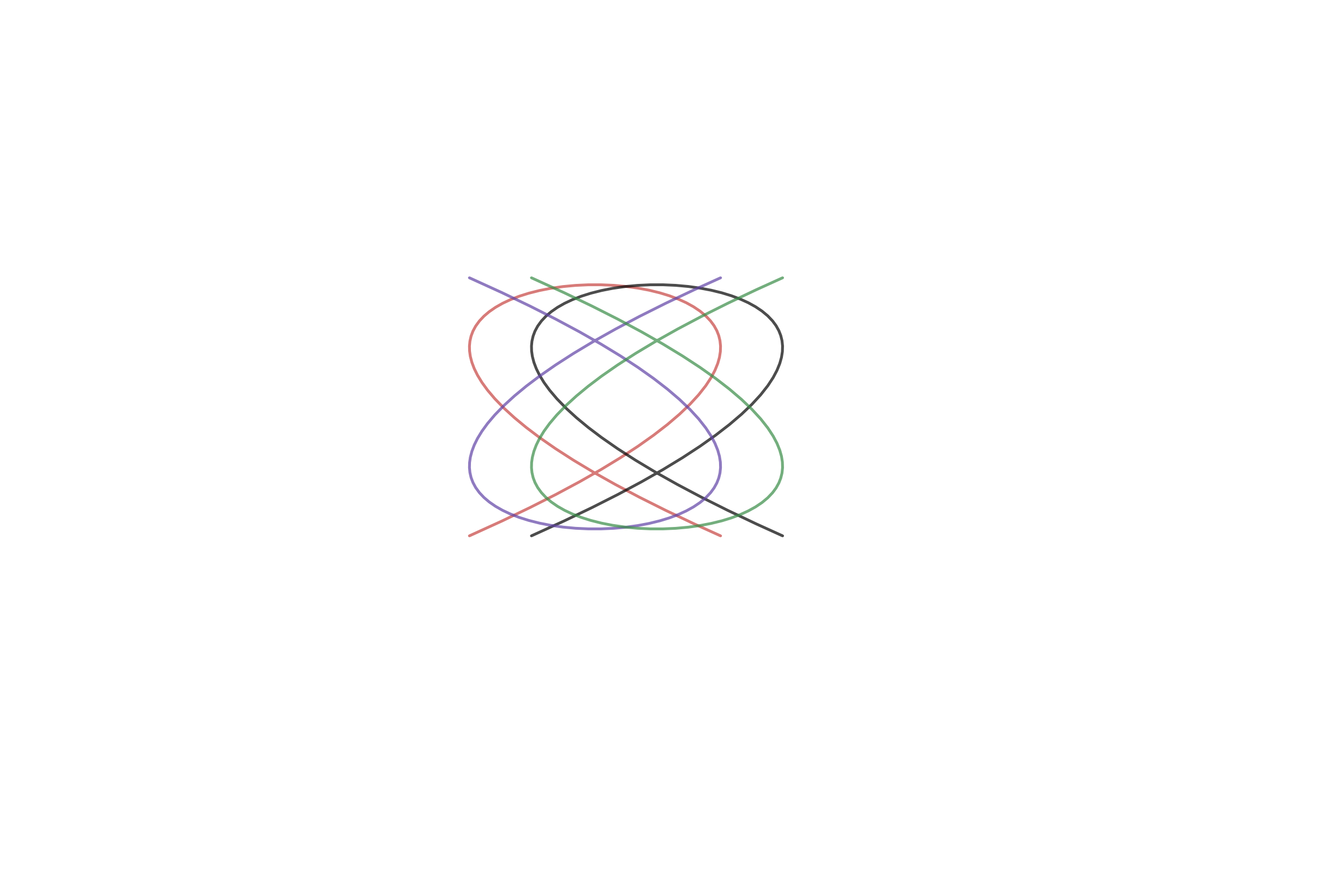}
\end{center}
\caption{An expressive curve whose four components are singular cubics related to each other
either by a horizontal translation or a dilation (with $c=-1$).
Each pair of components intersect at 3 or 6 hyperbolic nodes, respectively.}
\vspace{-.1in}
\label{fig:four-cubics}
\end{figure}

\begin{remark}
\label{rem:arrangements-poly}
More generally, consider a collection of expressive polynomial curves related to each other
by affine changes of variables (equivalently, affine transformations of the plane~$\AA^2$).
Suppose that for every pair of curves in this collection,
the number of hyperbolic nodes in their intersection attains the upper bound
for an appropriate version of
Proposition~\ref{pr:dilations-poly}, \ref{pr:shifts-poly}, \ref{pr:shifts-lissajous},
\ref{pr:rescale-lissajous}, or Remark~\ref{rem:generic-crossings}.
Then the union of the curves in the given collection is an expressive curve.
See Figures~\ref{fig:three-nodal-cubics}--\ref{fig:four-cubics}.
\end{remark}

\clearpage

\newpage

\section{Arrangements of trigonometric curves}
\label{sec:arrangements-trig}

In this section, we generalize Example~\ref{ex:circle-arrangements} to arrangements of curves
obtained from a given trigonometric curve by shifts, dilations, and/or rotations.

We continue to use the notation \eqref{eq:Cab}--\eqref{eq:Cabs=V(G)} for the shifted and dilated curves.

\begin{corollary}
\label{cor:shift-dilate-trig-2pts-transversal}
Let $C$ be a trigonometric curve of degree~$2d$,
with two local branches at infinity centered at distinct points $p,\overline{p}\in C\cap L^\infty$.
Suppose that $\mt(C,p)=d$, i.e., these branches are transversal to~$L^\infty$.
Then we have, for $a,b\in\CC$ and $c\in\CC^*$:
\begin{equation}
\label{eq:dilate-trig-transversal}
(C\cdot \Cabc)_p \ge d^2.
\end{equation}
If $C$ and $\Cabc$ intersect in $N$ points in the affine plane~$\AA^2$,
\redsf{counting with multiplicity,}
then $N\le 2d^2$.
\end{corollary}

\begin{proof}
We apply Proposition~\ref{pr:gen-dilation}, with $S^-(\Gamma(C,p))=\frac12 d^2$,
to obtain~\eqref{eq:dilate-trig-transversal}.
It follows that $N\le (2d)^2-2\cdot \frac12 d^2=2d^2$.
\end{proof}

\begin{figure}[ht]
\begin{center}
\includegraphics[scale=0.27, trim=16cm 23.5cm 16cm 24cm, clip]{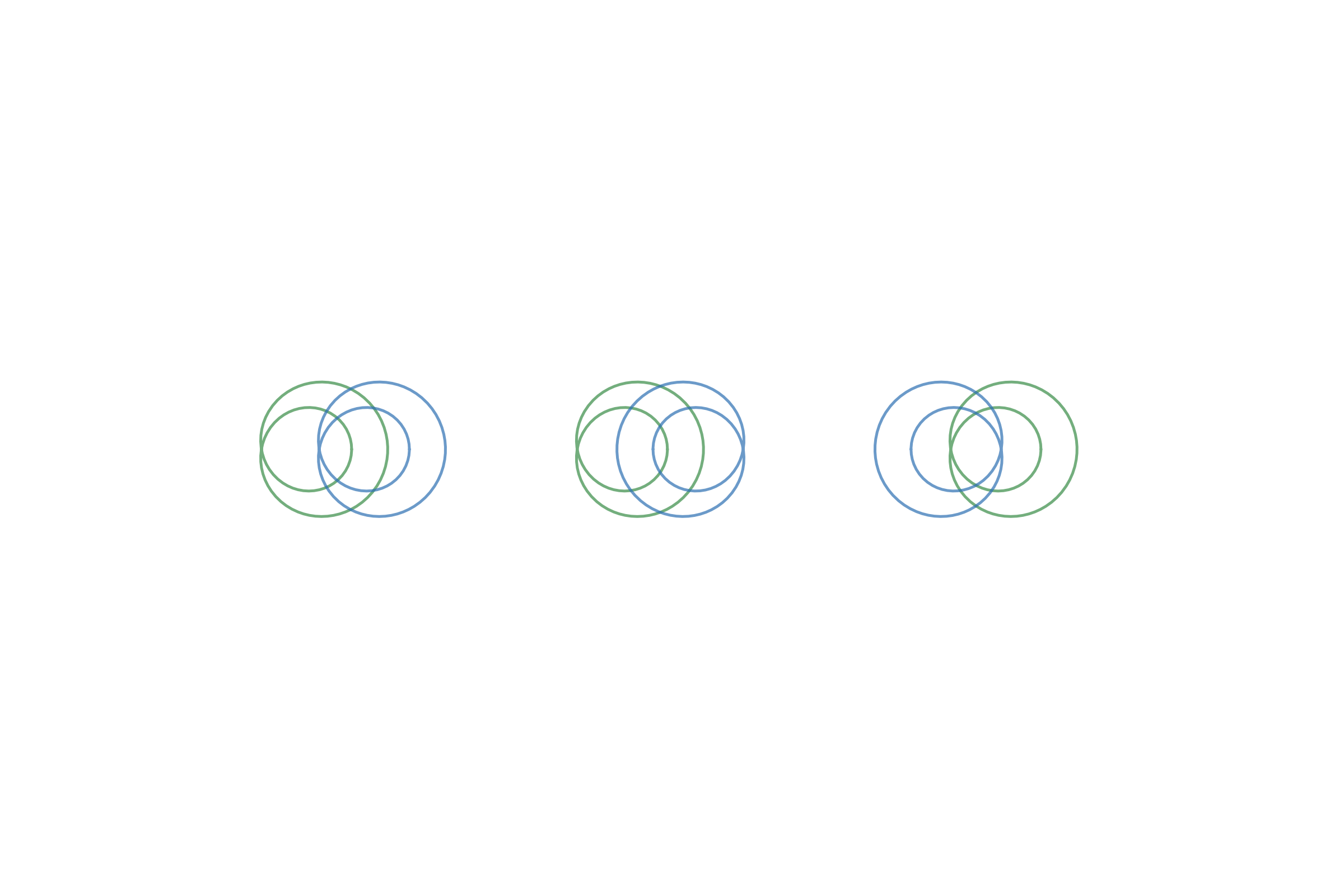}
\vspace{-.1in}
\end{center}
\caption{A union of two lima\c cons related to each other by a shift or dilation
is expressive if they intersect at 8 hyperbolic nodes.}
\vspace{-.1in}
\label{fig:two-limacons-3}
\end{figure}

\begin{example}
\label{ex:limacon-shifts}
Let $C$ be an epitrochoid with parameters $(2,-1)$, i.e., a lima\c con.
It has two conjugate points at infinity, each an ordinary cusp transversal to~$L^\infty$.
This is a quartic trigonometric curve, so by~Corollary~\ref{cor:shift-dilate-trig-2pts-transversal}
(with $d=2$), any two shifts/dilations of~$C$ intersect in at most 8~points in the affine plane.
Thus, if they intersect at 8 hyperbolic nodes, then $C\cup \Cabc$ is expressive by Theorem~\ref{th:reg-expressive}.
More generally, an arrangement of lima\c cons related to each other by shifts and dilations
gives an expressive curve if any two of these lima\c cons intersect at 8 hyperbolic nodes.
See Figure~\ref{fig:two-limacons-3}.
\end{example}

\begin{corollary}
\label{cor:shift-dilate-trig-2pts-tangent}
Let $C$ be a trigonometric curve of degree~$2d$,
with two local branches at infinity centered at distinct points $p,\overline{p}\in C\cap L^\infty$.
Suppose that $m=\mt(C,p)<d$, i.e., $C$ is tangent to~$L^\infty$.
Then we have, for $a,b\in\CC$ and $c\in\CC^*$:
\begin{align}
\label{eq:dilate-trig}
(C\cdot \Cabc)_p &\ge dm, \\
\label{eq:shift-trig}
(C\cdot \Cab)_p &\ge (d-1)(m+1)+\gcd(d,m).
\end{align}
If $C$ and $\Cabc$ (resp., $\Cab$)
intersect in~$N$ (resp.,~$M$) points in the affine plane~$\AA^2$, then
\begin{align}
\label{eq:}
 N &\le 4d^2-2dm; \\
 M &\le 4d^2-2(d-1)(m+1)-2\gcd(d,m).
\end{align}
\end{corollary}

\begin{proof}
The proof is analogous to the proof of Corollary~\ref{cor:shift-poly}.
We apply Propositions \ref{pr:gen-dilation} and~\ref{pr:gen-shift} to
the case under consideration, taking into account that
\begin{align*}
S^-(\Gamma(C,p))&=\tfrac12 dm, \\
(C\cdot L^\infty)_p&=d,\\
r(Q)=\eta(Q)&=\gcd(d,m).
\qedhere
\end{align*}
\end{proof}

\begin{figure}[ht]
\begin{center}
\includegraphics[scale=0.28, trim=16cm 25.5cm 24cm 14cm, clip]{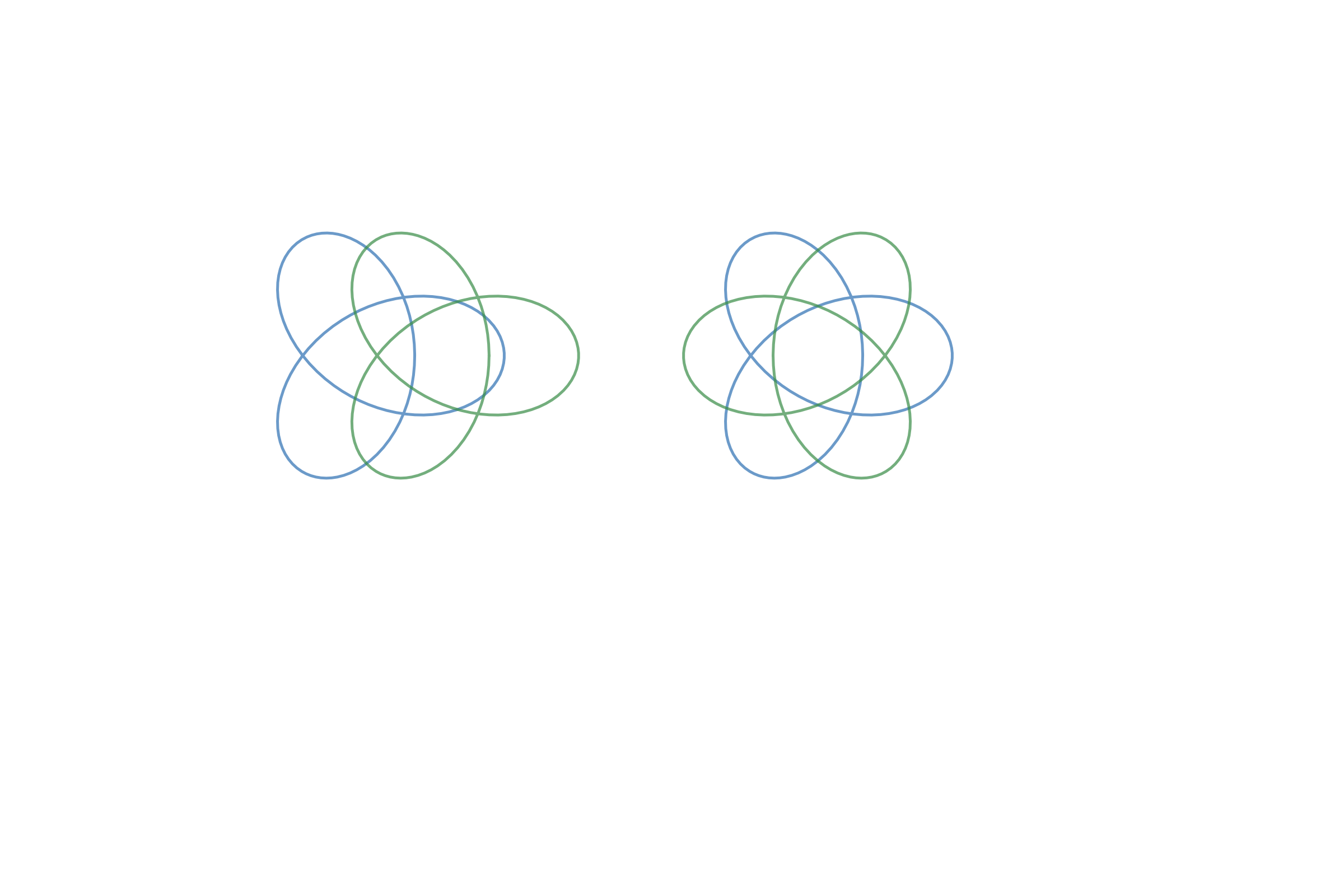}
\vspace{-.1in}
\end{center}
\caption{A union of two 3-petal hypotrochoids related to each other by a shift (resp., dilation)
is expressive if they intersect at~10 (resp.,~12) hyperbolic nodes.}
\vspace{-.1in}
\label{fig:two-hypotrochoids}
\end{figure}

\begin{example}
\label{ex:hypotrochoid-shifts}
Let $C$ be a hypotrochoid with parameters $(2,1)$, cf.\ Examples~\ref{ex:hypotrochoids}
and~\ref{ex:hypotrochoids-again}.
It has two conjugate points at infinity; at each of them, $C$~is smooth and has a simple
(order~2) tangency to~$L^\infty$.
By 
\eqref{eq:dilate-trig}
(with $d=2$ and $m=1$), any dilation of~$C$ intersects~$C$ in the affine plane~$\AA^2$ in at most 12~points.
Similarly, by~\eqref{eq:shift-trig}, any shift of~$C$ intersects~$C$ in~$\AA^2$ in at most 10~points.
When these bounds are attained, and all intersections are hyperbolic nodes,
the union of the two curves is expressive.
See Figure~\ref{fig:two-hypotrochoids}.
\end{example}

\begin{corollary}
\label{cor:shift-dilate-trig-1pt}
Let $C=Z(F)$ be a trigonometric projective curve of degree~$2d$,
with two local complex conjugate branches $Q,\overline Q$ centered at the same point
$p\in C\cap L^\infty$.
Denote $\mt(C,p)=2\mt Q=2m$. Then we have, for $a,b\in\CC$ and $c\in\CC^*$:
\begin{align}
\label{eq:4dm}
(C\cdot \Cabc)_p &\ge 4dm, \\
\label{eq:shift-trig-1pt}
(C\cdot \Cab)_p &\ge
\begin{cases}
4dm-2m+2d+2\gcd(d,m)-2 &\text{if}\ (Q\cdot\overline Q)_p=dm,\\
4dm-2m+2d+2\gcd(d,m)    &\text{if}\ (Q\cdot\overline Q)_p>dm.
\end{cases}
\end{align}
\end{corollary}

\begin{proof}
Once again, we apply Propositions~\ref{pr:gen-dilation} and~\ref{pr:gen-shift}, with
\begin{align*}
S^-(\Gamma(C,p))&=2dm, \\
(C\cdot L^\infty)_p&=2d,\\
r(Q)&=\gcd(d,m), \\
\eta(Q)&=\begin{cases}\gcd(d,m),\ &\text{if}\ (Q\cdot\overline Q)_p=dm,\\
2\gcd(d,m),\ &\text{if}\ (Q\cdot\overline Q)_p>dm.\end{cases}
\end{align*}
and similarly for $\overline Q$.
\end{proof}

\pagebreak[3]

\begin{figure}[ht]
\begin{center}
\includegraphics[scale=0.18, trim=22cm 21cm 24cm 18cm, clip]{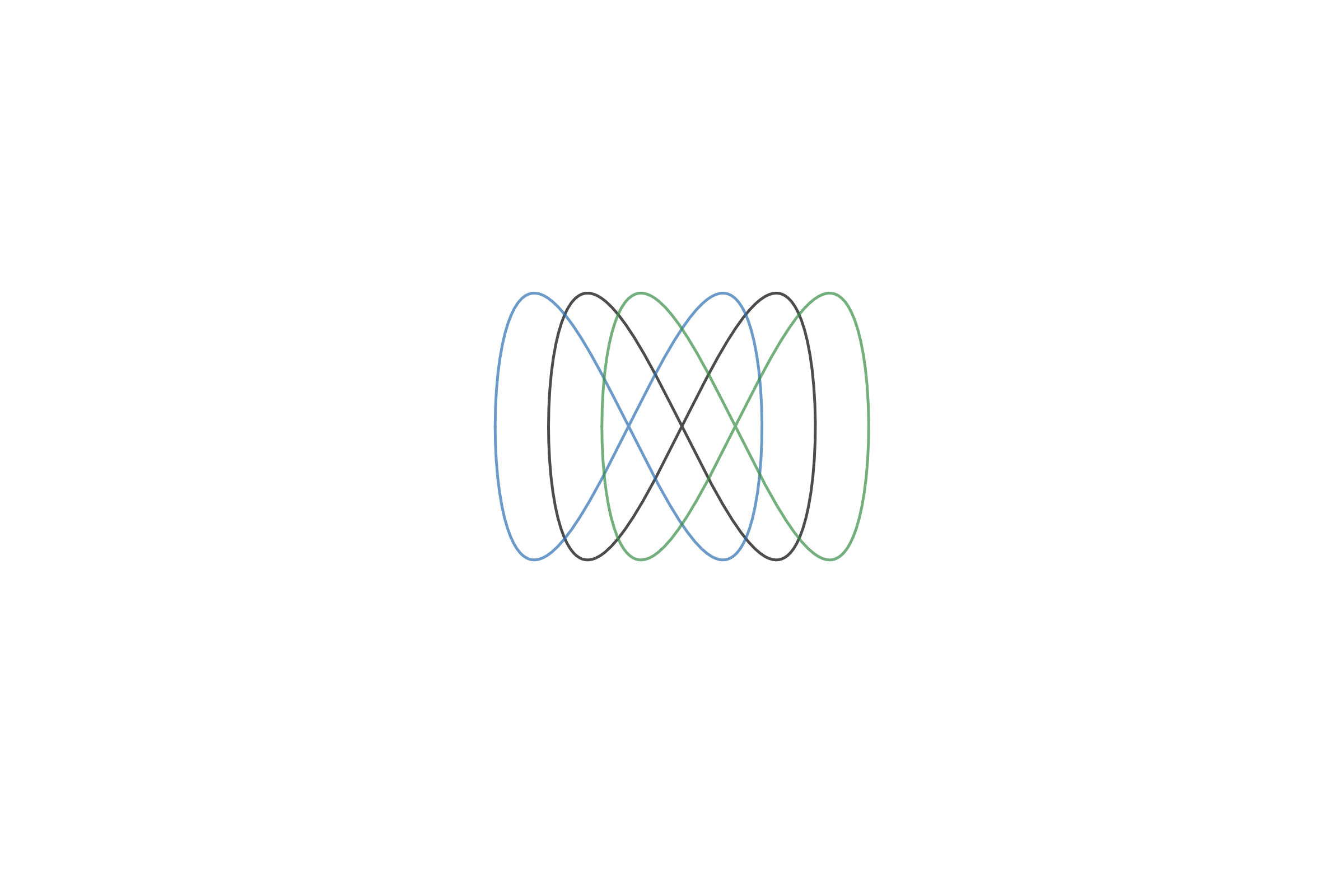}
\includegraphics[scale=0.19, trim=24cm 19cm 26cm 21.5cm, clip]{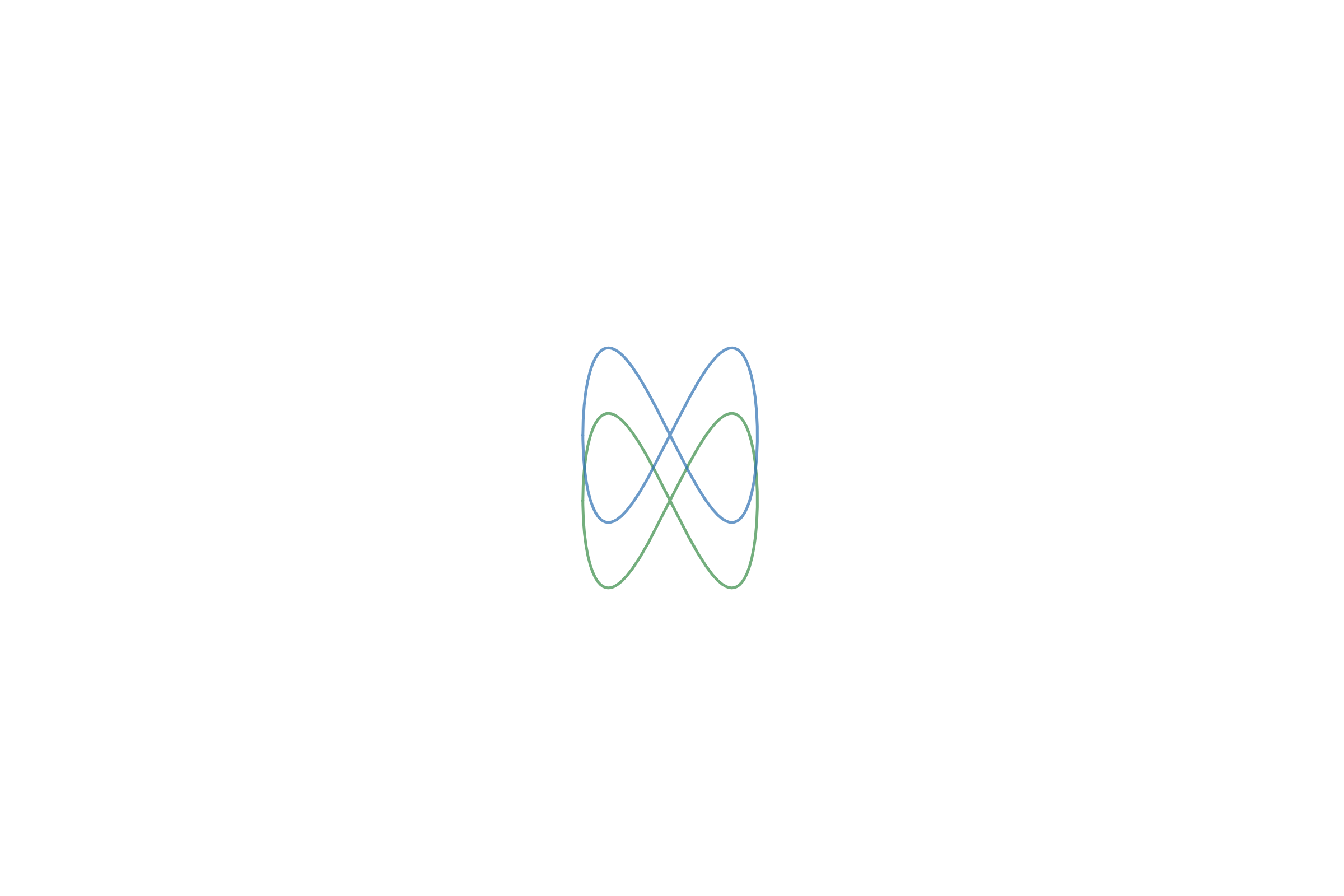}
\vspace{-.05in}
\end{center}
\caption{Left: an expressive curve whose three components are translations of the same lemniscate.
Each pair of components intersect at six hyperbolic~nodes.
Right: two lemniscates differing by a vertical shift, see Example~\ref{ex:lemniscate-vertical}.}
\vspace{-.1in}
\label{fig:three-lemniscates}
\end{figure}

\begin{example}
\label{ex:shifts-of-lemniscates}
Let $C$ be a lemniscate of Huygens
\begin{equation}
\label{eq:huygens-again}
y^2+4x^4-4x^2=0,
\end{equation}
see Example~\ref{ex:lemniscates}.
It has a single point~$p=(0,1,0)$ at infinity, with two conjugate local branches $Q$ and~$\overline Q$.
These branches are tangent to each other and to~$L^\infty$;
all these tangencies are of order~2.
Thus Corollary~\ref{cor:shift-dilate-trig-1pt} applies, with $d=2$, $m=1$, and $(Q\cdot \overline Q)=2=dm$.
The bound~\eqref{eq:shift-trig-1pt} yields $(C\cdot \Cab)_p\ge 10$,
implying that $C$ and a shifted curve~$\Cab$ intersect in~$\AA^2$ in at most 6~points.
Thus, any arrangement of shifts of~$C$ which intersect pairwise transversally in 6~real points
produces an expressive curve
(assuming all these double points are distinct).
See Figure~\ref{fig:three-lemniscates}.
\end{example}

For special choices of shifts and dilations,
the bounds in the corollaries above~can be strengthened,
leading to examples of expressive curves whose components intersect~in fewer real points
than one would ordinarily expect.
Here are two examples:

\begin{example}
\label{ex:lemniscate-vertical}
Let $C$ be the lemniscate~\eqref{eq:huygens-again}.
Since $C$ is a Lissajous-Chebyshev curve with parameters~$(4,2)$,
by Proposition~\ref{pr:shifts-lissajous}, its vertical shift $C_{0,b}$ intersects $C$ in~$\AA^2$
in at most 4~points.
Hence $C \cup C_{0,b}$ is expressive for $b\in(-2,2)$. See Figure~\ref{fig:three-lemniscates}.
\end{example}

\begin{figure}[ht]
\begin{center}
\includegraphics[scale=0.19, trim=20cm 23.5cm 22cm 21cm, clip]{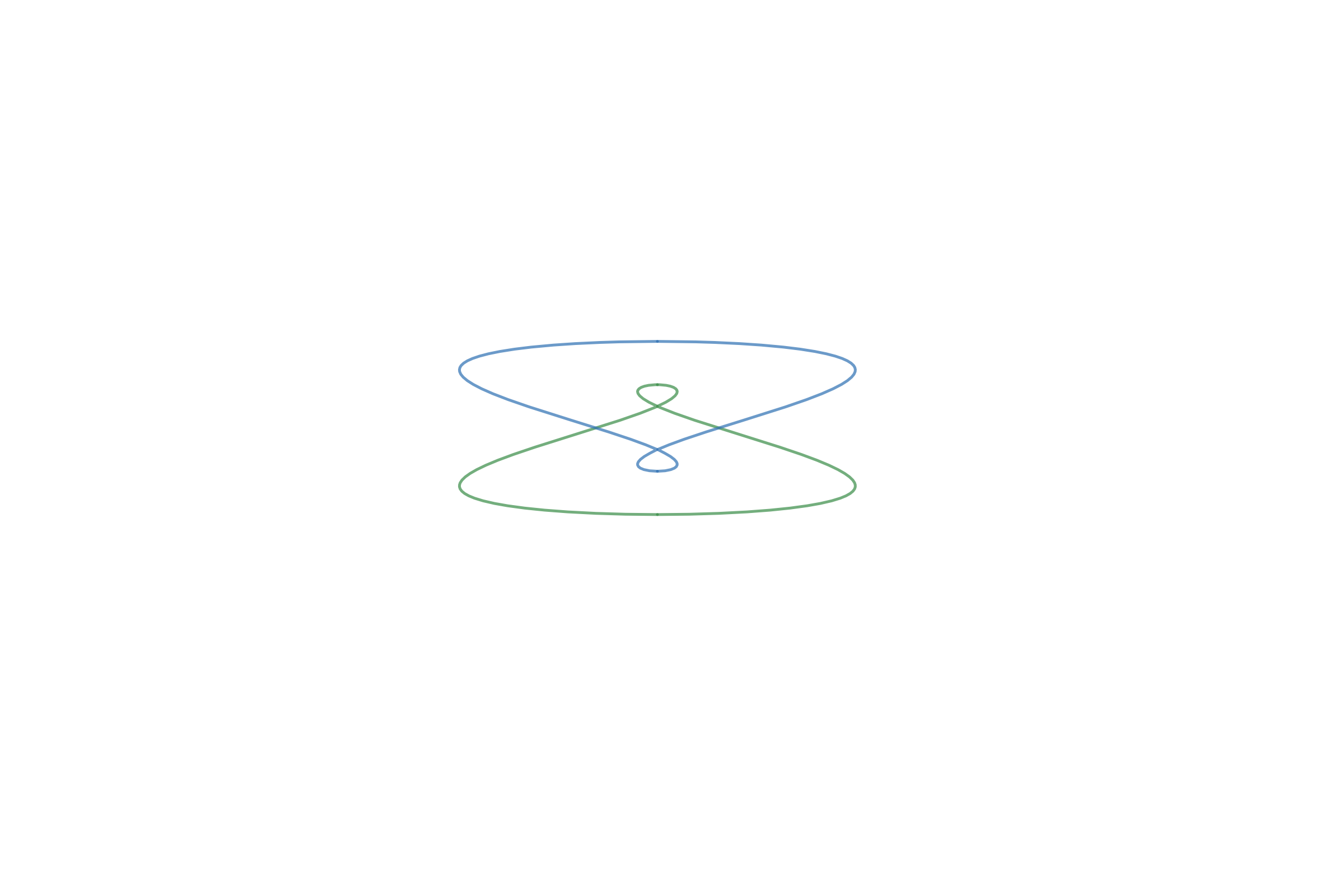}
\vspace{-.1in}
\end{center}
\caption{A two-component expressive curve $C\cup C_{-1}$
from Example~\ref{ex:two-exotic-lemniscates}.}
\vspace{-.1in}
\label{fig:two-exotic-lemniscates}
\end{figure}

\begin{example}
\label{ex:two-exotic-lemniscates}
Let $C=V(G(x,y))$ be the trigonometric curve  defined by the polynomial
\begin{equation}
\label{eq:exotic-lemniscate}
G(x,y)=x^2+y^4-11y^2+18y-8=x^2+(y+4)(y-1)^2(y-2).
\end{equation}
It is easy to see that $C$ is expressive,
and that $C$ intersects its dilation/reflection $C_{-1}=V(G(-x,-y))$ in two points in~$\AA^2$,
both of which are real hyperbolic nodes.
Hence the union $C\cup C_{-1}$ is expressive.
See Figure~\ref{fig:two-exotic-lemniscates}.
\end{example}

\newpage

\section{Alternative notions of expressivity}
\label{sec:alternative-expressivity}

In this section, we discuss two alternative notions of expressivity.
For the first notion, algebraic curves are treated as subsets of~$\RR^2$,
instead of the scheme-theoretic point of view that we adopted above in this paper.
For the second notion, bivariate polynomials are replaced by arbitrary smooth functions
of two real variables.

\medskip

Viewing real algebraic curves set-theoretically,
as ``topological curves'' in the real affine~plane, we arrive at the following definition:

\begin{definition}
\label{def:expressive-curve}
Let $\mathcal{C}\subset\RR^2$ be the set of real points of a real affine algebraic curve,
see Definition~\ref{def:curves}.
Assume that $\mathcal{C }$ is nonempty, with no isolated points.
We say that $\mathcal{C }$ is \emph{expressive}
if its (complex) Zariski closure~$C=\overline{\mathcal{C}}$ is an expressive plane algebraic curve.
Thus, a subset $\mathcal{C }\subset\RR^2$ is expressive if
\begin{itemize}[leftmargin=.2in]
\item
$\mathcal{C }$ is the set of real points of a real affine plane algebraic curve,
\item
$\mathcal{C }$ is nonempty, with no isolated points, and
\item
the minimal polynomial of $\mathcal{C }$ is expressive, see Definition~\ref{def:expressive-poly}.
\end{itemize}
\end{definition}

As always, one should be careful when passing from a real algebraic set to an algebraic curve,
or the associated polynomial.
A polynomial $G(x,y)\!\in\!\RR[x,y]$ can be expressive while the real algebraic set~$\VR(G)$ is not;
see Example~\ref{example1} below.
Conversely, $\VR(G)$ can be expressive while $G(x,y)$ is~not, see Example~\ref{ex:xy(x^2+y^2+1)}.
That's because $G$ may not be the minimal polynomial for~$\VR(G)$.

\begin{example}[cf.\ Examples~\ref{ex:(x^2+z^2)(yz^2-x^3+x^2y)},
\ref{ex:(x^2+z^2)(yz^2-x^3+yx^2)-regular}, and~\ref{ex:reducible-quintic}]
\label{example1}
The real polynomials 
\begin{align*}
G(x,y)&=(x^2+1)(x^2y-x^3+y), \\
\widetilde G(x,y)&=x^2y-x^3+y,
\end{align*}
define the same (connected) real algebraic set
\[
\mathcal{C }= \VR(G)=\VR(\widetilde G)
=\{(x,y)\in\RR^2\mid y=\tfrac{x^3}{x^2+1}\},
\]
cf.\ \eqref{eq:y=x^3/...}.
As we saw in Example~\ref{ex:reducible-quintic},
$G$ is expressive while $\widetilde G$ is not.
Consequently, the affine algebraic curve~$V(G)$ is expressive
while its real point set, the real topological curve~$\mathcal{C }$ is not---because
$\widetilde G$, rather than~$G$, is the minimal polynomial of~$\mathcal{C }$.
\end{example}

\begin{example}
\label{ex:xy(x^2+y^2+1)}
The real polynomials $G(x,y)$ and $\tilde G(x,y)$ given by
\begin{align*}
G(x,y)&=xy(x^2+y^2+1), \\
\widetilde G(x,y)&=xy,
\end{align*}
define the same real algebraic set $\mathcal{C }\!=\! \VR(G)\!=\!\VR(\widetilde G)$.
Clearly, $\widetilde G$ is expressive, and so is the affine curve $V(\widetilde G)$,
or the projective curve~$Z(xy)$.
Since $\widetilde G$ is the minimal polynomial of~$\mathcal{C }$,
this real algebraic set is expressive as well.
On~the other~hand,
\begin{align*}
G_x&= 3x^2y+y^3+y=y(3x^2+y^2+1),  \\
G_y&=x^3+3xy^2+x=x(x^2+3y^2+1),
\end{align*}
and we see that $G$ has 9 critical points:
\[
\text{$(0,0)$, $(0,\pm i)$, $(\pm i,0)$,
$(\tfrac12 i, \pm\tfrac12 i)$, $(-\tfrac12 i, \pm\tfrac12 i)$.}
\]
Since 8~of these points are not real, the polynomial~$G$ is not expressive;
nor are the curves $V(G)$ and $Z(xy(x^2+y^2+z^2))$.
\end{example}

\pagebreak[3]

The following criterion is a direct consequence of Theorem~\ref{th:reg-expressive}.

\begin{corollary}
Let $\mathcal{C }\subset\RR^2$ be the set of real points of a real affine algebraic~curve.
Assume that $\mathcal{C }$ is connected, and contains at least two (hence infinitely many) points.
Then the following are equivalent:
\begin{itemize}[leftmargin=.2in]
\item
the minimal polynomial of $\mathcal{C }$
(or the Zariski closure $C=\overline{\mathcal{C }}$, cf.\ Definition~\ref{def:expressive-curve})
is expressive and $L^\infty$-regular;
\item
each component of $C=\overline{\mathcal{C}}$ is trigonometric or polynomial,
and all singular points of $C$ in the complex affine plane~$\AA^2$ are real hyperbolic nodes.
\end{itemize}
\end{corollary}

We conclude this section by discussing the challenges involved in
extending the notion of expressivity to  arbitrary
(i.e., not necessarily polynomial)
smooth real functions of two real arguments.

\begin{remark}
\label{rem:expressive-analytic}
Let $G:\RR^2\to\RR$ be a smooth function.
Suppose that $G$ satisfies the following conditions:
\begin{itemize}[leftmargin=.2in]
\item
the set $\VR(G)=\{G(x,y)=0\}\subset\RR^2$ is connected;
\item
$\VR(G)$ is a union of finitely many immersed circles and open intervals
which intersect each other and themselves transversally, \linebreak[3]
in a finite number of points;
\item
the complement $\RR^2\setminus\VR(G)$ is a union of a finite number of disjoint open~sets;
\item
all critical points of $G$ in $\RR$ have a nondegenerate Hessian;
\item
these critical points are located as follows:
\begin{itemize}[leftmargin=.2in]
\item[$\circ$]
one critical point inside each bounded connected component of $\RR^2\setminus\VR(G)$;
\item[$\circ$]
no critical points inside each unbounded connected component of $\RR^2\setminus\VR(G)$;
\item[$\circ$]
a saddle at each double point of $\VR(G)$.
\end{itemize}
\end{itemize}
One may be tempted to call such a (topological, not necessarily algebraic) curve $\VR(G)$ expressive.
Unfortunately, this definition turns out to be problematic,
as one and the same curve $\mathcal{C}=\RR^2$
can be defined by two different smooth functions one of which satisfies the above-listed conditions
whereas the other does not.
An example is shown in Figure~\ref{fig:isolines},
with the polynomial $G$ given by
\[
G(x,y)=(\tfrac{x^{2}}{16}+y^{2}-1)((x-1)^{2}+(y-1)^{2}-1).
\]
As Figure~\ref{fig:isolines} demonstrates,
the function $G$ is not expressive in any reasonable sense.
At the same time, $\mathcal{C}$~can be transformed into an expressive curve
(a~union of two circles) by a diffeomorphism of~$\RR^2$,
and consequently can be represented as the vanishing set of a smooth function
satisfying the conditions listed above.
\end{remark}

\pagebreak

\section{Regular-expressive divides}
\label{sec:regular-expressive-divides}

It is natural to wonder which divides arise from expressive curves (perhaps satisfying additional technical conditions).

If $G(x,y)$ is an expressive polynomial, then all singular points of~$V(G)$ are (real) hyperbolic nodes, so the divide $D_G$ is well defined, see Definition~\ref{def:divide-polynomial}.
The class of divides arising in this way
(with or without the additional requirement of $L^\infty$-regularity)
is however too broad to be a natural object of study:
as Example~\ref{example1} demonstrates, a non-expressive polynomial may become expressive
upon multiplication by a polynomial with an empty set of real zeroes.

With this in mind, we propose the following definition.

\begin{definition}
\label{def:regular-expressive-divide}
A divide $D$ is called \emph{regular-expressive}
if there exists an $L^\infty$-regular expressive curve $C=V(G)$
with real irreducible components such that $D=D_G$.
\end{definition}

\begin{remark}
Some readers might prefer to just call such divides ``expressive'' rather than regular-expressive.
We decided against the shorter term, as it would misleadingly
omit a reference to the $L^\infty$-regularity requirement.
\end{remark}

By Theorem~\ref{th:reg-expressive}, a connected divide is regular-expressive
if and only if it arises from an algebraic curve
whose components are real and either polynomial or trigonometric,
and all of whose singular points in the affine plane are hyperbolic nodes.

We note that a regular-expressive divide is isotopic to
an expressive algebraic subset of~$\RR^2$, in the sense of
Definition~\ref{def:expressive-curve},
or more precisely to an intersection of such a subset with a sufficiently large disk
$\Disk_R$, see~\eqref{eq:Disk_R}.

Numerous examples of regular-expressive divides are scattered throughout this paper.

\begin{problem}
\label{problem:reg-expr-divide}
Find a criterion for deciding whether a given divide is regular-expressive.
\end{problem}

\newsavebox{\crossing}
\setlength{\unitlength}{0.8pt}
\savebox{\crossing}(10,10)[bl]{
\thicklines
\qbezier(5,5)(7,10)(10,10)
\qbezier(5,5)(3,0)(0,0)
\qbezier(5,5)(3,10)(0,10)
\qbezier(5,5)(7,0)(10,0)
}

\newsavebox{\opening}
\setlength{\unitlength}{0.8pt}
\savebox{\opening}(10,10)[bl]{
\thicklines
\qbezier(0,0)(-5,0)(-5,5)
\qbezier(0,10)(-5,10)(-5,5)
}

\newsavebox{\closing}
\setlength{\unitlength}{0.8pt}
\savebox{\closing}(10,10)[bl]{
\thicklines
\qbezier(0,0)(5,0)(5,5)
\qbezier(0,10)(5,10)(5,5)
}

\begin{remark}
Problem~\ref{problem:reg-expr-divide} appears to be very difficult.
It seems to be even harder to determine whether a given divide can be realized
by an expressive curve of a specified degree.
For example, the divide
\setlength{\unitlength}{0.8pt}
\begin{picture}(26,10)(-5,-5)
\put(0,0){\makebox(0,0){\usebox{\opening}}}
\put(0,0){\makebox(0,0){\usebox{\crossing}}}
\put(10,0){\makebox(0,0){\usebox{\crossing}}}
\put(20,0){\makebox(0,0){\usebox{\closing}}}
\end{picture}
can be realized by an expressive sextic (the $(2,6)$-Lissajous curve, see Figure~\ref{fig:lissajous-chebyshev})
but not by an expressive quadric---even though there exists a (non-expressive) quadric
realizing this divide.
\end{remark}

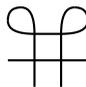
\begin{figure}[ht]
\begin{center}
\setlength{\unitlength}{1pt}
\begin{picture}(30,31)(0,0)
\thicklines
\put(0,10){\line(1,0){30}}
\put(10,20){\line(1,0){10}}
\put(10,0){\line(0,1){20}}
\put(20,0){\line(0,1){20}}
\qbezier(10,20)(0,20)(0,25)
\qbezier(5,30)(0,30)(0,25)
\qbezier(5,30)(10,30)(10,20)
\qbezier(20,20)(30,20)(30,25)
\qbezier(25,30)(30,30)(30,25)
\qbezier(25,30)(20,30)(20,20)
\end{picture}
\vspace{-.1in}
\end{center}
\caption{A connected divide which is not regular-expressive.} 
\vspace{-.2in}
\label{fig:nonexpressive-divide}
\end{figure}

Here is a non-obvious non-example:

\begin{proposition}
\label{pr:not-reg-expr-divide}
The divide shown in Figure~\ref{fig:nonexpressive-divide} is not regular-expressive.
\end{proposition}

\begin{proof}
Suppose on the contrary that the divide~$D$ in Figure~\ref{fig:nonexpressive-divide} is
regular-expressive.
By Theorem \ref{th:reg-expressive}, $D$~must come from a plane curve
$C$ consisting of two polynomial components: $C=K\cup L$.
One of them, say~$L$, is smooth.
By the Abhyankar-Moh theorem \cite[Theorem~1.6]{AM-1975},
there exists a real automorphism of the affine plane that transforms $L$ into a real straight line.
So without loss of generality, we can simply assume that $L$ is a line.
The other component $K$ has degree $d\ge 4$.
Since $C$ is expressive, the projective line $\hat L$
must intersect the projective closure~$\hat K$ of~$K$ at the unique point $p\in K\cap L^\infty$ with multiplicity $d-2$. The curve $\hat K$ has a unique tangent line~$L^\infty$ at~$p$.
Therefore any projective line $\hat L'\ne L^\infty$ passing through $p$ intersects $\hat K$ at $p$ with multiplicity $d-2$. Equivalently, every affine line $L'$ parallel to $L$ intersects $K$ in at most two points (counting multiplicities). However, shifting $L$ in a parallel way until it intersects $K$ at a node for the first time, we obtain a real line $L'$ parallel to $L$ and crossing $K$ in at least four points in the affine plane, a contradiction.
\end{proof}

We recall the following terminology, adapting it to our current needs.

\begin{definition}
\label{def:pseudoline}
A (simple) \emph{pseudoline arrangement} is a connected divide
whose branches are embedded intervals any two of which intersect at most once.

A pseudoline arrangement is called \emph{stretchable} if it is isotopic
to a configuration of straight lines, viewed within a sufficiently large disk.
\end{definition}

\begin{proposition}
\label{pr:non-stretchable-is-not-reg-expr}
Let D be a pseudoline arrangement in which any two pseudolines intersect.
Then the divide D is regular-expressive if and only if it is stretchable.
\end{proposition}

\begin{proof}
The ``if'' direction has already been established, see Example~\ref{ex:line-arrangements}.
Let us prove the converse.
Suppose that a pseudoline arrangement $D$ with $n$ pseudolines $D_1,...,D_n$
is the divide of an $L^\infty$-regular expressive curve~$C$ with real irreducible components.
We need to show that $D$ is stretchable.

By Theorem \ref{th:reg-expressive}, $C$~must consist of $n$ polynomial components $C_1,...,C_n$. Since $C_1$ is smooth, the Abhyankar-Moh theorem \cite[Theorem 1.6]{AM-1975}
implies that a suitable real automorphism of the affine plane
takes $C_1$ to a straight line (and leaves the other components polynomial).
So without loss of generality, we can assume that $C_1=\{x=0\}$.
Let $i\in\{2,...,n\}$ be such that $\deg C_i=d\ge 2$.
(If there is no such~$i$, we are done.)
The projective line $\hat C_1$ intersects the projective curve~$\hat C_i$
at the point $p=(0,1,0)=\hat C_i\cap L^\infty$ either with multiplicity $d-1$
(if $C_1$ and $C_i$ intersect in the affine plane) or with multiplicity~$d$
(if $C_1$ and $C_i$ are disjoint there).
Note that $\hat C_1$ cannot be tangent to $\hat C_i$ at $p$, since $\hat C_i$ is unibranch at $p$, and the infinite line $L^\infty$ is the tangent to $\hat C_i$ at $L^\infty$. It follows that $(\hat C_1\cdot\hat C_i)_p<(\hat C_i\cdot L^\infty)_p=d$ and therefore $(\hat C_1\cdot\hat C_i)_p=d-1$. Since $\hat C_1$ is transversal to $\hat C_i$ at~$p$, we have $\mt(\hat C_i,p)=d-1$.
It~follows that the (affine) equation of $C_i$ does not contain monomials with exponents of $y$ higher than~$1$.
That is, $C_i=\{yP_i(x)=Q_i(x)\}$ with $P_i,Q_i$ coprime real polynomials.
Moreover, the polynomiality of $C_i$ implies that $P_i(x)$ is a nonzero constant, say~$1$.
Finally, recall that any two components $C_i$ and~$C_j$, $2\le i<j\le n$,
intersect in at most one point in the affine plane, and if they do, the intersection is transversal.
So if $C_i=\{y=Q_i(x)\}$ and $C_j=\{y=Q_j(x)\}$, then $\deg(Q_i-Q_j)\le 1$.
It follows that a suitable automorphism of the affine plane $(x,y)\mapsto(x,y-R(x))$, with $R(x)\in\RR[x]$,
leaves $C_1$ invariant while simultaneously taking all other components $C_2,...,C_n$ to straight lines.
\end{proof}


\pagebreak

\section{Expressive curves vs.\ morsifications}
\label{sec:expressive-vs-algebraic}

In this section, we compare two classes of divides:
\begin{itemize}[leftmargin=.2in]
\item
the \emph{algebraic} divides, which arise from real morsifications of isolated singularities of real plane curves,
\redsf{see, e.g., \cite[Definition~2.2]{FPST}};
\item
the \emph{regular-expressive} divides, which arise from $L^\infty$-regular expressive curves with real components,
see Section~\ref{sec:regular-expressive-divides}.
\end{itemize}


There are plenty of divides (such as, e.g., generic real line arrangements)
which are both algebraic and regular-expressive.

\begin{proposition}
\label{pr:reg-expr-not-alg}
None of the following regular-expressive divides is algebraic:
\begin{itemize}[leftmargin=.3in]
\item[\rm{(a)}]
the divides shown in Figures \ref{fig:multi-limacons}-\ref{fig:5-petal};
\item[\rm{(b)}]
the divide shown in Figure \ref{fig:newton-degenerate};
\item [\rm{(c)}]
any divide containing two branches whose intersection is empty.
\end{itemize}
\end{proposition}

\begin{proof}
(a) Let $D$ be the divide shown in Figure~\ref{fig:5-petal} on the right.
Suppose that $D$ is algebraic.
Pick a point $p$ inside the ``most interior'' region of~$D$.
Any real line through $p$ intersects~$D$ (or any curve isotopic to it)
in at least $6$ points, counting multiplicities.
Hence the underlying singularity has multiplicity~$\ge 6$.
Such a singularity must have Milnor number $\ge(6-1)^2=25$.
On the other hand, $D$ exhibits only $21$ critical points
($10$ saddles at the hyperbolic nodes and $11$ extrema, one per region), a contradiction.

The same argument works for the other divides except for the left divides in Figures \ref{fig:3-petal} and \ref{fig:5-petal}, which require a slightly more complicated treatment.
We leave it to the reader as an exercise.

(b) Suppose on the contrary that our divide $D$ is algebraic.
The corresponding singularity must be unibranch, with Milnor number~$12$.
Let $T$ be the sole region of~$D$ whose closure has zero-dimensional intersection
with the union of closures of the unbounded connected components.
Every real straight line crossing~$T$ intersects $D$ in at least $4$ points (counting multiplicities).
It follows that the underlying singularity~has multiplicity~$\ge4$. \redsf{The simplest unibranch singularity 
of multiplicity~$4$ is $y^4-x^5=0$ (up to topological equivalence); its Milnor number is~$12$. 
Hence we cannot encounter any other (more complicated) singularity.} 
The link of the singularity $y^4-x^5=0$ is a $(4,5)$ torus knot.
Its Alexander polynomial is
(see, e.g., \cite[p.~131, formula~(1)]{SW})
\begin{equation}
\label{eq:Alexander-poly-1}
\frac{(1-t^{4\cdot5})(1-t)}{(1-t^4)(1-t^5)}=1-t+t^4-t^6+t^8-t^{11}+t^{12}.
\end{equation}
On the other hand, \redsf{as shown by N.~A'Campo \cite[Theorem 2]{acampo-1999} (reproduced in \cite[Theorem~7.6]{FPST}),
the link of a singularity is isotopic to the link of the divide of its morsification,
as defined by A'Campo; see, e.g., \cite[Definition~7.1]{FPST}.}
In the terminology of \cite[Section 11]{FPST}, the divide $D$ is scannable of multiplicity~$4$.
Its link can be computed as the closure of a 4-strand braid~$\beta$
constructed by the Couture-Perron algorithm \cite[Proposition 2.3]{CP}
(reproduced in \cite[Definition~11.3]{FPST}):
\[
\beta=\sigma_1\sigma_3\sigma_2\sigma_3\sigma_2\sigma_1\sigma_3\sigma_2\sigma_3\sigma_2\sigma_1\sigma_3\sigma_2\sigma_3\sigma_2.
\]
Direct computation shows that the Alexander polynomial of this knot is equal to
\[
(1-t+t^2)(1-t^2+t^4)(1-t^3+t^6)=1-t+t^4-t^5+t^6-t^7+t^8-t^{11}+t^{12},
\]
which is different from~\eqref{eq:Alexander-poly-1}.

(c) Follows from \cite[Proposition 1.8(ii)]{BK}.
\end{proof}

By Proposition~\ref{pr:reg-expr-not-alg}, a connected line arrangement containing parallel lines
is regular-expressive but not algebraic.
We next give an example of the opposite~kind.

\begin{figure}[ht]
\begin{center}
\hspace{.1in}
\setlength{\unitlength}{1pt}
\begin{picture}(190,160)(-10,-30)
\put(0,0){\line(1,1){160}}
\put(20,20){\line(3,1){160}}
\put(20,20){\line(-3,-1){30}}
\put(25,55){\line(9,-1){155}}
\put(25,55){\line(-9,1){35}}
\put(20,20){\line(5,3){160}}
\put(20,20){\line(-5,-3){30}}
\put(80,80){\line(2,-1){100}}
\put(80,80){\line(-2,1){90}}
\put(160,40){\line(-1,1){120}}
\put(160,40){\line(1,-1){20}}
\put(100,0){\line(0,1){160}}
\put(100,140){\line(-1,-3){46.67}}
\put(100,140){\line(1,3){6.67}}
\put(100,140){\line(1,-3){46.67}}
\put(100,140){\line(-1,3){6.67}}
\put(6,-5){\makebox(0,0){5}}
\put(-14,-2){\makebox(0,0){8}}
\put(-12,17){\makebox(0,0){9}}
\put(-9,66){\makebox(0,0){7}}
\put(-7,131){\makebox(0,0){6}}
\put(48,160){\makebox(0,0){4}}
\put(89,164){\makebox(0,0){3}}
\put(100,166){\makebox(0,0){2}}
\put(111,164){\makebox(0,0){1}}

\put(20,20){\red{\circle*{3}}}
\put(70,50){\red{\circle*{3}}}
\put(80,80){\red{\circle*{3}}}
\put(100,46.67){\red{\circle*{3}}}
\put(100,100){\red{\circle*{3}}}
\put(100,140){\red{\circle*{3}}}
\put(120,80){\red{\circle*{3}}}
\put(160,40){\red{\circle*{3}}}
\put(128,56){\red{\circle*{3}}}
\end{picture}
\qquad\quad
\includegraphics[scale=0.25]{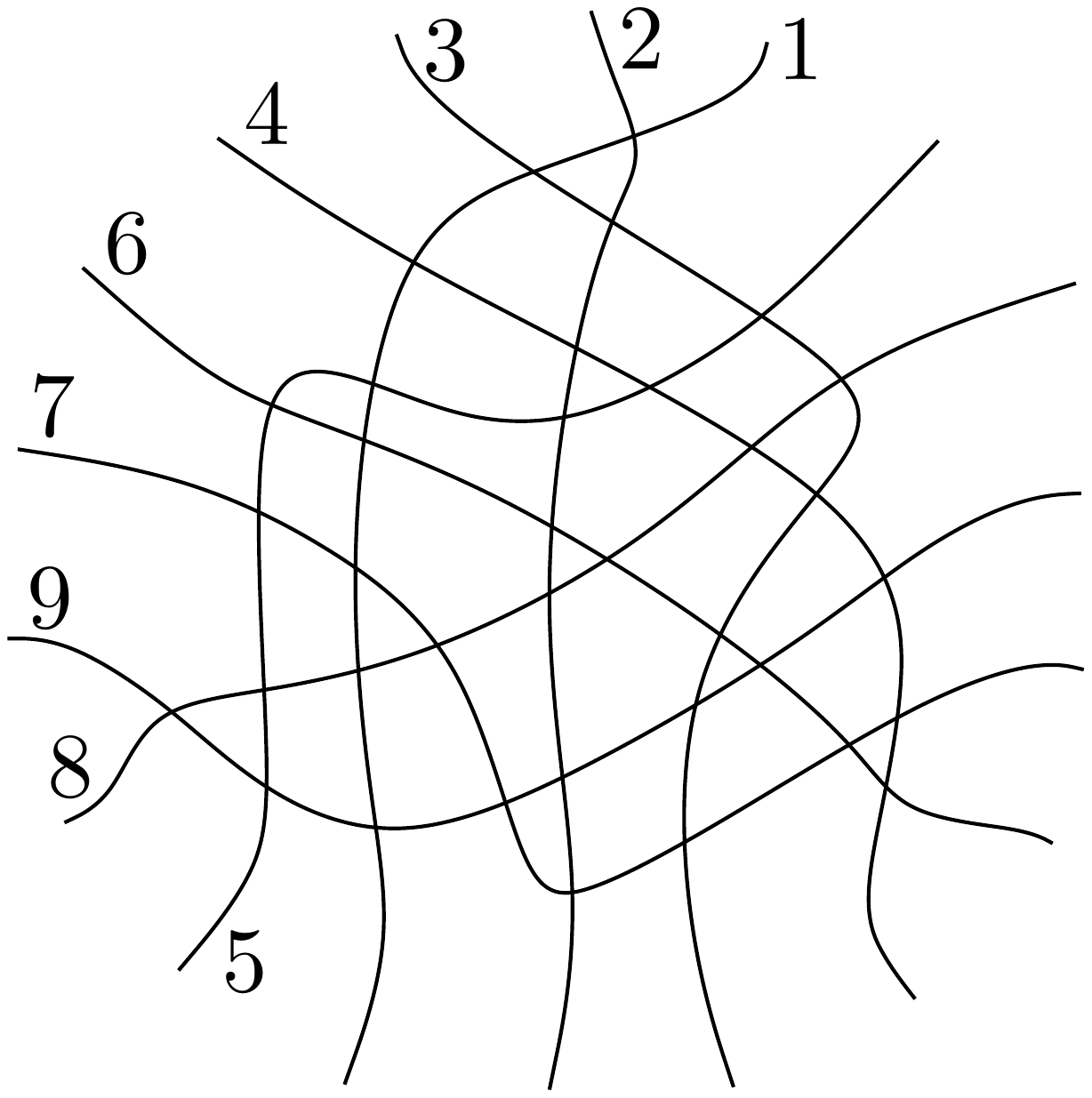}
\hspace{-1.5in}
\setlength{\unitlength}{1pt}
\begin{picture}(100,160)(0,0)
\put(-22,75){\red{\circle*{3}}}
\put(19.5,61){\red{\circle*{3}}}
\put(64.6,72.2){\red{\circle*{3}}}

\put(-2,89){\red{\circle*{3}}}
\put(44.5,83.6){\red{\circle*{3}}}

\put(-9.5,119.5){\red{\circle*{3}}}
\put(26.4,122){\red{\circle*{3}}}
\put(56.5,114.5){\red{\circle*{3}}}

\put(26.7,149.5){\red{\circle*{3}}}
\end{picture}
\end{center}
\caption{Left: a Pappus configuration of straight lines.
Right~\cite[Figure~14]{AKPV}:
a non-Pappus pseudoline arrangement obtained by locally deforming
the Pappus configuration around each of its 9 triple points.
(These points correspond to the 9 marked triangular regions on the right.)
}
\vspace{-.25in}
\label{fig:non-pappus-arr}
\end{figure}

\begin{proposition}
Let $D$ be the ``non-Pappus'' pseudoline arrangement
shown in Figure~\ref{fig:non-pappus-arr} on the right.
The divide $D$ is algebraic but not regular-expressive.
\end{proposition}

\begin{proof}
It is well known that the non-Pappus pseudoline arrangement~$D$ is non-stretchable,
see \cite[Section~3]{AKPV}, \cite[Section~5.3]{PL}.


We next show that the divide $D$ is algebraic.
For the benefit of the readers who may not be experts in singularity theory,
we begin by recalling some background.

Any isolated curve singularity possesses
a versal deformation with finite-dimen\-sion\-al base;
see \cite[Section~II.1.3]{GLS} for a brief account of the theory of versal deformations.
In particular, a versal deformation induces any other deformation. 
Furthermore, if the singular point is real and a versal deformation is conjugation-invariant with a smooth base, then it induces any other conjugation-invariant deformation with a smooth base.
(Indeed, an analytic map $(\CC^n,0)\to(\CC^N,0)$ that takes real points to real points
is given by germs of analytic functions with real coefficients.)

If a singularity (at the origin) is given by $f(x,y)=0$, then the deformation
\[
f(x,y)+\sum_{i+j\le d}t_{ij}x^iy^j=0,\quad (t_{ij})\in\mathbf{B}^N_\eps,\ N=\textstyle\binom{d+2}{2},
\]
is versal if $d\ge\mu(f,0)-1$; here $\mathbf{B}^N_\eps\subset\CC^N$ denotes a sufficently small open disc centered at the origin.

If $C=V(F(x,y))$ is an affine curve and $p_1,...,p_r$ are some of its isolated singularities,
then the deformation
\[
F(x,y)+\sum_{i+j\le d}t_{ij}x^iy^j=0,\quad (t_{ij})\in\mathbf{B}^N_\eps,\ N=\textstyle\binom{d+2}{2},
\]
where $d\ge\sum_i\mu(C,p_i)-1$,
is a joint versal deformation of the singular points $p_1,...,p_r$.
That is, it simultaneously induces arbitrary individual deformations at all these singular points;
see the versality criterion \cite[Theorem~II.1.16]{GLS} and \cite[Proposition~3.4.6]{GLS1}.

Let
\[
\prod_{i=1}^9(a_ix+b_iy+c_i)=0
\]
be a Pappus configuration without parallel lines.
The family
\[
F_t(x,y)=\prod_{i=1}^9(a_ix+b_iy+c_it)=0, \quad t\in[0,1],
\]
is a deformation of the ordinary $9$-fold singular point
$F_0(x,y)=0$, whose members $F_t=0$, $t\ne0$, are Pappus configurations
differing from each other by a homothety.
It is induced by a versal deformation
\begin{align}
\label{eq:versal-pappus}
&F_0(x,y)+\sum_{i+j\le d}t_{ij}x^iy^j=0,\quad (t_{ij})\in\mathbf{B}^N_\eps, \\
\intertext{where}
\notag
&\quad N=\textstyle\binom{d+2}{2}, \quad d=\mu(F_0,0)-1=63.
\end{align}
Since the total Milnor number of the Pappus configuration
(i.e., the sum of Milnor numbers at all singular points)
is less than $\mu(F_0,0)=64$,
the deformation~\eqref{eq:versal-pappus} is
joint versal for all the singularities of any given curve $F_t(x,y)=0$, $0<t\ll1$.

We now construct a morsification of the singularity
$F_0=0$ which is isotopic to our divide~$D$.
We begin by deforming it into a Pappus configuration $F_t=0$ with $0<t\ll1$.
Then, by variation of the monomials up to
degree $63$, we deform each triple point of the latter curve into appropriate
three double intersections (while keeping the nodes of the Pappus configuration).
Since $t$ can be chosen arbitrarily small, the curve isotopic to the constructed one appears arbitrarily close to the original germ $\{F_0=0\}$. In view of the fact that all the strata in the real part of the discriminant in the versal deformation base are semialgebraic sets, by the arc selection lemma \cite[Lemma 3.1]{Milnor}, there exists a real analytic deformation
\[
\widetilde F_t(x,y)=0, \quad 0\le t\ll1
\]
of $\widetilde F_0=F_0=0$, whose members $\widetilde F_t=0$, $t\ne0$, realize the desired morsification.
\end{proof}

\end{document}